\documentclass[a4paper,11pt]{book}

\usepackage[german,english]{babel}
\usepackage{amsmath,latexsym,amssymb,bezier,epic}
\newtheorem{thm}{Theorem}[chapter]
\newtheorem{prop}[thm]{Proposition}
\newtheorem{lemma}[thm]{Lemma}
\newtheorem{cor}[thm]{Corollary}
\newtheorem{main}{Main Theorem}

\newtheorem{conj}{Conjecture}

\newenvironment{defi}
{\stepcounter{thm}\begin{description}\item[Definition
\thechapter.\arabic{thm}:]}
{\end{description}}

\newenvironment{nota}
{\stepcounter{thm}\begin{description}\item[Notation
\thechapter.\arabic{thm}:]}
{\end{description}}

\newenvironment{remark}
{\stepcounter{thm}\begin{description}\item[Remark
\thechapter.\arabic{thm}:]}
{\end{description}}

\newenvironment{examp}
{\stepcounter{thm}\begin{description}\item[Example
\thechapter.\arabic{thm}:]}
{\end{description}}

\newcommand{\proof}[1][]{{\it Proof#1: }}
\newcommand{\qed}[1][3mm]{\hspace*{\fill} $\Box$ \vspace{#1}}

\renewcommand{\P}{{{\mathbf P}^1}}
\newcommand{\NN}{{\mathbf N}}
\newcommand{\ZZ}{{\mathbf Z}}
\newcommand{\CC}{{\mathbf C}}
\newcommand{\PP}{{\mathbf P}}

\newcommand{\RR}{{\mathbf R}}
\newcommand{\EE}{{\mathbf E}}
\newcommand{\LL}{{\mathbf L}}

\newcommand{\del}{\partial}
\newcommand{\ON}{\operatorname}

\renewcommand{\a}{\alpha}
\renewcommand{\b}{\beta}
\newcommand{\cg}{\gamma}
\newcommand{\G}{\Gamma}
\newcommand{\Cg}{\Gamma}
\renewcommand{\d}{\delta}

\newcommand{\e}{\varepsilon}

\newcommand{\varth}{\vartheta}
\renewcommand{\k}{\kappa}
\renewcommand{\l}{\lambda}

\newcommand{\x}{\xi}
\newcommand{\s}{\sigma}
\newcommand{\oo}{\omega}
\renewcommand{\o}{\omega}

\newcommand{\tto}{\longrightarrow}
\newcommand{\inj}{\hookrightarrow}
\newcommand{\surj}{\to\!\!\!\!\!\!\!\!\tto}

\newcommand{\inv}{^{^{-1}}}
\newcommand{\invv}{^{-2}}

\newcommand{\afami}{{\cal A}}
\newcommand{\bfami}{{\cal B}}
\newcommand{\cfami}{{\cal C}}
\newcommand{\dfami}{{\cal D}}
\newcommand{\efami}{{\cal E}}
\newcommand{\ffami}{{\cal F}}
\newcommand{\hfami}{{\cal H}}

\newcommand{\lfami}{{\cal L}}

\newcommand{\ofami}{{\cal O}}

\newcommand{\sfami}{{\cal S}}

\newcommand{\wfami}{{\cal W}}
\newcommand{\xfami}{{\cal X}}

\newcommand{\atil}{\tilde{a}}
\newcommand{\ftil}{\tilde{f}}
\newcommand{\ptil}{\tilde{p}}

\newcommand{\cutoff}[1]{}
\newcommand{\labell}[1]{\label{#1}}
\newcommand{\zspace}[1][-4mm]{\rule[#1]{0mm}{2mm}}

\newcommand{\br}{\ON{Br}}
\newcommand{\pbr}{\ON{PBr}}
\newcommand{\psl}{\ON{PSL}}

\newcommand{\sbr}{\ON{Br}^{\!\!^{s}}}

\newcommand{\lspan}{\langle}
\newcommand{\rspan}{\rangle}

\newcommand{\ethe}{{e_\varth}}
\newcommand{\etheO}{e_{\varth_0}}
\newcommand{\etheI}{e_{\varth_1}}

\newcommand{\vartho}{\varth^o}
\newcommand{\ups}[1]{^{(#1)}}
\newcommand{\jith}{\varth^o}   
\newcommand{\idx}{I_n}

\newcommand{\cable}{\delta_\phi}
\newcommand{\curve}{\cfami}

\def\pr{\operatorname{pr}}
\def\Aut{\operatorname{Aut}}
\def\Hom{\operatorname{Hom}}
\def\id{\operatorname{id}}
\def\conju{\operatorname{conj}}
\def\im{\operatorname{im}}

\newcommand{\permgroup}{\Sigma}

\newcommand{\fund}[1]{\pi_1(#1)}

\newcommand{\piele}{g}  
\newcommand{\point}{p}

\newcommand{\smin}{\setminus}

\newcommand{\diff}{\ON{Diff}}

\newcommand{\mfrak}{{\sf m}}
\newcommand{\mideal}{\mfrak}

\newcommand{\bif}{\ON{Bif}}

\newcommand{\jac}{{J}}

\newcommand{\riso}{\rho}

\newcommand{\pull}{\varphi}

\newcommand{\fk}{{\mathbf F}_k} 
\newcommand{\hirz}{{\mathbf F}}
\newcommand{\ph}{PH^0}
\newcommand{\isec}{C_{-k}} 
\newcommand{\low}{^{\phantom{1}}}
\newcommand{\pnull}{\l}
\newcommand{\qnull}{\xi}
\newcommand{\peins}{\l_1}
\newcommand{\qeins}{\xi_1}
\newcommand{\imag}{\sqrt{-1}}
\newcommand{\bdiff}{\varphi} 
\newcommand{\Phisl}{{\mit\Phi}}
\newcommand{\slz}{\ON{SL_2\ZZ}}

\begin{document}

\begin{titlepage}

\vspace*{20mm}

\begin{center}
\bf
\Huge
Braid Monodromy of Hypersurface Singularities\\[84mm]
\end{center}

\begin{center}
\large
dem Fachbereich Mathematik der Universit\"at Hannover\\
zur Erlangung der {\it venia legendi} f\"ur das Fachgebiet Mathematik\\
vorgelegte Habilitationsschrift\\[33mm]
\end{center}

\begin{center}
\large
von\\
{\bf Michael L\"onne}\\
2003
\end{center}

\end{titlepage}

\tableofcontents

\newpage


\section*{Introduction}
\addcontentsline{toc}{chapter}{introduction}

Complex geometry can certainly be seen as a major source for the
development and
refinement of topological concepts and topological methods.

To exemplify this claim, we like to give to instances, which also
will have impact
on the proper topic of this work.

First there is the paper of Lefschetz on the topology of complex projective
manifolds, which only later were adequately expressed in the language
of algebraic
topology. For example the Picard Lefschetz formula of ordinary double
points is due
to this paper.

Second we want to mention the theorem of van Kampen. It yields, in
quite general
situations, a presentation of the fundamental group of a union of spaces in
terms of
presentations of their fundamental groups. Originally conceived while
investigating
the fundamental group of plane curve complements, it is in its abstract form a
standard topic of basic algebraic topology and a backbone for geometric and
combinatorial group theory.\\[-2mm]

On the other hand new topological concepts are often tested in the
reals of complex
geometry. One may observe that many classifying spaces,
Eilenberg-MacLane space in
particular, have a natural complex structure and can thus be
considered to belong
to complex geometry.\\[4mm]

A prominent example for the fruitful interplay of geometric, topological and
combinatorial methods is singularity theory, into which the present
work has to be
subsumed.\\[-2mm]

Given a holomorphic function $f$ or a holomorphic function germ it is standard
procedure to consider a versal unfolding which is given by a function
$$F(x,z,u)=f(x)-z+\sum b_iu_i.$$

In case of a semi universal unfolding the unfolding dimension is given by the
Milnor number $\mu=\mu(f)$ and we get a diagram
$$
\begin{array}{cccccl}
z,u_1,...,u_{\mu-1} & \CC^\mu & \supset & \dfami & = & \{(z,u)| F(0,z,u)=
0=\nabla F(0,z,u)\}\\
\downarrow & \downarrow \\
u_1,...,u_{\mu-1} & \CC^{\mu-1} & \supset & \bfami & = & \{u|
F(\_,0,u) \text{ is
not Morse}\}
\end{array}
$$
The restriction $p|_\dfami$ of the projection to the discriminant is
a finite map,
such that the branch set coincides with the bifurcation set $\bfami$.

One contribution of the present work is to show, that a suitable
restriction of $p$
to a subset of $p\inv(\CC^{\mu-1}\smin\bfami)\smin\dfami$ is a fibre
bundle in a
natural way. Its fibres a diffeomorphic to the $\mu$-punctured disc and its
isomorphism type depends only on the right equivalence class of $f$.\\[4mm]

When the focus was on the case of simple hypersurface singularities,
this aspect
was not needed, since there is a lot of additional structure one may resort to.

In this case the fundamental groups of discriminant complements of functions of
type ADE are given by the Artin-Brieskorn groups of the same type.
Moreover these
groups have a natural presentation encoded by the Dynkin diagram of that type.

The complements of discriminants and of bifurcation sets were shown to be
Eilenberg-MacLane spaces and homogeneous spaces. Moreover they were related to
natural combinatorial structures via their Weyl groups.

More of this abundance of structure and relations will be used in
chapter four. But
sadly enough it only covers the simple singularities. We can observe
that partial
aspects can be generalized -- especially to parabolic and hyperbolic
singularities
-- but progress to arbitrary singularities has been sparse and slow.

On the other hand, parts of the theory prospered when they became the starting
point of their own theory. Artin Brieskorn groups have lead  to
generalized Artin
groups and the theory of Garside groups now subsumes them into a very
active field
of research.\\[3.5mm]

Having succeeded in describing the discriminant complement in the
case of simple
singularities, Brieskorn, in \cite{Br}, casts a light on some
problems, which he
intended for guidelines to the case of more general singularities. Among other
problems he asked for the fundamental group and suggests to obtain these groups
from a generic plane section using the theorem of Zariski and of van Kampen.
But up to now, only in the case of simply elliptic singularities
presentations of
the fundamental group have been given.\\[3.5mm]

Independently -- initiated by Moishezon two decades ago -- the study
of complements
of plane curves by the methods of Zariski and van Kampen has been
revived and has
found a lot of applications. Conceptionally recast as braid monodromy
theory it has
been successfully used for projective surfaces and symplectic
four-manifolds alike
by investigating branch curves of finite branched maps to $\PP^2$.

The theory of braid monodromy has been generalized to the complements
of hyperplane
arrangements and it has found an interesting new interpretation in the theory
polynomial coverings by Hansen.\\[3.5mm]

The braid monodromy we develop in this work is based on this
interpretation. In its
context the fibre bundle obtained from $p|_\dfami$ naturally gives
rise to a braid
monodromy homomorphism, which then can be made a braid monodromy
invariant of the
unfolded function $f$.

As in the case of plane curves the method of van Kampen leads to an explicit
presentation of the fundamental group of the discriminant complement
$\CC^\mu\smin\dfami$ in terms of generators and relations.\\[3.5mm]

Having accomplished this aim of more theoretical nature, we address next the
problem to find the invariants and the group presentations for
$\pi_1(\CC^\mu\smin\dfami)$ in case of polynomial functions of the
kind given by
\pagebreak
$f(x)=\sum x_i^{l_1+1}$.

Pham investigated this class of function in the spirit of Lefschetz.
He computed
the homology of the regular fibre and then gave the global monodromy
transformation thus generalizing the Picard Lefschetz situation $l_i=1$.

Brieskorn exploited the same class of functions. He showed some of
their links to
be examples of exotic spheres. In his list \cite{Br} of problems he
asks for the
intersection lattice of $f$.

This problem has soon found a solution by a paper of Hefez and Lazzeri
\cite{HL}. Their article has quite an impact on the present work, we
owe them the
description of a Milnor fibre and the choice of a natural geometrically
distinguished path system.

We follow common convention by calling functions $f$ of this class
Brieskorn Pham
polynomials.\\[4mm]

We succeed to solve the Brieskorn problem of three decades ago in one
go for the
large and infinite class of Brieskorn Pham polynomials. Though
generally speaking we
follow the approach suggested by Brieskorn, our method to determine
the presentation
of the fundamental groups deviates in some essential aspects.
To have explicit formula for the bifurcation divisor, we are forced to consider
plane sections of $\CC^\mu$, which fall short of the genericity
conditions in even
several ways. Nevertheless by a substantial amount of additional arguments and
concepts, we finally get the desired results on the braid monodromy.\\[4mm]

The presentations of fundamental groups thus obtained depend on the
Brieskorn Pham
polynomial chosen. They are natural generalizations of the
presentations of Artin
Brieskorn groups associated to the simple singularities.
As in the case of simple singularities we can show, that they are
determined by a
intersection graph of $f$, given in \cite{HL}. Thus a further result has found
an adequate generalization.

Its interesting to note, that also triangles, i.e.\ 2-simplices of the Dynkin
diagram, make their contribution to the relations of the
presentation. Surely one
may expect, that the methods of combinatorial group theory will
eventually provide
a lot of additional properties of these groups.\\[4mm]

With a chapter on elliptic fibrations we want to point to the fact,
that also in
the realm of compact manifolds the concept of braid monodromy may
result in new and
fruitful observations. Elliptic surfaces are good candidates, since in families
almost always the fibration map deforms well, so we can make the singular value
divisor of such a family the object of our braid monodromy
considerations.\\[4mm]

Concerning future developments we may only speculate. Nevertheless
in the presence
of such a lot of open problems we venture to finish our last chapter by some
conjectures, the choice being led by personal interest and the newly gained
insight.
\pagebreak

We like to give a short outline of particular chapters.\\

The first two chapters are mainly of an introductory character. The
first reviews
braid monodromy. We start with braid monodromy of plane curves in the spirit
of Moishezon and proceed like Hansen to get braid monodromy of horizontal
divisors and of affine hypersurface germs. The result of van Kampen on
fundamental groups is developed in each set up. Interspersed we mention results
of Libgober on the complement of plane curve and applications by Moishezon and
Teicher to the theory of branched covers of the projective plane.

In the second we review basic notions of singularity theory. We introduce
discriminant divisors which we consider as a horizontal divisor over
truncated versal unfoldings. We close the chapter with the definition of our
new braid monodromy invariants for right equivalence classes of singular
functions and the implications for the fundamental group of the discriminant
complement.
\\

With the third chapter we enter our computations of the braid monodromy of
Brieskorn Pham polynomials. The equations of the discriminant and the
bifurcation set of their unfoldings by linear polynomials are the main topic of
this chapter. We then define a distinguished system of paths in
regular fibres of a certain kind.

In the forth chapter the special case of singularities of type $A_n$ is
solved and the results prepared for later use in an inductive argument.
\\

The fifth chapter the versal braid monodromy and provides the means to
compute the braid monodromy of Brieskorn Pham polynomials from the versal braid
monodromy of two one-parameter families of functions.

This is computed in the sixth chapter for one of the families in case of
Brieskorn Pham polynomials defined on the plane. We have to develop a big
machinery to distill from our geometric insight the concrete results we want
to prove.

In the seventh chapter we conclude the computation of the braid monodromy by an
inductive argument. Again we have to present more geometric notations and
results.
\\

The eighth chapter is devoted to the study of elliptic surfaces we mentioned
before. We relate each family of elliptic surfaces with a family of divisors in
Hirzebruch surfaces and can thus make use of a detailed study of plane
polynomial functions.
\\

In the final chapter we compute the fundamental group of discriminant
complements in case of Brieskorn Pham polynomials. We consider and prove a
close relationship to the Dynkin diagram found by Pham. Some immediate
corollaries to general function are presented and all these results are used as
motivation for the concluding conjectures.
\\

It is my pleasure to express my thanks to Prof.\ Ebeling, who introduced me to 
the beautiful topic of singularity theory, and
to my colleges in Hannover for their interest and many fruitful discussions.

While special thanks go to Andrea Honecker for the proofreading, I want to thank
my family and all my
friends for constant support.


\chapter{introduction to braid monodromy}

Given a singular curve $\curve$ in the affine plane $\CC^2$ it is natural to
ask for the topology of the complement $\CC^2\smin\curve$. The study of its
fundamental group $\pi_1(\CC^2 \setminus \curve)$
for various types of algebraic curves is a classical
subject going back to the work of Zariski. An algorithm for its computation was
given by van Kampen in
\cite{vK}. It was obtained again by Moishezon as an application of his notion
of braid monodromy, which he introduced in \cite{Mo2} and elaborated with Teicher
in subsequent papers, eg.\ \cite{Mo4,MT}. Libgober \cite{L1} finally proved that
the
$2$-complex associated to the braid monodromy even captures all homotopy
properties of the curve complement
$\CC^2\setminus\curve$.
\\

Before generalizing the considerations to complements of divisors in affine
space, we present the interpretation given in \cite{CS} of the process by which
the braid monodromy of a curve $\curve$ is defined.  It is close to the
approach in
\cite{L1}, but uses a self-contained argument based on Hansen's theory
of polynomial covering maps, \cite{H1}, \cite{H2}.

Given a simple
Weierstrass polynomial $f:X\times \CC\to \CC$ of degree $n$, we consider
the complement of its zero locus $Y=X\times \CC \setminus \{f(x,z)=0\}$.
In Theorem~\ref{braid/bdl}, we show that the projection
$p=\pr_1|_Y: Y\to X$ is a fiber bundle map, with structure group
the braid group $\br_n$, and monodromy the homomorphism from
$\pi_1(X)$ to $\br_n$ induced by the coefficient map of $f$.
\\

This result can be applied when a plane curve $\curve$ is defined by a
polynomial $f$, and $X=\CC\setminus \{y_1,\dots, y_s\}$ is the set of regular
values of a generic linear projection, so that by restriction to $X\times\CC$
the polynomial $f$ becomes simple Weierstrass of degree $n$.
The braid monodromy of $\curve$ is simply the
coefficient homomorphism, $a_*:F_s\to \br_n$.

Obviously $a_*$ depends on the choice of a generic projection, of loops
representing a basis of $F_s$, of an identification of mapping classes with braid
group generators, and of basepoints.
However, the braid-equivalence class of the monodromy -- the double
coset $[a_*]\in \br_s\backslash \Hom(F_s,\br_n)\slash \br_n$,
where $\br_s$ acts on the left by the Artin representation, and $\br_n$ acts on
the right by conjugation -- is uniquely determined by $\curve$.
\pagebreak

\begin{remark}
Recall that the braid monodromy depends not only on the number and types of the
singularities of a curve but is also sensitive to their relative positions as
is shown by the famous example of Zariski \cite{Z1}, \cite{Z2} consisting
of two sextics, both with six cusps, one with all cusps on a conic,
the other not.

It even captures more than the fundamental group of the curve complement as is
shown in \cite{CS}, and one may hope that it detects to some extend the
homeomorphism type of the complement or the ambient homeomorphism type of the
curve.
\end{remark}

When passing to higher dimensions we assume to be given a Weierstrass
polynomial $f:\CC^r\times\CC\to\CC$ defining a horizontal divisor $\dfami$
over $\CC^r$. If $X:=\CC^r\setminus\bfami$ is the set of regular values, the
complement of the bifurcation divisor $\bfami$ of the branched covering
$\dfami\to\CC^r$, then the restriction of $f$ to $X\times\CC$ is a simple
Weierstrass polynomial of degree $n$ equal to the degree of the covering.
The braid monodromy is again the coefficient homomorphism
$a_*:\pi_1(X)\to\br_n$.
Also the method to compute the fundamental group of the plane curve
complement extends to the given situation and provides the tool to get the
fundamental group $\pi_1(\CC^r\setminus\dfami)$.

We push the generalization even further to include the case of analytic
germs. With a generic choice of local coordinates the Weierstrass preparation
theorem can be applied and provides us with a Weierstrass polynomial
which is simple in the complement of the germ of a divisor. Again the
subsequent definitions generalize.

\section{Polynomial covers and $\br_n$-bundles}

We begin by reviewing polynomial covering maps.  These
were introduced by Hansen in \cite{H1}, and studied to some detail in
his book \cite{H2}. Together with the by now classical
book of Birman
\cite{Bi} it should serve as the basic
reference for this section.  We then consider the relation between bundles of
punctured discs, whose structure group is Artin's braid group $\br_n$, and
polynomial $n$-fold covers.

\subsection{Polynomial covers}

Let $X$ be a path-connected topological manifold.  A {\it Weierstrass
polynomial} of degree $n$ is a map $f: X\times \CC \to \CC$ given by
$$
f(x,z)=z^n+\sum_{i=1}^n a_i(x) z^{n-i},
$$
with continuous coefficient maps $a_i:X\to \CC$. If $f$ has no multiple roots
for any $x \in X$, then $f$ is called a {\it simple Weierstrass polynomial}.

Given such $f$, the restriction of the first-coordinate projection map $X
\times \CC \to X$ to the subspace
$$
E=E(f)=\{ (x,z)\in X\times \CC \mid f(x,z)=0\}
$$
defines an $n$-fold topological cover $\pi=\pi_f:E \to X$, the
{\it polynomial covering map} associated to $f$.
\pagebreak

Since $f$ has no multiple roots, the {\it coefficient map}
$$
a=(a_1,\dots ,a_n):X \to \CC^n
$$
takes values in the complement $B^n=\CC^n \setminus \Delta_n$ of the discriminant
set $\Delta_n$, which is a tautology by the definition of $\Delta_n$ as the set
of coefficient $n$-tuples such that the corresponding polynomial of degree $n$
has at least one multiple root.

Over $B^n$, there is a tautological $n$-fold
polynomial covering
\begin{equation}
\label{taut}
\pi_n:=\pi_{_{f_n}}: E(f_n) \to B^n,
\end{equation}
determined by the tautological
Weierstrass polynomial
$$
f_n:\CC^n\times\CC\tto\CC,\quad(x_1,...,x_n,z)\mapsto z^n + \sum_{i=1}^n x_i
z^{n-i}.
$$
The polynomial cover $\pi_f:E(f)\to X$ can then be identified with the
pull-back of $\pi_n:E(f_n) \to B^n$ along the coefficient map
$a: X\to B^n$.

This can be interpreted on the level of fundamental groups
as follows.  The fundamental group of the configuration space,
$B^n$, of $n$ unordered points in $\CC$ is the group $\br_n$
of braids on $n$ strands.
The map $a$ determines the {\it coefficient homomorphism}
$a_*: \pi_1(X) \to \br_n$, unique up to conjugacy.

Recall that the structure group of a topological $n$-sheeted cover is the
permutation group $\permgroup_n$ and the associated {\it cover monodromy} is a
homomorphism from the fundamental group of the base to $\permgroup_n$. So it is
immediate that the monodromy of the polynomial cover
$\pi:E\to X$ factors through the coefficient homomorphism
$a_*$ and the canonical surjection $\br_n\to \permgroup_n$. In fact this
condition is sufficient, i.e.\ one
may characterize polynomial covers by this factorization property of their
permutation monodromy map.

Assume now that the simple Weierstrass polynomial $f$
is {\it completely solvable}, i.e.\ $f$ factors as
$$
f(x,z)=\prod_{i=1}^n (z-b_i(x)),
$$
with continuous roots $b_i: X \to \CC$.  Since the Weierstrass
polynomial $f$ is simple, the {\it root map} $b=(b_1,\dots ,b_n):X\to
\CC^n$ takes values in the complement $P^n$ of the
braid arrangement  
\zspace[-2mm]
$\afami_n=\{(w_1,...,w_n)|w_i\neq w_j\,\forall i<j\}$ in $\CC^n$.
Over $P^n$, there is a canonical $n$-fold covering map, 
\zspace[-2mm]
$\tilde\pi_n=\pi_{\tilde f_n}: E(\tilde f_n) \to P^n$, determined by the
\zspace[-2mm]
Weierstrass polynomial $\tilde f_n(w,z)=(z-w_1)\cdots (z-w_n)$.
Evidently, the cover $\pi_f:E\to X$ is the pull-back of 
$\tilde \pi_n: E(\tilde f_n)\to P^n$ along the root map $b: X\to P^n$.

The fundamental group of the configuration space,
$P^n$, of $n$ ordered points in $\CC$ is the group,
$\pbr_n=\ker(\br_n\to\Sigma_n)$, of pure braids on $n$ strands.
The map $b$ determines
the {\it root homomorphism} $b_*: \pi_1(X) \to \pbr_n$,
unique up to conjugacy.
The polynomial covers which are trivial covers (in the usual sense)
are precisely those for which the coefficient homomorphism
$a_*$ has image in the subgroup $\pbr_n$ of $\br_n$.

\subsection{$\br_n$-Bundles}

The group $\br_n$ is isomorphic the mapping class group 
$\ON{Map}^n(D^2)$ of orientation-preserving diffeomorphisms of 
the disc $D^2$, permuting a collection of $n$ marked points.
A natural way to fix an isomorphism is the choice of a 'frame' as in \cite{MT}.
But it pays off immediately if we choose instead a different geometric object
(which incidentally is called  'bush' in \cite{MT}).

\begin{defi}
A {\it geometrically distinguished} system of paths in $D^2$ with respect to
marked points $y_i,1\leq i\leq n,$ is a finite sequence of paths $p_i:[0,1]\to
D^2$, such that:
\begin{enumerate}
\item
$y_0:=p_i(0)$ is the unique point in the intersection of two paths,
$p_i(1)=y_i$,
\item
on some loop around $y_0$ in $\RR^2$ the intersection points with the $p_i$ are
in bijection with the indices $i=1,...,n$, preserving the order.
\end{enumerate}

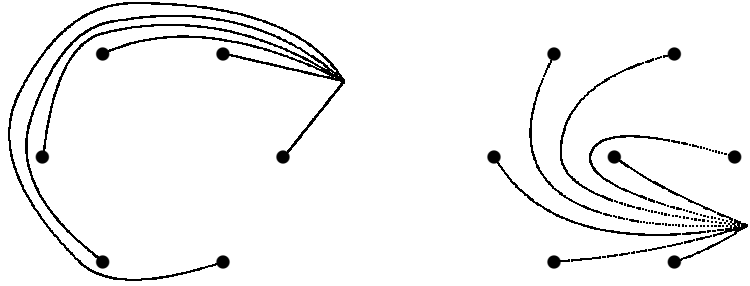
\begin{figure}[h]
\begin{center}
\setlength{\unitlength}{2mm}
\begin{picture}(60,20)(-9,1)

\put(-2,0){\setlength{\unitlength}{1.6mm}
\begin{picture}(30,20)(0,0)
\put(10,5){\circle*{1}}
\put(20,5){\circle*{1}}
\put(5,13.7){\circle*{1}}
\put(10,22.3){\circle*{1}}
\put(25,13.7){\circle*{1}}
\put(20,22.3){\circle*{1}}

\bezier{140}(30,20)(27.5,16.85)(25,13.7)
\bezier{140}(30,20)(25,21.15)(20,22.3)
\bezier{140}(30,20)(18,26)(10,22.3)
\bezier{140}(30,20)(20,26.5)(10,24)
\bezier{140}(10,24)(6,23)(5,13.7)
\bezier{140}(30,20)(22.25,27)(11,25)
\bezier{140}(11,25)(6.5,24.5)(4,17)
\bezier{140}(4,17)(2,11)(10,5)
\bezier{140}(30,20)(26,26.5)(12.5,26.5)
\bezier{140}(12.5,26.5)(6.75,26.5)(3,19)
\bezier{140}(3,19)(0,13)(8,5)
\bezier{140}(8,5)(11,2)(20,5)
\end{picture}
}

\put(28,0){\setlength{\unitlength}{1.6mm}
\begin{picture}(30,20)
\put(10,5){\circle*{1}}
\put(20,5){\circle*{1}}
\put(5,13.7){\circle*{1}}
\put(15,13.7){\circle*{1}}
\put(10,22.3){\circle*{1}}
\put(25,13.7){\circle*{1}}
\put(20,22.3){\circle*{1}}

\bezier{40}(26,8)(23,6)(20,5)
\bezier{80}(26,8)(18,5.5)(10,5)
\bezier{70}(26,8)(18.5,11)(15,13.7)
\bezier{140}(26,8)(10,5)(5,13.7)
\bezier{140}(26,8)(2,7)(10,22.3)
\bezier{80}(26,8)(10.5,9.5)(10.5,14)
\bezier{80}(10.5,14)(10.5,20)(20,22.3)
\bezier{80}(26,8)(12,11)(13,14)
\bezier{80}(13,14)(14,17)(25,13.7)
\end{picture}
}

\end{picture}
\end{center}

\caption{Examples of geometrically distinguished path systems}
\end{figure}

They give rise to a basis $t_1,...,t_n$ of free generators of
$\pi_1(D^2\setminus\{y_1,...,y_n\},y_0)$, which accordingly is called {\it
geometrically distinguished}, too.
\end{defi}

The free basis of the fundamental group is unambiguously represented by loops
which are are obtained by replacing each path $p_i$ by a sufficiently close
noose embedded into $D^2\setminus\{y_1,...,y_n\}$, based at $y_0$, and linking
the marked point $y_i$ once.

\begin{figure}[h]
\setlength{\unitlength}{.5mm}
\begin{picture}(60,33)

\put(180,14){\begin{picture}(5,2)
\put(5,5){\circle*{2}}

\bezier{300}(0,5)(0,7.5)(2,9)
\bezier{300}(2,9)(5,11)(8,9)
\bezier{300}(8,9)(10,7.5)(10,5)
\bezier{300}(10,5)(10,2.5)(8,1)
\bezier{300}(8,1)(5,-1)(2,1)
\bezier{300}(2,1)(0,2.5)(0,5)

\bezier{300}(8,1)(20,-5)(70,-15)

\end{picture}}

\put(30,14){\begin{picture}(5,0)
\put(5,5){\circle*{2}}

\bezier{310}(5,5)(30,-5)(70,-15)

\end{picture}}

\end{picture}
\caption{Noose associated to a path}
\end{figure}
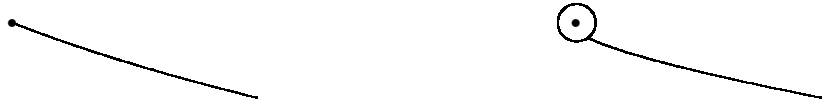

Upon identifying $\pi_1(D^2 \setminus \{ n \text{ points}\})$ in such a way
with the free group $F_n$, the action of $\br_n$ on $\pi_1$ yields the 
{\it Artin representation}, $\a_n:\br_n\to \Aut(F_n)$.  As shown by 
Artin, this representation is faithful.  Hence, we 
may identify a braid $\b\in \br_n$ with the 
corresponding braid automorphism, $\a_n(\b)\in \Aut(F_n)$.
\\

Now let $f:X\times \CC\to \CC$ be a simple Weierstrass polynomial.
Let $\pi_f: E(f)\to X$ be the corresponding polynomial
$n$-fold cover, and $a:X \to B^n$ the coefficient map.
Consider the complement
$$
Y=Y(f)=X\times \CC \setminus E(f),
$$
and let $p=p_f: Y(f)\to X$ be the restriction of
$\pr_1:X\times \CC \to X$ to $Y$.

The fibre over a single point we denote by
$\CC_n:=\CC \setminus \{ n \text{ points}\}$.

\begin{thm}
\labell{braid/bdl}
The map $p:Y\to X$ is a locally trivial bundle, with structure group
$\br_n$ and fiber $\CC_n$.  Upon
identifying $\pi_1(\CC_n)$ with $F_n$, the monodromy
of this bundle may be written as $\a_n\circ a_*$, where
$a_*:\pi_1(X)\to \br_n$ is the coefficient homomorphism.

Moreover, if $f$ is completely solvable, the structure group
reduces to $\pbr_n$, and the monodromy is $\a_n\circ b_*$,
where $b_*:\pi_1(X)\to \pbr_n$ is the root homomorphism.
\end{thm}

\proof  We first prove the theorem for the configuration
spaces, and their canonical Weierstrass polynomials.
Start with $X=P^n$, $f=\tilde f_n$, and the canonical cover $\tilde \pi_n:
E(\tilde f_n)\to P^n$. Clearly, $Y(\tilde f_n)=\CC^{n+1} \setminus E(\tilde
f_n)$ is equal to the configuration space $P^{n+1}$.
Let $p_{\ftil_n}:P^{n+1}\to P^n$ be the restriction of
$\pr_1:\CC^n\times \CC\to \CC^n$.  As shown by Faddell and Neuwirth~\cite{FN},
this is a bundle map, with fiber $\CC_n$, and monodromy the restriction of
the Artin representation to $\pbr_n$.

Next consider $X=B^n$, $f=f_n$, and the canonical cover
$\pi_n:E(f_n) \to B^n$. Forgetting the order of the points defines
a covering projection from the ordered to the unordered configuration
space, $\kappa_n:P^n \to B^n$.

In coordinates, $\kappa_n(w_1,\dots ,w_n)=(x_1,\dots ,x_n)$, where
$x_i=(-1)^i s_i(w_1,\dots ,w_n)$, and $s_i$ are the elementary
symmetric functions.  By Vieta's formulas, we have
$$
\ftil_n(w,z) = f_n(\kappa_n(w),z).
$$
Let $Y^{n+1}=Y(f_n)$ and $p_{f_n}: Y^{n+1} \to B^n$.
By the above formula, we see that
$\kappa_n\times \id :P^n\times \CC \to B^n\times \CC$ restricts
to a map $\bar\kappa_{n+1}:Y(\ftil_n) \to Y(f_n)$, which fits into
the fiber product diagram
$$
\begin{array}{ccl}
P^{n+1}       & \stackrel{p_{\ftil_n}}{\tto} & P^n \\
\quad\downarrow\bar\kappa_{n+1}   & &   \downarrow{\kappa_n}\\
Y^{n+1}      & \stackrel{p_{f_n}}{\tto} &     B^n
\end{array}
$$
where the vertical maps are principal $\Sigma_n$-bundles.

Since the bundle map $p_{\ftil_n}: P^{n+1}\to P^n$ is naturally equivariant with
respect to the $\Sigma_n$-actions, the map on quotients,
$p_{f_n}: Y^{n+1} \to B^n$, is also a bundle map, with fiber $\CC_n$,
and monodromy action the Artin representation of $\br_n$.
This finishes the proof in the case of the canonical
Weierstrass polynomials over configuration spaces.

Now let $f:X\times \CC\to \CC$ be an arbitrary simple Weierstrass
polynomial.  We then have the following Cartesian square:
$$
\begin{array}{ccl}
Y      & \tto &   Y^{n+1} \\[1mm]
\quad\big\downarrow\, p\,   & &      \big\downarrow \,{p_{f_n}} \\[1mm]
X      & \tto &   B^n
\end{array}
\pagebreak
$$
In other words, $p: Y\to X$ is the pull back of the bundle
$p_{f_n}: Y^{n+1}\to B^n$ along the coefficient map $a$.
Thus, $p$ is a bundle map, with fiber $\CC_n$, and monodromy
representation $\a_n\circ a_*$.  When $f$ is completely
solvable, the bundle $p: Y\to X$ is the pull back of
$p_{\ftil_n}: P^{n+1}\to P^n$ along the root map $b$.
Since $a_*=b_*$ then, the monodromy is as claimed.
\qed[3mm]

Let us summarize the above discussion of braid bundles over
configuration spaces.  From the Faddell-Neuwirth
theorem \cite{FN}, it follows that $P^n$ is a $K(\pi_1,1)$ space.
Since $B^n$ is covered by $P^n$, it is also an $K(\pi_1,1)$ space.

So in both cases the groups are discrete and any of their representations as
groups of diffeomorphisms determines a locally trivial bundle over the
corresponding configuration space $P^n$, resp.\ $B^n$.

\begin{examp}
We considered two bundles over $P^n$ obtained in this way:
\begin{itemize}
\item[(i)] $\tilde\pi_n:E(\ftil_n) \to P^n$, by the trivial representation of
$\pbr_n$ on $\{1,\dots ,n\}$;
\item[(ii)] $p_{\ftil_n}:P^{n+1} \to P^n$, by the geometric Artin
representation of $\pbr_n$ on $\CC_n$.
\end{itemize}
Over $B^n$ we have even seen three instances of bundles of this kind:
\begin{itemize}
\item[(iii)] $\kappa_n:P^n \to B^n$, by the canonical surjection
$\br_n\to \Sigma_n$, acting by left translations on the discrete set $\Sigma_n$,
\item[(iv)] $\pi_n:E(f_n) \to B^n$, by the above, followed by the
permutation representation of $\Sigma_n$ on $\{1,\dots ,n\}$;
\item[(v)] $p_n:Y^{n+1} \to B^n$, by the (geometric) Artin
representation of $\br_n$ on $\CC_n$.
\end{itemize}
\end{examp}

Note finally, that $\pi_1(Y^{n+1})$ is isomorphic to
$\br_{n,1}=F_n\times_{\a_n} \br_n$, the group of braids on $n+1$
strands that fix the endpoint of the last strand, and that
$Y^{n+1}$ is a $K(\br_{n,1},1)$ space.


\section{The braid monodromy of a plane algebraic curve}

We are now ready to define the braid monodromy of an algebraic
curve in the complex plane.  The construction, based on
classical work of Zariski and van~Kampen, is due to
Moishezon \cite{Mo2}.
Though we want to include the results of Libgober \cite{L1}, \cite{L2},
\cite{L3}, we prefer to interpret the construction in the context
established in the previous section.

\subsection{The construction}

Let $\curve$ be a reduced curve in the affine plane $\CC^2$ defined by a
polynomial $f$.  Let $\pi:\CC^2\to \CC$ be a linear projection, such that
no fibre of $\pi$ is a component of $\curve$,
and let $\sfami =\{y_1,\dots, y_s\}$ be the set of points in $\CC$ for which
the corresponding fiber of $\pi$
is tangent to $\curve$ or contains a singular point of $\curve$.

In case we assume $\pi$ to be {\it generic} with respect to $\curve$, we mean
that, for each $k$, the line $\lfami_k=\pi^{-1}(y_k)$ passes
transversally through regular points of $\cfami$ except for either a single
regular point $v_k$ at which it is a simple tangent or for one singular point
$v_k$ at which it is not contained in the tangent cone.
\pagebreak

Let $\lfami$ denote the union
of the lines $\lfami_k$, and let $y_0$ be a basepoint in $\CC\setminus \sfami$.
The definition of the braid monodromy of $\curve$ depends on two observations:
\begin{enumerate}
\item[(i)] The restriction of the projection map,
$p:\CC^2 \setminus \curve \cup \lfami \to \CC \setminus \sfami$,
is a locally trivial bundle.
\end{enumerate}

Identify the fiber $p^{-1}(y_0)$ with $\CC_n$ for the appropriate $n$ and fix a
basepoint $v_0\in\CC_n$. The monodromy of $\curve$ is, by definition, the
holonomy of the bundle,
$\rho: \pi_1(\CC \setminus \sfami,y_0) \to \Aut(\pi_1(\CC_n,v_0))$.
Upon identifying $\pi_1(\CC \setminus \sfami,y_0)$ with $F_s$,
and $\pi_1(\CC_n,v_0)$ with $F_n$ using geometrically distinguished systems of
paths with respect to the exceptional points, this can be written as $\rho: F_s
\to
\Aut(F_n)$.
\begin{enumerate}
\item[(ii)] The image of $\rho$ is contained in the
braid group $\br_n$ (viewed as a subgroup of $\Aut(F_n)$ via the Artin
embedding $\a_n$).
\end{enumerate}

\begin{defi}
The homomorphism $\a:F_s\to\br_n$ determined by $\a_n\circ\a=\rho$ is called
{\it braid monodromy homomorphism} of $\curve$ {\it with respect
to $\pi$}.

In case $\pi$ is generic, $\a:F_s\to\br_n$ is simply called {\it braid
monodromy homomorphism} of $\curve$.
\end{defi}

We shall present a self-contained proof of the two assertions,
and, in the process, identify the map $\a$.  The first
assertion is well-known, and can also be proved by standard techniques
(using blow-up and Ehresmann's criterion---see \cite{Di}, page~123),
but our approach sheds some light on the underlying
topology of the situation.

We may assume
-- after a linear change of variables in $\CC^2$ if necessary -- that
$\pi$ is the projection map $\pr_1$ onto the first coordinate.
In the chosen coordinates, the defining polynomial $f$ of $\curve$
may be written as
$f(x,z)=z^n+\sum_{i=1}^n a_i(x) z^{n-i}$.
Since $\curve$ is reduced, for each $x\notin \sfami$, the equation $f(x,z)=0$
has $n$ distinct roots.  Thus $f$ is a simple Weierstrass polynomial
over $\CC\setminus \sfami$, and
\begin{equation}
\label{st}
\pi=\pi_{f}:\curve \setminus \lfami \cap \curve \to \CC\setminus \sfami
\end{equation}
is the associated polynomial $n$-fold cover.

Note that
$Y(f)=((\CC \setminus \sfami)\times \CC) \setminus (\curve \setminus
\lfami \cap \curve)$ is equal to $\CC^2 \setminus \curve \cup \lfami$.  By
theorem~\ref{braid/bdl}, the restriction of $\pr_1$ to $Y(f)$,
\begin{equation}
\label{nd}
p:\CC^2 \setminus \curve \cup \lfami \to \CC \setminus \sfami,
\end{equation}
is a bundle map, with structure group $\br_n$, fiber $\CC_n$, and
monodromy homomorphism
\begin{equation}
\label{rd}
\a=a_*:\pi_1(\CC\setminus \sfami) \to \br_n.
\end{equation}
This proves assertions (i) and (ii), and implies immediately:

\begin{prop}
\labell{bm/ch}
The braid monodromy of a plane algebraic
curve coincides with the coefficient homomorphism of the
associated polynomial cover.
\end{prop}

\subsection{Braid equivalence}

The braid monodromy of a plane algebraic curve is
not unique, but rather, depends on the choices made in defining it.
This indeterminacy was studied to its full extend by Libgober in \cite{L2},
\cite{L3}. To make the analysis more precise, we first need a definition.

\begin{defi}
Two homomorphisms $\a:F_s\to \br_n$ and
$\a':F_s\to \br_n$ are {\it braid equivalent} if there exist
an automorphisms $\psi\in \Aut(F_s)$ in the image of the Artin representation 
and $\phi\in \br_n$ such that
$\a'(\psi(g)) = \phi^{-1}\cdot \a(g)\cdot \phi$,
for all $g\in F_s$.  In other words, the following diagram commutes
$$
\begin{array}{ccl}
F_s      &  \stackrel{\a}{\tto} &     \br_n \\
\,\,\,\downarrow{\psi}  & &   \downarrow\conju_{\phi} \\
F_s   &     \stackrel{\a'}{\tto}   &  \br_n.
\end{array}
$$
\end{defi}

\begin{thm}  
\labell{well}
The braid monodromy of a plane algebraic
curve $\curve$ is well-defined up to braid-equivalence, and so are the braid
monodromies with respect to a fixed linear projection.
\end{thm}

\proof
First fix the projection.
The identification $\pi_1(\CC\setminus \sfami)=F_s$ depends on the
choice of a distinguished system of paths, and any two such choices
yield monodromies which differ by a braid automorphism of $F_s$,
see \cite{L2}.  Furthermore, there is the choice of basepoints,
and any two such choices yield monodromies differing by a conjugation
in $\br_n$.

Finally, one must analyze the effect of a change in the choice of
generic projection.   Let $\pi$ and $\pi'$ be two such projections,
with critical sets $\sfami$ and $\sfami'$, and braid monodromies
$a_*:\pi_1(\CC\setminus \sfami)\to \br_n$ and  $a_*':\pi_1(\CC\setminus
\sfami')\to \br_n$. Libgober \cite{L3} shows that there is a homeomorphism
$h:\CC\to\CC$, isotopic to the identity, and taking $\sfami$ to $\sfami'$,
for which the isomorphism
$h_*:\pi_1(\CC \setminus \sfami) \to \pi_1(\CC \setminus \sfami')$
induced by the restriction of $h$ satisfies $a_*'\circ h_* =a_*$.
From the construction, we see that $h$ can be taken to be the
identity outside a ball of large radius (containing $\sfami\cup \sfami'$).
Thus, once the identifications of source and target with $F_s$ using
distinguished systems of paths are made,
$h_*$ can be written as a braid automorphism of $F_s$:  $h_*=
\psi$, since $\br_n$ acts transitively on the isotopy classes of distinguished
path systems. We obtain $a_*'\circ \psi= a_*$,
completing the proof.
\qed[3mm]

Thus, we may regard the braid monodromy of $\curve$ as a braid-equivalence
class, i.e., as a double coset
$[a_*]\in \br_s\backslash \Hom(F_s,\br_n)\slash \br_n$, uniquely determined
by $\curve$.  In fact, it follows from \cite{L3} that $[a_*]$ depends only
on the equisingular isotopy class of the curve.


\section{The fundamental group of a plane algebraic curve}

We now give the braid monodromy presentation
of the fundamental group of the complement of a plane algebraic
curve $\curve$. This presentation first appeared in the classical work
of van~Kampen and Zariski \cite{vK}, \cite{Z2}, and has been much
studied since, e.g.\ by Moishezon, Teicher \cite{MT,Te1},
Libgober \cite{L1,L2},
Rudolph
\cite{Ru} and many more.

\subsection{Braid monodromy presentation}
\labell{bmp}

We want to find a presentation of $\pi_1(\CC^2\setminus\curve)$. As a first
approximation we give a presentation of $\pi_1(\CC^2\setminus\curve\cup\lfami)$
which can be derived with the help of the previous discussions. The essential
step in then to extract enough information out of the embedding
$\CC^2\setminus\curve\cup\lfami\to\CC^2\setminus\curve$ to determine the
presentation of $\pi_1(\CC^2\setminus\curve)$ sought for.

The first necessary observation is the following:

\begin{lemma}
\labell{sect}
Given a Weierstrass polynomial $f:X\times\CC\to\CC$ there is a topological
section $s:X\to X\times\CC$ to the projection $X\times\CC\stackrel{\ON{pr}}\to
X$ with image in the complement of the zero set of $f$.
\end{lemma}

\proof
It is a well known fact, that the zeroes of a monic polynomial are bounded by
the sum of the absolute values of all coefficients.
Since the coefficients are continuous functions on $X$, this bound -- considered
as a complex valued function on $X$ -- defines a continuous section with image
disjoint from
$f\inv(0)$.
\qed[3mm]

Since the $s$-punctured complex line $\CC\setminus\sfami\cong\CC_s$ is a
$K(F_s,1)$, the long exact homotopy sequence of the bundle
$p:\CC^2 \setminus \curve \cup \lfami \to \CC \setminus \sfami$
yields a short exact sequence which is split due to the preceding lemma:
$$
1\to \pi_1(\CC_n) \to \pi_1(\CC^2 \setminus \curve \cup \lfami)
\stackrel{s_*}{\stackrel{\textstyle{\longleftarrow}}{\tto}}
{}_{\hspace*{-5.7mm}p_*\hspace*{3mm}}
\pi_1(\CC \setminus \sfami)\to 1.
$$
Moreover the action is given by
the braid monodromy homomorphism $a_*$ of (\ref{rd}),
so $\pi_1(\CC^2\setminus\cfami\cup\lfami)$ is the semi-direct product via
$\a_*$, a presentation of which can be derived from presentations of
$\pi_1(\CC_n)$ and $\pi_1(\CC\setminus\sfami)$.

As remarked before, these groups are free, but in the given geometric set up
the identifications with an abstractly presented free groups
$F_n=\lspan t_1,...,t_n|\,\rspan$ and $F_s=\lspan \cg_1,...,\cg_s|\,\rspan$
determined by geometrically distinguished bases are privileged.
\\

{\sloppypar
Having chosen such geometrically distinguished bases for $\CC_n$ and
$\CC\setminus\sfami$, which amounts to a distinguished choice of
isomorphisms $\pi_1(\CC_n)\cong\lspan t_1,...,t_n|\,\rspan$ and
\mbox{$\pi_1(\CC\smin\sfami)\cong\lspan \cg_1,...,\cg_s|\,\rspan$}, the split
sequence above naturally induces an isomorphism
$$
\pi_1(\CC^2 \setminus \curve \cup \lfami)=
\langle t_1,\dots t_n,\cg_1\dots ,\cg_s \mid
\cg_k^{-1}\cdot t_i\cdot \cg_k = a_*(\cg_k)(t_i)\rangle.
$$
}

To proceed we are in need of a result relating the fundamental groups given the
injective map $\CC^2\setminus\curve\cup\lfami\to\CC^2\setminus\curve$, where
$\lfami$ is a divisor with no component in common with $\curve$.
This is very much in the spirit of the Zariski van Kampen results although we owe
it mostly to \cite{Be}, to which we have added only a distinctive topological
flavour.

\begin{defi}
Let $D$ be a reduced divisor in affine space $\CC^n$. An element $\piele$ of
$\fund{\CC^n\setminus D}$ is called a {\it simple geometric element} if there is
an embedded oriented disc in $\CC^n$, which intersects $D$ transversally in a
unique point, such that the orientations of $D$ and the disc give the
orientation of
$\CC^n$ and such that the oriented boundary is freely homotopic to $\piele$.\\
\end{defi}

\begin{defi}
Let $D$ in $\CC^n$ be a reduced divisor and $\point$ a point of $\CC^n$ in the
complement of $D$. Let $D_0$ be an irreducible component of $D$. Then an element
$\piele$ of
$\fund{\CC^n\setminus D,\point}$ is called the {\it simple geometric element
associated to $D_0$} if $\piele$ is freely homotopic to the oriented boundary
of a disc intersecting $D$ transversally in a unique point of $D_0$.
\end{defi}

There are of course simple geometric elements in abundance, examples of which
are provided by the nooses associated to a system of paths.

\begin{lemma}
\labell{conjugate}
The simple geometric elements associated to the same component are
conjugate and any element conjugate to a simple geometric one is itself a
simple geometric element associated to the same component.
\end{lemma}

\proof
The open part $D_{reg}$ of the given component consisting of points which are
regular and not contained in any other component is path connected.
So if transversal discs to points in $D_{reg}$ are given, these points are
connected by a path embedded in $D_{reg}$. Along this path $D_{reg}$ is a
submanifold of real codimension two, so a normal disc bundle exists, which
shows that the boundaries of both discs are freely homotopic. Hence the
geometric elements are freely homotopic as well, which implies the first claim.
The second is obvious.
\qed[0mm]

\begin{lemma}
\labell{lbessis}
Let $U$ be a smooth connected complex variety. 
Let $D_1$ and $D_2$ be divisors with no irreducible component in common.
Then the naturally induced map $\fund{U-D_1\cup D_2}\to\fund{U-D_1}$ is surjective
and
\begin{enumerate}
\item
For any simple geometric element in $\fund{U-D_1}$ associated to an irreducible
component $J$ of $D_1$, there is a lift in $\fund{U-D_1\cup D_2}$ which is a simple
geometric element associated to $J$.
\item
The simple geometric elements of $\fund{U-D_1\cup D_2}$ associated to $D_2$
generate the kernel.
\end{enumerate}
\end{lemma}

\proof
Since a path is of real dimension 1 and the divisor $D_2$ is of real codimension
two, any path in the complement of $D_1$ is isotopic in $U-D_1$ to a path
disjoint from $D_2$, so surjectivity holds as claimed.

To address the claim $i)$, let any simple geometric element in $\pi_1(U-D_1)$ be
given which is associated to an irreducible component $J$ of $D_1$. Choose any
simple geometric element in $\pi_1(U-D)$ associated to $J$. Its image in
$\pi_1(U-D_1)$ is still a simple geometric element associated to $J$, hence
freely homotopic to the given one. By \ref{conjugate} they are even conjugate.

Since surjectivity is already established we may find a preimage of the
conjugating element by which we can conjugate the chosen geometric element to
get a preimage of the given element. So by \ref{conjugate} we have found a
preimage which itself is a geometric element.

We prove the last claim by induction on the number of irreducible components of
$D_2$:
If we have just a single component $J$ we know the claim to be true in case
$\dim U=1$ and $J$ consists only of an isolated point.
If $\dim U\geq 2$ any element in the kernel is represented by a path isotopic
in $U-D_1$ to the constant path. By a general position check we may assume that
the isotopy is transversal to $J$.

Hence a suitable modification can be found which is supported in the complement
of $J$ and exhibits the path to be isotopic to a concatenation of boundaries of
discs normal to $J$ and segments connecting these. So its class in
$\pi(U-D_1\cup J)$ is represented by a product of simple geometric elements and
their inverses.

Suppose finally $D_2=D_2'\cup J$, so by induction hypothesis the simple
geometric elements associated to $D_2'$ generate the kernel of $\fund{U-D_1\cup
D_2'}\to\fund{U-D_1}$ and the simple geometric elements associated to $J$
generate the kernel of
$$
\fund{U-D_1\cup D_2}\to\fund{U-D_1\cup D_2'}.
$$

By $i)$ the simple geometric elements associated to $D_2'$ lift to simple
geometric elements associated to $D_2'$, hence we can conclude that the kernel
of the composed map
$$
\pi_1(U-D_1\cup D_2)\to\pi_1(U-D_1\cup D_2')\to\pi_1(U-D_1)
$$
is generated by simple geometric elements associated to $D_2$ as claimed.
\qed[3mm]

We can now apply this technical lemma to the curves $\curve$ and $\lfami$ in
$\CC^2$.

\begin{lemma}
\labell{pres}
The fundamental group $\pi_1(\CC^2 \setminus \curve)$ of the
complement of the curve is the quotient of $\pi_1(\CC^2 \setminus \curve \cup
\lfami)$ by the normal closure of $F_s$ considered as a subgroup by the 
presentation above, thus it is presented as
$$
\pi_1(\CC^2 \setminus \curve)=
\langle t_1,\dots ,t_n\mid t_i = a_*(\cg_k)(t_i)\rangle.
$$
\end{lemma}

\begin{remark}
Some if the given relations may be trivial, e.g.\ if $a_*(\cg_k)=\s_1$ then
$t_i=a_*(\cg_k)(t_i)=t_i$ for $i>2$.
\end{remark}

\proof[ of lemma \ref{pres}]
First by \ref{lbessis} the map $i_*:\pi_1(\CC^2 \setminus \curve \cup
\lfami)\to\pi_1(\CC^2 \setminus \curve)$ is surjective. By the construction
$F_s$ is considered a subgroup via
$$
\begin{array}{ccrcl}
F_s & \cong & \pi_1(\CC \setminus \sfami,y_0) &
\stackrel{s_*}{\tto} &
\pi_1(\CC^2 \setminus \curve \cup \lfami,s(y_0))\\
&&& \cong & \pi_1(\CC^2 \setminus \curve \cup \lfami).
\end{array}
$$
Since $s(\CC)$ is disjoint from $\curve$ and a section to the projection
$\CC^2\to\CC$, it is contractible and all elements of $\pi_1(\CC\setminus
\sfami)$ are therefore mapped to the trivial class.
So as claimed the normal closure of $F_s$ is contained in the kernel of $i_*$.

We are left to prove that it actually coincides with this kernel. As we know by
\ref{lbessis} the kernel is generated by simple geometric elements associated
to the irreducible components of $\lfami$.
Pick $\lfami_k$ among the irreducible components, then a simple geometric
element in $\pi_1(\CC\setminus\sfami)$ associated to $y_k$ is mapped by $s_*$
to a simple geometric element in $\pi_1(\CC^2 \setminus \curve \cup \lfami)$
associated to $\lfami_k$. So by \ref{conjugate} each simple geometric element
associated to an irreducible component of $\lfami$ is conjugated to an element
in the image of $s_*$ and thus in the normal closure of $F_s$ as claimed.

Finally we have to derive the given presentation.
We start with a presentation of $\pi_1(\CC^2\smin\curve\cup\lfami)$. Since
$\cg_1,...,\cg_s$ generate the subgroup $F_s$ we get a presentation for
$\pi_1(\CC^r\smin\curve)$ by adding $\cg_1,...,\cg_s$ to the set of relations.
In the third step we get rid of generators $\cg_i$ and relations $\cg_i$
simultaneously and must replace $\cg_i$ by the identity in the remaining
relations. In fact we get the presentation as claimed.
\qed[3mm]

\begin{remark}
\labell{quot/triv}
The group $G(a_*)$ defined by presentation
(\ref{pres}) is the quotient of $F_n$ by the normal subgroup generated
by $\{\gamma(t)\cdot t^{-1} \mid \gamma\in \im(a_*),\ t\in F_n\}$.
In other words,
$G(a_*)$ is the maximal quotient of $F_n$
on which the representation $a_*:F_s \to \br_n$ acts trivially.
\end{remark}

For use in a later chapter we prove a slight variation of this claim
and compute the fundamental group of  the complement of $\cfami$
and a single component of $\lfami$.

\begin{lemma}
\labell{pres-vert}
Suppose $\lfami_1$ is a single component of $\lfami$, then the fundamental
group $\pi_1(\CC^2\smin\cfami\cup\lfami_1)$ is generated by $i_*(\pi_1(\CC_n))$
and any simple geometric element associated to $\lfami_1$.
\end{lemma}

\proof
By the same argument as before, we get a presentation
$$
\pi_1(\CC^2 \setminus \curve)=
\langle t_1,\dots ,t_n\mid \cg_1\inv\cdot t_i\cg_1 = a_*(\cg_1)(t_i),t_i =
a_*(\cg_k)(t_i), k>1\rangle.
$$
A simple geometric element associated to $\lfami_1$ is then in the conjugation
class of $\cg_1$, thus represented by a word $w\cg_1w\inv$. Due to the
relations $\cg_1\inv\cdot t_i\cdot\cg_1 = a_*(\cg_1)(t_i)$ we may assume that $w$
does not contain the letter $\cg_1$, since they can be moved to the end of $w$
and then cancel through $\cg_1$. But then it is obvious that the fundamental
group is generated by the $t_i$ and $w\cg_1w\inv$.
\qed

\begin{remark}
\sloppypar
Suppose $a_*':F_s \to \br_n$ is related to $a_*$
by a commutative diagram
$$
\begin{array}{ccccl}
F_s      &  \stackrel{a_*}{\tto} &     \br_n &
  \stackrel{\a_n}{\tto} &     \Aut(F_n) \\
\psi\downarrow\quad  & &   
 & &  
\,\,\downarrow\conju_{\phi} \\
F_s   &     \stackrel{a'_*}{\tto}   &  \br_n
&     \stackrel{\a_n}{\tto}   &  \Aut(F_n)
\end{array}
$$
with $\psi\in\Aut(F_s)$, $\phi\in\Aut(F_n)$ such that the restriction of $\phi$
to $\a_n(\br_n)$ is an isomorphism of $\a_n(\br_n)$.
Then $G(a_*)$ is
isomorphic to $G(a_*')$.  Indeed, this condition
can be written as
$\phi(a_*(g)(t)\cdot t^{-1})=a_*'(\psi(g))(\phi(t))\cdot
\phi(t)^{-1}$, \mbox{$\forall g\in F_s$}, \mbox{$\forall t\in F_n$}.
Thus $\phi\in \Aut(F_n)$ induces an isomorphism
$\bar\phi:G(a_*)\to G(a_*')$.

Since $a_*,a_*'$ need not to be braid equivalent we see that the fundamental
groups of complements can be isomorphic, in fact the curve complements can be
homotopy equivalent, for curves which have different braid monodromies. An
example of this kind was given in \cite{CS}.
\end{remark}

\subsection{braid monodromy generators}

We now make the presentation (\ref{pres}) more precise.
To this end recall that the braid group $\br_n$ can be presented by generators
$\sigma_1,\dots ,\sigma_{n-1}$ subjected to relations
$\sigma_i\sigma_{i+1}\sigma_i=\sigma_{i+1}\sigma_i\sigma_{i+1}
\ (1\le i<n-1),\ \sigma_i\sigma_j=\sigma_j\sigma_i
\ (|i-j|>1)$, see \cite{Bi}, \cite{H2}. 
The Artin representation
$\a_n:\br_n\to \Aut(F_n)$ is then given in terms of $F_n=\lspan
t_1,...,t_n|\,\rspan$ by:
$$
\sigma_i(t_j)=\begin{cases}  t_it_{i+1}t_i^{-1}&\text{ if } j=i,\\
t_i&\text{ if } j=i+1,\\
t_j &\text{ otherwise.}
\end{cases}
$$
So far we chose systems of paths in $\CC\smin\sfami$ and the fibre $\CC_n$ over
$y_0$ only. But of course we may do so in any fibre
$p\inv(y),y\in\CC\smin\sfami$ -- though not in a coherent way, as that would
amount to a global trivialization of the bundle.

As we have represented the elements $\cg_i$ by nooses we may divide them into
small loops $\o_i$ based at the neck, and long ropes $\eta_i$.

Like $\cg_i$ induce braid automorphism upon the choice of a system of paths
for $\CC_n$, so do the $\o_i$ upon a choice for the fibre over the neck of
the noose.

The identification of $\eta_i$ with an element in $\br_n$ depends on both
choices. Note that restricted to $\eta_i$ the bundle is trivializable, hence
one choice is sufficient as a global choice. So we can compare both choices.

In fact the different choices up to isotopy form a simply transitive orbit
under the Artin representation, hence an ordered pair determines uniquely the
corresponding transition braid.

These constructions fit naturally together to yield
$$
a_*(\cg_k)=\beta_k^{-1}\a_k \beta_k,
$$
where $\a_k\in\br_n$ is the monodromy along $\o_k$ and $\beta_k\in \br_n$ is the
transition along $\eta_k$.

So as one would like to
express these braids in terms of the standard generators $\sigma_i$ of $\br_n$,
one may try to accomplish this in two steps.

\begin{description}
\item[Step 1]
The structure of the isolated singularities in the fibre at $y_k$ determines
the local braid $\a_k\in\br_n$, upon a choice of a geometrically
distinguished system of paths for the fibre at $y_k^0$, the neck of the noose
and base point of $\o_k$.

This braid may be obtained from the Puiseux series expansion of the defining
polynomial of $\curve$ at each intersection point with $\pi\inv(y_k^0)$.
The actual algorithm is implicit in the work of Brieskorn, Kn\"orrer
\cite{BK} and Eisenbud, Neumann \cite{EN}.
\end{description}

\begin{description}
\item[Step 2]
A transition braid $\b_k$ depends on the relation between the choices of
systems of paths in the fibres to be compared along the path.

The pull back of the bundle along the path $\eta_k$ is a trivial $\br_n$-bundle,
hence the mapping class groups of the two fibres at the endpoints are
identified by parallel transport. Since both are identified with the abstract
braid group we get an isomorphism of braid groups.

This is actually obtained by inner conjugation since the isomorphisms are all
defined in terms of distinguished systems of paths.

The important point to note is that in case the local braid $\a_k$ involves few
strings we don't need to understand the parallel transport resp.\ $\b_k$ to
full extend, but only $\b_k$ up to the stabilizer of the local braid under the
conjugation action, in order to determine $\b_k\inv\a_k\b_k$.
\end{description}

Let us try the first step on some specific singularities and exemplify the final
remark of the second step.
\\

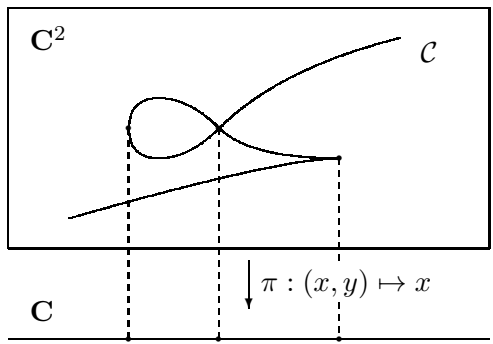
\begin{figure}[h]
\begin{center}
\setlength{\unitlength}{0.8mm}
\begin{picture}(80,55)(-40,-15)
\put(0,-2){\vector(0,-1){8}}
\put(2,-7){$\pi:(x,y)\mapsto x$}
\put(-40,-15){\line(1,0){80}}
\put(-38,-12){ $\CC$}
\put(-40,0){\line(1,0){80}}
\put(-40,40){\line(1,0){80}}
\put(-40,0){\line(0,1){40}}
\put(40,0){\line(0,1){40}}
\put(-38,33){ $\CC^2$}
\put(27,31){ $\cfami$}
\multiput(-20,20)(0,-2){18}{\line(0,-1){1}}
\multiput(-5,20)(0,-2){18}{\line(0,-1){1}}
\multiput(15,15)(0,-2){9}{\line(0,-1){1}}
\multiput(15,-9)(0,-2){3}{\line(0,-1){1}}
\put(-20,-15){\circle*{1}}
\put(-5,-15){\circle*{1}}
\put(15,-15){\circle*{1}}
\qbezier[200](25,35)(5,30)(-5,20)
\qbezier[60](-5,20)(-10,15)(-15,15)
\qbezier[40](-15,15)(-20,15)(-20,20)
\qbezier[40](-20,20)(-20,25)(-15,25)
\qbezier[60](-15,25)(-10,25)(-5,20)
\qbezier[110](-5,20)(0,15)(15,15)
\qbezier[320](15,15)(5,15)(-30,5)
\put(-20,20){\circle*{1}}
\put(-5,20){\circle*{1}}
\put(15,15){\circle*{1}}
\end{picture}
\end{center}

\caption{vertical tangent, node and cusp}
\end{figure}

\begin{examp}
Consider the plane curve $\curve: z^p - x^q = 0$.
The fundamental group of its complement was determined by Oka \cite{Ok}.
A look at Oka's computation reveals
a natural choice of system of paths such that the braid monodromy generator
is $(\sigma_1\cdots \sigma_{p-1})^q\in \br_p.$
For instance, to a simple tangency yields
$\sigma_1$, a node yields $\sigma_1^2$, and a cusp yields $\sigma_1^3$.
Hence the individual monodromies around the special points $y_k$ are conjugated
to powers of half-twists, the exponent being $1$ in the case of tangency
points, $2$ in the case of ordinary nodes, and $3$ in the case of ordinary cusps.

\end{examp}

\begin{examp}
If the monodromy $\a$ of $\o$ is conjugated to a power of a half twist as above,
then this twist is determined by an isotopy class of arcs between the two
punctures involved. Of course $\b\a\b\inv$ is a half twist again, so in fact we
need to determine the image of the arc under the parallel transport along
$\eta$ only.

\end{examp}

We close this section with two remarks touching only on the surface of recent
developments.

\begin{remark}
The program presented so far is well adapted to the study of branch curves of
generic projections of projective surfaces to complex projective plane. Such
curves have only singular points  which are ordinary cusps or nodes.

The methods have also been generalized to the symplectic set up.
\end{remark}

\begin{remark}
In the case of branch curves
Moishezon and Teicher,
\cite{Mo3,Te1}, developed a degeneration technique. They start with a
projective embedding of the complex surface
$X$, and deforms the image of this embedding to a singular configuration
$X_0$ consisting of a union of planes intersecting along lines. The
branch curve of a projection of $X_0$ to $\CC^2$ is therefore a union
of lines; the manner in which the smoothing of $X_0$ affects this curve
can be studied explicitly, by considering a certain number of standard local
models near the various points of $X_0$ where three or more planes
intersect.
\end{remark}


\section{braid monodromy of horizontal divisors}

We want now to generalize the considerations for plane algebraic curves in two
aspects. We go to higher dimensions, i.e.\ affine algebraic divisors in complex
affine space, and include the local case of analytic space germs in
$\CC^{r+1},0$. The latter is in fact a generalization of the local
considerations at singularities $v_k$ of a plane curve by means of Puiseux
series.

So we work with a polynomial $f$ or a holomorphic function germ
of vanishing order $n$ in $r+1$ variables. For simplicity we assume that the
projection is along the distinguished coordinate $z$ and that $f$ is monic of
degree $n$ in the polynomial ring over $\CC[x_1,...,x_r]$, resp.\ $z$-general in
$\CC[x,z]$. Then in the local analytic case we may use the Weierstrass
preparation theorem to end up in either case with the divisor given by the zero
set of a Weierstrass polynomial function $f:\CC^r\times\CC\to\CC$,
$$
f(x,z)\,=\,z^n+\sum_{i=0}^{n-1}a_i(x)z^i,
$$
with $a_i\in\CC[x]$ resp.\ $a_i\in\mideal_r:=\{a\in\ofami_r|a(0)=0\}$.

Let us resume the local analytic case later merely addressing the necessary
improvements and stick to the affine case now.

The set $\sfami\subset\CC^r$ of points for which the corresponding fibre is
tangent to the divisor $\hfami$ defined by $f$ or contains a singular point of
$\hfami$ is the zero locus of the resultant of $f$ and $\del_z f$ which
eliminates $z$ and coincides with the pull back of the discriminant $\Delta_n$
in $\CC^n$ along the coefficient map $a=(a_0,...,a_{n-1}):\CC^r\to\CC^n$.

Let $\lfami$ denote the linear extension of $\sfami$ to $\CC^r\times\CC$ by pull
back. In generalization of the plane curve case we get a locally trivial bundle
$$
\CC^{r+1}\setminus\hfami\cup\lfami\tto\CC^r\setminus\sfami
$$
with structure group $\br_n$, fibre $\CC_n$ and holonomy
\begin{equation}
\label{fifth}
\rho:\pi_1(\CC^r\setminus\sfami,y_0)\to\Aut(\pi_1(\CC_n,v_0))
\end{equation}
which is the coefficient homomorphism
$a_*:\pi_1(\CC^r\setminus\sfami,y_0)\to\br_n$,
if $\br_n$ is identified with its image under the Artin representation using a
suitable geometrically distinguished system of paths. So we define:

\begin{defi}
The homomorphism $\a:\pi_1(\CC^r\setminus\sfami)\to\br_n$ determined by the
composition
$\a_n\circ\a=\rho$ is called {\it braid monodromy homomorphism} of the
hypersurface $\hfami$ projected along $z$.
\end{defi}

It is only well defined up to a certain isomorphisms of the source and
conjugation in the target, but to proceed we won't have to specify them.

\subsection{braid monodromy presentation}

To get a presentation for $\pi_1(\CC^{r+1}\setminus\hfami)$ we argue as in
the curve case. So we start with the long exact homotopy sequence of the
bundle $\CC^{r+1}\smin\hfami\cup\lfami\to\CC^r\smin\sfami$.

This time our argument must be a bit more substantial to get a short exact
sequence out of it:

\begin{lemma}
\labell{del/triv}
The boundary map $\pi_2(\CC^r\smin\sfami)\to\pi_1(\CC_n)$ is trivial.
\end{lemma}

\proof
Of course the image of the boundary map is normal and abelian. We get the proof
now in cases $n=1$ and $n>1$ separately.
If $n=1$ then $f$ is a linear polynomial and $\sfami$ is empty so
$\pi_2(\CC^r\smin\sfami)=0$.
If $n>1$ the group $\pi_1(\CC_n)$ is free of rank at least two. But then only
the trivial subgroup is normal and abelian.
\qed[3mm]

Since the splitting argument \ref{sect} still applies, we can base the proof of
the following lemma on the short split exact sequence:
$$
1\to \pi_1(\CC_n) \to \pi_1(\CC^{r+1} \setminus \hfami \cup \lfami)
\stackrel{s_*}{\stackrel{\textstyle{\longleftarrow}}{\tto}}
{}_{\hspace*{-5.7mm}p_*\hspace*{3mm}}
\pi_1(\CC^r \setminus \sfami)\to 1,
$$
with action given by the braid monodromy homomorphism $a_*$ (\ref{fifth}).

\begin{lemma}
\labell{gpres}
The fundamental group $\pi_1(\CC^{r+1}\smin\hfami)$ of the complement of an
affine divisor $\hfami$ is the quotient of
$\pi_1(\CC^{r+1}\smin\hfami\cup\lfami)$ by the normal closure of
$\pi_1(\CC^r\smin\sfami)$ considered as a subgroup by the section $s_*$, thus
it is presented as
$$
\pi_1(\CC^{r+1}\smin\hfami)=\lspan t_1,...,t_n|\,t_i=a_*(\cg_k)(t_i)\rspan,
$$
where $\cg_1,...,\cg_m$ is any system of generators for
$\pi_1(\CC^r\smin\sfami)$.
\end{lemma}

\proof
By \ref{lbessis} the map
$i_*:\pi_1(\CC^{r+1}\smin\hfami\cup\lfami)\to\pi_1(\CC^{r+1}\smin\hfami)$ is
surjective. Since $s(\CC^r)$ is contractible the image by $s_*$ of
$\pi_1(\CC^r\smin\sfami)$ and its normal closure must be trivial in
$\pi_1(\CC^{r+1}\smin\hfami)$. As in the curve case \ref{pres} we can give an
argument relying on \ref{lbessis} that this is the kernel of $i_*$.

The presentation is then obtained using the split sequence, $s_*$, and $i_*$.
Let us suppose we are given an arbitrary presentation
$$
\pi_1(\CC^r\smin\sfami)=\lspan\cg_1,...,\cg_s|\,r_1,...,r_m\rspan,
$$
where $r_1,...,r_m$ are words in the listed generators.
Since the fibre of the bundle has fundamental group presented as
$\pi_1(\CC_n)=\lspan t_1,...,t_n|\,\rspan$, we get with the split exact
sequence:
$$
\pi_1(\CC^{r+1}\smin\hfami\cup\lfami)=\lspan
t_1,...,t_n,\cg_1,...,\cg_s|\,r_1,...,r_m,\cg_k
t_i\cg_k\inv=a_*(\cg_k)(t_i)\rspan
$$
As in the curve case the results on the kernel of $i_*$ imply, that we get a
presentation for $\pi_1(\CC^{r+1}\smin\hfami)$ if we add the relations
$\cg_1,...,\cg_s$.
So we may drop generators $\cg_k$ and replace $\cg_k$ in all relations by the
identity, to get a Tietze equivalent presentation. Because
$\cg_1,...,\cg_s,r_1,...,r_m$ become trivial relations, we may as well discard
them and end up with 
$$
\pi_1(\CC^{r+1}\smin\hfami)=\lspan t_1,...,t_n|\,t_i=a_*(\cg_k)(t_i)\rspan,
$$
as claimed, since in the beginning we may have chosen any set of generators
for $\pi_1(\CC^r\smin\sfami)$.
\qed[3mm]

In fact, to get presentations of fundamental groups as in lemma \ref{gpres}, we
may go back to the case of curves:

Consider a generic fibre $L$ of a generic fibration of $\CC^r\smin\sfami$. Then
$\pi_1(L\smin\sfami_L)$ surjects onto $\pi_1(\CC^r\smin\sfami)$.
So given generators $\cg_1,...,\cg_s$ for $\pi_1(L\smin\sfami_L)$ we get by the
last result a presentation
$$
\pi_1(\CC^{r+1}\smin\hfami)=\lspan t_1,...,t_n|\,t_i=a_*(\cg_k)(t_i)\rspan,
$$
If we define in the plane $p\inv(L)$ the curve $\curve_L:=\hfami\cap p\inv(L)$
then by our arguments the above presentation is also valid for
$\pi_1(p\inv(L)\smin\curve_L)$.
This is of course one of the well known Zariski van Kampen results.\\

\subsection{braid monodromy of local analytic divisors}

To get the analogous results in the case of local analytic divisors, we must go
through the same procedure, but we are forced to make some more choices in the
construction, which we later have to show to have no bearing on the definitions.

We already reached the stage at which we have to deal with a Weierstrass
polynomial $f(x,z)=z^n+\sum a_i(x)z^i$, $a_i(0)=0$.
It defines the given analytic germ $\hfami$ and for sufficiently small $\rho$
the coefficients $a_i$ are defined and bounded for $|x|<\rho$.
\\

Since the divisor was assumed to be reduced, the set of points $|x|<\rho$
such that $f(x,z)$ is a simple Weierstrass polynomial is the complement of a
proper analytic subset $\sfami_\rho$, which is defined by the
resultant of $f$ and
$\del_z f$ with respect to $z$. So we get an associated polynomial cover
over $X_\rho:=B_\rho\smin\sfami_\rho$:
$$
f:X_\rho\times\CC\to\CC,\,
(x,z)\mapsto z^n+\sum a_i(x)z^i.
$$
By the choice of $\rho$, the zeroes of $f$ are uniformly bounded for $x\in
B_\rho$, hence there is a section $s$ to the projection $\pi:X_\rho\times\CC\to
X_\rho$ which avoids the zero set of $f$.

As in the previous cases we have a locally trivial bundle
$$
Y_\rho\to X_\rho,\quad Y_\rho:=X_\rho\times\CC\smin f\inv(0),
$$
for which a holonomy $\pi_1(X_\rho,y_0)\to \Aut(\pi_1(\CC_n,v_0))$ is given
after a choice of base points, $y_0\in X_\rho$ and $v_0$ in the fibre
$\pi\inv(y_0)$.

For $\rho$ below some finite bound, the isomorphism type of $\pi_1(X_\rho,y_0)$
is well defined, hence we may define:

\begin{defi}
For $\rho$ sufficiently small,
the map $a_*:\pi_1(X_\rho,y_0)\to\br_n$ induced by
the coefficient map of $Y_\rho\to X_\rho$ is called {\it braid monodromy
homomorphism} of the germ $\hfami$ {\it projected along $z$}.
\end{defi}

It is well defined at least up to isomorphisms of the source and conjugation in
the target, which is all we need in the sequel.

Moreover the local analogue of lemma \ref{gpres} holds true:

\begin{lemma}
\labell{loc-pres}
Suppose $\hfami$ represents the germ of a divisor in $\CC^{r+1},0$, then
for $\e$ and $\rho=\rho(\e)$ sufficiently small, the isomorphism class
$\pi_1(B_\e\times B_\rho\smin\hfami)$ of fundamental groups
is  presented as
$$
\pi_1(B_\e\times B_\rho\smin\hfami)=\lspan
t_1,...,t_n|\,t_i=a_*(\cg_k)(t_i)\rspan,
$$
where $\cg_1,...,\cg_m$ is any system of generators for
$\pi_1(X_\rho,y_0)$.
\end{lemma}

\proof
Let $\e$ and $\rho$ be small enough to make sure that $B_\e\times
B_\rho\cap\hfami$ is a branched cover of degree $n$ over $B_\rho$ and branched
along an analytic subset with complement $X_\rho$.

Then $B_\e\times B_\rho\smin\hfami$ is a strong deformation retract of $Y_\rho$
and the claim follows along the line of argument of lemma \ref{gpres}.
\qed[0mm]

\begin{remark}
\labell{qhomogen}
If $\hfami$ is given by a quasi homogeneous polynomial $f$ with a good
$\CC^*$-action, then so is $\sfami_\rho$.
In that case the $\CC^*$-action can be used to show that $B_\e\times
B_\rho\smin\hfami$ is a strong deformation retract of $\CC^{r+1}\smin f\inv(0)$
by a map which respects the projection to $B_\e$ resp.\ $\CC^r$.

Hence all braid monodromy considerations for the hypersurface germ $\hfami$ are
equal to those for the affine hypersurface defined by $f$.
\end{remark}


\chapter{braid monodromy of singular functions}


In this chapter we review the notion of semi universal unfoldings of singular
function germs and show how a natural polynomial cover arises
in this set-up. The corresponding braid monodromy homomorphisms and braid
monodromy groups are then associated invariantly to equivalence classes of
function germs via their versal unfoldings.

In particular a presentation for the fundamental group of the discriminant
complement will be derived from the new invariants.

But to get that far, we have to make quick digression through the theory of
unfoldings of singular functions touching on such diverse notions as
discriminant sets, truncated unfoldings, and bifurcation sets.

\section{preliminaries on unfoldings}

The basic objects in singularity theory we start with are holomorphic function
germs on affine coordinate space, $f:\CC^n,0\to\CC$, which form the
$\CC$-algebra $\ofami_n$. It is a local algebra with maximal ideal
\begin{eqnarray*}
\mideal_n & = & \{f\in\ofami_n|\,f(0)=0\}
\end{eqnarray*}

Since a function on abstract affine space is identified with different elements
of $\ofami_n$ depending on a choice of coordinates, it is natural to consider
such elements to be equivalent, more precisely:

\begin{defi}
Given elements $f,g\in\ofami_n$ are called {\it right equivalent} or simply
{\it equivalent}, $f\sim g$, if there is a holomorphic map
$\riso:\CC^{n},0\to\CC^n,0$ such that $f(x)=g(\riso(x))$ and
$\riso$ is biholomorphic.
\\
If $X\subset\mfrak_{n}$ is an equivalence class and $f\in X$, then $f$ is
called a function of type $X$.
\end{defi}

\begin{examp}
The classes of simple singularities according to Arnold, cf. \cite{AGV1}, are
represented by
\begin{enumerate}
\item[$A_k$:]
$x_1^{k+1}+x_2^2+\cdots+x_n^2$,
\item[$D_k$:]
$x_1^2x_2+x_2^{k-1}+x_3^2+\cdots+x_n^2$,
\item[$E_6$:]
$x_1^3+x_2^4+x_3^2+\cdots+x_n^2$,
\item[$E_7$:]
$x_1^3+x_1x_2^3+x_3^2+\cdots+x_n^2$,
\item[$E_8$:]
$x_1^3+x_2^5+x_3^2+\cdots+x_n^2$,
\end{enumerate}
\end{examp}

\subsection{versal unfolding}

Proceeding deeper into the theory we introduce secondary objects to a given
function germ $f\in\ofami_n$:

\begin{defi}
A function germ $F$ on affine coordinate germ $\CC^n\times\CC^k,0$
is called {\it unfolding} of $f\in\ofami_n$, if
$F_0=f$ for $F_0:=F|_{\CC^n\times\{0\},0}$.
Then $k$ is called the {\it unfolding dimension} and $\CC^k,0$ the base
or the parameter space of the unfolding.
\end{defi}

This notion should be understood as a family of function germs in $\ofami_n$
parameterized over a space germ $\CC^k,0$. This interpretation is also at the
base of the equivalence notion induced on unfoldings.

\begin{defi}
Suppose $F,G\in\ofami_{n+k}$ are unfoldings of $f\in\ofami_n$.
Then $F,G$ are called {\it equivalent}, if there is a holomorphic map germ
$\riso:\CC^{n+k},0\to\CC^n,0$ such that $G(u,x)=F(u,\riso(u,x))$,
$\riso(0,x)=x$.
\end{defi}

And as with families one can consider the pull back of an unfolding along a map
to its base:
Suppose $F\in\ofami_{n+k}$ is an unfolding of $f\in\ofami_n$ and
$\pull:\CC^l,0\to\CC^k,0$ is a map germ, we call $G\in\ofami_{n+l}$ the {
unfolding} of $f$ {\it induced} from $F$ by $\pull$, if $G(v,x)=F(\pull(v),x)$.

Thus prepared we can now introduce the concept best suited for classification
in singularity theory.

\begin{defi}
If $F\in\ofami_{n+k}$ is an unfolding of $f\in\ofami_n$ then
$F$ is called {\it versal}, if each unfolding $G\in\ofami_{n+l}$ of $f$ is
equivalent to an unfolding induced from $F$.

A versal unfolding is called {\it miniversal} if it is versal and of minimal
dimension.
\end{defi}

Naturally one would ask to induce a given unfolding in a unique way but
experience taught to be content with uniqueness only of the differential of the
pull back map.

Accordingly the miniversal unfoldings which have got this property are also
called {\it semi universal}.

\begin{remark}
Versality is an open condition in the following sense. If a representative of a
versal unfolding $F$ is given and is defined for $u,0$, then $G$ given by
$G(v,x):=F(u+v,x)$ is an unfolding of the function $x\mapsto F(u,x)$.
\end{remark}

The case of functions with an isolated singularity -- which we are interested
in exclusively -- is characterized by the finite codimension in $\ofami_n$ of
the {\it Jacobian ideal} generated by the partial derivatives,
$J(f)=\lspan\del_1 f,...,\del_n f\rspan$.

In this case specific miniversal unfoldings and the mutual relations of versal
unfoldings can be described according to \cite{AGV1}, chapter 8:

\begin{prop}
\labell{model}
Given $f\in\mideal_{n}$ and $b_1,...,b_k\in\ofami_n$ such that the classes of
$b_1,...,b_k$ modulo $\jac(f)$ form a generating set (a basis) of
$\ofami_n/\jac(f)$. A versal (semi universal) unfolding of $f$ is then
given by $F(u_1,...,u_k,x):=f(x)+u_1b_1(x)+...+u_kb_k(x)$.
\end{prop}

\proof
By hypothesis the function $b_i$ and the ideal $\jac(f)$ span the tangent space
of $\ofami_n$ at $f$. The ideal $\jac(f)$ can be shown to be the tangent space
to the orbit of action by biholomorphic maps of the source of $f$. The point of
the proof is thus -- as indicated in \cite{AGV1} -- that the infinitesimal
transversality to the orbit of $f$ lifts to local transversality to all orbits
sufficiently close.

\begin{prop}
\labell{v-pull}
Suppose $F\in\ofami_{n+k}$ and $G\in\ofami_{n+l}$ are versal unfoldings of
$f\in\ofami_n$ and $k\geq l$. Then $F$ is equivalent to some unfolding
induced from $G$ by a map germ $\pull:\CC^k,0\to\CC^l,0$ of full rank.
\end{prop}

\proof
This result can be derived from the previous by a careful analysis of all the
definitions involved.

\subsection{discriminant set}
\newcommand{\Ftil}{\tilde F}

The discriminant set is the germ of subsets in the base of a versal unfolding
$F\in\ofami_{n+k}$ of a singular function $f\in\ofami$ given by those
parameters $u$ for which $0$ is a critical value of the function $F_u:x\mapsto
F(u,x)$.

To give a precise meaning to this description, we examine the situation for a
representative of the germ $F$. W.l.o.g.\ assume $f:\CC^n,0\to\CC,0\,$
to be a germ with an isolated critical point at $0$. Choose
sufficiently small neighbourhoods of zero $M=\{x|\,\|x\|\leq
\rho\}\subset\CC^n$
and $U=\{u|\,\|u\|\leq\d\}\subset\CC^k$, for which a representative $\Ftil$ is
defined.

By the curve selection lemma we may assume $\rho$ and $\d=\d(\rho)$
sufficiently small such that the level set $\{x|\,\Ftil(u,x)=0\}$ is
non-singular on the boundary $\del M$ of the ball $M$ and is
transverse to $\del
M$ for every $u\in U$, cf.\ \cite{Mil}.

If we distinguish $u\in U$ into singular and non-singular parameters according
to $V_u:=\{x|\,F(u,x)=0\}\cap M$ being singular or not, we define the
discriminant set by means of any representative of $F$:

\begin{defi}
The {\it discriminant set} in the base of the versal unfolding
$F\in\ofami_{n+k}$ of a singularity $f\in\ofami_n$ is represented by the
singular parameters in $U$ for a representative of $F$.
\end{defi}

Of course any other representative of $F$ coincides with the chosen one for
sufficiently small neighbourhoods and so does the corresponding set of singular
parameters.

\begin{examp}
Let $f(x)=x^3,F(u,x)=x^3+u_1x+u_0$. The discriminant set for $f$ is precisely
the set of pairs $u_1,u_0$ such that the polynomial $x^3+u_1x+u_0$ has multiple
roots. Hence the discriminant is the hypersurface germ cut out by the equation
$27u_0^2+4u_1^3=0$.
\end{examp}

This examples exposes a common feature of all discriminant sets:

\begin{prop}
\labell{d-irr}
The discriminant set in the base of a versal unfolding $F$ of a singular
function $f$ is an irreducible hypersurface germ.
\end{prop}

\proof
The proof uses suitable representatives again and exploits the fact, that the
discriminant is the image of the critical set, which is irreducible and
analytic, under a proper finite map, cf.\ \cite{ebe}
\qed[3mm]

\subsection{truncated versal unfolding}

Let $f:\CC^n,0\to\CC,0$ be a function germ. Instead of unfoldings
$F\in\ofami_{n+k}$, which can be considered as families of elements in
$\ofami_n$, we now confine ourselves to unfoldings $F$ such that
$F(u,0)\equiv0$, which can be considered as families of elements in
$\mideal_n$.

The definitions of equivalence, induced unfoldings, and versality carry over
without modifications to the present case.

\begin{defi}
If $F\in\mideal_{n+k}$ is an unfolding of $f\in\mideal_n$ then
$F$ is called a {\it truncated versal unfolding}, if $F(u,0)\equiv0$ and each
unfolding $G\in\mideal_{n+l}$ of
$f$ with $G(v,0)\equiv0$ is equivalent to an unfolding induced from $F$.
\end{defi}

The relation to versal unfoldings in the ordinary sense is readily given:

\begin{prop}
Suppose $F\in\ofami_{n+k}$ and $F^\#\in\mideal_{n+k-1}$ are unfoldings of
$f\in\mideal_k$ related by $F(u,x)=F^\#(u^\#,x)+u_0$. Then $F^\#$ is
a truncated
versal unfolding of $f$ if and only if $F$ is a versal unfolding of $f$.
\end{prop}

This stems from the fact that a miniversal truncated unfolding is given by
$F^\#(u,x)=f(x)+u_1b_1+...+u_{k-1}b_{k-1}$ if for $\mideal_n\to\mideal_n/J(f)$
the $b_i$ map to a basis of the quotient, cf.\ prop.\ \ref{model}.

Given a truncated versal unfolding $F^\#$ we call the versal unfolding $F$ of
the proposition the {\it completed versal unfolding}  corresponding to $F^\#$.

\subsection{bifurcation set}

A function is said to be a {\it Morse function}, if it has only non-degenerate
critical points and their values are distinct. For the definition of the
bifurcation set as the set of parameters in the base of a truncated versal
unfolding $F^\#$, for which the corresponding function
is not Morse, we have once again to resort to a local representative:

\begin{defi}
The {\it bifurcation set} in the base of a truncated versal unfolding
$F^\#$ of a
singular function $f$ is represented by the set of parameters for which the
corresponding function restricted to a sufficiently small ball $M$ is not a
Morse function.
\end{defi}

Again one should convince oneself that this definition does not depend on the
choice of representatives and sufficiently small neighbourhoods.

\begin{examp}
By $F^\#(u,x)=x^4+u_2x^2+u_1x$ a
truncated versal unfolding of $f(x)=x^4$ is given.
This polynomial has degenerate critical points if and only if its derivative
has multiple roots, i.e.\ $27u_1^2+8u_2^3=0$.
Moreover one can verify that critical values coincide only if $F_u$ is an even
function, i.e.\ $u_1=0$. The bifurcation set is thus the union of a cusp
corresponding to functions with degenerate critical points and a line
corresponding generically to functions with common values of distinct critical
points.
\end{examp}

The example is but an instance of general facts which are more involved than in
the case of discriminants and which are summarized below.

\begin{prop}
In the base of a truncated versal unfolding the
set of parameters such that the corresponding function has a degenerate
critical points (not of type $A_1$) defines an analytic hypersurface
germ. It is
empty only for $f$ of type $A_1$ and irreducible otherwise.
\end{prop}

\proof
Similar to the proof of prop.\ \ref{d-irr} one has to use the fact that the
determinant of the Hessian is transversal to the critical set and thus
determines a smooth analytic germ of dimension $k-1$ which can be shown to
project properly and finitely to the set under scrutiny, cf.\
\cite{Wi},\cite{poe}.
\qed[3mm]

The complement of this {\it strict bifurcation variety} consists of parameters
corresponding to
functions with non-degenerate critical points only, at least two of
which have a
common value. Its closure goes by several names, e.g.\ {\it Maxwell stratum}
\cite{AGV2}, or {\it conflict variety} \cite{Wi}, and its
decomposition into irreducible components was determined by Wirthm\"uller
\cite{Wi}, cf.\ also \cite{poe}:

\begin{prop}
The conflict variety in the base of a truncated versal unfolding is an
analytic hypersurface germ, it is
\begin{enumerate}
\item
empty, if $f$ is of type $A_1,A_2$,
\item
the union of three irreducible components, if $f$ is of type $D_4$,
\item
the union of two irreducible components, if $f$ is of type $D_\mu,\mu\geq5$,
\item
irreducible in all remaining cases.
\end{enumerate}
\end{prop}

\section{discriminant braid monodromy}

We finally draw closer to our proper objective.
For any given $f\in\mideal_n$ with isolated singularity consider a truncated
versal unfolding $F^\#$ and the corresponding completed versal unfolding $F$.

The discriminant is an analytic hypersurface germ, in the base $\CC^k,0$ of
$F$, hence given by a reduced holomorphic function $\Delta$.
Since $f$ is assumed to have an isolated singularity only, its singular value
$0$ is isolated, too. Therefore $\Delta$ is $u_0$-regular and we may define:

\begin{defi}
The {\it braid monodromy homomorphism} associated to the truncated versal
unfolding $F^\#$ is the braid monodromy homomorphism of the discriminant in the
base of the corresponding completed versal unfolding $F$ projected along $u_0$.
\end{defi}

It is this notion which is the central one of this paper and we will devote the
rest of the section to some of its general properties, before we start the
investigation for specific classes of singularities.

\subsection{basic properties}

First we identify the range of the braid monodromy homomorphism. We
actually link it to the {\it Milnor number} $\mu$ of the singular function $f$,
which can be defined as the multiplicity of the critical value.

\begin{lemma}
\sloppypar
The function $\Delta$ defining the discriminant divisor is $u_0$-regular of
\mbox{order} equal to the Milnor number of
$f$.
\end{lemma}

\proof
This follows from the fact that in the base of a versal unfolding of the
given kind the line $u_1=...=u_{k-1}=0$ is not in the tangent cone of the
discriminant. We have just to add the information that the Milnor number gives
the multiplicity of the discriminant.
\qed[3mm]

\begin{cor}
The range of the braid monodromy homomorphism is $\br_{\mu(f)}$.
\end{cor}

Next we identify the source with the help of the bifurcation set $\sfami_\rho$
in the base of $F^\#$ restricted to a ball of sufficiently small radius $\rho$
centered at the origin.

\begin{prop}
The braid monodromy homomorphism of a truncated versal unfolding $F^\#$ is
defined on the isomorphism class $\pi_1(B_\rho\smin\sfami_\rho)$ for $\rho$
sufficiently small.
\end{prop}

\proof
Note first that by the Weierstrass preparation theorem $\Delta$ can be assumed
to be a Weierstrass polynomial.
By the construction of the braid monodromy homomorphism of projected
hypersurface germs, it suffices to show, that the the discriminant set
coincides with the set of parameters for which $\Delta$ fails to be a simple
Weierstrass polynomial, at least for $\rho$ sufficiently small.

Now the multiplicity of a critical value is one, if and only if it is the value
of a single non-degenerate critical point. Therefore the number of critical
values drops for sufficiently small $\rho$, if and only if the corresponding
function is not a Morse function.
\qed[0mm]

\subsection{invariance properties}

Most important are the invariance properties of the braid monodromy
homomorphism. We have already noted, that the assignment of a braid monodromy
homomorphism to an affine hypersurface germ projected along a suitable
coordinate is not well-defined. But if we consider equivalence classes up to
isomorphisms of the source and conjugation in the range, then the
induced map is
well defined.
\\

In fact with this interpretation we get the invariance on the class of all
truncated versal unfoldings of equivalent functions with isolated singularity.

\begin{prop}
Given braid monodromy homomorphisms $\a(F),\a(G)$ associated to
truncated versal
unfoldings
$F^\#$ and $G^\#$ of equivalent functions $f,g$, then
there exists a commutative diagram
\begin{eqnarray*}
\pi_1(B^F_\rho\smin\sfami^F_\rho) & \stackrel{\a(F)}\to & \br_{\mu(f)}\\[2mm]
\stackrel{\cong}{}\,\downarrow\quad\quad && \,\downarrow\, conj_\phi\\
\pi_1(B^G_\rho\smin\sfami^G_\rho) & \stackrel{\a(G)}\to & \br_{\mu(g)},
\end{eqnarray*}
where $\mu(f)=\mu(g)$ and $conj_\phi$ is conjugation by $\phi\in\br_\mu$.
\end{prop}

\proof
Suppose w.l.o.g.\ that $k\geq l$ for the unfolding dimensions of $F$ and $G$.
By hypothesis and prop.\ \ref{v-pull} there is a map germ of full rank
$\pull^\#:\CC^{k-1},0\to\CC^{l-1},0$ and a map germ
$\riso^\#:\CC^{k-1}\times\CC^n,0\to\CC^n,0$ with $\riso^\#$ restricted to
$0\times\CC^n,0$ biholomorphic and
\begin{eqnarray*}
F^\#(u,x) & = &
G^\#(\pull^\#(u),\riso^\#(u,x)).
\end{eqnarray*}
The map germs $\pull:u\mapsto (u_0,\pull^\#(u^\#))$ and $\riso:u,x\mapsto
\riso^\#(u^\#,x)$ then have analogous properties with
\begin{eqnarray*}
F(u,x) & = &
G(\pull(u),\riso(u,x)).
\end{eqnarray*}
Hence we arrive at a commutative square of germs
$$
\begin{array}{lcl}
\CC^{k},0\, & \stackrel\pull\to & \CC^{l},0\\[2mm]
\,\big\downarrow\pi_{u_0}\, & & \,\big\downarrow\pi_{u_0}\\[2mm]
\CC^{k-1},0 & \stackrel{\pull^\#}\to & \CC^{l-1},0,
\end{array}
$$
which induces a commutative square for sufficiently small representatives.

Since the pull back of the discriminant for $G$ along $\pull$ yields
the discriminant for $F$, the bottom map $\pull^\#$ induces an isomorphism
of fundamental groups of local bifurcation complements.

Moreover the associated polynomial cover for $F$ is the pull back by
$\pull,\pull^\#$ of the polynomial cover associated to $G$. Hence the
holonomy is
the same.

The assertion of the proposition then follows, since conjugation on the right
corresponds to different choices of geometrically distinguished systems of
paths used for an identification of the holonomy group with the abstract group
$\br_\mu$.
\qed[0mm]

\subsection{invariants}

Though keeping in mind the invariance along equivalence classes of singular
functions, we will for linguistic reasons define an invariant of singular
functions:

\begin{defi}
The class of braid monodromy homomorphisms associated to a truncated versal
unfolding of a singular function $f$ with isolated singularity is called the
{\it braid monodromy homomorphism} of $f$.
\end{defi}

An invariant which is easier to handle is obtained by considering the image of
the braid monodromy homomorphism only, up to conjugation of course.

\begin{defi}
The braid monodromy group associated to $f$ is the conjugacy class of
subgroups of $\br_{\mu(f)}$ given as the image of a braid monodromy
homomorphism.
\end{defi}

Supposing we have generators for the braid monodromy group associated to $f$,
we get a presentation of the local fundamental group of the discriminant
complement:

\begin{lemma}
\labell{cpl-pres}
Suppose the braid monodromy group associated to a singular function
$f$ is given
as the conjugacy class determined by a subgroup of $\br_\mu$
generated by braids $\b_1,...,\b_n$. Then the isomorphism class
$\pi_1(B_\e\times B_\rho\smin\hfami_\Delta)$ for $\e,\rho(\e)$ sufficiently
small is represented by the finitely presented group
$$
\begin{array}{ccc}
&
\lspan t_1,...,t_\mu |\, t_i\inv\b_j(t_i) \rspan.
\end{array}
$$
\end{lemma}

\proof
If braids $\b_j$ generate the image of the braid monodromy homomorphism, their
preimages $\cg_j$ and generators of the kernel generate the source.
Hence the claim follows from \ref{loc-pres}, since the trivial braid yields the
trivial relation only.
\qed

In fact we should also include an algebraic observation, which reduces the
number of relations dramatically in case the generators are conjugated to
braids non-trivial only on a few strands.\\

First we note that the choice of the
generators $\b_j$ and of generators of the free group does not matter.

\begin{lemma}
Suppose $B$ is a subgroup of $\br_n$ which acts on generators $t_1,...,t_n$ of
a free group $F_n$ by the Artin representation. Then the normal closure of the
subgroup of $F_n$ generated by
$$
w_i\inv\b_j(w_i)
$$
is independent of the choice of a finite set of
generators $\{\b_j\}$ of $B$ and a finite set of generators $\{w_i\}$ of $F_n$.
\end{lemma}

Now we can use this to reduce the number of relations in case a braid generator
$\b$ is conjugated to a twist.

\begin{lemma}
\labell{few-rela}
Suppose $\b=\b_0\s_1^l\b_0\inv\in\br_n$, then the normal subgroup generated by
$t_i\inv\s(t_i), i=1,...,n $, is equal to the normal subgroup generated by
$$
\b_0(t_1\inv)\b_0(\s_1^l(t_1)),\b_0(t_2\inv)\b_0(\s_1^l(t_2)).
$$
\end{lemma}

\proof
By the previous lemma the first normal subgroup is equal to the normal subgroup
generated by $\b_0(t_i\inv)\b_0(\s_1^l(t_i))$. Since these elements are trivial
except for $i=\in\{1,2\}$, the claim follows.
\qed


\chapter{Hefez Lazzeri unfoldings}

The singular functions to which this study is devoted are the Brieskorn-Pham
polynomials in arbitrary dimensions.

\begin{defi}
A polynomial $f\in\CC[x_1,...,x_n]$ is called a {\it Brieskorn Pham
polynomial}, if there are positive integers $l_1,...,l_n$ and
$$
f(x_1,...,x_n)=x_1^{l_1+1}+\cdots+x_n^{l_n+1}.
$$
\end{defi}

They exhibit a lot of symmetry, most apparent the invariance of the
polynomial under multiplication of the coordinate $x_i$ with $(l_i+1)^{st}$
roots of $1$.
But more important to us is the invariance of its singular values under
multiplication with $l_i^{th}$ unit roots, for this invariance persists to
the singular values of all functions obtained by linear perturbation terms.
It is due to these symmetries that singular values and the bifurcation can be
given by explicit polynomials, whereas in general such description seems
quite unattainable.

In \cite{HL}, Hefez and Lazzeri exploited the unfolding over the linear
perturbation terms to some extend, computing the intersection lattice for
Brieskorn polynomials.
In order to extend their exploits, we review their results on the
discriminant of $f$ unfolded over the space of linear monomials and determine
the bifurcation divisor, too. Moreover we show how this sheds light on various
geometric aspects of the unfolding and of the corresponding perturbed
functions.

\section{discriminant and bifurcation hypersurface}

Our current objective is to compute the corresponding discriminant and
bifurcation divisor, the discriminant for the unfolding of
$f\in\CC[x_1,...,x_n]$ given by
\begin{eqnarray*}
F(x,\a,z) & := &
f(x)+z-\sum_{i=1}^{n}\a_i(l_i+1)x_i,
\end{eqnarray*}
and the bifurcation divisor for the unfolding
\begin{eqnarray*}
F^\#(x,\a) & := &
f(x)-\sum_{i=1}^{n}\a_i(l_i+1)x_i.
\end{eqnarray*}

The partial derivatives are given by
\begin{eqnarray*}
\del_i F(x,\a,z)\,=\,\del_i F^\#(x,\a) & = &
(l_i+1)(x_i^{l_i}-\a_i),\quad i=1,...,n.
\end{eqnarray*}

Of course we shouldn't neglect the obvious quasi-homogeneity apparent in this
situation:

\begin{lemma}
The functions $F,\del_iF$ are weighted homogeneous of degree
one, resp.\ $\frac{l_i}{l_i+1}$, if $wt(x_i)=\frac1{l_i+1}$, $wt(z)=1$ and
$wt(\a_i)=\frac{l_i}{l_i+1}$.
\end{lemma}

It implies in particular that the critical set, the discriminant divisor and
the bifurcation divisor are quasi-homogeneous in their ambient affine
spaces with a good $\CC^*$-action.\\

We adopt now the following notation referring to roots of unity:
\begin{eqnarray*}
\xi_i & : & \text{the primitive $l_i^{th}$ root of 1 of least angle in
}]0,2\pi[.
\end{eqnarray*}

Then we can state the two central results of this section:

\begin{lemma}
\labell{laz-dis}
The polynomial defining the critical value divisor is given by the expansion
of the formal product
$$
\prod_{1\leq\nu_1\leq l_1}\cdots
\prod_{1\leq\nu_n\leq l_n}
\left(-z+\sum_{i=1}^nl_i\xi^{\nu_i}_i\a_i^{\frac{l_i+1}{l_i}}\right).
$$
\end{lemma}

\proof
The discriminant is the set of points $(\a_1,...,\a_n,z)$ such that
$F$ and its partial derivatives $\del_{x_i}F$ vanish simultaneously at a
point $(x,\a,z)$. Since the discriminant is known to be an algebraic
hypersurface we are thus looking for a reduced monic polynomial
$p_\Delta\in\CC[\a_1,...,\a_n,z]$ with zero set equal to the discriminant set.

Therefore we try to eliminate the variables $x_i$ from the system of equations
$F=\del_iF=0$. In a first step we replace $x_i^{l_i}$ by $\a_i$ in the equation
$F=0$ and get instead
\begin{eqnarray*}
z & = & \sum_{i=1}^n\a_il_ix_i.
\end{eqnarray*}

Now it is helpful to consider the Galois extension of $\CC(\a_1,...,\a_n)$
defined by the polynomials $x_i^{l_i}-\a_i$ which we denote by
$\CC(a_1,...,a_n)$ with $a_i^{l_i}=\a_i$.

Then the system of equations is easily solvable in
$\CC(a_1,...,a_n)[x_1,...,x_n,z]$ with $x_i=\xi_i^{\nu_i}a_i$ to be
eliminated to get
\begin{eqnarray*}
z\,=\, \sum_{i=1}^n\a_il_ix_i & = & \sum_{i=1}^nl_i\xi_i^{\nu_i}a_i^{l_i+1},
\quad 0\leq \nu_i <l_i.
\end{eqnarray*}

So the corresponding discriminant is simply given by the polynomial
$$
\prod_{1\leq\nu_i\leq l_i}
\left(-z+\sum_{i=1}^nl_i\xi^{\nu_i}_ia_i^{l_i+1}\right).
$$
It is in fact a polynomial in $\CC[\a_1,...,\a_n,z]$ since it is invariant
under the action of the Galois group of the extension which acts by
multiplication of $\xi_i$ on $a_i$.
\qed[3mm]

In particular the Milnor number $\mu$ of $f$, which coincides with the
$z$-degree of the discriminant, is thus shown to be
$$
\mu=\prod_{i=1}^n l_i.
$$

Similarly we get a description of the bifurcation divisor:

\begin{lemma}
\labell{laz-bif}
The polynomial defining the bifurcation divisor in the base of $F^\#$ is given
by the expansion of the formal product
$$
\prod_{\begin{array}{c}\\[-7.5mm] \scriptscriptstyle
1\leq\nu_i\leq l_i \\[-2.5mm] \scriptscriptstyle
1\leq \eta_i\leq l_i\\[-2.5mm] \scriptscriptstyle
\pmb\nu\neq\pmb\eta
\end{array}}
\left(\sum_{i=1}^n l_i\a_i^{\frac{l_i+1}{l_i}}(\xi^{\nu_i}_i
-\xi^{\eta_i}_i)\right).
$$
\end{lemma}

\proof
Of course the polynomial is uniquely defined up to a constant as the
discriminant with respect to the variable $z$ of the polynomial defining the
discriminant divisor, given in lemma \ref{laz-dis}.

Passing to the Galois extension of $\CC(\a_1,...,\a_n)$ once again, we have
to compute the discriminant of a polynomial which is a product of linear
factors in $\CC(a_1,...,a_n)[z]$.
So up to a constant the discriminant is the product of the squares of the
mutual differences between the zeroes of distinct factors.
\begin{eqnarray*}
discr_z(p_\Delta) & = &
\prod_{{{1\leq\nu_i,\eta_i\leq l_i}} \atop
{\pmb\nu\neq\pmb\eta}}
\left(\sum_{i=1}^n l_i a_i^{l_i+1}(\xi^{\nu_i}_i
-\xi^{\eta_i}_i)\right).
\end{eqnarray*}
Since this polynomial is Galois invariant, the claim follows as above.
\qed[3mm]

Due to the importance of the bifurcation divisor this polynomial deserves a
proper name, it will be denote by $p_\bfami$ and henceforth be called the
{\it Hefez-Lazzeri polynomial} of $f$.

We end this section with some corollaries, highlighting some of the nice
geometric properties of the Hefez-Lazzeri unfoldings of Brieskorn Pham
polynomials.

\begin{lemma}
The polynomial $p_\bfami$ vanishes to order exactly $(l_i^2-1)\mu/l_i$ along
the hyperplane $\a_i=0$.
\end{lemma}

\proof
We have to show that $\a_i$ exponentiated to the given order is a divisor of
$p_\bfami$ while higher powers are not.
 From the factorization of $p_\bfami$ in $\CC[a_1,...,a_n]$, we see that $a_i$
divides a factor if and only if $\nu_j=\eta_j$ for all $j\neq i$. An easy
check shows that $(l_i^2-1)\mu$ is the vanishing order of $a_i$. So the
claim follows, for $a_i^{l_i}=\a_i$.
\qed[0mm]

\begin{lemma}
\labell{fdiscr-ind}
The leading coefficient of $p_\bfami\in\CC[\a_1,...,\a_n]$ with respect to the
variable $\a_n$ is the $l_n^{th}$ power of the Hefez Lazzeri polynomial of the
Brieskorn Pham function
$g(x_1,...,x_{n-1})=x_1^{l_1+1}+\cdots+x_{n-1}^{l_{n-1}+1}$.
\end{lemma}

\proof
In the factorization of $p_\bfami$ in $\CC[a_1,...,a_n]$ we note that a
factor is either a polynomial in $a_n$ with constant leading coefficient or
does not contain $a_n$ at all.
We conclude that up to a constant the leading coefficient with respect to the
variable $a_n$ is the product of those factors not containing $a_n$. So we
get this coefficient by collecting the factors with $\nu_n=\eta_n$ from the
factorization of $p_\bfami$:
$$
\prod_{\nu_n=\eta_n\leq l_n}
\prod_{{{\nu_i,\eta_i\leq l_i} \atop { i< n}} \atop
{\pmb\nu\neq\pmb\eta}}
\left(\sum_{i=1}^n l_ia_i^{l_i+1}(\xi^{\nu_i}_i-\xi^{\eta_i}_i)
\right)
=
\prod_{{{\nu_i,\eta_i\leq l_i} \atop { i< n}} \atop
{\pmb\nu\neq\pmb\eta}}
\left(\sum_{i=1}^{n-1} l_ia_i^{l_i+1}(\xi^{\nu_i}_i-\xi^{\eta_i}_i)
\right)^{l_n}
$$
Compare this product with the proof of \ref{laz-bif} and the claim is
immediate.
\qed

\begin{lemma}
Given a Brieskorn Pham polynomial $f$, then the generic degenerations in its
Hefez-Lazzeri unfolding are:
\begin{enumerate}
\item
at a generic point of the coordinate hyperplane $\a_i=0$, the corresponding
function has $l_i'=\mu/l_i$ critical points of type $A_{l_i}$ with distinct
critical values,
\item
at any point on the coordinate hyperplanes, the corresponding function is of
Brieskorn Pham type,
\item
at a generic point of any other component of the Lazzeri discriminant, the
corresponding projection has critical points of type $A_1$ with at least two
of common value.
\end{enumerate}
\end{lemma}

\proof
The Hessian of the function $F^\#(x,\a)$ with respect to the variables $x$
only is given by the diagonal matrix with entries
$l_i(l_i+1)x^{l_i-1}\delta_{i,j}$, which is of full rank outside the
hyperplanes $x_i=0$.

 From the gradient of $F^\#$ with respect to $x$ we deduce, that there is a
critical point on the $i^{th}$ hyperplane $x_i=0$ if and only if the
parameter $\a$ is on the hyperplane $\a_i=0$.

Hence for each parameter in the complement of the hyperplanes $\a_i=0$ the
Hessian of the corresponding function is of full rank, so all its critical
points are of type $A_1$ only.

Thus the parameter belongs to the bifurcation divisor if and only if the
function maps at least two of its critical points to the same value.

The bifurcation locus outside the coordinate hyperplanes is therefore part of
the Maxwell stratum and corresponds to the transversal intersection of
several branches of the discriminant locus, the number of which is half the
local degree of the bifurcation locus plus one.

On the other hand a parameter on a coordinate hyperplane is considered
generic if it does not belong to any other hyperplane.
Assume in this case w.l.o.g.\ that $\a_n=0$ , so the corresponding function is
\begin{eqnarray*}
F^\#(x,\a) & = & f(x)+\sum_{i=1}^{n-1}\a_i(l_i+1)x_i.
\end{eqnarray*}
At each point of its critical locus, the Hessian is of corank one,
therefore the type of the critical point is easily seen to be $A_{l_n}$.

The generalization to arbitrary parameters on the coordinate hyperplanes can
be checked to yield the second claim.
\qed[0mm]

\begin{lemma}
Suppose the exponents are in increasing order, $l_1\leq ...\leq l_n$, then
the total degree of the Hefez Lazzeri polynomial $p_\bfami$ is
$$
\mu\left(\sum_{i=1}^n\bigg(\frac{l_i^2-1}{l_i}\prod\limits_{i<
j\leq n}l_j\bigg)\right).
$$
\end{lemma}

\proof
In case $n=1$ the degree is $l^2-1$ in accordance with the claim since
$$p_\bfami=
\prod_{{{1\leq\nu,\eta\leq l}} \atop
{\nu\neq\eta}}
\left(l\a^{\frac{l+1}{l}}(\xi^{\nu}
-\xi^{\eta})\right).
$$
In case $n>1$ we argue by induction:
In each formal factor of $p_\bfami$ only pure monomials occur, so if $\a_1$
occurs, the corresponding monomial is the term of highest degree. The number
of such factors is
$\mu^2(l_1-1)/(l_1)$. The product of the other factors is nothing else
but the leading coefficient with respect to the variable $\a_1$.

Hence by lemma \ref{fdiscr-ind} and the induction hypothesis the total degree
is
\begin{eqnarray*}
&&
\mu^2\frac{l_1^2-1}{l_1^2}+
\mu\left(\sum_{i=2}^n\bigg(\frac{l_i^2-1}{l_i}\prod\limits_{i< j\leq
n}l_j\bigg)\right).
\end{eqnarray*}
But this is the claim, since $\mu/l_1=l_2\cdots l_n$.
\qed[0mm]

\section{Hefez Lazzeri path system}
\label{hl-path}

In case of the Brieskorn Pham polynomials Hefez and Lazzeri described an
method to define a path system using induction on the dimension for a regular
fibre obtained by a linear perturbation under the sole assumption that the
parameters $\a_i$ are of quite distinctive magnitude.

So if $\a_n\ll...\ll\a_2\ll\a_1$ the critical values are distributed on circles
of radius $l_n\a_n^{\frac{l_n+1}{l_n}}$ centred at the critical values of the
polynomial
$$
x_1^{l_1+1}-\a_1(l_1+1)x_1+
\cdots+x_{n-1}^{l_{n-1}+1}-\a_{n-1}(l_{n-1}+1)x_{n-1}
$$
according to lemma \ref{laz-dis}.

Therefore we need to replace all paths by $l_n$ copies and refine the system on
these discs then.\\

Inductively we start with a path system of a first kind, as depicted in figure
$1$.

\setlength{\unitlength}{3mm}
\begin{picture}(30,20)(-21,-10)

\put(0,5){\circle*{.2}}
\put(0,-5){\circle*{.2}}
\put(5,0){\circle*{.2}}
\put(-5,0){\circle*{.2}}

\put(-10,-10){\setlength{\unitlength}{6mm}\begin{picture}(10,10)
\bezier{30}(0,5)(0,7.5)(2,9)
\bezier{30}(2,9)(5,11)(8,9)
\bezier{30}(8,9)(10,7.5)(10,5)
\bezier{30}(10,5)(10,2.5)(8,1)
\bezier{30}(8,1)(5,-1)(2,1)
\bezier{30}(2,1)(0,2.5)(0,5)
\end{picture}}

\bezier{130}(5,0)(6,2)(10,0)

\bezier{30}(0,5)(.5,5)(1,5)
\bezier{300}(1,5)(3,5)(10,0)

\bezier{390}(-5,0)(-5,6)(1,6)
\bezier{390}(1,6)(5,6)(10,0)

\bezier{400}(0,-5)(-6,-5)(-6,.8)
\bezier{400}(-6,.8)(-6,7)(1,7)
\bezier{400}(1,7)(6.5,7)(10,0)

\end{picture}
\begin{figure}[ht]
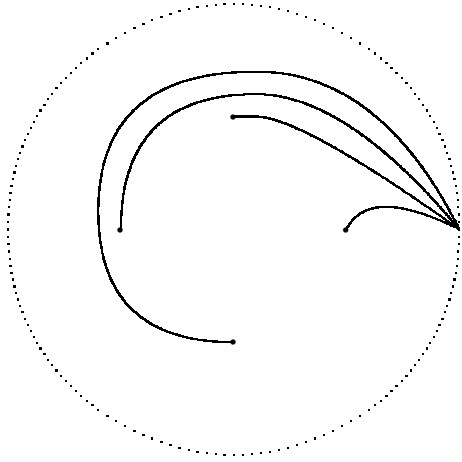

\begin{center}
\caption{path system of the first kind, $l=4$}
\end{center}
\end{figure}

For the second step we replace it by the system of the second kind, which is
given in the second figure, except for the fact that we did not care to
distinguish the copies we take of each single path.

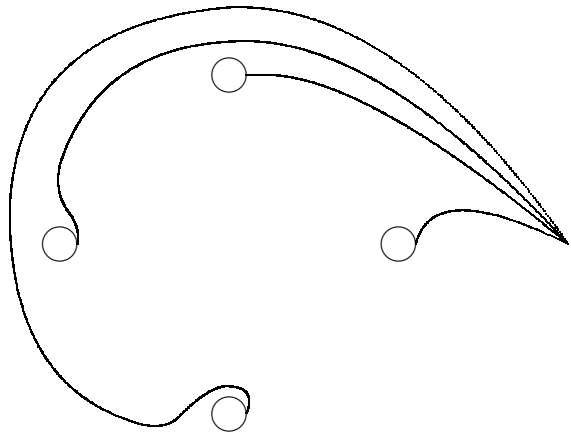
\begin{figure}[ht]
\begin{center}
\setlength{\unitlength}{.45mm}
\begin{picture}(180,150)(-70,-70)

\put(50,0){\circle{10}}
\put(0,50){\circle{10}}
\put(0,-50){\circle{10}}
\put(-50,0){\circle{10}}


\bezier{130}(55,0)(60,20)(100,0)


\bezier{30}(5,50)(8,50)(11,50)
\bezier{300}(11,50)(40,50)(100,0)

\bezier{300}(-50,24)(-52,16)(-48,10)
\bezier{300}(-48,10)(-44,5)(-45,0)
\bezier{300}(-50,24)(-38,60)(5,60)
\bezier{300}(5,60)(45,60)(100,0)

\bezier{80}(5,-50)(8,-42)(-1,-42)
\bezier{80}(-1,-42)(-6,-42)(-15,-51)

\bezier{80}(-15,-51)(-20,-56)(-30,-52)
\bezier{200}(-30,-52)(-65,-40)(-65,10)
\bezier{200}(-65,10)(-63,65)(00,70)
\bezier{200}(0,70)(50,72)(100,0)
\end{picture}
\caption{a path system of second kind in case $l=4$}
\end{center}
\end{figure}

Into each of the small discs we have to paste a path system of the first kind
to get a path system in case $n=2$.

Inductively then, we have to replace all systems of the first kind pasted in
step $n-1$ by path systems of the second kind and to paste systems of the first
kind into the new smaller discs left by the systems of the second kind.

These path systems are called {\it Hefez Lazzeri path systems}.

\chapter{singularities of type $A_n$}

The simplest singular functions for which we have to determine the braid
monodromy are those of type $A_n$.
So let us investigate the holomorphic function $f:x\mapsto x^{n+1}$.
A versal unfolding can be chosen to be quasi-homogeneous
$$
F(x,u,z)=x^{n+1}+z+\sum_{i=1}^{n-1}u_ix^i,
$$
while the associated truncated versal deformation is
$$
F^\#(x,u)=x^{n+1}+\sum_{i=1}^{n-1}u_ix^i.
$$
This case of functions of type $A_n$ is special in some aspects.
First, being quasi-homogeneous with good $\CC^*$-action, we may work in the
affine set up, \ref{qhomogen}, and avoid the cumbersome notation used for
germs. So we consider the projection $p$ from $\CC^n$ to $\CC^{n-1}$ of the
base spaces of the versal resp.\ truncated versal unfoldings above and get the
braid monodromy as the coefficient homomorphism for the bundle which we obtain
by restriction of $p$ to
$$
p:\CC^n\smin\Delta\cup p\inv(\bfami)\tto\CC^{n-1}\smin\bfami.
$$
Second, being of codimension one, there is another map besides $p$ to
which the theory of polynomial covers applies.
Consider the projection $\ptil$ from $\CC^{n+1}$ to $\CC^n$ from the domain of
the versal deformation $F$ to its base and the zero set
$F\inv(0)\subset\CC^{n+1}$.
Then there is a restriction, which again is a polynomial covering space, but
of degree $n+1$:
$$
\ptil:\CC^{n+1}\smin F\inv(0)\cup\ptil\inv(\Delta)\tto\CC^n\smin\Delta.
$$
Of course $\ptil$ is induced from the tautological bundle (\ref{taut}),
p.\pageref{taut}, we have seen before, but even more is true. It is a subbundle
which is a strong deformation retract.

When we now attack the braid monodromy of $f$ -- recall this to be the
braid monodromy of $p$ -- our argument is not straightforward, but relies
heavily on properties of $\ptil$ and citations from the literature.

But setting out with the fibration $p$, let us first choose a base point,
such that we can describe the corresponding fibre in more detail.
To rely on previous results we choose $u=(0,...,0,-\a(n+1))\in\CC^{n-1}$
from the base of the Hefez-Lazzeri unfolding
with corresponding function $f_\a=x^{n+1}-\a(n+1)x$.
The fibre of $p$ at $u$, which we henceforth refer to by $L_\a$, falls into
the topological class we denoted by $\CC_n$, but we easily compute from
$\del_x f_\a=0$ that the critical locus is at the $n^{th}$ roots of $\a$ and
that the fibre is punctured at the critical values which coincide with the
$n^{th}$ roots of $\a^{n+1}$ multiplied by $n$, cf.\ lemma \ref{laz-dis}. Thus
we get even the true geometric picture for $L_\a$:

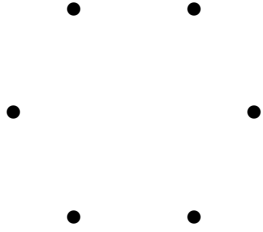
\begin{figure}[h]
\begin{center}
\setlength{\unitlength}{1.6mm}
\begin{picture}(30,20)(0,2)
\linethickness{1pt}
\put(10,5){\circle*{1}}
\put(20,5){\circle*{1}}
\put(5,13.7){\circle*{1}}
\put(10,22.3){\circle*{1}}
\put(25,13.7){\circle*{1}}
\put(20,22.3){\circle*{1}}
\end{picture}
\end{center}
\end{figure}

This punctured line $L_\a$ is a subset in the base of the bundle $\ptil$, and
as we said, an understanding of the bundle $\ptil$ will eventually lead to
the braid monodromy of $p$.

Looking for a distinguished regular base point of the restricted bundle
$\ptil_\a:=\ptil|_{L_\a}$ the origin is the obvious choice.
The corresponding fibre is an affine line punctured at the zeroes of
$x^{n+1}-\a(n+1)x$, i.e..\ at $0$ and the $n^{th}$ roots of $\a(n+1)$.

\begin{figure}[h]
\begin{center}
\setlength{\unitlength}{2.2mm}
\begin{picture}(30,20)(0,3.5)
\linethickness{1pt}
\put(10,5){\circle*{1}}
\put(20,5){\circle*{1}}
\put(5,13.7){\circle*{1}}
\put(15,13.7){\circle*{1}}
\put(10,22.3){\circle*{1}}
\put(25,13.7){\circle*{1}}
\put(20,22.3){\circle*{1}}
\end{picture}
\end{center}
\end{figure}

But we need more information to compute the monodromy.
Assuming $\a$ to be real, we consider the punctures in the fibres of $\ptil$
over points $\rho\geq0$ of the positive real axis in the base $L_\a\cong\CC$.
We make the following observations:
\begin{enumerate}
\item
By the sign rule of Descartes the number of positive real zeroes of the
polynomial $x^{n+1}-\a(n+1)x+\rho$ is either two or none.
\item
Let $\rho_\a:=n\sqrt[n]{\a^{n+1}}$ be a positive real by choosing
$\sqrt[n]{\a^{n+1}}$ the only positive real among the $n^{th}$ roots of
$\a^{n+1}$. Then $\rho_\a$ is the base of the only singular fibre of $\ptil$
over the positive real axis, with a single ordinary double point
$x=\sqrt[n]{\a}$.
\item
Both facts together imply that the number of positive real zeroes along the
real axis $\rho\in\RR^{\geq0}$ is constantly two over $[0,\rho_\a[$, a zero
of multiplicity two at $\rho_\a$ and vanishes over $]\rho_\a,\infty[$.
\end{enumerate}

We elaborate on these observations to get the proofs of the following two
lemmas concerning the monodromy homomorphism to $\permgroup_{n+1}$ for the
cover associated to $\ptil$, resp.\ the braid monodromy of
$\ptil$ itself.
Though the first would certainly suffice in the course of the
argument, we couldn't help including the second, and if only to show the
particular flavour of arguments applied in braid monodromy considerations.

\begin{lemma}
\labell{a/mono}
\sloppypar
Given the simple branched cover of $\CC$ by $\ptil$ restricted to $F\inv(0)$
over the complex line $(0,-\a(n+1),z)$, $\a\neq0$ fix, in the base of the
versal unfolding $F$.
Suppose $\a_n$ is any $n^{th}$ root of $\a$, then the Hurwitz monodromy of
the cover associates to the radial path from $0$ around $n\a\a_n$ the
transposition of the points $0$ and $\a_n\sqrt[n]{n+1}$.
\end{lemma}

\proof
We just considered the radial path in the special case of
$\a,\a_n\in\RR^{>0}$. In that case the local monodromy is the
transposition of the merging points which originate in the said points $0$ and
$\a_n\sqrt[n]{n+1}$ of the fibre at the origin.

The setting is quasi-homogeneous with respect to weights $wt(x)=1$,
$wt(\a)=n$, $wt(z)=n+1$.
That's why multiplication by a suitable factor maps the special case
bijectively to any other and so the proof is obtained.
\qed

\begin{lemma}
Consider the cover $\ptil$ restricted to the complex line
$(0,-\a(n+1),z)$, $\a\neq0$ fix, in the base of the versal unfolding $F$.
Suppose $\a_n$ is any $n^{th}$ root of $\a$, then the braid monodromy
associates to the radial path from $0$ around $n\a\a_n$ in the base the twist
on the radial arc from $0$ to $\a_n\sqrt[n]{n+1}$ in the fibre.
\end{lemma}

\proof
By quasi-homogeneity it suffices to prove the special case of $\a\in\RR^{>0}$,
$\a_n\in\RR^{>0}$.
In that case the local monodromy is a twist since the degeneration is an
ordinary double point.
In a fibre sufficiently close to the degenerate fibre this is a twist on the
straight arc between the merging points.

We choose the nearby fibre at a positive real base parameter and get an arc
supported on the positive real axis in the fibre.

In the description of the general strategy for the computation of braid
monodromies we have mused for some time on the fact that it suffices to
understand the parallel transport of this arc to the reference fibre, which
is the fibre at the origin in the case at hand.

So we have to prove that the arc in a close by fibre can be transported to
the radial arc from $0$ to $\a_n\sqrt[n]{n+1}$ in the fibre at the origin.
In fact we can do so by choosing a suitable vector field.
Since the punctures have to be transported along integral curves, the vector
field is determined at the punctures, but otherwise we may smoothly
interpolate as we need.
At the two punctures moving on the positive real axis the vector field is
parallel to the axis. Moreover no other puncture hits the positive real axis,
as we observed above, due to the sign rule of Descartes.
So we interpolate by a vector field on the ray parallel to the axis, such
that the arc in the nearby fibre is stretched to the claimed interval.
\qed

As was noted previously, \ref{quot/triv}, the fundamental group
$\pi_1(\CC^n-\Delta)$ is a quotient of $\pi_1(\CC_n)$ such that the braid
automorphisms of
$\pi_1(\CC_n)$ descend to trivial automorphisms of $\pi_1(\CC^n-\Delta)$.

So the braid monodromy homomorphism of $p$ with range considered as
automorphism of $\pi_1(\CC_n)$ may map only to elements $a_*(\b)$ which
stabilize the coefficient homomorphism $\atil_*$ for $\ptil$ followed by the
natural homomorphism $\pi:\br_{n+1}\to\sfami_{n+1}$:

\begin{eqnarray*}
\pi_1(\CC_n) & \stackrel{\pi\circ \atil_*}{\tto} & \sfami_{n+1}\\
a_*(\b)\big\downarrow & & \big\downarrow\id\\
\pi_1(\CC_n) & \stackrel{\pi\circ \atil_*}{\tto} & \sfami_{n+1}
\end{eqnarray*}

To make the argument more explicit we introduce geometric bases in the
two punctured lines.
If we do so for the fibre of $p$ punctured at the $n^{th}$ roots of
$\a^{n+1}$ dilated by $n$ in the way as depicted on the left and for the fibre
of
$\ptil$ punctured at zero and the $n^{th}$ roots of $\a(n+1)$ as shown on the
right

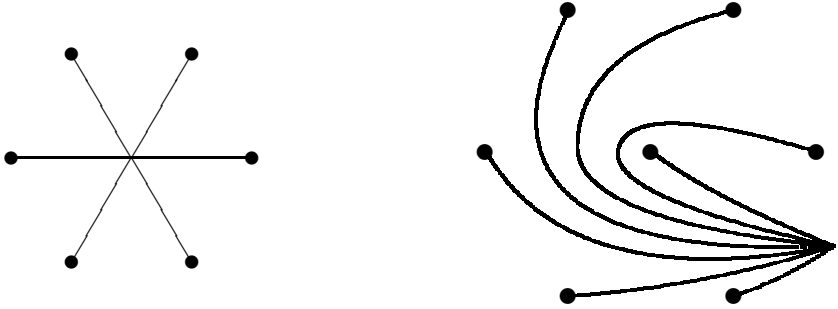
\begin{figure}[h]
\begin{center}
\setlength{\unitlength}{2mm}
\begin{picture}(60,25)

\put(-2,3.5){\setlength{\unitlength}{1.6mm}
\begin{picture}(30,20)(0,2)
\put(10,5){\circle*{1}}
\put(20,5){\circle*{1}}
\put(5,13.7){\circle*{1}}
\put(10,22.3){\circle*{1}}
\put(25,13.7){\circle*{1}}
\put(20,22.3){\circle*{1}}

\drawline(15,13.7)(10,5)
\drawline(15,13.7)(20,5)
\drawline(15,13.7)(10,22.3)
\drawline(15,13.7)(25,13.7)
\drawline(15,13.7)(20,22.3)
\drawline(15,13.7)(5,13.7)
\end{picture}
}

\put(28,2){\setlength{\unitlength}{2.2mm}
\begin{picture}(30,20)(0,3.5)
\linethickness{1pt}
\put(10,5){\circle*{1}}
\put(20,5){\circle*{1}}
\put(5,13.7){\circle*{1}}
\put(15,13.7){\circle*{1}}
\put(10,22.3){\circle*{1}}
\put(25,13.7){\circle*{1}}
\put(20,22.3){\circle*{1}}

\bezier{40}(26,8)(23,6)(20,5)
\bezier{80}(26,8)(18,5.5)(10,5)
\bezier{120}(26,8)(18.5,11)(15,13.7)
\bezier{180}(26,8)(10,5)(5,13.7)
\bezier{190}(26,8)(2,7)(10,22.3)
\bezier{120}(26,8)(10.5,9.5)(10.5,14)
\bezier{120}(10.5,14)(10.5,20)(20,22.3)
\bezier{120}(26,8)(12,11)(13,14)
\bezier{120}(13,14)(14,17)(25,13.7)
\end{picture}
}

\end{picture}
\end{center}
\end{figure}
then the mapping class group is identified with
$\br_{n+1}=\lspan\s_i|...\rspan$
in such a way that the twists on the radial arcs considered in the lemma
above correspond to $\s_{1,2},\s_{1,3},...,\s_{1,l+1}$.

The coefficient homomorphism $a_*$ is
given on the generator system $\lspan t'_1,...,t'_n|\rspan$ of $F_n$ thus
determined by
$$
t_i' \mapsto \s_{1,i+1}.
$$
If we work instead with a geometric basis given by the next figure, the map
$a_*$ is given on the generator set $t_1,...,t_n$ thus defined by
$$
F_n\to\br_{n+1},\quad t_i\mapsto \s_i \quad\big( \,
\stackrel{\pi}{\mapsto}\,(\,i\,\,i\!+\!\!1\,)\,\big)
$$
as can be checked inductively.
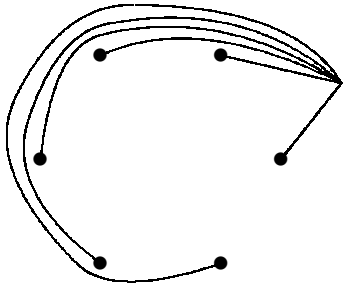
\begin{center}
\setlength{\unitlength}{2mm}
\begin{picture}(60,20)

\put(16,.5){\setlength{\unitlength}{1.6mm}
\begin{picture}(30,20)(0,2)
\put(10,5){\circle*{1}}
\put(20,5){\circle*{1}}
\put(5,13.7){\circle*{1}}
\put(10,22.3){\circle*{1}}
\put(25,13.7){\circle*{1}}
\put(20,22.3){\circle*{1}}

\bezier{140}(30,20)(27.5,16.85)(25,13.7)
\bezier{140}(30,20)(25,21.15)(20,22.3)
\bezier{140}(30,20)(18,26)(10,22.3)
\bezier{140}(30,20)(20,26.5)(10,24)
\bezier{140}(10,24)(6,23)(5,13.7)
\bezier{140}(30,20)(22.25,27)(11,25)
\bezier{140}(11,25)(6.5,24.5)(4,17)
\bezier{140}(4,17)(2,11)(10,5)
\bezier{140}(30,20)(26,26.5)(12.5,26.5)
\bezier{140}(12.5,26.5)(6.75,26.5)(3,19)
\bezier{140}(3,19)(0,13)(8,5)
\bezier{140}(8,5)(11,2)(20,5)
\end{picture}
}

\end{picture}
\end{center}
Not least do we prefer the second choice since the corresponding
group of stabilizing automorphisms in the image of the Artin
homomorphism $\br_n\to\Aut(F_n)$ has been investigated in several aspects which
-- when brought together -- are sufficient to yield our goal.

Let us first cite from the results of Catanese/Wajnryb \cite{CW}, Kluitmann
\cite{klui}, and D\"orner \cite{Doe}, who gave the stabilizer group for the
naturally induced homomorphism to the symmetric group $\sfami_{n+1}$ in terms
of generators as follows.

\begin{lemma}
\label{perm/stab}
Suppose $h:\lspan t_1,...,t_n|\,\rspan\to\sfami_{n+1}$ is given on
generators by
$t_i\mapsto (\,i\,\,i\!+\!\!1\,)$ and $\a:\br_n\to\Aut(\lspan
t_1,...,t_n|\,\rspan)$ is the Artin homomorphism, then the set of braids
$\b\in\br_n$ such that
$h\circ\a(\b)=h$ is a subgroup, called the braid stabilizer group of $h$,
which is generated by
$$
\d_2^3,...,\d_n^{n+1},
$$
or equivalently by
$$
\s_i^3,\s_{i,j}^2,|i-j|\geq2.
$$
\end{lemma}

\proof
We refer to \cite{CW} and give only the relations between the two sets of
generators which can be proved inductively as in \cite{klui}, \cite{Doe}:
\begin{eqnarray*}
\d_i^{i+1} & = &
\big(\s_{1}^3(\s_{1,3}^2\cdots\s_{1,i}^2)\big)\big(\s_2^3(\s_{2,4}^2\big)
\cdots\s_{2,i}\big)\cdots\big(\s_{i-1}^3\big)
\end{eqnarray*}
respectively for the other direction:
\begin{eqnarray*}
\s_{1}^3 & = & \d_2^3\\
\s_{1,3}^2 & = & \d_2^{-3}(\d_3^4)^2\d_2^{-3}\d_3^{-4}\\
\s_{1,i+1}^2 & = &
\d_i^{-i-1}\d_{i+1}^{i+2}\d_{i+1}^{i+2}\d_i^{-i-1}\d_{i+1}^{-i-2}\\
\s_{i}^3 & = & (\d_n^{n+1})^i\s_1^3(\d_n^{n+1})^{-i}\\
\s_{i+1,j}^2 & = & (\d_n^{n+1})^i\s^2_{1,j-i}(\d_n^{n+1})^{-i}
\end{eqnarray*}
\qed[7mm]

The stabilizer of the braid monodromy homomorphism with respect to the same
geometric bases can be deduced now. Of course it must be a subgroup of the
stabilizer group just computed and a direct check shows that it is even the
whole group.

\begin{lemma}
\label{braid/stab}
Suppose $h:\lspan t_1,...,t_n|\,\rspan\to\br_{n+1}$ is given on generators by
$t_i\mapsto \s_i$ and $\a:\br_n\to\Aut(\lspan
t_1,...,t_n|\,\rspan)$ is the Artin homomorphism, then the set of braids
$\b\in\br_n$ such that
$h\circ\a(\b)=h$ is a subgroup, called the braid stabilizer group of $h$,
which is generated by
$$
\d_1^3,...,\d_n^{n+1},
$$
or equivalently by
$$
\s_i^3,\s_{i,j}^2,|i-j|\geq2.
$$
\end{lemma}

\proof
We show the stabilizing property for the second set of generators. To do so
we recall the automorphisms associated to $\s_i\in\br_n$:
$$
\sigma_i(t_k)=\begin{cases}  t_it_{i+1}t_i^{-1}&\text{ if } k=i,\\
t_i&\text{ if } k=i+1,\\
t_k &\text{ otherwise.}
\end{cases}
$$
and compute then the automorphisms associated to $\s_i^3,\s_{i,j}^2$:
\begin{eqnarray*}
\sigma^3_i(t_k) & = & \left\{\begin{array}{ll}
t_it_{i+1}t_it_{i+1}t_i^{-1}t_{i+1}\inv t_i\inv&\text{ if } k=i,\\
t_it_{i+1}t_it_{i+1}\inv t_i\inv&\text{ if } k=i+1,\\
t_k &\text{ otherwise.}\end{array}\right.\\
\sigma_{i,j}(t_k) & = & \left\{\begin{array}{ll}
t_it_jt_i\inv&\text{ if } k=i,\\
t_i&\text{ if } k=j,\\
t_it_j\inv t_kt_jt_i\inv&\text{ if } i<k<j,\\
t_k &\text{ otherwise.}\end{array}\right.\\
\sigma_{i,j}^2(t_k) & = & \left\{\begin{array}{ll}
t_it_jt_it_j\inv t_i\inv&\text{ if } k=i,\\
t_it_jt_i\inv&\text{ if } k=j,\\
t_it_jt_i\inv t_j\inv t_k(t_it_jt_i\inv t_j\inv)\inv&\text{ if } i<k<j,\\
t_k &\text{ otherwise.}\end{array}\right.
\end{eqnarray*}
So we are left checking for $i<n-1$ and $i<k<j$:
\begin{eqnarray*}
\s_i & = &
\s_i\s_{i+1}\s_i\s_{i+1}\s_i^{-1}\s_{i+1}\inv \s_i\inv\\
\s_{i+1} & = &
\s_i\s_{i+1}\s_i\s_{i+1}\inv \s_i\inv\\
\s_i & = &
\s_i\s_j\s_i\s_j\inv \s_i\inv\\
\s_j & = &
\s_i\s_j\s_i\inv\\
\s_k & = &
\s_i\s_j\s_i\inv \s_j\inv \s_k(\s_i\s_j\s_i\inv \s_j\inv)\inv
\end{eqnarray*}
which follow easily from the braid relations.
\qed

The counterpart of these results is provided by the paper of Looijenga,
\cite{loo}. He already investigated -- in different terminology -- the
coefficient homomorphism of the projection $p$ and he observed:

\begin{lemma}
The image of the coefficient homomorphism of the projection $p$ coincides
with the stabilizer of the monodromy homomorphism of the finite cover
associated to $\ptil$.
\end{lemma}

\proof
We have seen before that there is an inclusion of these groups. So the
seemingly weaker result of Looijenga that both groups are conjugate subgroups
of $\br_n$ immediately implies the stronger form.
\qed

We must admit, that we fall short of finding the braid monodromy for the
projection $p$, but what we have determined is the braid monodromy group of
$p$ with respect to the base choices specified above.

\begin{lemma}
\labell{brmon-An}
The braid monodromy group of the singular function $x\mapsto x^{n+1}$ is
identified by the geometric bases above with the subgroup of $\br_{n}$
generated by
$$
\s_i^3,\s_{i,j}^2,|i-j|\geq2.
$$
\end{lemma}

\proof
The braid monodromy group is the image of the coefficient homomorphism. By
the result of Looijenga this coincides with the stabilizer, so the claim
follows.
\qed[3mm]

For later application we translate this result back into a statement on
mapping classes of the reference fibre at a positive real $\a$.

\begin{lemma}
\label{An-mono}
A set of mapping classes which generate the braid monodromy in the fibre at
$(0,...,0,-\a(n+1))$ punctured at the roots $n\a^{\frac{n+1}n}$
is given by
\begin{enumerate}
\item
the $\frac32$-twists on the straight arcs joining consecutive punctures,
\item
the full twists on arcs joining non-consecutive punctures in the
complement of the inscribed polygon and the open cone defined by the $1$
and the last puncture.
\end{enumerate}
\end{lemma}

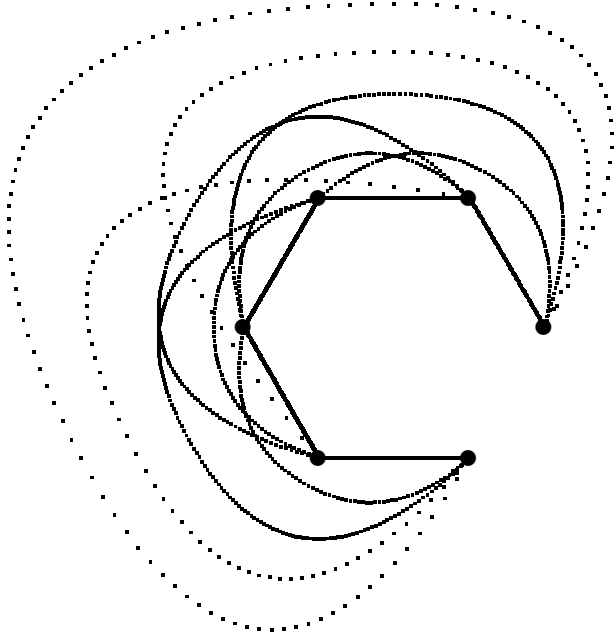
\begin{figure}[h]
\begin{center}
\setlength{\unitlength}{2mm}
\begin{picture}(30,40)(-6,-6)
\linethickness{1pt}
\put(10,5){\circle*{1}}
\put(20,5){\circle*{1}}
\put(25,13.7){\circle*{1}}
\put(20,22.3){\circle*{1}}
\put(5,13.7){\circle*{1}}
\put(10,22.3){\circle*{1}}
\bezier{190}(10,5)(15,5)(20,5)
\bezier{190}(10,22.3)(7.5,18)(5,13.7)
\bezier{190}(10,5)(7.5,9.3)(5,13.7)
\bezier{100}(10,22.3)(15,22.3)(20,22.3)
\bezier{100}(20,22.3)(22.5,18)(25,13.7)
\bezier{40}(3,13.7)(3,20.1)(10,22.3)
\bezier{40}(3,13.7)(3,7.3)(10,5)
\bezier{40}(9,3.3)(3.5,6.5)(5,13.7)
\bezier{40}(9,3.3)(14.6,.1)(20,5)
\bezier{40}(9,24)(3.5,20.9)(5,13.7)
\bezier{40}(9,24)(14.6,27.2)(20,22.3)
\bezier{40}(21,24)(26.5,20.9)(25,13.7)
\bezier{40}(21,24)(15.4,27.2)(10,22.3)
\bezier{80}(10,22.3)(-5.5,18.1)(1.5,5.9)
\bezier{80}(1.5,5.9)(8.5,-6.2)(20,5)
\bezier{80}(10,5)(-5.5,9)(1.5,21.5)
\bezier{80}(1.5,21.5)(8.5,33.5)(20,22.3)
\bezier{80}(15,29.2)(1.1,29.2)(5,13.7)
\bezier{80}(15,29.2)(29.9,29.2)(25,13.7)
\bezier{30}(25,13.7)(34,32)(15,32)
\bezier{40}(15,32)(-13,32)(10,5)
\bezier{30}(20,5)(7,-10.6)(-1.5,4.1)
\bezier{40}(-1.5,4.1)(-15.5,28.3)(20,22.3)
\bezier{30}(25,13.7)(37.5,35.2)(15,35.2)
\bezier{30}(20,5)(7.5,-16.5)(-4.4,2.4)
\bezier{50}(-4.4,2.4)(-23.3,35.2)(15,35.2)
\end{picture}
\end{center}
\end{figure}

\proof
We have to show that the generators in lemma \ref{brmon-An} and the arcs
described in the assertion are related as claimed by the geometric basis
chosen above.

More precisely, the generators with exponent $3$ correspond to
$\frac32$-twists on arcs obtained by joining consecutive paths of the
geometric bases up to isotopy, the generators with exponent $2$ correspond to
full twists on arcs obtained by joining non-consecutive paths.

The arcs thus obtained can be characterized as in the assertion due to the
fact that the complement of the inscribed polygon and the given cone is
simply connected and contains the geometric basis.
\qed[2mm]

Having thus computed our first braid monodromy group, we can deduce by
\ref{cpl-pres} a presentation for the fundamental group of the discriminant
complement of the $A_n$ singularity, which of course
is well known since long:

\begin{cor}
The fundamental group of the complement of the discriminant for
$f(x)=x^{n+1}$ is isomorphic to the braid group on
$n+1$ strands given
by the presentation
$$
\langle t_1,...,t_n|t_it_{i+1}t_i=t_{i+1}t_it_{i+1},
t_it_j=t_jt_i\text{ if }|i-j|>1\rangle.
$$
\end{cor}


The results of this chapter should also be regarded as a tool to compute
braid monodromy groups of complicated singularities. In our present set up,
we start with a generic one parameter family of functions induced from the base
of the truncated versal unfolding, compute the local monodromy and its
parallel transport to a common reference fibre.

In the next chapter we prove that we may consider instead special families to
compute the braid monodromy group. We will then have to compute not the
monodromy but the monodromy group of each degeneration -- what we can do now,
if the degeneration is of type $A_n$. Parallel transport has still to be
performed, but in a much simpler family.

In fact it will be families induced from the Hefez-Lazzeri unfolding.


\chapter{results of Zariski type}
\label{Zar}


Having defined the braid monodromy group of a singular function we only
succeeded to compute it for functions of type $A_n$ by means of strong results
cited from the literature.
To proceed we have to develop powerful methods for the computations of braid
monodromy groups.

Generally speaking, in this chapter we will link the braid monodromy of a
singular function to the braid monodromy of adjacent functions. We will
actually determine the braid monodromy of a function from a \mbox{\it tame
$\ell$-perturbation}, a suitably defined unfolding over a two-dimensional base,
using the degeneration properties over the conflict divisor and braid
monodromy of adjacent singularities in this family only.
In fact this method can be applied in such a way that -- in principal -- the
braid monodromy group of a Brieskorn-Pham polynomial can be computed from its
Hefez-Lazzeri unfolding and the monodromy groups of adjacent Brieskorn Pham
polynomials.

The actual execution of this computation and the set up of the necessary
induction are topics of subsequent chapters.

\section{generalization of Morsification}

A Morsification of a singular function $f\in\ofami_n$ is usually defined as a
map representing an unfolding of $f$
\begin{eqnarray*}
\CC^n\times\CC & \tto & \CC\\
x,\l & \mapsto & f_\l(x)
\end{eqnarray*}
such that for generic $\l$ the function $f_\l$ is a Morse function.

Given any versal unfolding of $f$ represented by a map
\begin{eqnarray*}
\CC^n\times\CC^k & \tto & \CC \\
x,u & \mapsto & f(x)-u_0+\sum u_ib_i
\end{eqnarray*}
with $b_i\in\mideal_n$, cf.\ \ref{model},
then a Morsification as above can also be understood as an unfolding
\begin{eqnarray*}
\CC^n\times\CC^2 & \tto & \CC \\
x,\l,u & \mapsto & f_\l(x)-u
\end{eqnarray*}
which is induced by a map $\CC^2\to\CC^k$ such that the restriction to a line
with $\l$ equal to a generic constant maps onto a line transversal to the
discriminant.

In fact we get a pencil of lines in the base of the versal unfolding
parameterized by $\l$, such that all lines sufficiently close to the
line $\l=0$
are transversal to the discriminant.\\

We want to have a notion which generalizes this property to the case of
truncated versal unfoldings.

However we do not have a preferred element like the constant
$1\in \ofami_n$ anymore, but have to choose among the elements of
$\mideal_n$ which are not in the Jacobi ideal $J_f$. In fact we will allow any
choice among the linear functions yielding Morsifications - as given by the
usual existence proof for Morsifications, \cite{ebe}.

As long as this choice is unspecified we denote it by $\ell$, otherwise its
place in the subsequent definitions can be taken by the polynomial actually
chosen.

We consider now two parameter unfoldings of a singular function $f$ in the
maximal ideal $\mideal_n$,
\begin{eqnarray*}
\CC^n\times\CC^2 & \tto & \CC \\
x,\l,u & \mapsto & f_\l(x)+u\ell(x)
\end{eqnarray*}
which are of course induced from any truncated versal deformation of $f$.
So we may also consider the associated pencil of lines parameterized by $\l$ in
the base of the truncated versal unfolding.

\begin{defi}
A two parameter unfolding as above is called \mbox{\it
$\ell$-perturbation}, if all lines of the pencil sufficiently close to the
origin meet the bifurcation set in isolated points only.

By assumption on $\ell$ the line through the origin is not contained in the
bifurcation set.
\end{defi}

\begin{defi}
A two parameter unfolding as above is called \mbox{\it
$\ell$-generification}, if
all lines of the pencil sufficiently close to the origin are transversal to the
bifurcation set. (In particular they meet the bifurcation set in generic
points only corresponding to functions which have non-degenerate
critical points
only with distinct critical values except for a unique critical point of type
$A_2$ or a pair with conflicting values).
\end{defi}

Any $\ell$-generification or $\ell$-perturbation of a function of Milnor number
$\mu$ determines a
$\CC_\mu$-bundle over a multiply punctured disc, well defined up to fibration
isomorphism, since the lines of the pencil close to the origin are
transverse to
the bifurcation set in a uniform way.\\

As a Morsification may serve for the computation of monodromy groups so a
generification of a function $f$ can replace its truncated versal unfolding:

\begin{lemma}
\labell{generif}
The braid monodromy group of a $\ell$-generification is equal to the braid
monodromy group.
\end{lemma}

\proof
It just suffices to point to the analogous argument in the case of a
Morsification. The important point to note is, that the lines of the pencil are
generic with respect to the bifurcation set and therefore the induced map on
fundamental groups surjects.
\qed[3mm]

\section{versal braid monodromy group}
\label{versal-b}

Given a one-parameter family of monic polynomials we have formerly
divided the computation of the braid monodromy group into two steps.
First we assign the local monodromy generator to a local Milnor fibre, second
we use parallel transport to get mapping classes in just one regular fibre.
Upon the choice of a geometrically distinguished system of paths in that fibre,
the subgroup generated by the transported classes is identified with the braid
monodromy group.
\\

In a similar approach, we will assign a group instead of a generator to local
Milnor fibres close to each singular fibre, and we will then use parallel
transport of the group elements to get mapping classes again in a single
regular fibre.

It is the local assignment which we have to define carefully to get a sensible
additional notion of braid monodromy.

In fact it will only be defined for one parameter families of monic polynomials
associated to a family of functions on which -- for technical
simplicity only --
we impose the further restriction of {\it tameness}.

\begin{defi}
A one parameter family of functions is called {\it tame} if of each function
only non-degenerate critical points may have conflicting values.
\end{defi}

\begin{figure}[h]
\setlength{\unitlength}{2mm}
\begin{picture}(10,35)(-25,0)

\put(5,10){\begin{picture}(5,2)

\bezier{300}(0,5)(0,7.5)(2,9)
\bezier{300}(2,9)(5,11)(8,9)
\bezier{300}(8,9)(10,7.5)(10,5)
\bezier{300}(10,5)(10,2.5)(8,1)
\bezier{300}(8,1)(5,-1)(2,1)
\bezier{300}(2,1)(0,2.5)(0,5)

\end{picture}}

\put(5,22){\begin{picture}(5,0)

\bezier{300}(0,5)(0,7.5)(2,9)
\bezier{300}(2,9)(5,11)(8,9)
\bezier{300}(8,9)(10,7.5)(10,5)
\bezier{300}(10,5)(10,2.5)(8,1)
\bezier{300}(8,1)(5,-1)(2,1)
\bezier{300}(2,1)(0,2.5)(0,5)

\end{picture}}

\put(5,2){\setlength{\unitlength}{1mm}
\begin{picture}(5,0)(-3.75,0)

\bezier{300}(0,2.5)(0,3.75)(2,4.5)
\bezier{300}(2,4.5)(5,5.5)(8,4.5)
\bezier{300}(8,4.5)(10,3.75)(10,2.5)
\bezier{300}(10,2.5)(10,1.25)(8,.5)
\bezier{300}(8,.5)(5,-.5)(2,.5)
\bezier{300}(2,.5)(0,1.25)(0,2.5)

\put(5,2.5){\circle*{1}}

\end{picture}}

\drawline(10,8)(10,35)

\bezier{30}(7.5,7)(7.5,21.5)(7.5,36)
\bezier{30}(12.5,7)(12.5,21.5)(12.5,36)


\bezier{200}(10,27)(5,27)(1,30)
\bezier{200}(10,27)(5,27)(1,24)
\bezier{20}(10,27)(15,27)(19,30)
\bezier{20}(10,27)(15,27)(19,24)


\drawline(19,12)(1,18)
\drawline(19,18)(1,12)

\end{picture}
\caption{Milnor fibration}
\label{mil}
\end{figure}
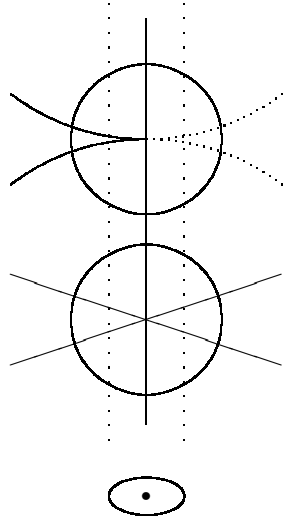

Given a tame family of functions then locally at a singular function $f$ the
associated family of monic polynomials $p_\l$ is parameterized by $\l$ in a
disc such that the coefficient map is holomorphic and such that the polynomial
is a simple Weierstrass polynomial for $\l\neq0$:
$$
p_\l:(\l,x)\mapsto x^n+\sum_{i=0}^{n-1}a_i(\l)x^i.
$$
Suppose now $p_0$ not to be simple with roots denoted by $v_j$. Then for $\e$
and $\d=\d(\e)$ sufficiently small, the local family
$$
Y:=\CC\times B_\d\smin p_\l\inv(0)
$$
is trivializable over $B_\d$ in the complement of $\cup_j B_\e(v_j)$, cf.\
fig.\ \ref{mil}.

We conclude that all mapping classes of the local family $Y$ can be given with
support on the intersection $\cup_j D_j$ of a local Milnor fibre with $\cup_j
B_\e(v_j)$.
Moreover we notice that the restriction to a disc $D_j$ yields the braid
monodromy transformations of the local discriminant divisor at $v_j$.\\

We finally assign a group of mapping classes to $p_\l$ by a choice of a group
of mapping classes supported on $D_j$ for each root $v_j$ of $p_0$.

Consider first the case that $v_j$ is a multiple root of $p_0$, which is the
image of a single critical point $c_j$ of $f$.
Then for the germ of $f$ at $c_j$ braid monodromy yields a well
defined group of
mapping classes supported on $D_j$.

In case $v_j$ is the image of several non-degenerate critical points of $f$, we
simply choose the group of mapping classes of $D_j$ which fix the punctures and
which therefore correspond to pure braids.

By tameness these are the only possible cases. The generalization to other
families is conceptionally straight forward, but notationally a mess, so we
decided to skip it here.

Anyway we should sum up our definition:

\begin{defi}
Given a tame family of functions we assign a group of mapping classes to a
Milnor fibre of each singular fibre of the associated family. The generators
are given by mapping classes in the Milnor fibres $D_j$ of the multiple roots
$v_j$:
\begin{enumerate}
\item
all mapping classes which fix the punctures in case $v_j$ is the image of
non-degenerate critical points only,
\item
the mapping classes provided by braid monodromy of a critical point $c_j$ in
case $c_j$ is the only critical point which maps to $v_j$.
\end{enumerate}
The group of mapping classes in a regular fibre obtained by parallel transport
is called the {\it versal braid monodromy group} of the tame family
of functions.
\end{defi}

Of course this determines a well defined conjugacy class of subgroups of the
braid group $\br_n$ upon the choice of a geometrically distinguished system of
paths in a regular fibre.\\

\begin{remark}
\label{perturb}
We note that a given $\ell$-perturbation of a singular function determines a
one-parameter family of functions up to smooth fibration isomorphism, hence in
case the family is tame we may also speak of the {\it versal braid
monodromy} of
the $\ell$-perturbation.

Moreover we may extend the notion of {\it versal braid monodromy} to the case
where several tame families are given which have a regular fibre in
common: Then
it denotes the group of mapping classes of the common fibre which is
generated by
the subgroups which are the versal braid monodromies of the separate families.
\end{remark}

By the definition we have to consider all possible transports to the reference
fibre, but in fact we can restrict the computation to the transport along the
paths of a distinguished system.

\begin{lemma}
\labell{vers-geom}
The versal braid monodromy of a family is obtained if the locally assigned
groups are transported along the paths of a geometrically distinguished system.
\end{lemma}

\proof
The key observation is, that the local braid monodromy transformation belongs
to the locally assigned group. We can thus conclude as in the classical case.
\qed[0mm]

\section{comparison of braid monodromies}

In this section we relate the versal braid monodromy of a tame
$\ell$-perturbation to the braid monodromy of a versal family, which will
eventually justify the name.

\begin{prop}
\labell{vers-comp}
The braid monodromy group of a function $f$ is equal to the versal braid
monodromy group of any of its $\ell$-perturbation which is tame.
\end{prop}

\proof
Suppose the $\ell$-perturbation is represented by a map
\begin{eqnarray*}
\CC^n\times\CC^2 & \tto & \CC\\
x,\l,u & \mapsto & f_\l(x)+u\ell(x).
\end{eqnarray*}
We extend it to an unfolding with base of dimension $3$
\begin{eqnarray*}
\CC^n\times\CC^3 & \tto & \CC\\
x,\l,u_1,u & \mapsto & f_\l(x)+u_1 b_1(x)+u\ell(x),
\end{eqnarray*}
such that the restriction to $\l=0$ is a $\ell$-generification.
We can check that the bifurcation set is of codimension one and its singular
locus of codimension two in this base.
This is immediate from our assumptions, since a non empty part of a
generification is induced from the complement of the bifurcation set
and another
non empty part from the complement of its singular locus.

The conclusion still holds for the dominantly induced unfolding
\begin{eqnarray*}
\CC^n\times\CC^3 & \tto & \CC\\
x,\l,\l_1,u & \mapsto & f_\l(x)+\l\l_1 b_1(x)+u\ell(x),
\end{eqnarray*}
Restricted to $\l_1=0$ this is the $\ell$-perturbation we started with.
Restricted to some fixed $\l_1\neq0$ and sufficiently small it is still a
$\ell$-perturbation of $f$, but its base only meets the singular locus of the
bifurcation set in isolated points. Hence we get even a
$\ell$-generification of $f$.\\

For the rest of the proof we fix $\l$ at a sufficiently small non-zero value
and consider various restrictions of the family
\begin{eqnarray*}
F:\CC^n\times\CC^2 & \tto & \CC\\
x,\l_1,u & \mapsto & f_\l(x)+\l\l_1 b_1(x)+u\ell(x).
\end{eqnarray*}
For example $F|_{\l_1=0}$ is the tame family of which we want to understand the
versal braid monodromy and for sufficiently small $\eta$ we get a tame family
$F|_{\l_1=\eta}$ which has braid monodromy equal to that of $f$ by lemma
\ref{generif}.\\

For each critical parameter $y_i$ in the bifurcation set on the line $\l_1=0$
we choose a local ball $U_i$ in the base of $F$ centered at $y_i$. We
fix $\eta$
sufficiently small and the tubular neighbourhood $N_\eta$ of the line $\l_1=0$
of radius
$\eta$, such that the bifurcation set of $F|_{N_\eta}$ is contained
in the union
of the
$U_i$ and its singular locus is a subset of the $y_i$.

The braid monodromy of $F|_{N_\eta}$ is then equal to the braid
monodromy of $f$,
since it contains the family $F|_{\l_1=\eta}$.
On the other hand it is generated by the braid monodromies of the
$F|_{U_i}$ and parallel transport over the complement of the $U_i$.

This should be compared to the fact that the versal braid monodromy of
$F|_{\l_1=0}$ is generated by the versal braid monodromies of $F|_{E_i}$ --
where $E_i$ denotes the intersection of $\l_1=0$ with $U_i$ -- and parallel
transport over the complement of the $E_i$.
\\

We have therefore accomplished a major reduction step in the proof:
\begin{quote}
It suffices to prove that the versal braid monodromy of
$F|_{E_i}$ is equal to the braid monodromy of $F|_{U_i}$ for each $i$,
since the complement of the $U_i$ in $N_\eta$ retracts onto the
complement of the
$E_i$ on $\l_1=0$.
\end{quote}

We move thus our attention to the discriminant family of $F$ restricted to a
single ball $U$ of the base. The restriction to $E$ yields a
discriminant family with a single singular fibre for which we gave a local
description in section \ref{versal-b} already.

In fact this description extends to $U$ if it is chosen appropriately: The
complement
$Y$ of the discriminant in $\CC\times U$ is trivializable over $U$ in the
complement of balls $B_\e(v_j)$ centered at the roots $v_j$ on the fibre over
$y$.

The braid monodromy of $F|_U$ and the versal braid monodromy of $F|_E$ can thus
be considered as a group of mapping classes which are supported on the
intersection $\cup_j D_j$ of a local Milnor fibre with $\cup_j B_\e(v_j)$.\\

According to the decomposition of the discriminant into connected components
$\dfami_j$ over $U$, the bifurcation divisor decomposes as
$\bfami=\cup_j\bfami_j$, such that each divisor $\bfami_j$ is the branch locus
of the finite map of $\dfami_j$ onto $U$.

As the $B_\e(v_j)$ are disjoint the braid monodromy transformations along
simple geometric elements associated to different parts of $\bfami$ commute.

In particular the braid monodromy transformation along a simple geometric
element based at the chosen Milnor fibre and associated to $\bfami_j$ can be
chosen with support in $D_j$.
\pagebreak

Now let us consider first a multiple root $v_j$ of the critical fibre at $y$
which is the image of several non-degenerate critical points of the
corresponding
function. Then the corresponding local discriminant divisor $\dfami_j$ in
$B_\e(v_j)$ has irreducible components in bijection to the preimages.
Hence all mapping classes in the braid monodromy of $F|_U$ restrict to
mapping classes of $D_j$ which fix the punctures pointwise.

On the other hand $E':=U\cap\{\l_1=\eta\}$ is transversal to the
bifurcation set, so the divisorial discriminant components in $B_\e(v_j)$ meet
pairwise, transversally, and over distinct points of the bifurcation set
$\bfami_j\cap E'$. This implies that the braid monodromy of
$F|_{E'}$ contains all pure mapping classes of $D_j$, i.e.\ the group of
mapping classes which are supported on $D_j$ and fix the punctures pointwise.
Hence this braid monodromy contains all mapping classes we assign to
$v_j$ to get
the versal braid monodromy group of $F|_E$
\\
\begin{figure}[h]
\setlength{\unitlength}{2mm}
\begin{picture}(50,29)(-10,3)


\drawline(10,25)(10,9)
\drawline(26,33)(10,25)
\drawline(26,17)(10,9)
\drawline(26,33)(26,17)

\bezier{120}(12,16)(18,21)(24,26)
\bezier{120}(12,23)(18,21)(24,19)
\bezier{120}(12,12)(15,18.5)(18,21)
\bezier{120}(24,22)(21,23.5)(18,21)

\drawline(30,25)(30,9)
\drawline(46,33)(30,25)
\drawline(46,17)(30,9)
\drawline(46,33)(46,17)

\bezier{120}(32,16)(38,21)(44,26)
\bezier{120}(32,23)(38,21)(44,19)
\bezier{190}(44,24)(36,26)(32,12)

\end{picture}
\caption{generification of smooth branches}
\label{genfic-red}
\end{figure}
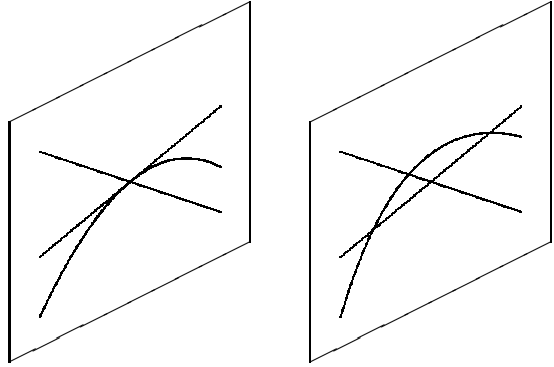

Similarly we argue in case the multiple root $v_j$ is the image of a unique
critical point $c_j$.
Then $B_\e(v_j)$ can be considered as a discriminant family induced from the
base of a versal truncated unfolding of the function at $c_j$.

It is in fact a generification, since its bifurcation set $\bfami_j$ is met by
$E$ in a single point only and transversally by $E'$.

Hence the braid monodromy of $F|_U$ contains the braid monodromy of the
function at $c_j$ considered as mapping classes on $D_j$ extended by the
identity to the Milnor fibre of $F|_U$, which is just what we assigned to
$v_j$ to get the versal braid monodromy group.
\\
\setlength{\unitlength}{2mm}
\begin{picture}(50,35)(-10,0)


\drawline(10,25)(10,9)
\drawline(26,33)(10,25)
\drawline(26,17)(10,9)
\drawline(26,33)(26,17)

\bezier{120}(12,23)(16,21)(18,21)
\bezier{120}(18,21)(20,23)(24,30)
\bezier{16}(18,21)(15,18)(12,12)
\bezier{10}(18,21)(14,19)(12,19)
\bezier{10}(24,17)(22,21)(18,21)
\bezier{10}(24,21)(22,21)(18,21)

\drawline(30,25)(30,9)
\drawline(46,33)(30,25)
\drawline(46,17)(30,9)
\drawline(46,33)(46,17)

\bezier{120}(32,23)(36,21)(40,20)
\bezier{120}(36,19)(40,23)(44,30)
\bezier{70}(36,19)(37.5,20.5)(40,20)
\bezier{16}(36,19)(34,16)(32,12)
\bezier{10}(36,19)(34,18)(32,19)
\bezier{10}(44,17)(42,19)(40,20)
\bezier{10}(44,21)(42,19)(40,20)


\end{picture}
\begin{figure}[h]
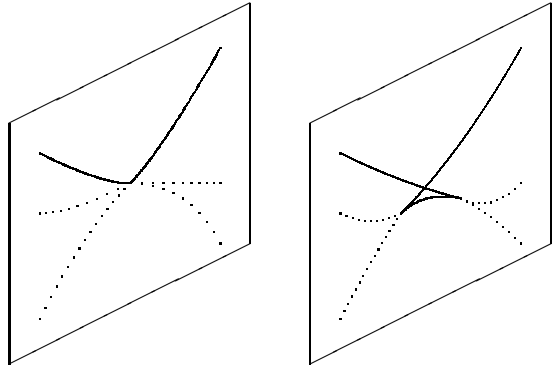

\caption{generification of irreducible branch}
\label{genfic-irr}
\end{figure}

So we have shown, that the versal braid monodromy of $F|_E$ is contained in the
braid monodromy of $F|_U$.
\\

For the reverse inclusion it suffices to argue that all braid monodromy
transformations along simple geometric elements belong to the versal braid
monodromy of $F|_E$. Suppose the element is associated to a component of
$\bfami_j$, then the corresponding monodromy transformation is supported on the
intersection with $B_\e(v_j)$, and therefore it is equal to its restriction to
$D_j$ extended by the identity.

If $v_j$ is of the first kind considered above, then we noted that the
restriction to $D_j$ of any monodromy transformation belongs to the pure
mapping classes of $D_j$. Their extension by the identity are thus elements of
the versal braid monodromy group of $F|_E$ assigned to $v_j$.

If $v_j$ is of the second kind, then the restriction to $D_j$ must be an
element of the braid monodromy of the discriminant family given by the
restriction to $B_\e(v_j)$. Again the extension by the identity is an element
of the versal braid monodromy group of $F|_E$ assigned to $v_j$.
\qed

This result can be generalized to arbitrary perturbations, but we need
only this form and have therefore preferred to avoid the bulk of technicalities
involved in the general case.

\section{Hefez-Lazzeri base}

In the case of a Brieskorn-Pham polynomial we want finally reduce the
computation of the braid monodromy group to the computation of versal braid
monodromy groups of families which are monomial perturbations induced from the
Hefez-Lazzeri unfolding.

\begin{eqnarray*}
\CC^n\times\CC^n & \tto & \CC\\
x,u & \mapsto & \sum_i(x_i^{l_i+1}+u_ix_i)
\end{eqnarray*}

\begin{nota}
We introduce the following shorthand notation for families parameterized by
$\a$ with fixed real constants $\e_i>0$:\\
By $f_\a(x_1,x_2)$ denote the $\a$-family
$$
x_1^{l_1+1}-\a(l_1+1)x_1+x_2^{l_2+1}-\e_2(l_2+1)x_2.
$$
By $g_\a(x_1,x_2)$ denote the $\a$-family
$$
x_1^{l_1+1}-(l_1+1)x_1+x_2^{l_2+1}-\a\e_2(l_2+1)x_2.
$$
By $f_\a(x_1,...,x_n)$ denote the $\a$-family
$$
x_1^{l_1+1}-\a(l_1+1)x_1+x_2^{l_2+1}-\e_2(l_2+1)x_2+\cdots
+x_n^{l_n+1}-\e_n(l_n+1)x_n.
$$
By $g_\a(x_1,...,x_n)$ denote the $\a$-family
$$
x_1^{l_1+1}-(l_1+1)x_1+x_2^{l_2+1}-\a\e_2(l_2+1)x_2+\cdots
+x_n^{l_n+1}-\a\e_n(l_n+1)x_n.
$$
\end{nota}

As a first step we have to show that at least for suitable choices of constants
the versal braid monodromy is defined for the given families.

\begin{lemma}
\labell{gtame}
The families $g_\a(x_1,...,x_n)$ are tame.
\end{lemma}

\proof
A degenerate critical point may only occur if the Hessian determinant
vanishes. But for any family induced from the Hefez Lazzeri unfolding this
determinant is constantly equal to
$$
x_1^{l_1-1}\cdot...\cdot x_n^{l_n-1}\prod_il_i(l_i+1).
$$
Hence at least one coordinate $x_i$ of a degenerate critical point must vanish.
The vanishing of the gradient then implies via
$$
\del_if(x,u)=(l_i+1)x_i^{l_i}+u_i
$$
that the corresponding parameter vanishes too.
So in the case of the induced family $g_\a$ we conclude that $\a=0$ if a
degenerate critical point occurs.

In that case the critical points are determined by
$$
x_i^{l_i}=1 \wedge x_2=...=x_n=0
$$
and we get a bijection between the set of critical points and their values:
$$
\{(\xi_1^k,0,...,0)\}  \stackrel{1:1}{\longleftrightarrow}
\{-l_1\x_1^k\},\qquad \xi_1^{l_1}=1.
$$
Hence the families $g_\a$ are indeed tame.
\qed

\begin{lemma}
\labell{ftame}
The family $f_\a(x_1,...,x_n)$ is tame if the function $f_0|_{x_1=0}$
is a Morse function.
\end{lemma}

\proof
By the preceding proof we need to worry only about the function for $\a=0$,
since otherwise all critical points are non-degenerate.

On the other hand any critical point of $f_0$ is situated on the hyperplane
$x_1=0$ due to
$$
\del_1f_0(x)=(l_1+1)x_1^{l_1}=0.
$$
So each must also be a critical point of $f_0|_{x_1=0}$.
If now $f_0|_{x_1=0}$ is a Morse function then it maps the set of critical
points bijectively onto the set of critical values, hence no pair of critical
points of $f_0$ may map to the same critical value.
\qed

Indeed we can deduce from the cases considered in the preceding proofs the
following criterion for the tameness of a family induced from the Hefez Lazzeri
base.

\begin{lemma}
\labell{htame}
A one-parameter family of functions induced from the Hefez Lazzeri
base is tame,
if any function is one of the following list:
\begin{enumerate}
\item
a function induced from the complement of the coordinate hyperplanes,
\item
a function induced from a coordinate axis,
\item
a function induced from a point on just one coordinate hyperplane $\a_i=0$,
such that its restriction to $x_i=0$ is a Morse function.
\end{enumerate}
\end{lemma}

We call a positive real constant $\e_2$ resp.\ a tuple  $\e_2,...,\e_n$ of
positive real constants admissible, if the fibre corresponding to $g_1$ in the
discriminant family associated to $g_\a$ is regular.
The condition is met if and only if $g_1$ is a Morse function.

\begin{lemma}
\labell{admiss}
If positive real constants $\e_2,...,\e_n$ are chosen generically, then
they are admissible and $g_1(x_1,...,x_n)$ is a Morse function.
\end{lemma}

\proof
There is a Zariski open set of complex constants, such that $g_1$ is a Morse
function. Furthermore
we note that $g_1$ is not a Morse function for all choices of complex constants
if and only if the defining polynomial of the bifurcation divisor given in
\ref{laz-bif} vanishes on the hyperplane $\a_1=1$.
Since that is not the case, the Zariski set above is non-empty.

But then it must contain a dense subset of all tuples of positive real
constants, too.
\qed

We can now consider the case $n=2$. Then $f_0|_{x_1=0}$ is a Morse function,
thus for an admissible choice of $\e_2>0$ the versal braid monodromy of the
families $f_\a$ and $g_\a$ is defined and we can establish the
following result:

\begin{prop}
\labell{ZA}
For an admissible choice $\e_2>0$
the braid monodromy of the function $f=x_1^{l_1+1}+x_2^{l_2+1}$ is given by the
versal braid monodromy of the families $f_{\a,\e}(x_1,x_2)$ and
$g_\a(x_1,x_2)$, where the parameter of the second family may be restricted to
the unit disc.
\end{prop}

\proof
By the result of the last section it suffices to show that the versal
braid monodromy of the two families is equal to the versal braid monodromy of
some $\ell$-perturbation of the function $f$.

So let us first pick an $\ell$-perturbation which we want to compare to the
families $f_\a$ and $g_\a$:

Since the bifurcation set in the Hefez Lazzeri base is a divisor, there is a
transverse line and the corresponding linear polynomial $\ell$ may even be
assumed to be different from $x_1$ by genericity.
Hence the $\ell$-perturbation
\begin{eqnarray*}
x,\a,\l & \mapsto & f(x)+\a x_1+\l\ell
\end{eqnarray*}
is just induced from the Hefez Lazzeri unfolding by a change of parameters.
Denote by $h_\l$ the family of functions which is induced from an
affine line $L_h$ in the Hefez Lazzeri base parallel and sufficiently close
to the line given by $\ell$. Moreover from the first two criteria of lemma
\ref{htame} it is immediate that this family is tame.

Note that incidentally the families $f_1$ and $g_1$ are induced from affine
lines $L_f$ and $L_g$ in the Hefez Lazzeri base parallel to the $u_2$-axis and
the $u_1$-axis respectively.\\

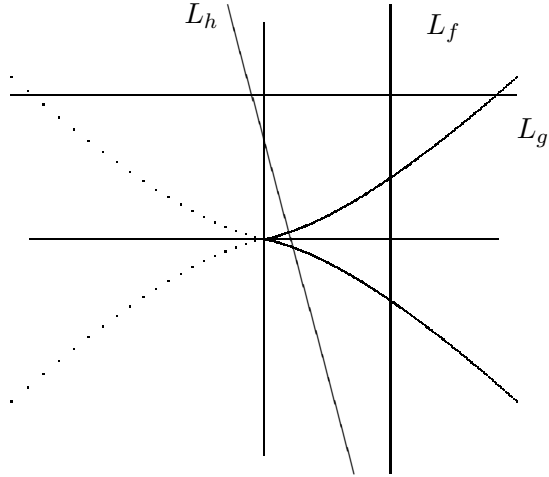
\begin{figure}[h]
\setlength{\unitlength}{2.4mm}
\begin{picture}(10,25)(-20,0)

\drawline(10,0)(10,24)
\drawline(-3,12)(23,12)


\bezier{20}(10,12)(5,13)(-4,21)
\bezier{20}(10,12)(5,11)(-4,3)
\bezier{200}(10,12)(15,13)(24,21)
\bezier{200}(10,12)(15,11)(24,3)

\drawline(8,25)(15,-1)
\put(5.6,24){$L_h$}
\drawline(17,25)(17,-1)
\put(19,23.4){$L_f$}
\drawline(-4,20)(24,20)
\put(24,17.6){$L_g$}

\end{picture}
\caption{Hefez Lazzeri base}
\label{Hef}
\end{figure}

Next we make some observations regarding the fundamental group of the
bifurcation complement:

Let $L'$ denote the intersection of a line $L$ with the complement of $\bfami$.
We first note that $\bfami$ is a curve without vertical components with respect
to the linear projection parallel to $L_h$. So as we remarked before, lemma
\ref{pres} and its proof imply
$$
\pi_1(L_h')\surj\pi_1(\CC^2\smin\bfami).
$$
Moreover we extended this result to handle the projection parallel to $L_f$.
Then the $u_1$-axis is the unique vertical component of $\bfami$ and by lemma
\ref{pres-vert} we get
$$
\pi_1(L_f')*\ZZ\surj\pi_1(\CC^2\smin\bfami),
$$
where the free generator is given by a geometric element supported on $L_g$ and
associated to $\a_1=0$. Therefore the braid monodromies of the families $h_\l$
resp.\ $f_\a$ and $g_\a|$, the restriction of $g_\a$ to any base
containing the unit disc, are both given by the braid monodromy of the Hefez
Lazzeri unfolding.\\

A further ingredient of the argument is the geometry of the bifurcation divisor
with respect to the lines $L_f$ and $L_g$:

The Hefez Lazzeri base is of dimension two and the bifurcation divisor is
quasi homogeneous with respect to a good $\CC^*$-action.
The line $L_f$ has a common point with every $\CC^*$ orbit except the
$u_2$-axis, so we may conclude that each component of the bifurcation divisor,
which obviously is the closure of a $\CC^*$-orbit, meets $L_f$ in at least one
point or $L_g$ at its intersection with $u_1=0$.\\

Finally we bring back our attention to the versal braid monodromies associated
to the families $h_\l$ resp.\ $f_\a$ and $g_\a|$:

To compare them we chose a path $p_1$ in the Hefez Lazzeri base which connects
the base points of the respective reference fibres and which is disjoint from
the bifurcation divisor.

As we noted in the proof of lemma \ref{vers-geom}
the versal braid monodromy of a tame family contains its braid monodromy.
Hence it suffices by the surjectivity result above to show that every
generator of the versal braid monodromy on one hand is equal to a generator on
the other hand transported along $p_1$ up to the braid monodromy of the
families.

The set of generators of the versal braid monodromy of the family $h_\l$ can be
chosen among parallel transports of generators in the local groups for $h_\l$.
Suppose $\b$ is such a generator associated to a degeneration of $h_\l$, i.e.\
to a multiple root $v_j$ in a singular fibre of the corresponding discriminant
family. Let $\bfami_0$ denote the component of the bifurcation divisor in the
Hefez Lazzeri base, from which this singular fibre is induced.

Then the same component meets the base of $f_\a$ or $g_\a|$ as we noted above.
So we may find a path $q_2$ in the smooth locus of $\bfami_0$ which connects
the singular fibres.
By equisingularity along the smooth part of $\bfami_0$
or by the explicit equations we get a Whitney stratification for the
discriminant family over a neighbourhood of $q_2$,
where the components of the smooth locus of the discriminant are the strata of
codimension one, and the components of its singular locus are the strata of
codimension two.

A lift $\hat q_2$ of $q_2$ to the codimension two stratum which contains
$v_j$ ends at some multiple root $v_j$ in the singular fibre at the end of
$q_2$.
Topological triviality along $\hat q_2$ implies that the local
families over $L_f$ or $L_g$ and $L_h$ restricted to
neighbourhoods of $v_j$ resp.\ $v_j'$ are topologically identified under
parallel transport. Hence $\b$ transported along a path $p_2$ in the
complement of the bifurcation divisor but sufficiently close to $q_2$ yields a
mapping class
$\b'$ which belongs to the local group of $f_\a$ or $g_\a|$.\\

Therefore the parallel transport of $\b$ along any path $p_h$ in $L_h'$ and the
parallel transport of $\b'$ along a path $p_{fg}$ in $L_h'\cup L_g'$ yield
mapping classes which are equal up to parallel transport along the closed path
obtained as the union of $p_1,p_2,p_h$ and $p_{fg}$.\\

But then by the surjectivity on fundamental groups the same class is
represented by a path in the base of $f_\a$ and $g_\a$, hence any parallel
transport of $\b$ in the base of $h_\l$ is equal to some parallel transport of
$\b'$ in the base of $f_\a$ and $g_\a$.\\

The argument can also be read with $h_\l$ taking the role of $f_\a$ and $g_\a|$
and vice versa, and then yields the reverse implication.
\qed

We have to extend the result to the higher dimensional case. In principal the
argument is the same but now some particular steps become more involved since
the components of the bifurcation divisor are no longer smooth outside the
origin. So to work on the smooth locus we have to add some more genericity
assumptions

\begin{lemma}
\labell{fgeneric}
For a generic choice of admissible constants $\e_2,...,\e_n$, the base of
$f_\a$ meets the reduced bifurcation divisor in regular points only.
\end{lemma}

\proof
The claim can be deduced by a dimension count. The set of singular points of
the reduced bifurcation divisor is of codimension two in the Hefez
Lazzeri base.

On the other hand the complex constants $\e_2,...,\e_n$ parameterize the
parallels to the $u_1$-axis in the base.
A generic line of this family does not meet the singular set above, so there is
a non-empty Zariski open set of constants, such that the corresponding line is
disjoint from the singular set.

But then even most tuples of positive real constants must belong to this set.
\qed

\begin{prop}
\labell{ZB}
For a generic choice of admissible constants $\e_2,...,\e_n>0$ the braid
monodromy of the function $x_1^{l_1+1}+\cdots+x_n^{l_n+1}$ is given by the
versal braid monodromy of the families $f_{\a}(x_1,...,x_n)$ and
$g_\a(x_1,...,x_n)$ where the parameter of the second family may be restricted
to the unit disc.
\end{prop}

\proof
Due to lemma \ref{vers-comp} again it suffices to show that the versal braid
monodromy of the families $f_\a$ and $g_\a|$ is equal to the versal braid
monodromy of some tame $\ell$-perturbation of the function
$f=x_1^{l_1+1}+\cdots+x_n^{l_n+1}$.

We assume that admissible constants $\e_i>0$ are chosen in such way, that the
base $L_f$ of $f_\a$ meets the reduced bifurcation divisor only in regular
points, cf.\ lemma \ref{fgeneric}.
By assumption $f_0|_{x_1=0}=g_1|_{x_1=0}$ is a Morse function. But then so is
$g_\a|_{x_1=0}$ for almost all values of the parameter $\a$.\\

We define $\ell=cx_1+\e_2(l_2+1)x_2+...+\e_n(l_n+1)x_n$ for a generic constant
$c\neq0$ and an unfolding with two-dimensional base $E$
\begin{eqnarray*}
x,u,\l & \mapsto & f(x)+ ux_1+\l\ell.
\end{eqnarray*}
Then let $h_\l$ be the family induced for $u=\e_1$ constant and sufficiently
close to 0.
Furthermore note that $f_\a$ is induced from $E$ for $\l=-1$ and
$u=-\a(l_1+1)-c$ while $g_\a$ is induced for $\l=-\a,u=\a c$.

Since the only functions in $h_\l$ induced from the coordinate hyperplanes of
the Hefez Lazzeri base are $h_0$ and $h_{-\e_1/c}$ we conclude that the family
$h_\l$ is tame by \ref{htame}:
\begin{enumerate}
\item
$E$ considered as a plane in the Hefez Lazzeri base is not contained in the
coordinate hyperplanes, hence each function $h_\l$, $\l\neq0,-\e_1/c$ is of the
first kind considered in \ref{htame},
\item
$h_0(x)=f(x)+\e_1x_1$ is induced from a coordinate axis, so it is of the second
kind considered in \ref{htame},
\item
$h_{-\e_1/c}|_{x_1=0}=g_{-\e_1/c}|_{x_1=0}$ is a Morse function for $\e_1$
sufficiently small, so $h_{-\e_1/c}$ is of the third kind considered in
\ref{htame}.
\end{enumerate}

Next we make some observations on the bifurcation locus restricted to $E$
considered as a plane in the Hefez Lazzeri base.

By assumption $g_1$ is a Morse function, hence the bifurcation locus is of
codimension one in $E$.
The restriction of the Hefez Lazzeri base to $E$ is obtained by imposing some
linear relations on the base coordinates. From the defining equation
\ref{laz-bif} of the bifurcation divisor $\bfami$ we can read off the fact, that
each component of the restriction $\bfami_E=\bfami\cap E$ contains the origin.
Therefore each component different from the $u_1$-axis meets $L_f$.

This situation is very similar to that encountered in the case $n=2$. Since
$L_h$ is a general line in $E$ with respect to $\bfami_E$ and since by
assumption $L_f$ meets
$\bfami_E$ in its regular locus only, we may conclude as before:
\begin{eqnarray*}
\pi_1(L_h') & \surj & \pi_1(E\smin\bfami_E),\\
\pi_1(L_f')*\ZZ & \surj & \pi_1(E\smin\bfami_E),
\end{eqnarray*}
where the free generator maps to any geometric element supported on $L_g$
associated to the $u_1$-axis.

Since $\pi_1(E\smin\bfami_E) \cong \pi_1(\CC^n\smin\bfami)$, the braid
monodromies of the families $h_\l$ resp.\ $f_\a$ and $g_\a|$ are both given by
the braid monodromy of the Hefez Lazzeri unfolding.\\

The argument then proceeds as in the $n=2$ case.
\qed


\chapter{braid monodromy of plane curve families}
\label{Plane}



This chapter can be considered the central one. Admittedly it is devoted to a
partial computation only of the braid monodromy group of a plane Brieskorn Pham
polynomial. Nevertheless it will be the essential ingredient for the
computation of braid monodromy groups associated to arbitrary Brieskorn Pham
polynomials which will be the topic of the next chapter.

Given a plane Brieskorn Pham polynomial $x_1^{l_1+1}+x_2^{l_2+1}$, we
consider for generic small $\e_2>0$ the one parameter family of functions
defined in the last chapter:
$$
f_\a(x_1,x_2)=x_1^{l_1+1}-\a(l_1+1)x_1+x_2^{l_2+1}-\e_2(l_2+1)x_2.
$$
It corresponds to a line in the Hefez Lazzeri base which in turn can be
considered as a plane in the base of a suitable truncated miniversal unfolding.
We have already obtained the necessary
formulae for the singular value divisor in the third chapter which
are of course
needed for further calculations.

Using the preceding chapter we will actually compute the versal braid monodromy
of the family $f_\a$.
So remember that we have two tasks. For each singular function we have to
determine the local group which is a group of mapping classes in a regular
fibre close to the corresponding fibre in the discriminant family.
We have then to find the group of mapping classes in a distinguished reference
fibre generated by translates of the elements of the local groups obtained by
parallel transport.\\

Let us muse a moment on what we have to do: We are given a one
dimensional base of
a family with a lot of
degeneration points. Close to these points we may choose local Milnor fibres
and are then given
mapping classes in these Milnor fibres. Next we have to choose a
geometrically distinguished system
of paths joining them to a global reference fibre.
Finally we transport the local mapping classes to the reference
fibre and determine the
subgroup of braids they generate under the identification of the mapping
class group with the braid
group given by the Hefez Lazzeri system of paths in the fibre.

All steps will be addressed in this chapter, though not in the given order.
Moreover we need to modify this approach, so let us point out some
details:\\

Instead of parallel transport in the discriminant family we will consider
parallel transport in a
closely related family from the first section on. Only in section \ref{six/eit}
we will return to the proper
discriminant family and exploit the relation of both families to transfer
the result to the family
where we actually need them.\\

The computation of the local monodromy groups is postponed to section
\ref{six/sevn}, so we can concentrate on
the parallel transport in the first sections.
Note that we won't transport general classes but exploit the bijection between
twist classes and isotopy classes of arcs which join the two punctures twisted
around each other. In fact it suffices to consider the parallel transport
of a specific sort of arcs only since we will prove in section \ref{six/sevn},
that all local mapping classes
form a group generated by twists on such arcs.\\

The most difficult task is the parallel transport of mapping classes. As
said we will in fact
transport the arcs such that the corresponding twists generate the mapping
class groups.

Since all the paths we use decompose into radial and circular segments, we
study the differential flow over such segments and determine the parallel
transport of arcs along radial segments in the first section.
It will be possible in subsequent sections
to determine the parallel transport of arcs, since we arrange our arcs to be
determined by geometric data which are preserved and some geometric datum which
changes in a way we can actually measure.

\section{parallel transport in the model family}
\label{six/one}

In this section we consider the parallel transport in a punctured disc
bundle associated to the
discriminant bundle of the function families $f_\a$. Recall that we have to
transport the local
groups to a global reference fibre along a geometrically distinguished
system of paths. Such a
system can be chosen to consist of paths which are obtained by unions of
radial segments with
circular segments on circles centred at the origin.

So we find suitable vector fields and differentiable flows along such
segments. Thereby we get
representatives of mapping classes which in turn will be used to get
representatives of parallel
transport of braid mapping classes.
\\

Let us consider then the punctured disc bundle
associated to the line arrangement
$$
\prod_{\xi_1^{l_1},\xi_2^{l_2}=1}
(z-\l\xi_1-\eta_2\xi_2)\quad = \quad 0.
$$
which we call the {\it model discriminant family} associated to $l_1,l_2$,
$\eta_2\ll1$.

In section \ref{hl-path} we have indexed the paths of the Hefez Lazzeri path
system for the fibre at $\l=1$ by elements of the set
$\{i_1i_2\,|\,1\leq i_1\leq l_1$, $1\leq i_2\leq l_2\}$.
We assign indices accordingly to the punctures in the fibre at $\l=1$, to the
lines of the arrangement and hence to any puncture in any fibre.

Our aim is to describe the parallel transport along radial paths and along
circle segments with radius $1$ or close to $0$.
We will find appropriate diffeomorphisms and obtain transported arcs.

\begin{nota}
We introduce polar coordinates $\l=t\ethe$, $\ethe:=e^{i\varth}$ of
unit absolute value and $t\in\RR^{\geq0}$.
\end{nota}

\begin{defi}
A parameter $t\ethe$ is called {\it critical}, if there is a pair
$i_1i_2,j_1j_2$ of indices such that the corresponding lines meet at
$t\ethe$.\\
The pair may be specified and $t\ethe$ called {\it critical for the pair}
$i_1i_2,j_1j_2$.
\end{defi}

Let us first outline our general approach. For a family we first give a vector
field on its total space. Next we check that the punctures form integral
curves, so the corresponding flow preserves the punctures. Then we obtain some
of the properties of the induced diffeomorphisms, to get finally the parallel
transport of some geometric objects.\\
As most important technical tool we employ {\it bump functions}:

\begin{nota}
We introduce smooth functions $\chi,\chi_\e:\CC\to\RR$ for any real $\e>0$:
\begin{eqnarray*}
   \chi &:& 0\leq\chi(z)=\chi(|z|)\leq 1,\chi(z)=0\text{ if }|z|\geq1,
\chi(z)=1\text{ if }|z|\leq\frac12,\\
   \chi_\e &:& \chi_\e(z)=\chi(z/\e),
\end{eqnarray*}
with support contained in the unit disc, resp.\ the disc of radius $\e$.
\end{nota}

For a start we consider the local situation, i.e.\ the family $z-t=0$,
$|z|\leq1$, $t\in[-\frac12,\frac12]$.

\begin{lemma}
\labell{vect/disc}
On the unit disc $|z|\leq1$ there is the vector field
$$
v(z)=\chi(z),
$$
which induces a family $\phi_t$, $t\in[-\frac12,\frac12]$ of diffeomorphisms of
the disc such that
\begin{enumerate}
\item
$\phi_t$ preserves horizontal lines,
\item
$Re\, \phi_{t_2}(z)\geq Re\, \phi_{t_1}(z)$ if $t_2\geq t_1$, i.e.\ $\phi_t$
propagates on horizontal lines,
\item
$\phi_t(0)=t$, i.e.\ $z=t$ is an integral curve.
\end{enumerate}
\end{lemma}

\proof
The first and second property are consequences of the vector field directing
parallel to the positive real numbers. The last is due to the vector
field being
the derivative of the function
$z\mapsto Re\, z$ in the disc of radius $\frac12$.
\qed

\begin{examp}
As an illustration we sketch the images of some vertical lines under
$\phi_0,\phi_\frac14$ and $\phi_\frac12$:

\setlength{\unitlength}{3.5mm}
\begin{picture}(34,12)(-2,-2)
\put(0,0){\begin{picture}(10,10)(-5,-5)

\bezier{300}(-5,0)(-5,2.5)(-3,4)
\bezier{300}(-3,4)(0,6)(3,4)
\bezier{300}(3,4)(5,2.5)(5,0)
\bezier{300}(5,0)(5,-2.5)(3,-4)
\bezier{300}(3,-4)(0,-6)(-3,-4)
\bezier{300}(-3,-4)(-5,-2.5)(-5,0)


\drawline(0,5)(0,-5)
\drawline(3,4)(3,-4)
\drawline(-3,4)(-3,-4)
\drawline(-1.5,4.73)(-1.5,-4.73)
\drawline(1.5,4.73)(1.5,-4,73)

\put(-.5,-6.2){$\phi_0$}

\end{picture}}

\put(12,0){\begin{picture}(10,10)(-5,-5)

\bezier{300}(-5,0)(-5,2.5)(-3,4)
\bezier{300}(-3,4)(0,6)(3,4)
\bezier{300}(3,4)(5,2.5)(5,0)
\bezier{300}(5,0)(5,-2.5)(3,-4)
\bezier{300}(3,-4)(0,-6)(-3,-4)
\bezier{300}(-3,-4)(-5,-2.5)(-5,0)


\bezier{100}(0,5)(0,4)(.75,2.5)
\bezier{100}(.75,2.5)(1.25,1.5)(1.25,0)
\bezier{100}(1.25,0)(1.25,-1.5)(.75,-2.5)
\bezier{100}(.75,-2.5)(0,-4)(0,-5)

\bezier{100}(3,4)(3,3.2)(3.3,2)
\bezier{100}(3.3,2)(3.5,1.2)(3.5,0)
\bezier{100}(3.5,0)(3.5,-1.2)(3.3,-2)
\bezier{100}(3.3,-2)(3,-3.2)(3,-4)

\bezier{100}(-3,4)(-3,3.2)(-2.7,2)
\bezier{100}(-2.7,2)(-2.5,1.2)(-2.5,0)
\bezier{100}(-2.5,0)(-2.5,-1.2)(-2.7,-2)
\bezier{100}(-2.7,-2)(-3,-3.2)(-3,-4)

\bezier{100}(-1.5,4.73)(-1.5,3.84)(-1.05,2.4)
\bezier{100}(-1.05,2.4)(-.75,1.44)(-.75,0)
\bezier{100}(-.75,0)(-.75,-1.44)(-1.05,-2.4)
\bezier{100}(-1.05,-2.4)(-1.5,-3.84)(-1.5,-4.73)

\bezier{100}(1.5,4.73)(1.5,3.84)(2.02,2.4)
\bezier{100}(2.02,2.4)(2.37,1.44)(2.37,0)
\bezier{100}(2.02,-2.4)(2.37,-1.44)(2.37,0)
\bezier{100}(1.5,-4.73)(1.5,-3.84)(2.02,-2.4)

\put(-.5,-6.2){$\phi_\frac14$}

\end{picture}}

\put(24,0){\begin{picture}(10,10)(-5,-5)

\bezier{300}(-5,0)(-5,2.5)(-3,4)
\bezier{300}(-3,4)(0,6)(3,4)
\bezier{300}(3,4)(5,2.5)(5,0)
\bezier{300}(5,0)(5,-2.5)(3,-4)
\bezier{300}(3,-4)(0,-6)(-3,-4)
\bezier{300}(-3,-4)(-5,-2.5)(-5,0)


\bezier{100}(0,5)(0,4)(1.5,2.5)
\bezier{100}(1.5,2.5)(2.5,1.5)(2.5,0)
\bezier{100}(2.5,0)(2.5,-1.5)(1.5,-2.5)
\bezier{100}(1.5,-2.5)(0,-4)(0,-5)

\bezier{100}(3,4)(3,3.2)(3.6,2)
\bezier{100}(3.6,2)(4,1.2)(4,0)
\bezier{100}(4,0)(4,-1.2)(3.6,-2)
\bezier{100}(3.6,-2)(3,-3.2)(3,-4)

\bezier{100}(-3,4)(-3,3.2)(-2.4,2)
\bezier{100}(-2.4,2)(-2,1.2)(-2,0)
\bezier{100}(-2,0)(-2,-1.2)(-2.4,-2)
\bezier{100}(-2.4,-2)(-3,-3.2)(-3,-4)

\bezier{100}(-1.5,4.73)(-1.5,3.84)(-.6,2.4)
\bezier{100}(-.6,2.4)(0,1.44)(0,0)
\bezier{100}(0,0)(0,-1.44)(-.6,-2.4)
\bezier{100}(-.6,-2.4)(-1.5,-3.84)(-1.5,-4.73)

\bezier{100}(1.5,4.73)(1.5,3.84)(2.55,2.4)
\bezier{100}(2.55,2.4)(3.25,1.44)(3.25,0)
\bezier{100}(2.55,-2.4)(3.25,-1.44)(3.25,0)
\bezier{100}(1.5,-4.73)(1.5,-3.84)(2.55,-2.4)

\put(-.5,-6.2){$\phi_\frac12$}

\end{picture}}

\end{picture}
\end{examp}

So we may infer the following result.

\begin{lemma}
The family of diffeomorphisms represents the diffeotopy class
associated to the family of punctured discs given by the function $f_t(z)=z-t$.
\end{lemma}

\subsubsection*{radial families}

Now we investigate the model discriminant family restricted to a radial path
$t\etheO$, $t\in[t_0,1]$. We will
consider the case only when this restriction has constant number of punctures,
in which case we call it a {\it regular family}.
Since all punctures depend affine linearly on the parameter without ever
meeting, the local situation is modeled on the case considered first.\\
When passing to a global view, we want to understand the corresponding parallel
transport diffeomorphism mapping the initial fibre to the terminal fibre.
Considered as an endomorphism of the plane it is seen to be
supported on the set of points which are close enough to some puncture at
some parameter, i.e.\ close enough to the union of their {\it traces}:

\begin{defi}
The {\it trace} of index $i_1i_2$ in a family is the set of points $z$
in the plane $\CC$
such that $z$ is a puncture of index $i_1i_2$ for some parameter of the family
base.
\end{defi}

This we can make explicit with a quick check:

\begin{lemma}
\labell{vect/radi}
Let $\e>0$ be bounded from above by half the minimal distance between punctures
in the fibres of the regular family over $t\etheO$, $t\in[t_0,1]$.
Then the punctures form integral curves for the vector field
\begin{eqnarray*}
v_\e(z,t) & = &
\sum_{\xi_1^{l_1},\xi_2^{l_2}=1}
\chi_\e(z-t\etheO\xi_1-\eta_2\xi_2)\etheO\xi_1,
\end{eqnarray*}
and the corresponding diffeomorphisms are supported on the $\e$-neighbourhood
of the union of all traces.
\end{lemma}

Hence parallel transport only affects small neighbourhoods of the punctures.
Any arc will be changed only due to the movements of its endpoints
and of the critical values which come close enough, to distances less
than $\e$ in fact.
So we can imagine what happens to a given arc in the fibre at $t_0$:

Let the arc be a piece of rope.
As the parameter $t$ increases additional rope is laid out on the
traces of both the critical values which form the ends of the arc.
A critical value about to cross the arc will push it ahead and lay out a double
rope behind forming a loop around its trace.

Likewise any time a critical value crosses a trace along which a multiple rope
has previously laid down, it picks this rope up and pushes a multiple loop into
it along its own trace.

So in the end the rope is lain down in an arbitrary small neighbourhood of the
union of all traces, in fact the union can be restricted to that part of each
trace traced after the corresponding critical value
picked up rope for the first time.\\

We want to apply parallel transport to a very restricted set of arcs:

\begin{defi}
Given a critical parameter $t_0\etheO$ for the index pair $i_1i_2, j_1j_2$, an
arc between the corresponding critical points in the fibre at $t_1\etheI$ is
called {\it local v-arc} if
\begin{enumerate}
\item
it is supported on the corresponding traces,
\item
the difference $t_1-t_0$ is positive and small compared to the distances of
critical parameters.
\end{enumerate}

In case of $j_1-i_1=l_1/2$, $i_2\neq j_2$, the traces of the corresponding
critical points in a fibre $t_1\etheI$ meet only if
$\varth_1=\varth_0$. In this
case we allow
$\varth_1\neq\varth_0$ nevertheless and concede that the local v-arc are
supported  on the traces except for a small part to join them. 
\end{defi}

\begin{defi}
Parallel transport of a local v-arc in a radial family by
the differentiable flow to radius $t=1$ yields an arc called {\it tangled v-arc}.
\end{defi}

\begin{defi}
An arc in the fibre at $t_1\etheI$ is called {\it local w-arc} if
\begin{enumerate}
\item
it connects punctures of indices $i_1i_2,j_1j_2$, $i_1^+<j_1,i_2=j_2$, by four
line segments, $i_1^+=i_1+1$,
\item
two segments are supported on the traces of the two punctures,
\item
the central pair forms a sharp wedge over the trace of the puncture of index
$i_1^+i_2$,
\item
its length and $t_1$ are small compared to the distance of critical
parameters.
\end{enumerate}
\end{defi}

\begin{defi}
Parallel transport of a local w-arc in a radial family by
the differentiable flow to radius $t=1$ yields an arc called {\it tangled w-arc}.
\end{defi}

To describe the local situation at a crossing of two or more critical
points, we
consider {\it tangled tails} of punctures. These one should imagine just as a
piece of rope laid out by a critical point on its trace and tangled by
subsequent critical points. Looking locally at the tail implies that it may
decompose into several pieces.

\begin{examp}
Imagine a crossing of just two traces, then the  tangled tails look
locally like

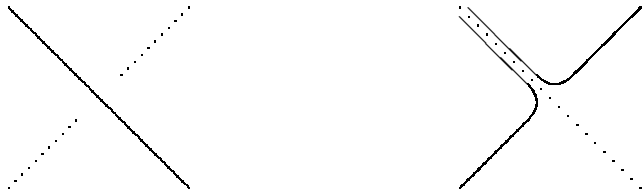
\setlength{\unitlength}{3mm}
\begin{picture}(40,10)

\put(10,5){\begin{picture}(5,5)

\bezier{140}(-4,4)(0,0)(4,-4)
\bezier{10}(-4,-4)(-2.5,-2.5)(-1,-1)
\bezier{10}(4,4)(2.5,2.5)(1,1)


\end{picture}}

\put(30,5){\begin{picture}(5,5)

\bezier{20}(-4,4)(0,0)(4,-4)
\bezier{80}(-4,-4)(-2.5,-2.5)(-1,-1)
\bezier{70}(4,4)(2.5,2.5)(1,1)

\put(-4,3.6){\line(1,-1){3}}
\put(-3.6,4){\line(1,-1){3}}

\bezier{35}(1,1)(.2,.2)(-.6,1)
\bezier{35}(-1,-1)(-.2,-.2)(-1,.6)

\end{picture}}

\end{picture}

(The critical points pass from bottom to top, the first from left to right, the
second from right to left.)
\end{examp}

By construction a local v-arc is approximately supported on tails hence so is
the transported arc throughout the radial family. In fact more is true. At each
crossing of critical points, to which the transported arc comes close, it is
approximately supported on the tails of the crossing punctures:

\begin{lemma}
\labell{loc/arc}
Locally at a crossing $P$ all local components of a tangled v-arc
can be assumed to be arbitrarily close approximations to one of the tangled
tails of the punctures passing through $P$.
\end{lemma}

\proof
All local components are laid out by a critical point which pushes them through
$P$. Hence the smaller $\e$ is, the better the approximation will be.
\qed

\begin{examp}
For the family with $l_1=3,l_2=2$ and $\varth=\frac\pi{12}$ we have
the sketches of a local v-arc at
$.12\ethe$ and its parallel transports at $.56\ethe,.78\ethe$ and $\ethe$
together with the traces of all critical values.

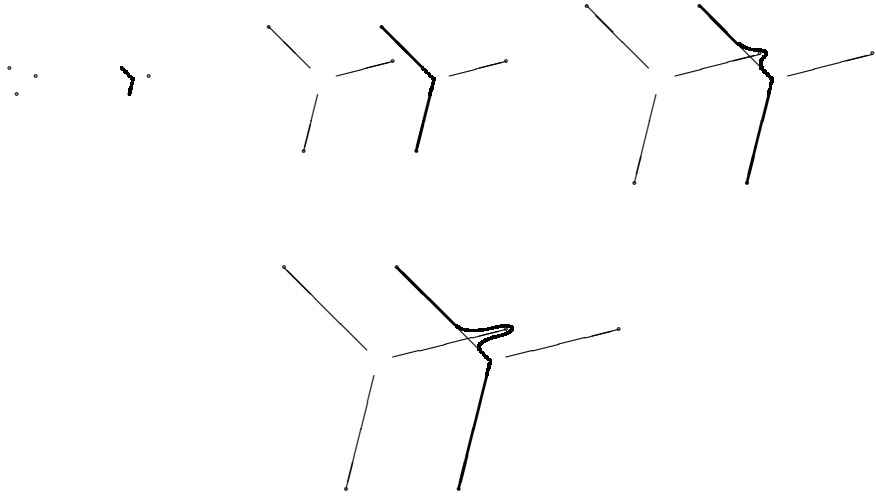
\setlength{\unitlength}{2.5mm}

\begin{picture}(25,25)(-4,0)

\put(-4,15){
\begin{picture}(12,12)(-6,-6)

\put(-3.2,-.8){\circle{.2}}
\put(-3.6,.6){\circle{.2}}
\put(-2.2,.2){\circle{.2}}

\put(2.8,-.8){\circle{.2}}
\put(2.4,.6){\circle{.2}}
\put(3.8,.2){\circle{.2}}

\thicklines
\drawline(2.8,-.8)(3,0)
\drawline(3,0)(2.4,.6)
\thinlines

\end{picture}}

\put(12,15){

\begin{picture}(12,12)(-6,-6)



\put(-3.95,-3.8){\circle{.2}}
\put(-5.8,2.8){\circle{.2}}
\put(.8,.95){\circle{.2}}

\put(2.05,-3.8){\circle{.2}}
\put(.2,2.8){\circle{.2}}
\put(6.8,.95){\circle{.2}}

\put(-3.2,-.8){\line(-1,-4){.75}}
\put(-3.6,.6){\line(-1,1){2.2}}
\put(-2.2,.2){\line(4,1){3}}

\thicklines
\put(2.8,-.8){\line(-1,-4){.75}}
\put(2.78,-.8){\line(-1,-4){.75}}
\put(2.4,.6){\line(-1,1){2.2}}
\put(2.37,.6){\line(-1,1){2.2}}
\thinlines
\put(3.8,.2){\line(4,1){3}}

\thicklines
\drawline(2.8,-.8)(3,0)
\drawline(3,0)(2.4,.6)
\thinlines

\end{picture}}

\put(30,15){
\begin{picture}(12,12)(-6,-6)



\put(-4.375,-5.5){\circle{.2}}
\put(-6.9,3.9){\circle{.2}}
\put(2.3,1.375){\circle{.2}}

\put(1.625,-5.5){\circle{.2}}
\put(-.9,3.9){\circle{.2}}
\put(8.3,1.375){\circle{.2}}

\put(-3.2,-.8){\line(-1,-4){1.175}}
\put(-3.6,.6){\line(-1,1){3.3}}
\put(-2.2,.2){\line(4,1){4.5}}

\thicklines
\put(2.8,-.8){\line(-1,-4){1.175}}
\put(2.78,-.8){\line(-1,-4){1.175}}
\thinlines
\put(2.4,.6){\line(-1,1){3.3}}
\put(3.8,.2){\line(4,1){4.5}}

\thicklines
\put(1.1,1.9){\line(-1,1){2}}
\put(1.13,1.9){\line(-1,1){2}}
\bezier{100}(1.1,1.9)(1.6,1.55)(2.15,1.55)
\bezier{20}(2.15,1.55)(2.85,1.55)(2.4,1.1)
\bezier{100}(2.4,1.1)(2.15,.85)(2.4,0.6)
\thinlines

\thicklines
\drawline(2.8,-.8)(3,0)
\drawline(3,0)(2.4,.6)
\thinlines

\end{picture}}

\put(15,0){
\begin{picture}(12,12)(-6,-6)



\put(-4.7,-6.8){\circle{.2}}
\put(-8,5){\circle{.2}}
\put(3.8,1.7){\circle{.2}}

\put(1.3,-6.8){\circle{.2}}
\put(-2,5){\circle{.2}}
\put(9.8,1.7){\circle{.2}}

\put(-3.2,-.8){\line(-1,-4){1.5}}
\put(-3.6,.6){\line(-1,1){4.4}}
\put(-2.2,.2){\line(4,1){6}}

\thicklines
\put(2.8,-.8){\line(-1,-4){1.5}}
\put(2.78,-.8){\line(-1,-4){1.5}}
\thinlines
\put(2.4,.6){\line(-1,1){4.4}}
\put(3.8,.2){\line(4,1){6}}

\thicklines
\put(1.1,1.9){\line(-1,1){3.1}}
\put(1.13,1.9){\line(-1,1){3.1}}
\bezier{100}(1.1,1.9)(1.6,1.4)(3.2,1.8)
\bezier{20}(3.2,1.8)(3.9,1.975)(4.1,1.775)
\bezier{20}(4.1,1.775)(4.1,1.525)(3.4,1.35)
\bezier{100}(3.4,1.35)(2,1)(2.4,0.6)
\thinlines

\thicklines
\drawline(2.8,-.8)(3,0)
\drawline(3,0)(2.4,.6)
\thinlines

\end{picture}}

\end{picture}

In this example only one additional critical value is
entangled.
\end{examp}

Finally we observe, that almost all radial families are regular.

\begin{defi}
The angle $\varth_1$ is called {\it regular} if for all $t_1>0$ the family over
the line segment from $t_1e_{\varth_1}$ to
$e_{\varth_1}$
\begin{enumerate}
\item
is a regular family, i.e.\ the segment does not pass a critical parameter,
\item
has no pair of distinct traces having more than one point in common,
\end{enumerate}
\end{defi}

Note that distinct traces have at most one point in common, if and only if no
trace contains a point $\eta_2\xi_2$, $\xi_2^{l_2}=1$.
So we can show:

\begin{lemma}
\labell{reg/angle}
Given a critical parameter $t_0\etheO$, then there is a regular angle
$\varth_1$ arbitrarily close to $\varth_0$.
\end{lemma}

\proof
Both regularity conditions are open.
Hence it suffices to prove that for each condition there is an angle
arbitrarily
close to $\varth_0$ such that the condition is met. For both
conditions this is easy to see, since in both cases a finite set is to be
avoided by the traces and the traces of one puncture along different radial
paths have no point in common.
\qed

\subsubsection*{circular families}

Similar to the case of radial families we can get hold of a
diffeomorphism which
represents the parallel transport over circular segments in the base. On
particular subsets the map is in fact quite easily described.

\begin{lemma}
\labell{vect/circ}
Given the vector field
\begin{eqnarray*}
v(z,\varth) & = &
i\bigg(z+\sum_{\xi_1^{l_1}=1}\chi_{2\eta_2}(z-\ethe\xi_1)(\ethe\xi_1-z)\bigg)
\end{eqnarray*}
then
\begin{enumerate}
\item
the punctures of the model family over the circle of radius $t=1$ form integral
curves,
\item
supposing $|z_0-\ethe\xi_1|\leq2\eta_2$ the flow of $v$
preserves the distance of $z_0(\varth)$ and $\ethe\xi_1$.
\item
supposing $|z_0-\xi_1|\geq2\eta_2$ for all $\xi_1, \xi_1^{l_1}=1$,
$z_0(\varth)=z_0\ethe$ is an integral curve.
\end{enumerate}
\end{lemma}

\proof
i) Each puncture forms a curve $\ethe\xi_1+\eta_2\xi_2$,
$\xi_1^{l_1},\xi_2^{l_2}=1$, for which we can check the integrality condition:
$$
\frac{d}{d\varth}(\ethe\xi_1+\eta_2\xi_2)\quad=\quad
i\ethe\xi_1\quad=\quad
i(\ethe\xi_1+\eta_2\xi_2-\eta_2\xi_2)\quad=\quad
v(\ethe\xi_1+\eta_2\xi_2,\varth).
$$
ii) We have to show that the following complex numbers
considered as real vectors are perpendicular for all $\varth$:
$$
(\ethe\xi_1-z_0(\varth))\cdot\frac{d}{d\varth}(\ethe\xi_1-z_0(\varth)).
$$
Both points move along integral curves, hence
$\frac{d}{d\varth}\ethe\xi_1=i\ethe\xi_1$ and

\begin{eqnarray*}
\frac{d}{d\varth}z_0(\varth) & = & v(z_0(\varth),\varth)\\
& = &
iz_0(\varth)+i\chi_{2\eta_2}(z_0(\varth)-\ethe\xi_1)(\ethe\xi_1-z_0(\varth)).
\end{eqnarray*}

Since $\chi$ is a real valued function, the second function
is a purely imaginary multiple of the first, hence they are orthogonal
at all $\varth$.\\
iii) Again we have only to check an integrality condition
$$
\frac{d}{d\varth}z_0\ethe\quad=\quad
iz_0\ethe\quad=\quad
v(z_0\ethe,\varth).
$$
\qed
\\
Let us rephrase the result of the lemma in more geometrical terms:
\begin{enumerate}
\item
the flow realises parallel transport in the model family over circle segments
of radius $t=1$,
\item
the $2\eta_2$-discs at points $\xi_1$, $\xi_1^{l_2}=1$ are mapped bijectively to
$2\eta_2$-discs of the transported points preserving the distance,
\item
points outside these discs are mapped by a rigid rotation around the origin.
\end{enumerate}

\begin{lemma}
\labell{vect/core}
Given the vector field for $\e<<\eta_2$
\begin{eqnarray*}
v(z) & = &
i\sum_{\xi_2^{l_2}=1}\chi_{4\e}(z-\eta_2\xi_2)(z-\eta_2\xi_2)
\end{eqnarray*}
then
\begin{enumerate}
\item
the punctures of the model family over the circle of radius $t=\e$
form integral
curves,
\item
suppose $|z_0-\eta_2\xi_2|\leq2\e$, $\xi_2^{l_2}=1$ then the curves
$z_0(\varth)=(z_0-\eta_2\xi_2)\ethe+\e_1\xi_2$ are integral for the flow of $v$,
\item
suppose $|z_0-\eta_2\xi_2|\geq4\e$ for all $\xi_2, \xi_2^{\l_2}=1$, then
$z_0(\varth)=z_0$ is an integral curve.
\end{enumerate}
\end{lemma}

\proof
i) Since each puncture is on a curve $\e\ethe\xi_1+\eta_2\xi_2$, the
assertion follows from case ii).

ii) We check the integrality condition:
$$
\frac{d}{d\varth}((z_0-\eta_2\xi_2)\ethe+\eta_2\xi_2)\quad=\quad
i(z_0-\eta_2\xi_2)\ethe\quad=\quad
v((z_0-\eta_2\xi_2)\ethe+\eta_2\xi_2),\varth).
$$

iii) Since the vector field vanishes at these points constant curves are
integral curves.
\qed
\\
Again we restate these results in geometrical terms:
\begin{enumerate}
\item
the flow realises parallel transport in the model family over circle segments
of radius $t=\e$,
\item
the $2\e$-discs at points $\eta_2\xi_2$, $\xi_2^{l_2}=1$, are rotated rigidly
under parallel transport,
\item
points outside $4\e$-discs of these points stay fix.
\end{enumerate}


\section{from tangled v-arcs to isosceles arcs}
\label{six/two}

In this section we consider two different kinds of mapping classes in a
fibre of large radius.
Both kinds are twists on embedded arcs. So we may equally well investigate
these arcs.
Arcs of the first kind are called tangled v-arcs, they are obtained from
local v-arcs by
parallel transport along a radial path using the differentiable flow of the
preceding section.

Arcs of the second kind are called {\it isosceles arcs}.
They are supported on traces of two punctures and form the two sides of an
approximate isosceles triangle. Again the degenerate case requires extra care.
If two traces are parallel but close, an arc
which is supported on these traces except for a small join between
them is called
a {\it straight isosceles arc}.

An isosceles arc is said to correspond to a tangled v-arc if it
connects the same
punctures.
In general these two arcs are not isotopic.
But we will define a group of mapping classes such that they belong to one
orbit.
In fact we will give some arcs, such that the
group generated by the full twists on these arcs will do.
They will be called {\it bisceles arcs} for the reason that they are supported
on segments of two traces not necessarily of similar length.

Note that by this definition all isosceles arcs are subsumed under
the notion of
bisceles arcs except for the straight isosceles arcs.\\

We want to encode the isotopy class of a tangled v-arc into a
planar diagram in the fibre at $\ethe$.
This diagram will consist of all the traces each of which is directed
from its {\it source point} -- which is one of $\xi_2$, $\xi_2^{l_2}=1$ -- to its
puncture.

Apart from the source points, there are only ordinary crossings,
which are given
by the mutual transversal intersection of several traces.

Crossings which are sufficiently close to the tangled v-arc are called {\it
vertices} of the diagram. The segments of traces close to the tangled v-arc are
called {\it essential traces}, they connect a vertex to a puncture.

At each vertex we put an order on the essential traces. The first or {\it
dominant trace} is the one which passed last, which is incidentally
the one such
that the puncture end is closest.
The other follow according to increasing distance to their puncture end.
The order can be made explicit by labels assigned to the essential traces at
each vertex. We can also make the dominant trace pass over by
replacing the other
traces by broken lines. Finally the lines are labeled at their ends
by the index
of the corresponding puncture.

We define the {\it essential diagram} to be obtained by discarding all lines
except the essential traces and we notice that the tangled v-arc is still
determined by this datum.\\

\begin{examp}
  From the tangled v-arc of the previous example, we get the following
diagram for $l_1=3,l_2=2,\varth=\frac\pi{12}$, in which we have discarded all
traces which do not pass a vertex.

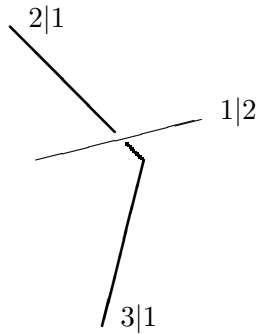
\begin{figure}[h]
\setlength{\unitlength}{2.4mm}
\begin{picture}(20,20)(-28,-10)

\put(7.2,2.3){$1|2$}
\put(1.7,-9.2){$3|1$}
\put(-3.4,7.4){$2|1$}

\put(-3,0){\line(4,1){9.2}} 
\thicklines
\put(3,0){\line(-1,-4){2.3}} \drawline(3,0)(.7,-9.2)
\put(1.4,1.6){\line(-1,1){5.8}} \drawline(1.4,1.6)(-4.4,7.4)
\put(3,0){\line(-1,1){.9}} \drawline(3,0)(2.1,.9)
\thinlines

\end{picture}
\end{figure}
\end{examp}

\begin{defi}
No essential diagram contains a directed cycle, hence the {\it
height} function on
vertices is well-defined by

\begin{eqnarray*}
\ON{ht}(P) & = & \max_{P'< P}(\ON{ht}(P'),0)+1.
\end{eqnarray*}

where the maximum is taken over all vertices $P'$ between $P$ and a puncture on
an essential trace. Each such vertex is called
{\it subordinate} to $P$.
\end{defi}

Given an essential diagram we consider {\it simple transformations}
at vertices.
We may change the crossing order at a vertex $P$ if and only if all
traces through $P$ are dominant at each subordinate vertex.
Note that on transformed diagrams we have to make the order explicit, since it
can no longer be read off the distances to the punctures.

The first observation is that we can change  an essential diagram by simple
transformations only to get a diagram in which the traces of the
v-arc punctures
are dominant at all vertices they cross.

\begin{lemma}
\labell{apex/change}
Given any vertex there is a composition of simple transformations which
changes the crossing order at this vertex but nowhere else.
\end{lemma}

\proof
If the vertex is of height one we can change it by a simple transformation.
If not, a simple transformation can only be performed if the essential traces
are dominant on subordinate vertices. But then we can argue inductively on
the height of the vertex. All subordinate vertices are of less height, so by
induction we may assume the existence of a composite transformation which
makes the traces under consideration dominant there.

Then we can perform the simple transformation to change the local
order. Finally
we invoke the inverse of the composite transformation to put all other
transformed vertices back to their initial state.
\qed

In particular, a series of simple transformations can be found such that the
traces of the v-arc punctures become dominant.

\begin{examp}
We illustrate this procedure in the following sequence of diagrams. In each step
we perform a simple transformation on some vertices which do not share
subordinate vertices.

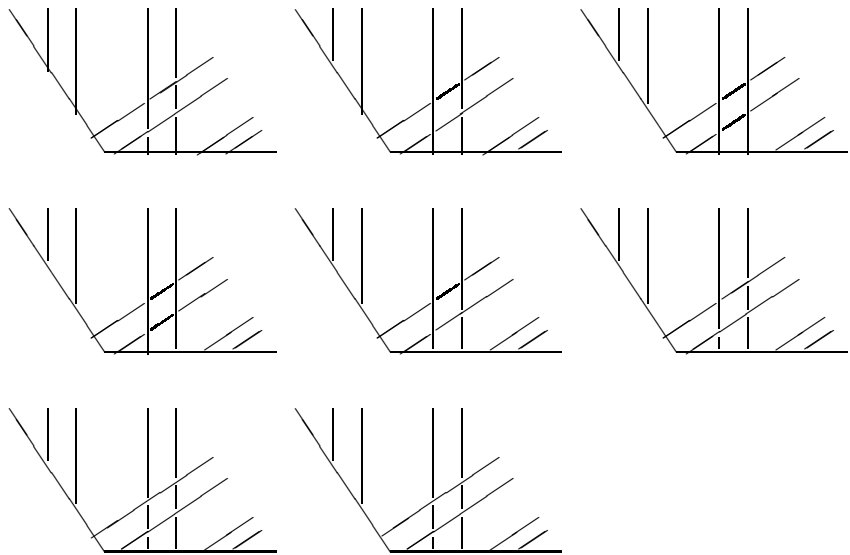
\setlength{\unitlength}{3.8mm}

\begin{picture}(30,21)

\put(4,15){
\begin{picture}(2,2)

\put(0,0){\line(1,0){6}}
\put(0,0){\line(-2,3){3.33}}

\put(5.5,.75){\line(-3,-2){1.25}}    
\put(5.2,1.2){\line(-3,-2){1.95}}    
\put(4.3,2.55){\line(-3,-2){3.95}}   
\put(3.8,3.3){\line(-3,-2){2.2}}    
\put(1.4,1.7){\line(-3,-2){1.85}}    

\put(2.5,5){\line(0,-1){2.4}}  
\put(1.5,5){\line(0,-1){4.2}}  
\put(2.5,2.3){\line(0,-1){.8}}  
\put(1.5,.5){\line(0,-1){.6}}  
\put(2.5,1.2){\line(0,-1){1.3}}  
\put(-1,5){\line(0,-1){3.7}}  
\put(-2,5){\line(0,-1){2.2}}  

\end{picture}}

\put(14,15){
\begin{picture}(2,2)

\put(0,0){\line(1,0){6}}
\put(0,0){\line(-2,3){3.33}}

\put(5.5,.75){\line(-3,-2){1}}    
\put(5.2,1.2){\line(-3,-2){1.95}}    
\put(4.3,2.55){\line(-3,-2){2.7}}   
\put(1.4,.616){\line(-3,-2){1.05}}   
\put(3.8,3.3){\line(-3,-2){1.2}}    
\put(2.4,2.36){\line(-3,-2){.8}}    
\drawline(2.4,2.36)(1.65,1.86)
\put(1.4,1.7){\line(-3,-2){1.85}}    

\put(2.5,5){\line(0,-1){3.5}}  
\put(1.5,5){\line(0,-1){5.1}}  
\put(2.5,1.2){\line(0,-1){1.3}}  
\put(-1,5){\line(0,-1){3.7}}  
\put(-2,5){\line(0,-1){1.8}}  

\end{picture}}

\put(24,15){
\begin{picture}(2,2)

\put(0,0){\line(1,0){6}}
\put(0,0){\line(-2,3){3.33}}

\put(5.5,.75){\line(-3,-2){1}}    
\put(5.2,1.2){\line(-3,-2){1.7}}    
\put(4.3,2.55){\line(-3,-2){1.7}}   
\drawline(2.4,1.283)(1.65,.783)
\put(2.4,1.283){\line(-3,-2){.8}}   
\put(1.4,.616){\line(-3,-2){1.05}}   
\put(3.8,3.3){\line(-3,-2){1.2}}    
\drawline(2.4,2.36)(1.65,1.86)
\put(2.4,2.36){\line(-3,-2){.8}}    
\put(1.4,1.7){\line(-3,-2){1.85}}    

\put(2.5,5){\line(0,-1){5.1}}  
\put(1.5,5){\line(0,-1){5.1}}  
\put(-1,5){\line(0,-1){3.3}}  
\put(-2,5){\line(0,-1){1.8}}  

\end{picture}}

\put(4,8){
\begin{picture}(2,2)

\put(0,0){\line(1,0){6}}
\put(0,0){\line(-2,3){3.33}}

\put(5.5,.75){\line(-3,-2){1}}    
\put(5.2,1.2){\line(-3,-2){1.7}}    
\put(4.3,2.55){\line(-3,-2){1.7}}   
\drawline(2.4,1.283)(1.65,.783)
\put(2.4,1.283){\line(-3,-2){.8}}   
\put(1.4,.616){\line(-3,-2){1.05}}   
\put(3.8,3.3){\line(-3,-2){1.2}}    
\drawline(2.4,2.36)(1.65,1.86)
\put(2.4,2.36){\line(-3,-2){.8}}    
\put(1.4,1.7){\line(-3,-2){1.85}}    

\put(2.5,5){\line(0,-1){4.9}}  
\put(1.5,5){\line(0,-1){5.1}}  
\put(-1,5){\line(0,-1){3.3}}  
\put(-2,5){\line(0,-1){1.8}}  

\end{picture}}

\put(14,8){
\begin{picture}(2,2)

\put(0,0){\line(1,0){6}}
\put(0,0){\line(-2,3){3.33}}

\put(5.5,.75){\line(-3,-2){1}}    
\put(5.2,1.2){\line(-3,-2){1.7}}    
\put(4.3,2.55){\line(-3,-2){2.7}}   
\put(1.4,.616){\line(-3,-2){1.05}}   
\put(3.8,3.3){\line(-3,-2){1.2}}    
\drawline(2.4,2.36)(1.65,1.86)
\put(2.4,2.36){\line(-3,-2){.8}}    
\put(1.4,1.7){\line(-3,-2){1.85}}    

\put(2.5,5){\line(0,-1){3.5}}  
\put(1.5,5){\line(0,-1){4.9}}  
\put(2.5,1.2){\line(0,-1){1.1}}  
\put(-1,5){\line(0,-1){3.3}}  
\put(-2,5){\line(0,-1){1.8}}  

\end{picture}}

\put(24,8){
\begin{picture}(2,2)

\put(0,0){\line(1,0){6}}
\put(0,0){\line(-2,3){3.33}}

\put(5.5,.75){\line(-3,-2){1}}    
\put(5.2,1.2){\line(-3,-2){1.7}}    
\put(4.3,2.55){\line(-3,-2){3.95}}   
\put(3.8,3.3){\line(-3,-2){2.2}}    
\put(1.4,1.7){\line(-3,-2){1.85}}    

\put(2.5,5){\line(0,-1){2.4}}  
\put(1.5,5){\line(0,-1){4.2}}  
\put(2.5,2.3){\line(0,-1){.8}}  
\put(1.5,.5){\line(0,-1){.4}}  
\put(2.5,1.2){\line(0,-1){1.1}}  
\put(-1,5){\line(0,-1){3.3}}  
\put(-2,5){\line(0,-1){1.8}}  

\end{picture}}

\put(4,1){
\begin{picture}(2,2)

\put(0,0){\line(1,0){6}}
\put(0,0){\line(-2,3){3.33}}

\put(5.5,.75){\line(-3,-2){1}}    
\put(5.2,1.2){\line(-3,-2){1.7}}    
\put(4.3,2.55){\line(-3,-2){3.7}}   
\put(3.8,3.3){\line(-3,-2){4.25}}    

\put(1.5,5){\line(0,-1){3.1}}  
\put(2.5,5){\line(0,-1){2.4}}  
\put(1.5,1.6){\line(0,-1){.8}}  
\put(2.5,2.3){\line(0,-1){.8}}  
\put(1.5,.5){\line(0,-1){.4}}  
\put(2.5,1.2){\line(0,-1){1.1}}  
\put(-1,5){\line(0,-1){3.3}}  
\put(-2,5){\line(0,-1){1.8}}  

\end{picture}}

\put(14,1){
\begin{picture}(2,2)

\put(0,0){\line(1,0){6}}
\put(0,0){\line(-2,3){3.33}}

\put(5.5,.75){\line(-3,-2){1}}    
\put(5.2,1.2){\line(-3,-2){1.7}}    
\put(4.3,2.55){\line(-3,-2){3.7}}   
\put(3.8,3.3){\line(-3,-2){4.1}}    

\put(1.5,5){\line(0,-1){3.1}}  
\put(2.5,5){\line(0,-1){2.4}}  
\put(1.5,1.6){\line(0,-1){.8}}  
\put(2.5,2.3){\line(0,-1){.8}}  
\put(1.5,.5){\line(0,-1){.4}}  
\put(2.5,1.2){\line(0,-1){1.1}}  
\put(-1,5){\line(0,-1){3.3}}  
\put(-2,5){\line(0,-1){1.8}}  

\end{picture}}

\end{picture}
\end{examp}

The important step is to see, that for any simple transformation at a
vertex $P$
there is a choice of a mapping class such that
\begin{enumerate}
\item
the mapping class is given by a product of full twists on bisceles
arcs supported
on the essential traces through $P$,
\item
a diagram transformed by a sequence of simple transformations encodes
the isotopy
class of the tangled v-arc transformed by the composition of the chosen mapping
classes.
\end{enumerate}

For the induction in the proof of the following lemma we need also a relation
between tails at a vertex.

\begin{defi}
At a vertex a tail {\it dominates} another one, if it is isotopic to its trace
up to an isotopy fixing the endpoints of both tails but not necessarily the
punctures not involved.
\end{defi}

\begin{lemma}
\label{diffeo}
Given a diagram with orders at its vertices which are obtained by a
composition of
simple transformations from those of the essential diagram of a tangled v-arc.
Then there is a diffeomorphism such that
\begin{enumerate}
\item
  it represents a mapping class which is a product of full twists on
bisceles arcs
supported on essential traces,
\item
it is supported close to the essential traces,
\item
locally at every vertex the dominant trace is close to the image of the
corresponding tail.
\end{enumerate}
\end{lemma}

\proof
We assume in addition that each simple transformation reverses the order of
consecutive traces and start an induction on the number of such transformations
in the composite transformation.

So we consider a simple transformation. For simplicity we first assume that the
vertex at which the order is changed is met by only two essential traces. By
assumption these traces are dominant at subordinate vertices, hence
we can depict
the tangled tails of the two punctures involved as follows:\vspace{2mm}

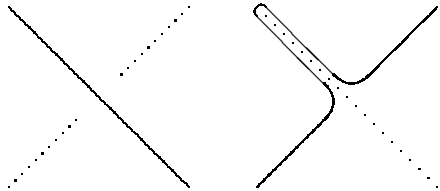
\setlength{\unitlength}{3mm}
\begin{picture}(40,10)

\put(16,5){\begin{picture}(5,5)

\bezier{140}(-4,4)(0,0)(4,-4)
\bezier{10}(-4,-4)(-2.5,-2.5)(-1,-1)
\bezier{10}(4,4)(2.5,2.5)(1,1)


\end{picture}}

\put(27,5){\begin{picture}(5,5)

\bezier{20}(-4,4)(0,0)(4,-4)
\bezier{80}(-4,-4)(-2.5,-2.5)(-1,-1)
\bezier{70}(4,4)(2.5,2.5)(1,1)

\put(-4,3.6){\line(1,-1){3}}
\put(-3.6,4){\line(1,-1){3}}

\bezier{35}(1,1)(.2,.2)(-.6,1)
\bezier{35}(-1,-1)(-.2,-.2)(-1,.6)

\bezier{8}(-4,3.6)(-4.2,3.8)(-4,4)
\bezier{8}(-3.6,4)(-3.8,4.2)(-4,4)

\end{picture}}

\end{picture}

\begin{quote}
(The critical points pass from bottom to top, the first from left to right, the
second from right to left.)
\end{quote}
\pagebreak

Now a full twist on the bisceles arc with the appropriate choice of orientation
can be performed close to these traces to yield:\vspace{2mm}

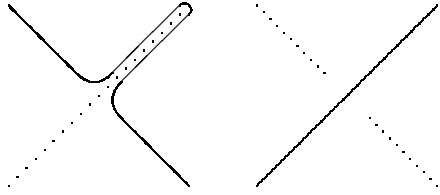
\setlength{\unitlength}{3mm}
\begin{picture}(40,10)

\put(27,5){\begin{picture}(5,5)

\bezier{140}(4,4)(0,0)(-4,-4)
\bezier{10}(4,-4)(2.5,-2.5)(1,-1)
\bezier{10}(-4,4)(-2.5,2.5)(-1,1)


\end{picture}}

\put(16,5){\begin{picture}(5,5)

\bezier{20}(4,4)(0,0)(-4,-4)
\bezier{70}(4,-4)(2.5,-2.5)(1,-1)
\bezier{70}(-4,4)(-2.5,2.5)(-1,1)

\put(4,3.6){\line(-1,-1){3}}
\put(3.6,4){\line(-1,-1){3}}

\bezier{35}(-1,1)(-.2,.2)(.6,1)
\bezier{35}(1,-1)(.2,-.2)(1,.6)

\bezier{8}(4,3.6)(4.2,3.8)(4,4)
\bezier{8}(3.6,4)(3.8,4.2)(4,4)

\end{picture}}

\end{picture}

Hence our claim is true in this case.

The same applies if there are more essential traces and we want to reverse the
order of the first two, since the corresponding tails are not effected by tails
of lower order.
\\

The situation changes drastically if our simple transformation reverses the
order of traces none of which is dominant. Then the picture is modified by the
essential traces of larger order pushing loops into the depicted tails.

But on the same time they push loops into the bisceles arc and hence into the
support of the diffeomorphism we want to perform.
Hence we need only to show that this pushed diffeomorphism will do.
\\

Of course it has the second property. It also has the first property since the
full twist on the modified bisceles arc is isotopic to the full twist on the
bisceles arc conjugated by full twists on bisceles arcs with apex in the same
vertex.

The third property is given, since the dominant traces and the corresponding
tails are locally not changed except for the explicit case considered first,
where the property can be simply checked.

Moreover for the induction process we should notice that any of our
diffeomorphisms preserves domination of a tail over another one, except that it
exchanges the role of the tails corresponding to the traces of which the order
has been reversed.
\\

To proceed our induction the first two properties are no obstacle. But we have
to prove that the third property is preserved when performing an additional
transformation.

If the additional transformation does not affect a dominant trace, then neither
does the diffeomorphism we perform. Since it also preserves the corresponding
tail, we are done in this case.

So let us assume the additional transformation affects a dominant trace. Then
the diffeomorphism we choose also affects both the trace and the tail. What we
have to show is that the image of the tail which was second before and is first
now has the claimed property.

By assumption this tail is only tangled along the essential traces through the
vertex under consideration. Moreover we may assume that it dominates all tails
through this vertex apart from the dominant one. Hence it is only
tangled by the
dominant trace and our diffeomorphism can be chosen to map it close
to its trace as in the case depicted above.
\qed

\begin{lemma}
Given a tangled v-arc there is a mapping class given by a composition of full
twists on bisceles arcs supported on essential traces which maps the tangled
v-arc to the isotopy class of the corresponding isosceles arc.
\end{lemma}

\proof
By lemma \ref{apex/change} there is a composition of simple transformation
changing vertex orders of the essential diagram of the given tangled v-arc in
such a way that the traces of both puncture ends are dominant at each vertex.

Then by lemma \ref{diffeo} there is a diffeomorphism representing a mapping
class as in the claim, which maps the tangled tails in such a way that locally
at each vertex the dominant trace is close to its tail.

Thus the images of the tangled tails of both puncture ends may no longer
deviate from the traces at any vertex. So they are isotopic to the traces and
we conclude that the image arc is isotopic to the corresponding isosceles arc.
\qed

We did not bother to adjust our arguments explicitly for $j_1-i_1=l_1/2$,
since we can choose $0<\e\ll t_0|\varth_0-\varth_1|$ small in
comparison with the minimal diameter of local neighbourhoods of vertices.\\

We close this section with two observation, which will be used later:

\begin{remark}
\label{shorter}
All bisceles arcs supported on essential traces are -- apart from the obvious
one -- not isosceles arcs, since one critical point has to pass after
the other.

For the same reason, the {\it length} of each bisceles arc supported on
essential traces of a tangled v-arc is bounded by the length of the
corresponding
isosceles arc.\\ The {\it length} is defined to be the maximum of the lengths
of the two sides.
\end{remark}


\section{from bisceles arcs to coiled isosceles arcs}
\label{six/tri}

We stay in the same fibre as before, so we work in the same group of
mapping classes.
And we are still interested into orbits of subgroups generated by full twists
on bisceles arcs.

We have accomplished so far, that we can express a tangled v-arc by means of
an isosceles arc and twists on bisceles arcs.
Now in a similar way we want to relate bisceles arcs and straight
isosceles arcs
to a third kind of arcs called {\it coiled isosceles arcs}.
With straight isosceles arcs we will deal only at the end of the section.

Again a bisceles arc and the associated coiled isosceles arc connect the same
pair of punctures and -- though not isotopic in general -- belong to one orbit
of a group generated by twists on specific bisceles arcs.

To make these statements precise, we first need to introduce some
more geometric
notions.

\begin{defi}
The {\it central core} is the disc of radius $\eta_2$ at the origin with
all source points distributed on its boundary circle.
\end{defi}

\begin{defi}
The {\it peripheral cores} are the discs of radius $\eta_2$ centred at the points
$\xi_1\ethe,\xi_1^{l_1}=1$. All critical points for $\l=\ethe$ are distributed
on their boundaries, the {\it peripheral circles}.
\end{defi}

By looking at the following sketches we notice that a
bisceles arc can take
essentially two different positions relative to a peripheral
core which contains one of its punctures.

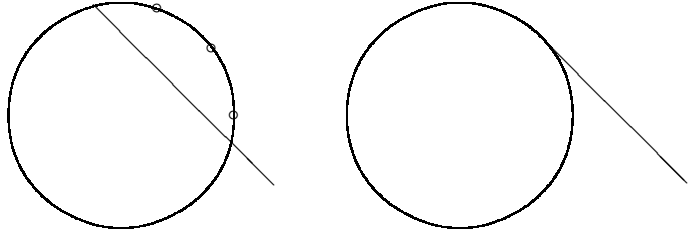
\setlength{\unitlength}{3mm}
\begin{picture}(40,10)

\put(10,0){\begin{picture}(10,10)

\bezier{300}(0,5)(0,7.5)(2,9)
\bezier{300}(2,9)(5,11)(8,9)
\bezier{300}(8,9)(10,7.5)(10,5)
\bezier{300}(10,5)(10,2.5)(8,1)
\bezier{300}(8,1)(5,-1)(2,1)
\bezier{300}(2,1)(0,2.5)(0,5)

\put(3.8,9.9){\line(1,-1){8}}

\put(10,5){\circle{.3}}
\put(9,8){\circle{.3}}
\put(6.6,9.75){\circle{.3}}

\end{picture}}

\put(25,0){\begin{picture}(10,10)

\bezier{300}(0,5)(0,7.5)(2,9)
\bezier{300}(2,9)(5,11)(8,9)
\bezier{300}(8,9)(10,7.5)(10,5)
\bezier{300}(10,5)(10,2.5)(8,1)
\bezier{300}(8,1)(5,-1)(2,1)
\bezier{300}(2,1)(0,2.5)(0,5)

\put(9.1,8){\line(1,-1){6}}

\end{picture}}

\end{picture}

\begin{defi}
A bisceles arc is called {\it unobstructed} if it is isotopic to some arc
supported outside the peripheral cores. It is called {\it obstructed} otherwise.

A bisceles arc of index pair $i_1i_2,j_1j_2$ is said to be {\it
obstructed on the
i-side}, if punctures of index $i_1i_2'$ are obstacles to unobstructedness.
\end{defi}

If a bisceles arc is obstructed then at least one side cuts
through the corresponding peripheral circle and thus divides the set of
critical points on the circle into two subsets.

\begin{defi}
If a bisceles arc is obstructed, then a set of critical points is called {\it
obstructing set}, if the bisceles arc is unobstructed in the complement of the
other punctures, i.e.\ isotopic to some arc supported outside the peripheral
cores.
\end{defi}

Since we may not isotopy arcs through punctures, we have to resort to changing
the isotopy class by means of full twists on some suitable bisceles arcs.
This has to be done in such a way, that up to isotopy the terminal part of
the obstructed bisceles arc is simply replaced by a
spiral segment coiled around the peripheral core.
\\

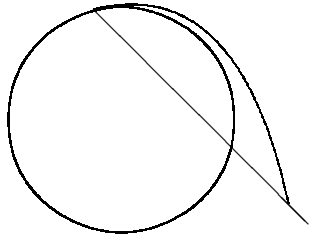
\setlength{\unitlength}{3mm}
\begin{picture}(20,10)(-10,-1)

\bezier{300}(0,5)(0,7.5)(2,9)
\bezier{300}(2,9)(5,11)(8,9)
\bezier{300}(8,9)(10,7.5)(10,5)
\bezier{300}(10,5)(10,2.5)(8,1)
\bezier{300}(8,1)(5,-1)(2,1)
\bezier{300}(2,1)(0,2.5)(0,5)

\put(3.8,9.9){\line(1,-1){9.5}}
\bezier{200}(3.8,9.9)(10.2,11.5)(12.4,1.3)

\end{picture}

To do so properly we choose a suitable obstructing set and
employ twists on arcs which
are supported on pairs of parallels to the sides of the bisceles arc and
which connect a point of the obstructing set to another one or to a
puncture of the
bisceles arc.\\

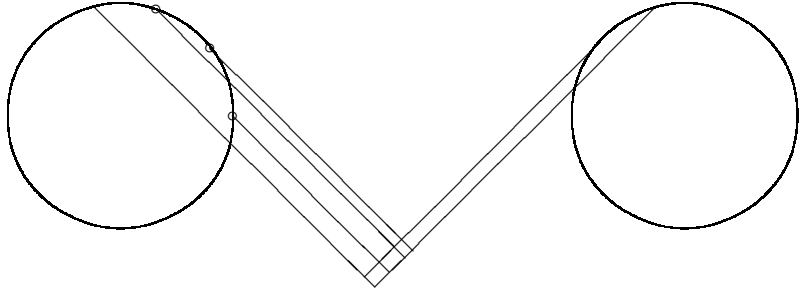
\setlength{\unitlength}{3mm}
\begin{picture}(40,15)

\put(5,5){\begin{picture}(10,10)

\bezier{300}(0,5)(0,7.5)(2,9)
\bezier{300}(2,9)(5,11)(8,9)
\bezier{300}(8,9)(10,7.5)(10,5)
\bezier{300}(10,5)(10,2.5)(8,1)
\bezier{300}(8,1)(5,-1)(2,1)
\bezier{300}(2,1)(0,2.5)(0,5)

\put(3.8,9.9){\line(1,-1){12.5}}

\put(10,5){\circle{.3}}
\put(9,8){\circle{.3}}
\put(6.6,9.75){\circle{.3}}

\put(10,5){\line(1,-1){6.94}}
\put(9,8){\line(1,-1){9}}
\put(6.6,9.75){\line(1,-1){11.04}}

\end{picture}}

\put(30,5){\begin{picture}(10,10)

\bezier{300}(0,5)(0,7.5)(2,9)
\bezier{300}(2,9)(5,11)(8,9)
\bezier{300}(8,9)(10,7.5)(10,5)
\bezier{300}(10,5)(10,2.5)(8,1)
\bezier{300}(8,1)(5,-1)(2,1)
\bezier{300}(2,1)(0,2.5)(0,5)

\put(1,8){\line(-1,-1){10.16}}
\put(3.8,9.9){\line(-1,-1){12.5}}

\end{picture}}

\end{picture}

By construction a bisceles arc bounds a well defined convex cone
which we call the {\it inner cone} of the bisceles arc.

Thus given an obstructed bisceles arc, the critical points on its
peripheral circles in the inner cone form a natural obstructing set and the
parallels for this obstructing set are naturally called either {\it inner
parallels} or {\it obstructing parallels} of the bisceles arc.
\\

Next we choose a topological disc, which contains the obstructed
bisceles arc and
its inner parallels, but no further critical point.
There is a natural way to  identify the mapping class group of
this disc with an abstract braid group:

Number all traces from left to right -- supposing the cone opens upwards as in
the sketch above.
Let $k'$ be the number of traces parallel to the first an let $k$ be the total
number of traces.
If $\s_{i,j}$, $1\leq i\leq k'<j\leq k$ is the class of the half twist on the
parallel supported on the $i^{th}$ and $j^{th}$ trace, then we put
\begin{eqnarray*}
\s_{i,j} & = & \s_{j,k}\s_{i,k}\s_{j,k}\inv \quad\text{ if }\quad 1\leq i<j\leq
k',\\
\s_{i,j} & = & \s_{1,i}\s_{1,j}\s_{1,i}\inv \quad\text{ if }\quad k'<i<j\leq k.
\end{eqnarray*}
Then considering the elements $\s_{i,i+1}$ as the Artin
generators of an abstract braid group yields the isomorphism, since it can be
checked that arcs for the $\s_{i,i+1}$ can be chosen in such a way
that they are
disjoint outside the punctures.

Under this identification the full twists on obstructing parallels are given by
$$
\s_{i,j}^2,\quad i\leq k'<j\leq k,\, (i,j)\neq(1,k).
$$
We can now prove the result concerning the new kind of arcs we want
to consider:

\begin{defi}
Any arc supported on two radial rays and two spiral segments in the
$\eta_2$-neighbourhoods of peripheral cores is called a {\it coiled
isosceles arc}.

Given a bisceles arc it is called the {\it associated coiled isosceles arc}, if
both are isotopic to each other up to full twists on inner parallels.
\end{defi}

\begin{examp}
\labell{ex/a}
Naturally we imagine a coiled arc to spiral monotonously towards the peripheral
cores. For $l_2=6$ and $l_2=4$ the given arc is a coiled isosceles
arc.\\

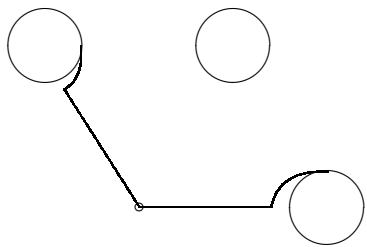
\setlength{\unitlength}{5mm}
\begin{picture}(20,6)(-12,-2)

\put(5,0){\circle{2}}
\put(2.5,4.3){\circle{2}}
\put(-2.5,4.3){\circle{2}}

\put(0,0){\circle{.2}}

\bezier{130}(0,0)(2,0)(3.5,0)
\bezier{180}(0,0)(-1,1.57)(-2,3.14)

\bezier{130}(-2,3.14)(-1.5,3.4)(-1.55,4.3)
\bezier{130}(3.5,0)(3.7,1)(5,.95)

\bezier{2}(-5,0)(-5,1)(-4.5,2)

\end{picture}
\end{examp}

\begin{remark}
By this definition an unobstructed bisceles arc is its own associated
coiled isosceles arc.
\end{remark}

\begin{lemma}
\label{coil-ex}
Given a bisceles arc, there is an associated coiled isosceles arc unique up to
isotopy.
\end{lemma}

\proof
Due to the remark above in case of unobstructed bisceles arcs there is nothing
to prove, because there are no inner parallels.

Otherwise, given a bisceles arc connecting punctures of indices
$i_1i_2,j_1j_2$,
an associated coiled isosceles arc -- if it exists -- must be
isotopic to an arc
supported in the topological disc considered above. But up to isotopy
there is a
unique arc in this disc which is supported in the complement of the peripheral
cores and which connects the same pair of punctures. Hence the uniqueness claim
is proved.

Then we consider the half twists corresponding to the bisceles arc and the arc
just considered. They are identified with $\s_{1,k}$ and
$\check\s_{i,k}$ (as
defined on page \pageref{sig}). Since by \ref{braid/x} they belong to an orbit
under conjugation by full twists on the inner parallels, so do the
corresponding
arcs and existence of an associated coiled isosceles arc is shown.
\qed

 From the simple observation that any side of a bisceles arc may only
cut through
either the central core or a peripheral core we can conclude -- relying on some
results from the appendix on geometry -- that obstructing parallels are in fact
bisceles arcs.

\begin{lemma}
\labell{apex/out}
If a bisceles arc is obstructed then its apex is outside the central
core.
\end{lemma}

\proof
Suppose the apex of a bisceles arc is inside the central core, then both its
sides pass through the central core. But then they are both disjoint to
the peripheral cores, so the bisceles arc is unobstructed.
\qed

\begin{lemma}
\label{bisc}
Each obstructing parallel is a bisceles arc.
\end{lemma}

\proof
Since an obstructing parallel is supported on the lines given by a pair of non
parallel sides, a parallel is not a bisceles arc if and only if the
intersection
point fails to belong to both sides.

Let us assume now that an obstructing parallel to the bisceles arc
is not a bisceles arc itself. Hence by the previous considerations
there must be
a source point in the closed cone defined by the bisceles arc which does not
contain its source points.
With lemma \ref{geo/x} we conclude that the apex of the bisceles arc
is contained
inside the circle through the three source points under consideration, that is
in the central core. This is a contradiction by lemma \ref{apex/out}.
\qed

We prefer to rephrase lemma \ref{coil-ex} using lemma \ref{bisc}.

\begin{prop}
\label{coil}
The class of a bisceles arc and the associated coiled isosceles arc belong
to the same orbit for the action of full twists on bisceles arcs
which are inner
parallels.
\end{prop}

For the closing remark we come back to the topic of straight isosceles.

\begin{remark}
A straight isosceles arc only occurs for $j_1-i_1=l_1/2$ and by a
short check we
see, that the corresponding traces are directing in opposite ways. So they come
close only if they pass the central core. Immediately we deduce, that
a straight
isosceles arc is isotopic to its associated coiled isosceles arc.
\end{remark}


\section{from coiled isosceles arcs to coiled twists}
\label{six/four}

The aim of this section is to identify the isotopy class of the transported arc
at $\l=1$ in terms of the Hefez Lazzeri system of paths.
In fact this system yields a well-defined identification of the mapping class
group of the corresponding fibre with the abstract braid group, so we finally
can even identify the twists on transported arcs with abstract braids.\\

We will see that a coiled isosceles arc transported along a circular segment at
radius $t=1$ is a coiled isosceles arc again, so we have to introduce notations
and definitions in such a way, that we get hold of those geometric properties
which eventually determine the braid associated to a coiled isosceles arc.\\

The fibre at $\l=1$, equal to the fibre at $\a=1$, is a punctured disc for
which Hefez and Lazzeri have given a strongly distinguished system of paths, of
which we should remind ourselves, \ref{hl-path}.

If paths $\oo_{i_1i_2}$, $1\leq i_1\leq l_1,1\leq i_2\leq l_2$, ordered
lexicographically, form the Hefez Lazzeri system, then up to isotopy they can be
obtained from two figures like the following in case $l_1=l_2=8$.

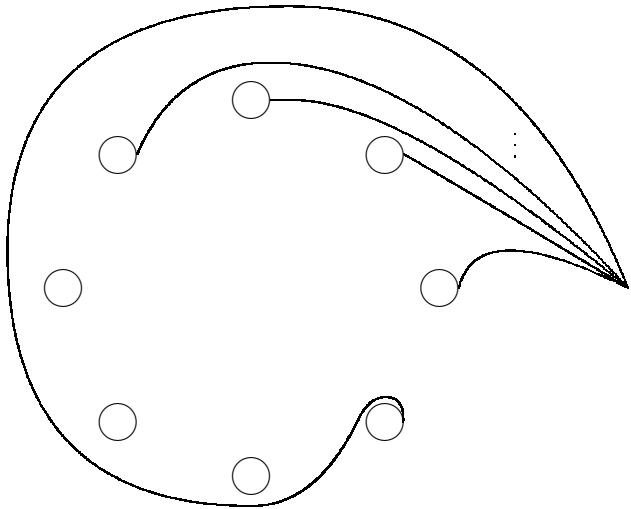
\setlength{\unitlength}{.5mm}
\begin{picture}(180,150)(-100,-70)

\put(50,0){\circle{10}}
\put(0,50){\circle{10}}
\put(0,-50){\circle{10}}
\put(-50,0){\circle{10}}
\put(-35.5,35.55){\circle{10}}
\put(35.55,35.55){\circle{10}}
\put(35.55,-35.55){\circle{10}}
\put(-35.55,-35.55){\circle{10}}

\bezier{2}(70,35)(70,38)(70,41)

\bezier{130}(55,0)(60,20)(100,0)

\bezier{240}(100,0)(70.5,18)(40.55,35.55)

\bezier{30}(5,50)(8,50)(11,50)
\bezier{300}(11,50)(40,50)(100,0)

\bezier{300}(-30.55,35.55)(-20,60)(5,60)
\bezier{300}(5,60)(45,60)(100,0)

\bezier{80}(40.3,-35.55)(40.5,-29)(35.55,-29)
\bezier{80}(35.55,-29)(31,-29)(28,-35.55)

\bezier{400}(28,-35.55)(17,-58)(0,-58)
\bezier{400}(-0,-58)(-65,-58)(-65,10)
\bezier{400}(-65,10)(-65,75)(10,75)
\bezier{400}(10,75)(70,75)(100,0)
\end{picture}

In the first figure a path has to be selected according to $i_1$. It terminates
at a disc (an $\eta_2$-neighbourhood of a peripheral core), which should be
replaced by the second figure.

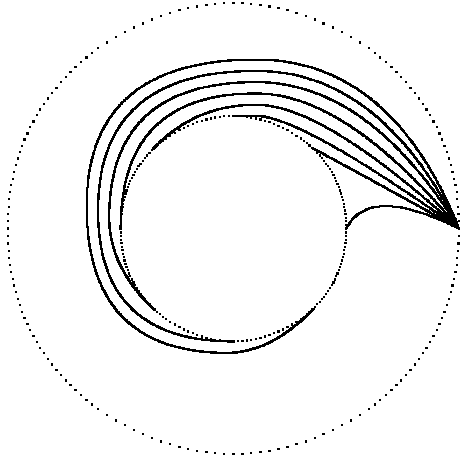
\begin{center}
\setlength{\unitlength}{3mm}
\begin{picture}(30,20)(-15,-10)

\put(-5,-5){\begin{picture}(10,10)
\bezier{30}(0,5)(0,7.5)(2,9)
\bezier{30}(2,9)(5,11)(8,9)
\bezier{30}(8,9)(10,7.5)(10,5)
\bezier{30}(10,5)(10,2.5)(8,1)
\bezier{30}(8,1)(5,-1)(2,1)
\bezier{30}(2,1)(0,2.5)(0,5)
\end{picture}}

\put(-10,-10){\setlength{\unitlength}{6mm}\begin{picture}(10,10)
\bezier{30}(0,5)(0,7.5)(2,9)
\bezier{30}(2,9)(5,11)(8,9)
\bezier{30}(8,9)(10,7.5)(10,5)
\bezier{30}(10,5)(10,2.5)(8,1)
\bezier{30}(8,1)(5,-1)(2,1)
\bezier{30}(2,1)(0,2.5)(0,5)
\end{picture}}

\bezier{130}(5,0)(6,2)(10,0)

\bezier{240}(10,0)(6.8,1.8)(3.55,3.55)

\bezier{30}(0,5)(.5,5)(1,5)
\bezier{300}(1,5)(3,5)(10,0)

\bezier{300}(-3.55,3.55)(-1.7,5.5)(1,5.5)
\bezier{300}(1,5.5)(4,5.5)(10,0)

\bezier{390}(-5,0)(-5,6)(1,6)
\bezier{390}(1,6)(5,6)(10,0)

\bezier{400}(-3.55,-3.55)(-5.5,-1.7)(-5.5,0.5)
\bezier{400}(-5.5,.5)(-5.5,6.5)(1,6.5)
\bezier{400}(1,6.5)(6,6.5)(10,0)

\bezier{400}(0,-5)(-6,-5)(-6,.8)
\bezier{400}(-6,.8)(-6,7)(1,7)
\bezier{400}(1,7)(6.5,7)(10,0)

\bezier{400}(3.55,-3.55)(1.7,-5.5)(-.3,-5.5)
\bezier{400}(-.3,-5.5)(-6.5,-5.5)(-6.5,1)
\bezier{400}(-6.5,1)(-6.5,7.5)(1,7.5)
\bezier{400}(1,7.5)(7,7.5)(10,0)

\end{picture}
\end{center}

The path selected in the second figure according to $i_2$ can be joint to the
first, so they represent the isotopy class of $\oo_{i_1i_2}$.

\begin{nota}
Denote by $\oo\ups{i_1}$ the positive loop around all $\oo_{i_1i_2}$, $1\leq
i_2\leq l_2$.
\end{nota}

Instead of a path segment as given in the second figure we may also join a path
which spirals around the core $n$ full times and then down to a puncture

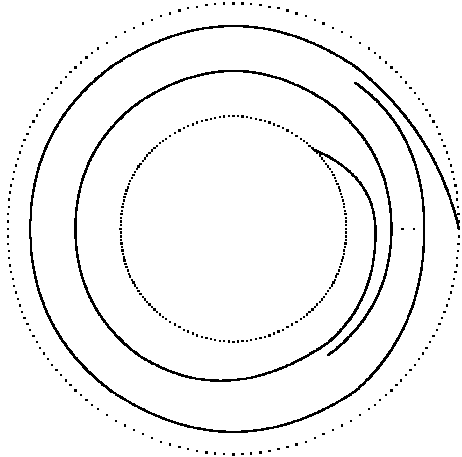
\begin{center}
\setlength{\unitlength}{3mm}
\begin{picture}(30,20)(-15,-10)

\put(-5,-5){\begin{picture}(10,10)
\bezier{30}(0,5)(0,7.5)(2,9)
\bezier{30}(2,9)(5,11)(8,9)
\bezier{30}(8,9)(10,7.5)(10,5)
\bezier{30}(10,5)(10,2.5)(8,1)
\bezier{30}(8,1)(5,-1)(2,1)
\bezier{30}(2,1)(0,2.5)(0,5)
\end{picture}}

\put(-10,-10){\setlength{\unitlength}{6mm}\begin{picture}(10,10)
\bezier{30}(0,5)(0,7.5)(2,9)
\bezier{30}(2,9)(5,11)(8,9)
\bezier{30}(8,9)(10,7.5)(10,5)
\bezier{30}(10,5)(10,2.5)(8,1)
\bezier{30}(8,1)(5,-1)(2,1)
\bezier{30}(2,1)(0,2.5)(0,5)
\end{picture}}

\put(-9,-9){\setlength{\unitlength}{5.4mm}\begin{picture}(10,10)
\bezier{300}(0,5)(0,7.5)(2,9)
\bezier{300}(2,9)(5,11)(8,9)
\bezier{300}(8,8.6)(9.7,7.4)(9.7,5)
\bezier{300}(8,9)(10,7.4)(10.55,5)
\bezier{300}(9.7,5)(9.7,2.5)(8,1)
\bezier{300}(8,1)(5,-1)(2,1)
\bezier{300}(2,1)(0,2.5)(0,5)
\end{picture}}

\put(-7,-7){\setlength{\unitlength}{4.2mm}\begin{picture}(10,10)
\bezier{230}(0,5)(0,7.5)(2,9)
\bezier{230}(2,9)(5,11)(8,9)
\bezier{230}(8,1)(10,2.5)(10,5)
\bezier{230}(8,9)(10,7.5)(10,5)
\bezier{230}(9.5,4.9)(9.5,2.7)(8,1.4)
\bezier{230}(8,1.4)(4.8,-.8)(2,1)
\bezier{230}(2,1)(0,2.5)(0,5)

\bezier{130}(7.52,7.52)(9.5,6.7)(9.5,4.9)
\end{picture}}

\bezier{2}(7,0)(7.5,0)(8,0)

\end{picture}
\end{center}

Such a path is naturally selected by an index $i_2'=i_2+nl_2$ and the notation
$\oo_{i_1i_2}$ is naturally extended to indices $i_1i_2$ with $i_2$
an arbitrary
integer.

\begin{remark}
For $i_1\ne j_1$ the paths $\oo_{i_1i_2},\oo_{j_1j_2}$ do not intersect,
whatever the integers $i_2,j_2$ are.\\
Moreover the loops $\oo\ups{k_1}$ can be chosen disjoint
from both.
\end{remark}

Now we can introduce twist braids corresponding to arcs which are determined by
suitable joins of paths.

\begin{nota}
$\s_{i_1i_2,j_1j_2}$ is the $\frac12$-twist on the union of $\oo_{i_1i_2}$ with
$\oo_{j_1j_2}$.
\end{nota}

\begin{nota}
$\tau_{i_1i_2,j_1j_2}$ is the $\frac12$-twist on the union of
$\oo_{i_1i_2}$ and
$\oo_{j_1j_2}$ with the $\oo\ups{k_1}$, $i_1\leq k_1<j_1$ in between.
\end{nota}

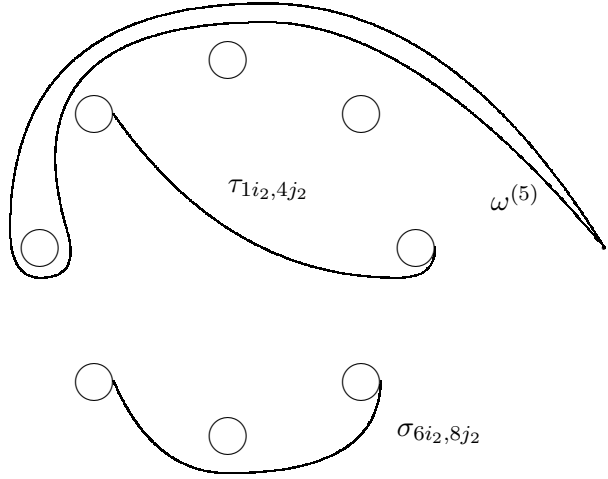
\begin{examp}\hspace*{\fill}\\

\setlength{\unitlength}{.5mm}
\begin{picture}(180,130)(-100,-70)

\put(70,10){$\oo\ups{5}$}
\put(0,15){$\tau_{1i_2,4j_2}$}
\put(45,-50){$\s_{6i_2,8j_2}$}

\put(50,0){\circle{10}}
\put(0,50){\circle{10}}
\put(0,-50){\circle{10}}
\put(-50,0){\circle{10}}
\put(-35.5,35.55){\circle{10}}
\put(35.55,35.55){\circle{10}}
\put(35.55,-35.55){\circle{10}}
\put(-35.55,-35.55){\circle{10}}

\put(100,0){\circle{1}}

\bezier{400}(-30.5,35.5)(0,-8)(45,-8)
\bezier{40}(45,-8)(55,-8)(55,0)

\bezier{240}(-30.5,-35.5)(-20,-60)(0,-60)
\bezier{239}(0,-60)(40,-60)(40.5,-35.5)

\bezier{70}(-43,5)(-39,-8)(-50,-8)
\bezier{390}(-43,5)(-59,60)(10,60)
\bezier{390}(10,60)(50,60)(100,0)

\bezier{70}(-58,5)(-58,-8)(-50,-8)
\bezier{400}(-58,5)(-58,65)(10,65)
\bezier{400}(10,65)(60,65)(100,0)

\end{picture}

\end{examp}

So far we have dwelled on the topology of the fibre at $\l=1$. Now we extract
the characteristic properties of the coiled isosceles.

\begin{defi}
The {\it winding angle} of a directed arc $\Cg$ in the plane of complex numbers
with respect to a disjoint point $z_0$ is -- in generalization of the winding
number of a closed curve -- defined by
\begin{eqnarray*}
\varth_\G  & := & \int_\G\frac{-\mbox{\sf
i}dz}{z-z_0},\quad(\mbox{\sf i}^2=-1).
\end{eqnarray*}
\end{defi}

\begin{nota}
We introduce notation for some characteristic angles:
\begin{enumerate}
\item
$\varth_1,\varth_2$, the angles between consecutive $l_1^{th}$,
resp.\ $l_2^{th}$
roots of unity,
\item
$\vartho:=(j_1-i_1)\varth_1$, the angle at the apex of the coiled isosceles arc
with index pair $i_1i_2,j_1j_2$, note that $0<\vartho<2\pi$,
\item
$\varth_i,\varth_j$, the winding angle of the i-side, resp.\ j-side starting at
the apex, with respect to the center of the core of the corresponding
peripheral
circle,
\item
$\vartho_j:=\varth_j+\vartho$, a useful shorthand.
\end{enumerate}
\end{nota}

The winding angle of a spiral is positive if it turns positively when
approaching the peripheral core.

\begin{examp}
In the example considered before, suppose the horizontal line supports the
$i$-side then
$\varth_i=-\frac\pi2,\varth_j=\frac\pi3$,
otherwise $\varth_i=\frac\pi3,\varth_j=\frac{\pi}6$.\\

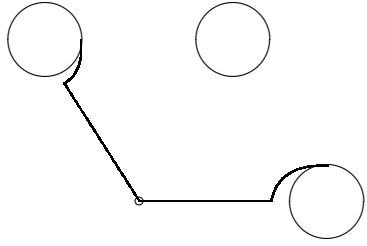
\setlength{\unitlength}{5mm}
\begin{picture}(20,6)(-12,-2)

\put(5,0){\circle{2}}
\put(2.5,4.3){\circle{2}}
\put(-2.5,4.3){\circle{2}}

\put(0,0){\circle{.2}}

\bezier{130}(0,0)(2,0)(3.5,0)
\bezier{180}(0,0)(-1,1.57)(-2,3.14)

\bezier{130}(-2,3.14)(-1.5,3.4)(-1.55,4.3)
\bezier{130}(3.5,0)(3.7,1)(5,.95)

\bezier{2}(-5,0)(-5,1)(-4.5,2)

\end{picture}
\end{examp}

We want now to pin down some geometric properties shared by the
coiled isosceles
arcs associated to bisceles arcs or straight isosceles arcs.

\begin{lemma}
The winding angle of a side of a coiled isosceles arc is in the open interval
$]-\frac{3\pi}2,\frac{3\pi}2[$.
\end{lemma}

\proof
The side of the bisceles arc is parallel to a side of the associated coiled
isosceles arc.
If the endpoint is on the half of the peripheral circle facing the origin, then
the side is unobstructed and the winding angle is therefore in the range
$[-\frac\pi2,\frac\pi2]$.

Otherwise it may be obstructed and there are two ways to make it unobstructed
depending on the other side. But in any case the absolute value of the winding
angle does not exceed $\frac{3\pi}2$.
\qed

\begin{lemma}
\labell{angle/coil}
The following inclusions hold:
\begin{eqnarray*}
\text{if }\jith \leq\pi: &\varth_i  \in & ]-\frac{3\pi}2,\frac\pi2],\\
&\varth_j  \in & [ -\frac\pi2,\frac{3\pi}2[;\\
\text{if }\jith >\pi: &\varth_i  \in & [-\frac\pi2,\frac{3\pi}2[,\\
& \varth_j  \in & ]-\frac{3\pi}2,\frac\pi2].
\end{eqnarray*}
\end{lemma}

\proof
If the endpoint of the $i$-side is on the half circle facing the origin, then
its winding angle is in $[-\frac\pi2,\frac\pi2]$.

If the endpoint is on the opposite half circle, then the winding angle is
in either $]-\frac{3\pi}2,-\frac\pi2[$ or $]\frac\pi2,\frac{3\pi}2[$ and the sign
depends on the second endpoint. The sign is that of
$\pi-\jith$ for the $i$-side and the opposite for the $j$-side.
\qed

Moreover the considerations of this proof immediately yield the observation:

\begin{lemma}
\labell{obstr}
The $j$-side of the bisceles arc is disjoint from the central core if
and only if
\begin{eqnarray*}
\text{either }\quad\jith \leq\pi &\text{ and}\quad\varth_j  \in &
]\frac\pi2,\frac{3\pi}2[,\\
\text{or }\quad\jith >\pi &\text{ and}\quad\varth_j  \in &
]-\frac{3\pi}2, -\frac\pi2[.
\end{eqnarray*}
\end{lemma}

The next thing we have to exploit is the fact that at a parameter $\l=\ethe$
bisceles do not exists for all index pairs.

\begin{lemma}
\labell{sine}
A bisceles arc with index pair $i_1i_2,j_1j_2$ exists at $\l=\ethe$ only if
\begin{enumerate}
\item
$i_2=j_2$,
\item
$\jith \leq\pi$ and $\sin \varth_i\leq \sin\vartho_j$,
$\vartho_j<\frac{3\pi}2$,
\item[or]
\item
$\jith >\pi$ and $\sin \varth_i\geq \sin\vartho_j$, $\vartho_j>\frac\pi2$.
\end{enumerate}
\end{lemma}

\proof
In the first case the claim is obvious since the traces have their source
points in common. So from now on we assume that the source points are
distinct.\\
Let us consider the case $\jith \leq\pi$ next. Then the possible
traces for the index $i_1i_2$ are sketched near to the central core as well as
the direction of possible traces with index $j_1j_2$. The second inequality
is now read off easily, since $ \sin \varth_i$ is the vertical component
of the $i$-side and $\sin\vartho_j$ the maximal vertical component of
the $j$-side.
\\

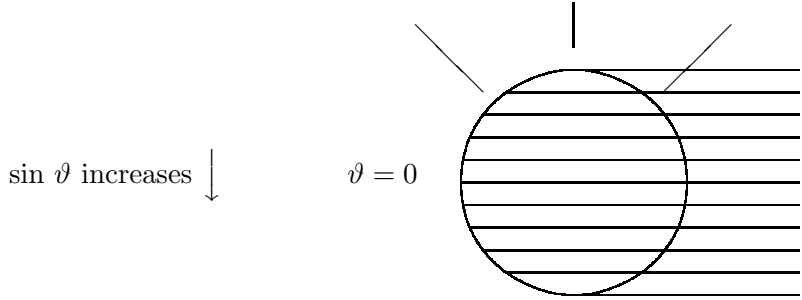
\setlength{\unitlength}{3mm}
\begin{picture}(30,17)(-20,-2)

\bezier{300}(0,5)(0,7.5)(2,9)
\bezier{300}(2,9)(5,11)(8,9)
\bezier{300}(8,9)(10,7.5)(10,5)
\bezier{300}(10,5)(10,2.5)(8,1)
\bezier{300}(8,1)(5,-1)(2,1)
\bezier{300}(2,1)(0,2.5)(0,5)

\drawline(0,5)(15,5)
\drawline(0.1,6)(15,6)
\drawline(0.4,7)(15,7)
\drawline(1,8)(15,8)
\drawline(2,9)(15,9)
\drawline(0.1,4)(15,4)
\drawline(0.4,3)(15,3)
\drawline(1,2)(15,2)
\drawline(2,1)(15,1)
\drawline(5,10)(15,10)
\drawline(5,0)(15,0)

\drawline(5,11)(5,13)
\drawline(1,9)(-2,12)
\drawline(9,9)(12,12)

\put(-20,5){sin $\varth$ increases $\Big\downarrow$}

\put(-5,5){$\varth=0$}

\end{picture}

Suppose now $\vartho_j$ exceeds $\frac{3\pi}2$, then $\varth_j$ exceeds
$\frac{\pi}{2}$ and by \ref{obstr} the $j$-side does not pass the central core.
To be cut properly by the trace of the $i$-side its source point must be on the
right hand half of the circle.
But the horizontal component of $\vartho_j$ is $\cos \vartho_j$ which is not
positive for $\vartho_j\in[\frac{3\pi}2,\jith +\frac{3\pi}2]$.

The final case $\jith >\pi$ can be handled in strict analogy.
\qed

Since these better bounds hold obviously in the case of straight isosceles we
get an improvement on the assertion of \ref{angle/coil}:

\begin{lemma}
\labell{angle/range}
The following inclusions hold:
\begin{eqnarray*}
\text{if }\quad\jith \leq\pi: &\varth_i  \in & ]-\frac{3\pi}2,\frac\pi2],\\
&\varth_j  \in & [-\frac\pi2,\frac{3\pi}2-\jith[;\\
\text{if }\quad\jith >\pi: &\varth_i  \in & [-\frac\pi2,\frac{3\pi}2[,\\
& \varth_j  \in & ]\frac\pi2-\jith,\frac\pi2].
\end{eqnarray*}
\end{lemma}

Now we combine the results to obtain a relation between the winding angles.

\begin{lemma}
\labell{angle/rel}
The winding angles are subject to
$$
\varth_i\leq\vartho_j\leq \varth_i+2\pi.
$$
\end{lemma}

\proof
Suppose $\jith >\pi$. Then by lemma \ref{angle/range}
$\vartho_j\in]\frac\pi2,\frac\pi2+\vartho]$ so
$$
i)\quad\varth_i\leq\vartho_j\qquad \text{or} \qquad
ii)\quad\varth_i,\vartho_j\in]\frac\pi2,\frac{3\pi}2[.
$$
Also the conditions in the
latter case imply $\varth_i\leq\vartho_j$, since $\sin \varth_i\geq
\sin\vartho_j$ by lemma \ref{sine} and the sine function is decreasing in
$]\frac\pi2,\frac{3\pi}2[$.

On the other hand by lemma \ref{angle/range}
$$
i)\quad\varth_i+2\pi\geq\vartho_j\qquad \text{or} \qquad
ii)\quad\varth_i+2\pi,\vartho_j\in]\frac{3\pi}2,\frac{5\pi}2[.
$$
and again the second case is a subcase of the first, since the sine function
is increasing on $]\frac{3\pi}2,\frac{5\pi}2[$.

The case $\jith \leq\pi$ is done analogously.
\qed

Next we investigate the impact of parallel transport.

\begin{lemma}
\labell{circ/trans}
Under parallel transport along a circle segment of winding angle $\varth$
at radius $t=1$ a coiled isosceles arc is mapped up to isotopy to a coiled
isosceles arc with winding angles changed \mbox{by $-\varth$}.
\end{lemma}

\proof
The line segments of the isosceles arc belong to the part which is
rotated rigidly
by the flow of the vector field in \ref{vect/circ}. The total rotation is of
angle $\varth$.\\
On the other hand the spirals are wound resp.\ unwound, since the endpoints are
fixed relative to their peripheral centres, while the points on the boundary of
the $2\eta_2$-discs are relatively rotated in opposite direction, hence the
amount and sign of the change in the winding angles.
\qed

\begin{remark}
If we introduce $\varth_{i}':=\varth_i-\varth$, (similarly
${\varth_{j}}:=\varth_j-\varth$), then $\varth_i$ is the $i$-side winding angle
of the isosceles arc transported from angular parameter $\varth$ to $\l=1$.
\end{remark}

\begin{lemma}
\labell{crit-twist}
Suppose a coiled isosceles arc is associated to a bisceles arc or to an
isosceles arc with index pair
$i_1i_2,j_1j_2$, $i_1<j_1$, $i_2\neq j_2$, then the full twist on any of its
parallel transports to $\l=1$ along a circle
segment of radius $t=1$ is identified with one of the abstract braid elements
$$
\tau_{i_1i_2',j_1j_2'}^2,\,1\leq i_1<j_1\le l_1,\,1\leq j_2'-i_2'< l_2.
$$
\end{lemma}

\proof
The transported coiled isosceles arc at $\l=1$ can be represented in a unique
way by the join of loops $\oo\ups{k_1}, i_1\leq k_1<j_1$ and two paths
$\oo_{i_1i_2'}$, $\oo_{j_1j_2'}$ with $i_2',j_2'$ suitable chosen. Hence the
corresponding half twist is identified with the abstract braid element
$\tau_{i_1i_2',j_1j_2'}$.

We note further that
\begin{eqnarray*}
(j_2'-1)\varth_2 & = & \varth_j'-\pi+(j_1-1)\varth_1,\\
(i_2'-1)\varth_2 & = & \varth_i'-\pi+(i_1-1)\varth_1.
\end{eqnarray*}
Computing the difference using $\varth_j'-\varth_i'=\varth_j-\varth_i$ and
$(j_1-i_1)\varth_1=\vartho$ we get:
\begin{eqnarray*}
(j_2'-i_2')\varth_2 & = & \varth_j-\varth_i+\vartho\,\;=\,\;\vartho_j-\varth_i.
\end{eqnarray*}
In case $j_2'-i_2'\leq0$ this implies $\vartho_j-\varth_i\leq0$, in case
$j_2'-i_2'\geq l_2$ we conclude $2\pi\leq\vartho_j-\varth_i$, so both these
cases contradict the assertion of lemma \ref{angle/rel}, since neither
$\vartho_j=\varth_i$ nor $\vartho_j=\varth_i+2\pi$ is possible under the
assumption $i_2\neq j_2$.
Therefore we get $1\leq j_2'-i_2'< l_2$, as claimed.
\qed

\begin{examp}
Suppose the example from page \pageref{ex/a} has been transported by an
angle $\varth=\frac{5\pi}3$ along the circle arc of radius $t=1$. Following
the recipe above we get:

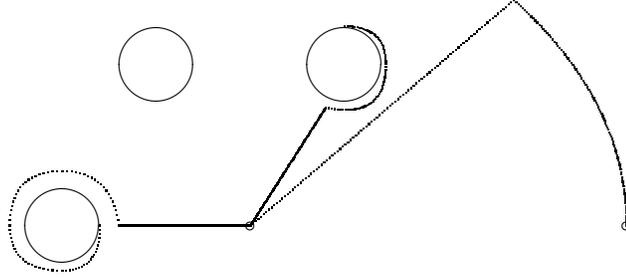
\setlength{\unitlength}{5mm}
\begin{picture}(20,6)(-12,-1)

\put(-5,0){\circle{2}}
\put(2.5,4.3){\circle{2}}
\put(-2.5,4.3){\circle{2}}

\put(0,0){\circle{.2}}
\put(10,0){\circle{.2}}

\bezier{100}(0,0)(3.5,3)(7,6)
\bezier{100}(7,6)(10,3.3)(10,0)


\bezier{130}(0,0)(-2,0)(-3.5,0)
\bezier{180}(0,0)(1,1.57)(2,3.14)

\bezier{20}(-3.5,0)(-3.7,1.45)(-5,1.45)
\bezier{20}(-5,1.45)(-6.4,1.45)(-6.4,0)
\bezier{20}(-6.4,0)(-6.4,-1.2)(-5,-1.2)
\bezier{20}(-5,-1.2)(-4,-1.2)(-4,0)

\bezier{30}(2.5,5.3)(3.6,5.3)(3.6,4.3)
\bezier{39}(3.6,4.3)(3.6,2.8)(2,3.14)

\end{picture}

Hence assuming $l_1=6,l_2=4$ the associated abstract braid is $\tau_{22,45}$.
\end{examp}

Let us call a coiled isosceles arc in the fibre at $\l=1$ {\it associated} to a
local v-arc, if it is obtained from the local v-arc by parallel transport along
a radial segment, a transformation by full twists to get the associated coiled
isosceles arc and parallel transport along a circle segment at $t=1$.

We note then the following converse to \ref{crit-twist}:

\begin{lemma}
\labell{surj}
Each element $\tau^2_{i_1i_2,j_1j_2}$, $1\leq i_1<j_1\leq l_1$, $1\leq
j_2-i_2<l_2,$ is the full twist on a coiled isosceles arc associated to a local
v-arc.
\end{lemma}

\proof
There is a local v-arc for each index pair $i_1i_2,j_1j_2$, $1\leq
i_1<j_1\leq l_1$, $1\leq i_2,j_2\leq l_2$. For each such local v-arc
there is an
associated coiled isosceles arc in the fibre at $\l=1$, which
determines some $\tau$ as above by \ref{crit-twist}. All
others are then obtained by changing the winding angle of the circular path
by suitable multiples of $2\pi$.
\qed

The case $i_2=j_2$ requires extra care. We have analogues to lemma
\ref{crit-twist} and lemma \ref{surj}.

\begin{lemma}
\labell{cri-twist}
Suppose a coiled isosceles arc is associated to a bisceles arc or to an
isosceles arc with index pair
$i_1i_2,j_1j_2$, $i_1<j_1$, $i_2=j_2$, then the full twist on any of its
parallel transports to $\l=1$ along a circle
segment of radius $t=1$ is identified with one of the abstract braid elements
$$
\tau_{i_1i_2',j_1j_2'}^2,\,1\leq i_1<j_1\le l_1,\,j_2'-i_2'\in\{0, l_2\}.
$$
\end{lemma}

\proof
With the same argument as in the proof of lemma \ref{crit-twist} we can exclude
the cases $j_2'-i_2'<0$ and $j_2'-i_2'>l_2$. Since $i_2,i_2'$ and
$j_2,j_2'$ may
only differ by multiples of $l_2$ we are left with the two possibilities of the
claim.
\qed

\begin{lemma}
\labell{sur}
Given an index pair $i_1i_2,j_1j_2$, $1\leq i_1<j_1\leq l_1$, $j_2-i_2=l_2$,
at least one of $\tau_{i_1i_2,j_1i_2},\tau_{i_1i_2,j_1j_2}$ is the
half twist on
a coiled isosceles arc associated to a local v-arc.
\end{lemma}

\proof
There is an local v-arc for the index pair $i_1i_2,j_1i_2$, $1\leq
i_1<j_1\leq l_1$, $1\leq i_2\leq l_2$.
We know that for each such local v-arc there is an
associated coiled isosceles arc in the fibre at $\l=1$, which
determines one of $\tau_{i_1i_2',j_1i_2'},\tau_{i_1i_2',j_1j_2'}$ as above by
\ref{cri-twist}. All others are then obtained by changing the winding angle of
the circular path by suitable multiples of $2\pi$.
\qed

We close this section be identifying the twists unambiguously under special
geometric assumptions.

\begin{lemma}
\labell{iso/zero}
Given any coiled isosceles arc with punctures of indices $i_1i_2,j_1j_2$,
$i_1<j_1$, $i_2=j_2$, facing the origin, then
the twist on any of its parallel transports to $\l=1$ along a circle segment of
radius $t=1$ is identified with one of the abstract braid elements
$$
\tau_{i_1i_2',j_1j_2'},\,1\leq i_1<j_1\le l_1,
$$
with $i_2'=j_2'$ if $\jith \leq\pi$ and $i_2'+l_2=j_2'$ if $\jith >\pi$.
\end{lemma}

\proof
We run through the same consideration as in \ref{crit-twist}. But in the final
step we are stuck since now $\varth_i=\vartho_j \mod 2\pi$, so
$\varth_i=\vartho_j$ and $\varth_i=\vartho_j-2\pi$ are possible
by \ref{angle/rel}.
Since both punctures face the origin by hypothesis,
$$
\varth_i\in[-\frac\pi2,\frac\pi2],
\varth_j\in[-\frac\pi2,\frac\pi2].
$$
In case $\varth_i=\vartho_j$ which corresponds to $i_2'=j_2'$ we must have
$\vartho=\varth_i-\varth_j\leq\pi$.
In case $\varth_i+2\pi=\vartho_j$ corresponding to $i_2'+l_2=j_2'$ we conclude
that $\vartho=2\pi+\varth_i-\varth_j\geq\pi$.
\qed

\begin{lemma}
\labell{iso/zero/ob}
Given any coiled isosceles arc associated to a
bisceles arc with punctures of indices $i_1i_2,j_1j_2$, $i_1<j_1$, $i_2=j_2$,
one of which exactly facing the origin, then the half twist on any of its
parallel transports to $\l=1$
along a circle segment of radius $t=1$ is identified with one of the abstract
braid elements
$$
\tau_{i_1i_2',j_1j_2'},\,1\leq i_1<j_1\le l_1,i_2'+l_2=j_2',
$$
under the assumption that $\jith \leq\pi$.
\end{lemma}

\proof
The claim is secured by similar considerations as in the proof of
\ref{iso/zero}. We know that $\varth_i=\vartho_j$ or $\varth_i=\vartho_j+2\pi$
and we imposed $\jith \leq\pi$.\\
If the puncture of index $i_1i_2$ faces the origin, then from
\ref{angle/coil} and \ref{obstr}:
$$
\varth_i\in[-\frac\pi2,\frac\pi2],
\varth_j\in]\frac\pi2,\frac{3\pi}2[.
$$
If the puncture of index $i_1i_2$ does not face the origin, then we get:
$$
\varth_i\in]-\frac{3\pi}2,-\frac\pi2[,
\varth_j\in[-\frac\pi2,\frac{\pi}2].
$$
In either case we can check that we are left with the possibility
$\varth_i=\vartho_j-2\pi$. Hence $j_2'=i_2'+l_2$ holds in the index pair of the
corresponding abstract braid.
\qed


\section{from local w-arcs to coiled twists}
\label{six/five}

Having understood the parallel transport of local v-arcs sufficiently well,
we can now consider the parallel transport of local w-arcs. They are only
considered close to the degeneration at $\l=0$ with $i_2=j_2$.\\

Let us first look at parallel transport along circular segments of very small
radius.

\begin{lemma}
Under parallel transport along circle segments of radius $\e<<\eta_2$
local w-arcs
are mapped to local w-arcs.
\end{lemma}

\proof
This is immediate, for the parallel transport can be realised by the
flow of the
vector field in \ref{vect/core}, which is rigid rotation for the support of the
local w-arcs.
\qed

Next local w-arcs are transported along a radial segment. We get then tangled
w-arcs in the fibre $\ethe$.
The simpler arcs, to which we want to compare them, are called {\it isosceles
w-arcs} and they relate to isosceles arcs as local w-arcs relate to
local v-arcs.

The isosceles w-arc of index pair $i_1i_2,j_1i_2$ can be best understood from
the isosceles arcs of index pairs $i_1i_2,i_1^+i_2$ and
$i_1^+i_2,j_1i_2$, which
are called its {\it constituents}. It is isotopic to the first
constituent acted on by
a positive half twist on the second. It can be chosen to be composed
of four line
segments, two of which are supported on the traces of the punctures
$i_1i_2$ and
$j_1i_2$, while the middle pair forms a sharp wedge over the trace of
the puncture
$i_1^+i_2$, cf.\ the example below.

An isosceles w-arc is called {\it corresponding} to a given tangled
w-arc, if both
connect the same pair of punctures.

The same methods as in the case of tangled v-arcs can now be employed to relate
tangled and isosceles w-arcs.

\begin{lemma}
Up to conjugation by full twists on bisceles arcs of shorter length a
tangled w-arc is isotopic to the corresponding isosceles w-arc.
\end{lemma}

We now make an observation which will help us to be concerned mostly with
isosceles w-arcs which are supported outside the peripheral circles except for
an arbitrarily small neighbourhood of the puncture $i_1^+i_2$. They shall be
referred to as {\it unobstructed isosceles w-arcs}.

\begin{lemma}
Any local w-arc for the index pair $i_1i_2,j_1i_2$ can be transported to
radius $t=1$ along a circle segment of radius $t=\e$ and a radial segment such
that the corresponding isosceles w-arc is unobstructed, except in case of
$j_1-i_1=(l_1+1)/2$.
\end{lemma}

\proof
In the cases under consideration either $j_1-i_1\leq l_1/2$ or
$j_1-i^+_1>l_1/2$.
We choose $\varth=(i_2-1)\varth_2-(j_1-1)\varth_1\mp\frac\pi2$
respectively, so we get
\begin{eqnarray*}
|(i_2-1)\varth_2-(k_1-1)\varth_1-\varth)|\geq\frac\pi2 & \text{ for }
k_1=i_1,i_1^+,j_1.
\end{eqnarray*}
Therefore at $\l=\ethe$ the punctures with index pairs $i_1i_2,i_1^+i_2,j_1i_2$
are all situated on the halfs of their peripheral circles facing the origin.
Accordingly the isosceles w-arc corresponding to the transported local w-arc is
unobstructed.
\qed

In the remaining case we can only arrange that the punctures of the $i$-side and
of the $j$-side face the origin.

\begin{lemma}
If $j_1-i_1=(l_1+1)/2$, then the local w-arc can be transported to
radius $t=1$ along a circle segment of radius $t=\e$ and a radial segment such
that
\begin{enumerate}
\item
only the wedge of the corresponding isosceles w-arc is obstructed,
\item
every critical point of its peripheral circle belongs to the inner cone
of either of the constituents or is of index $i_1^+i_2$.
\end{enumerate}
\end{lemma}

\proof
We choose $\varth=(i_2-1)\varth_2-(j_1-1)\varth_1+\frac\pi2$ and get
\begin{eqnarray*}
|(i_2-1)\varth_2-(k_1-1)\varth_1-\varth)|\geq\frac\pi2\phantom{.}
& \text{ for }
k_1=i_1,j_1,\\
|(i_2-1)\varth_2-(i^+_1-1)\varth_1-\varth)|<\frac\pi2.
\end{eqnarray*}
So the punctures of the $i$-side and of the $j$-side of the corresponding
isosceles w-arc face the origin as before, but the
wedge is obstructed. Since both inner angles are less than $\pi$, all
critical points on the corresponding peripheral circle belong to an
inner cone, except the puncture of index $i_1^+i_2$.
\qed

\begin{examp} An isosceles w-arc with obstructed wedge is
obtained in
case of $l_1=3,l_2=4$, $i_2=2$:\\

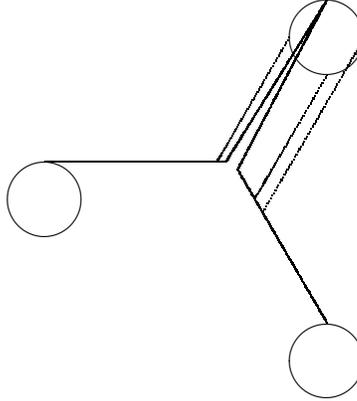
\setlength{\unitlength}{5mm}
\begin{picture}(20,9)(-16,-5)

\put(-5,0){\circle{2}}
\put(2.5,4.3){\circle{2}}
\put(2.5,-4.3){\circle{2}}

\bezier{200}(-5,1)(-2.5,1)(-.18,1)
\bezier{200}(0.125,0.785)(1.25,-1.15)(2.5,-3.3)

\bezier{200}(-.18,1)(1.16,3.15)(2.5,5.3)
\bezier{200}(0.125,0.785)(1.31,3.05)(2.5,5.3)

\bezier{60}(1.5,4.3)(.25,2.15)(-.42,1)
\bezier{70}(3.5,4.3)(2.25,2.15)(.8,-.34)
\bezier{60}(2.5,3.3)(1.25,1.15)(.58,0)


\end{picture}
\end{examp}

The final parallel transport of an isosceles w-arc along a circular segment at
radius $t=1$ can be understood using its constituents.

\begin{lemma}
\label{free}
Parallel transport along a circle segment of radius $t=1$ of an unobstructed
isosceles w-arc yields an arc isotopic to the parallel transport of its first
constituent acted upon by a positive half twist on the parallel transport of
its second constituent.
\end{lemma}

\proof
The relation between an isosceles w-arc and its constituents is preserved
under parallel transport.
\qed

\begin{lemma}
\label{wedge}
Up to full twists on bisceles arcs of shorter length an isosceles w-arc
obstructed on its wedge only is isotopic to the coiled isosceles arc associated
to its first constituent acted upon by a positive half twist on the coiled
isosceles arc associated to its second constituent.
\end{lemma}

\proof
The same full twists on inner parallels which map the constituents to
the isotopy
classes of their associated coiled isosceles arcs also maps the isosceles w-arc
to the isotopy class of the arc obtained from the associated coiled isosceles
arcs.
\qed

Let us finally summarize the results of this section:

\begin{lemma}
\labell{tri/twist}
Local w-arcs in a fibre close to the origin and twists among the following
elements with $1\leq i_1<j_1\leq l_1$, $i_2'=i_2-l_2$,
\begin{eqnarray*}
\tau_{i_1i_2',i^+_1i_2'}\inv\tau_{i^+_1i_2',j_1i_2'}^2\tau_{i_1i_2',i^+_1i_2'},
&& j_1-i_1\leq l_1/2,\\
\tau_{i_1i_2',i^+_1i_2'}\inv\tau_{i^+_1i_2',j_1i_2}^2\tau_{i_1i_2',i^+_1i_2'},
&& j_1-i_1\geq l_1/2+1,\\
\tau_{i_1i_2',i^+_1i_2}\inv\tau_{i^+_1i_2',j_1i_2}^2\tau_{i_1i_2',i^+_1i_2},
&& j_1-i_1=l_1/2+1/2.
\end{eqnarray*}
correspond in such a way that
\begin{enumerate}
\item
each local w-arc can be transported along a circle arc of radius $\e$ and a
radial segment to $t=1$, such that the twist on the corresponding isosceles
w-arc transports to $\l=1$ along the circle of radius $t=1$ to yield one of the
given twists up to conjugation by full twists on obstructing parallels to its
constituents.
\item
each given twist can be obtained from a local w-arc as in i).
\end{enumerate}
\end{lemma}

\proof
If $j_1-i_1\leq l_1/2$ then $(j_1-i_1^+)\varth_1,(i_1^+-i_1)\varth_1\leq\pi$,
hence by
\ref{iso/zero} and \ref{free} we can get an element of the first row,
since by the
braid relation it does not matter if we transform the first constituent by a
positive full twist on the second or if we transform the second by a negative
full twist on the first.

Similarly if $j_1-i_1^+>l_1/2$ then $(i_1^+-i_1)\varth_1\leq\pi$ but
$(j_1-i_1^+)\varth_1>\pi$, so we get a twist of the second row, again with
\ref{iso/zero} and \ref{free}.

In the final case we argue along the same line with \ref{iso/zero/ob} and
\ref{wedge}, so also
in case $j_1-i_1=l_1/2+1/2$ we get twists among the given ones.

As in the similar cases proved before, we get all twist this way as we can
transport around the circle at $t=1$ as many times as necessary.
\qed[2mm]


\section{the length of bisceles arcs}
\label{six/six}

In this section we want to compare the length of bisceles arcs to a real number
we assign to index pairs.

\begin{defi}
The modulus of a pair $i_1i_2,j_1j_2$ of indices is given by
$$
\eta_2\left|\frac{\sin(\pi\frac{i_2-j_2}{l_2})}
{\sin(\pi\frac{i_1-j_1}{l_1})}\right|.
$$
In this way a modulus is assigned to all objects with an index pair.
\end{defi}

Since modulus is in some way complementary to length, we introduce it also for
bisceles arcs.

\begin{defi}
The modulus of a bisceles arc is the shorter of the two distances
from the apex to
both source points.
\end{defi}

\begin{lemma}
The modulus $t_0$ of a critical parameter $t_0\etheO$ for the pair
$i_1i_2,j_1j_2$ coincides with the modulus for that index pair.
\end{lemma}

\proof
Given the traces at angle $\varth_0$ the pair corresponding to $i_1i_2,j_1j_2$
meet at an apex which forms an isosceles triangle with both source points on
the circle of radius $\eta_2$. So with $\delta=\pm\pi\frac{i_2-j_2}{l_2}$,
$\phi=\pm2\pi\frac{i_1-j_1}{l_1}$ the length of the sides equals the modulus as
can be seen from the following sketch.

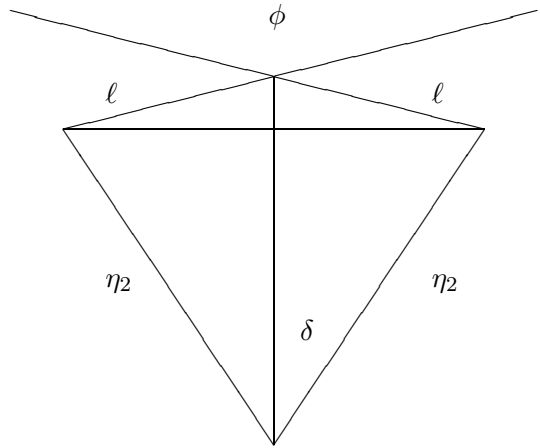
\begin{center}
\setlength{\unitlength}{7mm}
\begin{picture}(10,10)

\put(5,1){\line(2,3){4}}
\put(5,1){\line(-2,3){4}}
\put(1,7){\line(4,1){9}}
\put(9,7){\line(-4,1){9}}
\put(1,7){\line(1,0){8}}
\put(5,1){\line(0,1){7}}
\put(1.8,4){$\eta_2$}
\put(8,4){$\eta_2$}
\put(5.5,3){$\delta$}
\put(1.8,7.5){$\ell$}
\put(8,7.5){$\ell$}
\put(4.9,9){$\phi$}

\end{picture}
\end{center}
\qed

\begin{lemma}
\labell{mod/bnd}
The modulus of a bisceles arc bounds the modulus of the corresponding
index pair
from below. Equality holds only in the case that the bisceles arc is an
isosceles arc.
\end{lemma}

\proof
The apex of the bisceles arc which depends on the parameter angle $\varth$
determines
a triangle over the base
given by the two source points. The base and the angle over it
are independent of $\varth$, whereas the length of the shorter side is
the bisceles arc modulus. The modulus of the pair is the length of a
side if both
sides are equal which happens for a specific $\varth$.\\
The claim is now obvious from the following sketch, $m_b\geq
m_2=\min(m_1,m_2)$:

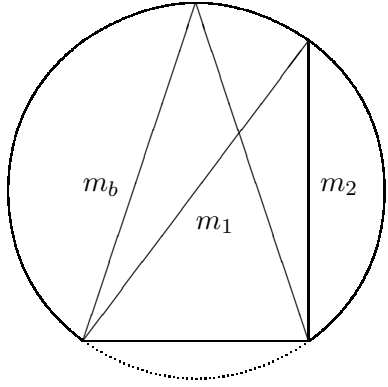
\begin{center}
\setlength{\unitlength}{5mm}
\begin{picture}(10,10)

\bezier{300}(0,5)(0,7.5)(2,9)
\bezier{300}(2,9)(5,11)(8,9)
\bezier{300}(8,9)(10,7.5)(10,5)
\bezier{300}(10,5)(10,2.5)(8,1)
\bezier{50}(8,1)(5,-1)(2,1)
\bezier{300}(2,1)(0,2.5)(0,5)

\put(2,1){\line(1,0){6}}
\put(2,1){\line(1,3){3}}
\put(8,1){\line(-1,3){3}}
\put(2,1){\line(3,4){6}}
\put(8,1){\line(0,1){8}}

\put(8.3,5){$m_2$}
\put(2,5){$m_b$}
\put(5,4){$m_1$}

\end{picture}
\end{center}
The algebraic argument reads as follows: By the cosine formula
$$
m_1^2+m_2^2-2m_1m_2\cos(apex)=2m_b^2-2m_b^2\cos(apex).
$$
We can get a lower bound for the l.h.s. assuming w.l.o.g.\ $m_1\leq m_2$:
\begin{eqnarray*}
m_1^2+m_2^2-2m_1m_2\cos(apex) & = &
(m_1-m_2)^2+m_1m_2(2-2\cos(apex))\\
& \geq & m_1^2(2-2\cos(apex))
\end{eqnarray*}
Then the conclusion $m_b\geq m_1$ is immediate.
\qed

Now we compare the modulus of arcs we encountered in preceding sections.

\begin{lemma}
\labell{ent/modu}
The modulus of a bisceles arc supported on the essential traces of a tangled
v-arc is strictly larger than the modulus of the corresponding isosceles arc.
\end{lemma}

\proof
This claim follows from lemma \ref{mod/bnd} above and the remark on page
\pageref{shorter}.
\qed

\begin{lemma}
\labell{ob/iso/modu}
An obstructing parallel to a bisceles arc is of strictly larger modulus.
\end{lemma}

\proof
Let us consider first the case that the obstructing parallel has a side
in common with the obstructed bisceles arc:\\
We have thus a triangle $ABC$ formed by the source points $A,C$ of
the traces of
obstructed bisceles arc and its apex $B$.
Similarly we have a triangle $AED$ formed by the source points $A,D$ of the
obstructing parallel and its apex $E$. We have $g_{AE}=g_{AB}$ and
$B\in\overline{AE}$. Moreover $g_{DE}\|g_{BC}$ and $D$ is separated from $A$ by
$g_{BC}$.
Denote by $F$ the intersection of $g_{DE}$ and $g_{AC}$. Then $g_{DE}$ is
divided into rays bounded by $E$ resp.\ $F$ and the finite segment
$\overline{EF}$.\\

Now $D$ may not be on the ray bounded by $E$, since then $B$ is in the interior
of $ACD$, hence in the central core contrary to the assumption on
obstructedness.
Neither may $D$ belong to $\overline{EF}$ since otherwise $\overline{BD}$ cuts
$\overline{BC}$ which is impossible since $\overline{BC}$ is on the obstructed
side of the bisceles arc and may hence not be cut by the chord
$\overline{BD}$ of
the central core.\\
So finally $F$ is on $\overline{DE}$ and we conclude with lemma
\ref{leng/inc}.\\

Suppose now that the obstructing parallel has no side in common with the
obstructed bisceles arc, then there is an intermediate obstructing
parallel which
has a side in common with each. So the full result is obtained in two steps as
above.

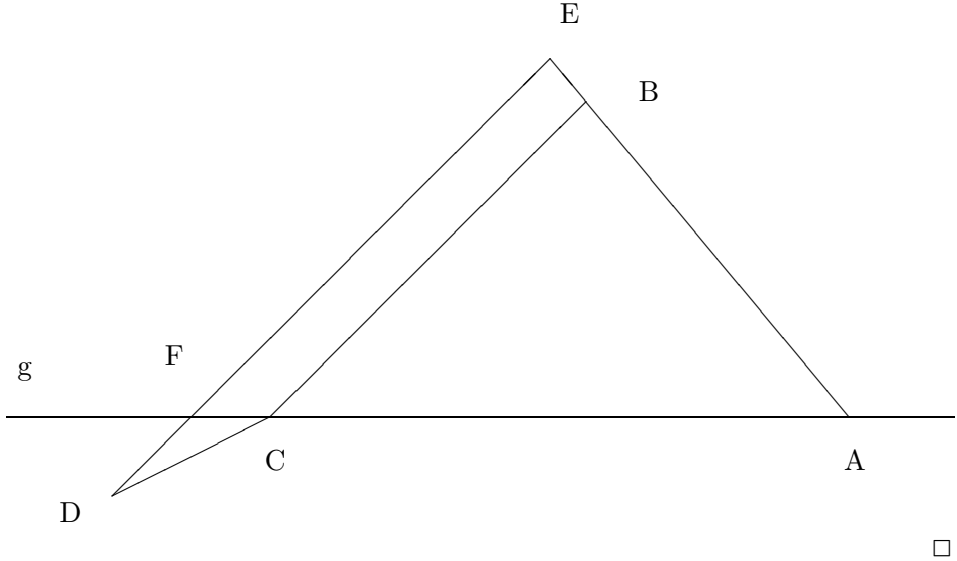
\setlength{\unitlength}{7mm}   
\begin{picture}(20,10)(0,1)

\put(1,3){\line(1,0){18}}
\put(6,3){\line(1,1){6}}
\put(17,3){\line(-5,6){5.68}}
\put(3,1.5){\line(1,1){8.32}}
\put(3,1.5){\line(2,1){3}}

\put(1.2,3.8){g}
\put(5.9,2){C}
\put(16.9,2){A}
\put(13,9){B}
\put(2,1){D}
\put(11.5,10.5){E}
\put(4,4){F}

\end{picture}
\qed

\begin{lemma}
\labell{ent/ob/modu}
The modulus of a bisceles arc supported on the essential traces of a tangled
v-arc or a tangled w-arc and the modulus of any of its obstructing parallels is
strictly larger than the modulus of the corresponding isosceles arc.
\end{lemma}

\proof
Thanks to \ref{ent/modu} we need only to argue for the obstructing parallels,
but their modulus is strictly bounded from below by the modulus of the
obstructed bisceles arc by lemma \ref{ob/iso/modu}.
\qed

\begin{lemma}
\labell{bis/modu/bnd}
The full twist on a bisceles arc which is not an isosceles arc
transported along
$t=1$ to $\l=1$ is in the group generated by all twists $\tau^2$ of modulus
larger than the modulus of the bisceles arc.
\end{lemma}

\proof
The obstructing parallels are bisceles arcs of strictly larger modulus.
Hence we may as well assume the bisceles arc to be unobstructed. Its parallel
transport is then isotopic to an arc defining some $\tau$ of larger modulus,
which is strictly larger in case the bisceles arc is no isosceles arc.
\qed


\section{the discriminant family}
\label{six/sevn}

In this section we will work with the discriminant family of the families of
function we consider. In order to compute the versal braid monodromy
in the next
section, we have to find the locally assigned groups. Moreover we need to
compare the parallel transport in the discriminant family to parallel transport
in the model family.\\

\begin{lemma}
\labell{cover}
The discriminant and the model discriminant family over the punctured parameter
bases have a common unramified cover.
\end{lemma}

\proof
The equation for the discriminant family has a formal factorisation
$$
\prod_{\xi_1^{l_1},\xi_2^{l_2}=1}
(z-\a^{\frac{l_1+1}{l_1}}\xi_1-\e_2^{\frac{l_2+1}{l_2}}\xi_2)\quad = \quad
\prod_{\xi_1^{l_1}=1}
((z-\a^{\frac{l_1+1}{l_1}}\xi_1)^{l_2}-\e_2^{l_2+1})\quad = \quad 0.
$$
as opposed to the equation for the model discriminant family:
$$
\prod_{\xi_1^{l_1},\xi_2^{l_2}=1}
(z-\l\xi_1-\eta_2\xi_2)\quad = \quad
\prod_{\xi_1^{l_1}=1}
((z-\l\xi_1)^{l_2}-\eta_2^{l_2})\quad = \quad 0.
$$
These equations coincide for $\eta_2^{l_2}=\e_2^{l_2+1}$ and
$\l^{l_1}=\a^{l_1+1}$. Hence the family parameterized by $\b$
$$
\prod_{\xi_1^{l_1}=1}
((z-\b^{l_1+1}\xi_1)^{l_2}-\e_2^{l_2+1})\quad = \quad 0.
$$
is isomorphic to the pull backs of the discriminant family and the model
discriminant family by
the covering map $\b\mapsto\a=\b^{l_1}$ resp.\ $\b\mapsto\l=\b^{l_1+1}$, if $
\e_21{l_2+1}=\eta_2^{l_2}$.

In this way we can understand polar coordinates of the bases of the two
discriminant families as different coordinates of the universal cover of the
bases punctured at the origin.\\

So with polar coordinates $r$ and $\theta$ in the base of the discriminant
family we can immediately compare parallel transport in the two families:

\begin{lemma}
Parallel transport in the discriminant family and in the model
discriminant family
coincides if
$r_0^{l_1+1}=t_0^{l_1}$ and $\theta(l_1+1)=\varth l_1$
\begin{enumerate}
\item
along radial paths $t\ethe$, $t\in[t_0,1]$ and $re_\theta, r\in[r_0,1]$,
\item
along circular paths of radius $1$ of winding angles $\varth$ and $\theta$
respectively.
\end{enumerate}
\end{lemma}

We can now define {\it standard paths} in the bases of both families by asking
them to be supported on radial segments and circular segments as in the lemma.

And for each standard path in one base we get another one in the other with the
same parallel transport.\\

\begin{examp}
A system of standard paths for the discriminant family associated to
$l_1=4,l_2=2$
is thus related to standard paths in the base of the model discriminant family:

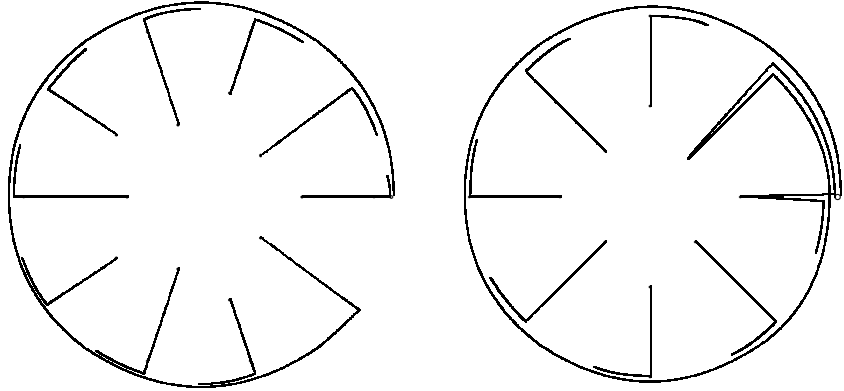
\begin{center}
\setlength{\unitlength}{6mm}
\begin{picture}(20,10)

\put(15,5){\setlength{\unitlength}{1.5mm}
\begin{picture}(40,40)

\put(-16,0){\circle*{.2}}
\put(0,-16){\circle*{.2}}
\put(0,16){\circle*{.2}}
\put(-8,0){\circle*{.2}}
\put(-4,4){\circle*{.2}}
\put(-4,-4){\circle*{.2}}
\put(0,-8){\circle*{.2}}
\put(0,8){\circle*{.2}}
\put(4,4){\circle*{.2}}
\put(4,-4){\circle*{.2}}
\put(8,0){\circle*{.2}}

\put(16.6,0){\circle{.45}}

\put(-16.5,-16.5){\setlength{\unitlength}{5mm}\begin{picture}(10,10)
\bezier{300}(0,5)(0,7.5)(2,9)
\bezier{300}(2,9)(5,11)(8,9)
\bezier{300}(8,9)(10,7.5)(10,5)
\bezier{270}(8,1.1)(9.7,2.4)(9.7,5)
\bezier{200}(8.2,8.2)(9.7,6.7)(9.7,5)
\bezier{300}(8,1.1)(5,-1)(2,1)
\bezier{300}(2,1)(0,2.5)(0,5)
\end{picture}}

\bezier{200}(8,0)(12,.1)(16.3,0.2)
\bezier{200}(8,0)(12,-.2)(15.3,-.4)
\bezier{200}(15.3,-.4)(15.2,-3.3)(14.6,-5)
\bezier{200}(-8,0)(-12,0)(-16,0)
\bezier{200}(-16,0)(-16,3.3)(-15.4,5)
\bezier{200}(0,8)(0,12)(0,16)
\bezier{200}(0,16)(3.3,16)(5,15.2)
\bezier{200}(0,-8)(0,-12)(0,-16)
\bezier{200}(0,-15.9)(-3.3,-15.9)(-5,-15.1)
\bezier{200}(4,-4)(6,-6)(11,-11)
\bezier{200}(11.1,-11.1)(9,-13.2)(7.2,-14)
\bezier{200}(-4,4)(-6,6)(-11,11)
\bezier{200}(-11.1,11.1)(-9,13.2)(-7.2,14)
\bezier{200}(4,4)(6,6)(10.8,10.8)
\bezier{200}(-11.1,-11.1)(-13.2,-9)(-14.2,-7.2)
\bezier{200}(-4,-4)(-6,-6)(-11,-11)
\bezier{200}(10.8,11.8)(16.3,6.7)(16.3,0)
\bezier{200}(10.8,11.8)(,)(4,4)

\end{picture}
}

\put(5,5){\setlength{\unitlength}{1.37mm}
\begin{picture}(40,40)

\put(-8,-6){\circle*{.2}}
\put(-8,6){\circle*{.2}}
\put(3,10){\circle*{.2}}
\put(3,-10){\circle*{.2}}
\put(10,0){\circle*{.2}}
\put(-7,0){\circle*{.2}}
\put(-2,-7){\circle*{.2}}
\put(-2,7){\circle*{.2}}
\put(6,-4){\circle*{.2}}
\put(6,4){\circle*{.2}}

\put(18.7,0){\circle{.3}}

\put(-18.5,-18.5){\setlength{\unitlength}{5.115mm}\begin{picture}(10,10)
\bezier{300}(0,5)(0,7.5)(2,9)
\bezier{300}(2,9)(5,11)(8,9)
\bezier{300}(8,9)(10,7.5)(10,5)
\bezier{70}(8,1)(8.5,1.4)(9.1,2)
\bezier{300}(8,1)(5,-1)(2,1)
\bezier{300}(2,1)(0,2.5)(0,5)
\end{picture}}

\bezier{120}(6,-4)(14,-10)(15.5,-11)
\bezier{120}(3,-10)(5,-16)(5.4,-17.2)
\bezier{120}(5.4,-17.2)(3,-18.2)(0,-18.2)
\bezier{120}(-2,-7)(-5,-16)(-5.4,-17.2)
\bezier{120}(-5.4,-17.2)(-8,-16.3)(-10,-15)
\bezier{120}(3,10)(5,16)(5.4,17.2)
\bezier{120}(5.4,17.2)(8,16.3)(10,15)
\bezier{120}(-2,7)(-5,16)(-5.4,17.2)
\bezier{120}(-5.4,17.2)(-3,18.2)(0,18.2)
\bezier{120}(-8,-6)(-14,-10)(-14.8,-10.5)
\bezier{120}(-14.8,-10.5)(-16.1,-9)(-17.2,-6)
\bezier{120}(6,4)(14,10)(14.8,10.5)
\bezier{120}(14.8,10.5)(16.1,9)(17.2,6)
\bezier{120}(-8,6)(-14,10)(-14.7,10.4)
\bezier{120}(-14.7,10.4)(-12.5,13.2)(-11,14.3)
\bezier{120}(-7,0)(-17,0)(-18,0)
\bezier{120}(-18,0)(-18,3)(-17.4,5)
\bezier{120}(10,0)(14,0)(18.5,0)
\bezier{60}(18.5,0)(18.5,1)(18.2,2)

\put(-17,0){\circle*{.2}}
\put(-14,10){\circle*{.2}}
\put(-14,-10){\circle*{.2}}
\put(-5,-16){\circle*{.2}}
\put(-5,16){\circle*{.2}}
\put(5,16){\circle*{.2}}
\put(5,-16){\circle*{.2}}
\put(-14,10){\circle*{.2}}
\put(-14,-10){\circle*{.2}}
\put(17,0){\circle*{.2}}
\put(-8,-6){\circle*{.2}}
\put(-8,6){\circle*{.2}}
\put(3,10){\circle*{.2}}
\put(3,-10){\circle*{.2}}
\put(10,0){\circle*{.2}}
\put(-7,0){\circle*{.2}}
\put(-2,-7){\circle*{.2}}
\put(-2,7){\circle*{.2}}
\put(6,-4){\circle*{.2}}
\put(6,4){\circle*{.2}}

\end{picture}
}
\end{picture}
\end{center}
\end{examp}

To get the versal braid monodromy of the discriminant family, we therefore need
to transfer the locally assigned groups from local Milnor fibres of the
discriminant family to local Milnor fibres of the model discriminant family and
transport them along all possible standard paths.\\

We assign a group to a local Milnor fibre in the model discriminant using the
fact that the fibre is isomorphic to a local Milnor fibre in the discriminant
family by way of the two finite covering maps.

\begin{lemma}
\labell{fake/ord/mondr}
The group assigned to a Milnor fibre at a regular parameter $t_1\etheI$,
sufficiently close to a singular parameter $t_0\etheO\neq0$ with
$t_1-t_0>0$, is generated by full twists on local v-arcs.
\end{lemma}

\proof
The singular fibre corresponds to a function with non-degenerate
critical points
only, cf.\ the proof of the lemmas \ref{gtame}, \ref{ftame}.
So by definition the locally assigned group is generated by mapping classes
fixing all punctures and supported on small discs each of which is a
Milnor fibre
for just one multiple puncture.\\

By close inspection we can see that the local v-arcs are supported on
such discs
and the full twists on local v-arcs generate the group of all mapping
classes of
each disc which preserve the punctures.
\qed

\begin{lemma}
\labell{fake/An/mondr}
The group assigned to a Milnor fibre at a regular parameter $t_1\etheI$,
sufficiently close to a singular parameter $\l=0$, is generated by
full twists on
local w-arcs and $\frac32$-twists on local v-arcs with index pair
$i_1i_2,i_1^+i_2$.
\end{lemma}

\proof
The singular fibre corresponds to a function which has $l_2$ critical points of
type $A_{l_1}$ with distinct critical values. So by definition the
group locally assigned to each disc, which is a local Milnor fibre of
a multiple
puncture, is generated by the mapping classes of the braid monodromy of the
singular function germ it corresponds to.

Each of the critical points of type $A_{l_1}$ is unfolded linearly,
so the local
Milnor fibre can be naturally identified with the Milnor fibre encountered in
lemma \ref{An-mono}. And in combination with lemma \ref{brmon-An} we conclude
that local generators are given by the $\frac32$-twists on v-arcs with index
pairs $i_1i_2,i_1^+i_2$ and full twists on arcs winding positively from a
puncture of index $i_1i_2$ to a puncture of index $j_1i_2$, $j_2>i_1^+=i_1+1$,
around all v-arcs.

By lemma \ref{iso/tri/gen} we can see that instead we can use the twists of the
claim to generate the same group.
\qed

To summarize the preceding discussion we should note:

\begin{remark}
The versal braid monodromy of the family of functions
$$
x_1^{l_1+1}-\a(l_1+1)x_1+x_2^{l_2+1}-\e_2(l_2+1)x_2
$$
is generated by the parallel transport of the appropriate twists as given by
lemma \ref{fake/ord/mondr} and lemma \ref{fake/An/mondr} along all
standard paths
in the model discriminant family.
\end{remark}



\section{conclusion}
\label{six/eit}

Finally we keep our promise and give braid elements which generate the versal
braid monodromy:

\begin{prop}
\labell{mainprop}
The versal braid monodromy the family of functions
$$
x_1^{l_1+1}-\a(l_1+1)x_1+x_2^{l_2+1}-\e_2(l_2+1)x_2
$$
is generated by twists ($i^+_1=i_1+1$, $i_2'=i_2-l_2$):
\begin{eqnarray*}
\tau^2_{i_1i_2,j_1j_2}, & & 1\leq i_1<j_1\leq l_1,1\leq j_2-i_2<l_2,\\
\tau^3_{i_1i_2',i^+_1i_2}, & & 1\leq i_1<i^+_1\leq l_1, 1\leq i_2\leq l_2,\\
\tau_{i_1i_2',i_1^+i_2}\inv\tau_{i^+_1i_2',j_1i_2}^2\tau_{i_1i_2',i_1^+i_2},
&  & 1<i_1^+<j_1\leq l_1, 1\leq i_2\leq l_2.
\end{eqnarray*}
\end{prop}

\proof
The versal braid monodromy of a one parameter family can by definition be
computed from their locally assigned groups of mapping classes and the parallel
transport of these groups along a distinguished system of paths in
the associated
discriminant family, cf.\ lemma \ref{vers-geom}\\

The locally assigned groups in the discriminant family were given
in lemma \ref{fake/ord/mondr} and lemma \ref{fake/An/mondr}
to be twists on local v-arcs and local w-arcs.

So by the closing remark of the last section parallel transport of local
v-arcs and local w-arcs along all possible standard paths in the base of the
model discriminant family generate the versal braid monodromy.\\

Note that the length of the circular part is not necessarily restricted to
$[0,2\pi[$.\\

We denote by $T$ the set of braid generators obtained by parallel transport and
identification using the Hefez Lazzeri path system in the fibre at $\l=1$.

$T$ is divided into subsets according to the index pair of the punctures
connected by the corresponding arc.\\

The given set $S$ of braid group elements is also divided into subsets
according to the modulus the index pairs of each element, which is unambiguous
since we note immediately that the modulus of all index pairs occurring in the
second and third row is zero.\\

Since the moduli of elements in $T$ and $S$ form a finite descending sequence
$m_1>...>m_n=0$, we can impose finite filtrations
$$
T_k:=\{\tau\in T| m(\tau)\geq m_k\},\quad
S_k:=\{\tau\in S| m(\tau)\geq m_k\}.
$$
To prove our claim, we are thus
left to check the hypotheses of lemma \ref{filt}:\\

Since $l_2>1$, the maximal modulus $m_1$ is positive. Hence $S_1$ only contains
twists on parallel transports of local v-arcs.
The local v-arcs of highest modulus get not tangled when transported along a
radial arc, since entangling bisceles arcs have to be of larger modulus
\ref{ent/modu}. The isosceles thus obtained are unobstructed, since
obstructing parallels would be of larger modulus, \ref{ob/iso/modu}.
By \ref{crit-twist} each element of $T_1$ is in $S_1$. Conversely by
\ref{surj}, each element in $S_1$ is an element in $T$ of equal modulus, hence
in $T_1$.\\

Given an element in $T$ of modulus $m_k>0$, which is the parallel
transport of an
local v-arc, then there is an element in $S$ obtained from the same local v-arc
transported along the same path, but conjugated by twists which are
the parallel
transports of bisceles arcs of strictly larger modulus, \ref{ent/modu},
\ref{ent/ob/modu}, \ref{ob/iso/modu}. So the second hypothesis of lemma
\ref{filt} holds for elements in $T_k-T_{k-1}$ of positive modulus.\\

Conversely each full twist in $S$ of positive modulus is obtained by parallel
transport from an local v-arc of equal modulus up to twists by entangled and
obstructing bisceles arcs, \ref{surj}. So due to \ref{ent/modu},
\ref{ent/ob/modu}, \ref{ob/iso/modu} again the third hypothesis
holds for the twists obtained from local v-arcs of positive modulus.\\

We are left with elements of modulus $m_n=0$ and consider local w-arcs first.
Though we have to transport along a standard path to get full twist elements in
$T_n-T_{n-1}$, we have to rely on a result which makes use of a different kind
of paths. We recall that in section \ref{six/five} we transported a local w-arc
along a circular arc of small radius, then along a radial segment and finally
along a circular segment of radius $t=1$.

Each standard path can be coupled with a path of the second kind in such a way,
that the closed path obtained as their join has winding number zero with
respect to the origin.
Hence parallel transport from $\l=1$ along this closed path amounts to
conjugation by a composition of full twists of positive modulus in $T_{n-1}$.

We deduce that the elements of $T$ obtained
by parallel transport along standard paths yield the elements given in lemma
\ref{tri/twist} up to conjugation by full twists in $T_{n-1}$ and $S_{n-1}$.
By lemma \ref{tri/conj} they are even conjugate to the full twists in
$S_n-S_{n-1}$.\\

This relation can obviously be reversed in the sense, that for each full twist
in $S_n$ of modulus zero we have an element in $T$ equal up to conjugation by
elements in $S_{n-1}$ and $T_{n-1}$.\\

Finally we have to address the $\frac32$-twists in $T$ a nd $S$.
A $\frac32$-twist in $T$ is obtained by the parallel transport of a local v-arc
with index pair $i_1i_2,i_1^+i_2$.
We conclude that up to full twists elements in $S_{n-1}$ and $T_{n-1}$ the
$\frac32$-twists in $T$ are among the elements given in lemma
\ref{cri-twist}. By
lemma \ref{iso/zero/conj} they are among the elements in $S$ even up to full
twists in $S_{n-1}$ and $T_{n-1}$.\\

Conversely the pairs of twists considered in lemma \ref{sur} correspond
bijectively to the
\mbox{$\frac32$-twists} of $S$ and are both equal up to conjugation by full
twists of positive modulus by lemma \ref{iso/zero/conj}. We deduce that also each
$\frac32$-twist of $S$ is an element of $T$ up to conjugation by full twists in
$S_{n-1}$ and $T_{n-1}$.
\qed

\section{appendix on plane elementary geometry}

\begin{lemma}
\labell{geo/x}
Given a proper triangle $ABC$ and lines through $A,C$ resp.\ $B,C$. If a fourth
point $D$ is in the opposite cone to $\overline{AB}$ at $C$, then $C$ is a
point in the disc with boundary through $ABD$.
\end{lemma}

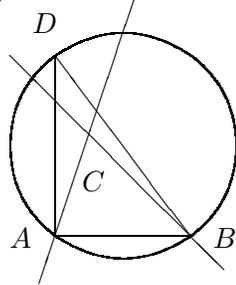
\setlength{\unitlength}{3mm}
\begin{picture}(30,12)(-10,-2)

\bezier{300}(0,5)(0,7.5)(2,9)
\bezier{300}(2,9)(5,11)(8,9)
\bezier{300}(8,9)(10,7.5)(10,5)
\bezier{300}(10,5)(10,2.5)(8,1)
\bezier{300}(8,1)(5,-1)(2,1)
\bezier{300}(2,1)(0,2.5)(0,5)

\put(2,1){\line(1,0){6}}
\put(8,1){\line(-3,4){6}}
\put(2,1){\line(0,1){8}}

\put(2,1){\line(1,3){3.5}}
\put(2,1){\line(-1,-3){.7}}
\put(8,1){\line(-1,1){8}}
\put(8,1){\line(1,-1){1.5}}

\put(0,.5){$A$}
\put(9,.5){$B$}
\put(3.2,3){$C$}
\put(1,10){$D$}

\end{picture}

\proof
Denote the line through a pair of points by $g$ indexed with the pair. Then the
assumptions can be stated as:
\begin{enumerate}
\item
$C,D$ are on the same side of the line $g_{AB}$,
\item
$A,D$ are on different sides of $g_{BC}$,
\item
$B,D$ are on different sides of $g_{AC}$.
\end{enumerate}
These imply
\begin{enumerate}
\item[iv)]
$A,C$ are on the same side of $g_{BD}$ by i) and ii),
\item[v)]
$B,C$ are on the same side of $g_{AD}$ by i) and iii).
\end{enumerate}
But i),iv) and v) together form the assertion of the lemma.
\qed

\begin{lemma}
\labell{ABCDEF}
Suppose in a quadrilateral $ABCD$ the points $B,D$ are on opposite
sides of the diagonal $d$ through $A$ and $C$. Let $E$ be the intersection of
the line through $A$ and $B$ and the parallel to $BC$ through $D$. If
$|AB|<|BC|$ then
$|AE|<|DE|$.
\end{lemma}

\proof
By hypothesis there is an intersection point $F$ on $\overline{DE}$ with the
diagonal $d$, so $|EF|<|DE|$. Hence the claim follows since by proportionality
$|AB|<|BC| \implies |AE|<|EF|$.

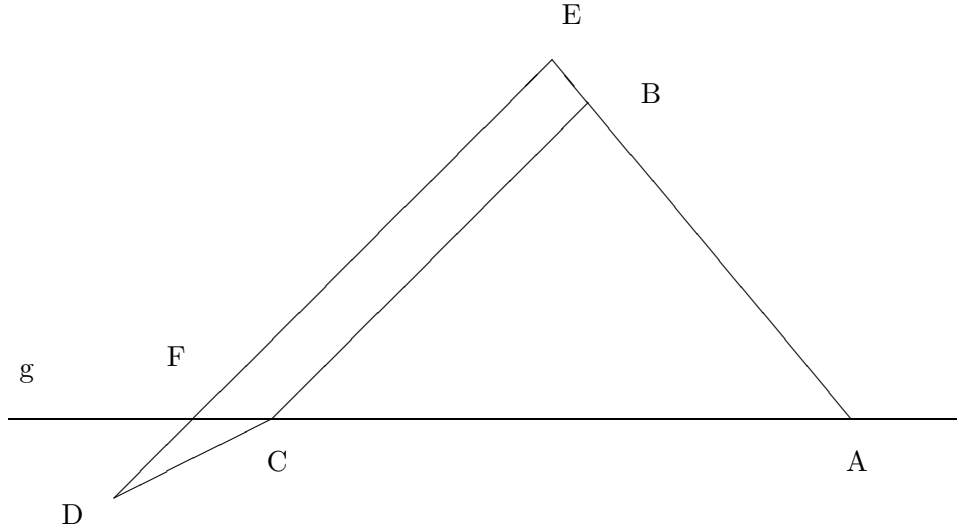
\setlength{\unitlength}{7mm}
\begin{picture}(20,11)

\put(1,3){\line(1,0){18}}
\put(6,3){\line(1,1){6}}
\put(17,3){\line(-5,6){5.68}}
\put(3,1.5){\line(1,1){8.32}}
\put(3,1.5){\line(2,1){3}}

\put(1.2,3.8){g}
\put(5.9,2){C}
\put(16.9,2){A}
\put(13,9){B}
\put(2,1){D}
\put(11.5,10.5){E}
\put(4,4){F}

\end{picture}
\qed

\begin{lemma}
\labell{leng/inc}
Suppose in a quadrilateral $ABCD$ the points $B,D$ are on opposite sides of the
diagonal $d_{AC}$ through $A$ and $C$ and the points $A,D$ are on
opposite sides
of the line through $B$ and $C$.
If $E$ is the intersection of the line
through $A$ and $B$ and the parallel to $BC$ through $D$ then
\begin{eqnarray*}
\min(|DE|,|AE|) & > & \min(|AB|,|BC|).
\end{eqnarray*}
\end{lemma}

\proof
By hypothesis $E$ and $A$ are on opposite sides of the line through $B$ and
$C$, hence $|AB|<|AE|$ and by proportionality $|EF|>|BC|$ which implies
$|DE|>|BC|$.

So with the result of lemma \ref{ABCDEF} we get in case $|AB|=\min(|AB|,|BC|)$
that $|DE|>|AE|>|AB|$. In case $|BC|=\min(|AB|,|BC|)$ we get the claim since
$|AE|>|AB|\geq |BC|$ and $|DE|>|BC|$.
\qed


\chapter{braid monodromy induction to higher dimension}


In the previous chapter we determined the braid monodromy of a one parameter
family of polynomials on the plane. Due to the results of
chapter \ref{Zar} this yields the main contribution to the braid monodromy
of a plane Brieskorn Pham polynomial.\\

In the present chapter we exploit the results of both chapters to set up an
induction for the computation of the braid monodromy for Brieskorn Pham
polynomials in arbitrary dimensions.
Since we have to relate the monodromies of families of different numbers of
variables and deal with an arbitrary number in the general case, we made some
efforts to chose our notation. We devote the first section to introduce this
notations, quite a few new definitions and to rephrase results of the preceding
chapters in a unified way.\\

In the second and third section the induction argument is given by way of
considering the versal braid monodromies of families of functions
\begin{eqnarray*}
f_\a & : &
x_1^{l_1+1}-\a(l_1+1)x_1+\sum_{i=2}^n\left(x_i^{l_i+1}-\e_i(l_i+1)x_i\right),\\
g_\a & : &
x_1^{l_1+1}-(l_1+1)x_1+\sum_{i=2}^n\left(x_i^{l_i+1}-\a\e_i(l_i+1)x_i\right),
\end{eqnarray*}
which were introduced in chapter \ref{Zar}.

To merge the various groups we have to choose the generating sets in many
different ways. Even though many of the braid computations have been put into
an appendix, the computational load is quite high.\\

The actual geometric argument which makes induction possible is presented in
the fourth and fifth section. We present a prove of a result which connects the
braid monodromy of a discriminant family to a family which is obtained by
replacing the divisor by a number of parallel copies. This $l$-companion family
is studied to the details to relate the associated versal monodromies.

Both sections should be regarded as a sort of appendix to which we
may resort on
need.

\section{preliminaries}

The choice of a Hefez Lazzeri base in the reference fibre provides a natural
bijection of the punctures with a multiindex set associated appropriately to
the given Brieskorn Pham polynomial.

\begin{nota}
Given a finite sequence $l_1,...,l_n$ of positive integers, define the
multiindex set $\idx=\idx(l_1,...,l_n)$ to be
\begin{eqnarray*}
\idx  & := & \{i_1...i_n\,|\,1\leq i_\nu\leq l_\nu,\,1\leq\nu\leq n\}
\end{eqnarray*}
equipped with the natural lexicographical order.
\end{nota}

While $i$ denotes an element $i_1...i_n$, we
will use $i^+$ for its immediate successor and $i':=i_1...i_{n-1}$ for the
naturally associated element in
$\idx(l_1,...,l_{n-1})$.\\

Whether the following property is given or not determines to some extend the
role of index pairs and index triples.

\begin{defi}
Multiindices $i,j\in \idx $ are called {\it correlated} or a {\it
correlated pair}, if $i<j$ and $j_\nu\in \{i_\nu,i_\nu+1\}$ , $1\leq\nu\leq n$.\\
Multiindices $i,j,k\in \idx $ are called {\it correlated} or a {\it correlated
triple}, if $i<j$, $i<k$, $j<k$ are correlated.\\
A quadruple of indices is called {\it correlated} if each pair is correlated.
\end{defi}

For the induction we need several homomorphisms between braid groups of
different numbers of strings.

\begin{defi}
Given a multiindex set $\{i_1...i_n\,|\,1\leq i_\nu\leq l_\nu,1\leq\nu\leq
n\}$\\
the {\it primary homomorphisms} are defined for $1\leq i_1\leq l_1$:
$$
\begin{array}{rccc}
\phi_{i_1}: & \br_{l_2...l_{n}} & \tto & \br_{l_1...l_n}\\
& \s_{i_2...i_n,j_2...j_n} & \mapsto &
\s_{i_1i_2...i_n,i_1j_2...j_n}
\end{array}
$$
the {\it secondary homomorphisms} are defined for $1\leq i_n\leq l_n$:
$$
\begin{array}{rccc}
\psi_{i_n}: & \br_{l_1...l_{n-1}} & \tto & \br_{l_1...l_n}\\
& \s_{i_1...i_{n-1},j_1...j_{n-1}} & \mapsto &
\s_{i_1...i_n,j_1...j_{n-1}i_n}
\end{array}
$$
the {\it $l_n$-band homomorphism}
$$
\begin{array}{rccl}
\eta_{l}: & \br_{l_1...l_{n-1}} & \tto & \br_{l_1...l_n}\\
\end{array}
$$
which assigns to a braid  a braid of $l_n$-times as many strand by
replacing each
strand with a ribbon of $l_n$ strands.

This definition should be understood from the following picture: To a twist,
i.e.\ an exchange of two strands, we associate a band or ribbon twist, i.e.\ an
exchange of two bands or ribbons into which $l$ strands each have
been assembled:

\setlength{\unitlength}{5mm}
\begin{picture}(20,15)(-6,-10)

\put(-.4,4.2){$i'1$}
\put(0.8,4.2){$i'2$}
\put(2.8,4.2){$i'l_n$}

\put(10.4,4.2){$j'1$}
\put(11.8,4.2){$j'2$}
\put(13.8,4.2){$j'l_n$}

\bezier{2}(1.5,2.5)(2,2.5)(2.5,2.5)
\bezier{2}(7.5,2.5)(8,2.5)(8.5,2.5)
\bezier{2}(12.5,2.5)(13,2.5)(13.5,2.5)

\drawline(0,-4)(.5,-3.7)
\drawline(1.2,-3.28)(1.5,-3.1)
\drawline(2,-2.8)(3.2,-2.08)
\drawline(3.9,-1.66)(11,2.6)

\drawline(1,-4.5)(1.5,-4.2)
\drawline(2,-3.9)(2.5,-3.6)
\drawline(3,-3.3)(4.1,-2.64)
\drawline(4.8,-2.22)(12,2.1)

\drawline(3,-5.5)(3.3,-5.32)
\drawline(3.8,-5.02)(4.2,-4.78)
\drawline(4.9,-4.36)(6,-3.7)
\drawline(6.8,-3.22)(14,1.1)

\drawline(0,-4)(0,-10)
\drawline(0,3)(0,-3)
\drawline(1,-4.5)(1,-10)
\drawline(1,3)(1,-2.5)
\drawline(3,-5.5)(3,-10)
\drawline(3,3)(3,-1.5)

\drawline(9,-8.6)(9,-10)
\drawline(9,-7.5)(9,-8.2)
\drawline(9,1.2)(9,.5)
\drawline(9,3)(9,1.6)
\drawline(9,.1)(9,-1.7)
\drawline(9,-5.1)(9,-7.1)
\drawline(9,-2.1)(9,-4.9)

\drawline(7,3)(7,.4)
\drawline(7,-7.4)(7,-10)
\drawline(7,-6.3)(7,-7.0)
\drawline(7,-.7)(7,-0)
\drawline(7,-4.1)(7,-5.9)
\drawline(7,-3.3)(7,-3.7)
\drawline(7,-1.1)(7,-2.9)

\drawline(11,-9.6)(0,-3)
\drawline(12,-9.1)(1,-2.5)
\drawline(14,-8.1)(3,-1.5)

\drawline(11,3)(11,2.6)
\drawline(11,-9.6)(11,-10)
\drawline(12,3)(12,2.1)
\drawline(12,-9.1)(12,-10)
\drawline(14,3)(14,1.1)
\drawline(14,-8.1)(14,-10)

\end{picture}
\begin{figure}[h]
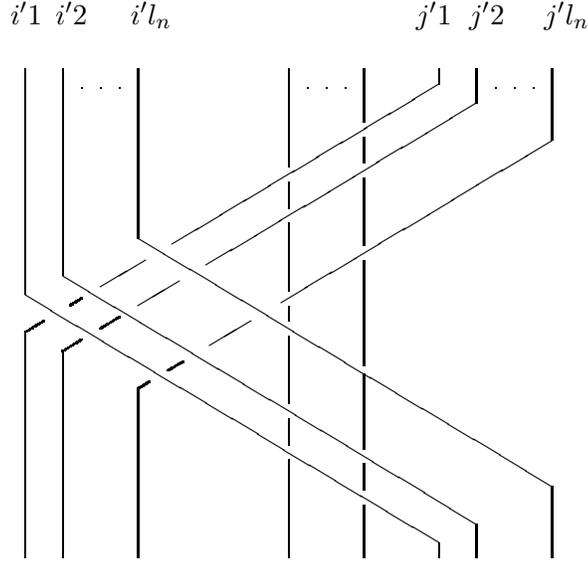

\caption{band twist, the image of a half twist under $\eta_{l_n}$}
\end{figure}
\end{defi}

\begin{defi}
The {\it level} of an index pair is the difference of the leading components.\\
The {\it level} of a braid $\s_{i,k}$ or $\tau_{i,k}$ is the level of its index
pair.
\end{defi}

\begin{defi}
The element $\delta_n:=\s_{1,2}\cdots\s_{n-1,n}$ is called the {\it fundamental
element} of the BKL presentation, cf.\ \cite{bkl} of the braid group $\br_n$,
$\delta^n$ is the full twist on the disc, i.e.\ a generator of the center of
$\br_n$.
\end{defi}

\begin{defi}
Given a pair of indices $i,j$ the associated {\it subcable twist} is defined as
the BKL fundamental word on the braid subgroup on the punctures with index $k,
i\leq k\leq j$, it can be given as
$$
\delta_{i,j}:=\prod_{i\leq k<j}\s_{k,k^+}.
$$
\end{defi}

\begin{defi}
The {\it $l_n$-cable twist} in $\br_{l_1...l_n}$ is defined to be the element
$$
\delta_{\phi,n}:=\prod_{i'\in I_{n-1}}\delta_{i'1,i'l_n}.
$$
\end{defi}

In the disguise of turning the peripheral circles we already considered this
cable twist, so in case $n=2$ conjugation
by the cable twist $\delta_{j_11,j_1l_2}$ yields maps
\begin{eqnarray*}
\tau_{i_1i_2,j_1j_2} & \mapsto & \tau_{i_1i_2,j_1j_2^+},\\
\tau_{j_1j_2,k_1k_2} & \mapsto & \tau_{j_1j_2^+,k_1k_2},\\
\s_{i_1i_2,j_1j_2} & \mapsto & \s_{i_1i_2,j_1j_2^+},\\
\s_{j_1j_2,k_1k_2} & \mapsto & \s_{j_1j_2^+,k_1k_2},
\end{eqnarray*}
where as usual we assume $i<j<k$.
\pagebreak

We extend the range of possible indices in higher dimensions accordingly.

\begin{defi}
Suppose $i'_n-i_n=:m_i,j'_n-j_n=:m_j$ then define
$$
\s_{i'i_n',j'j_n'}:=
\delta_{j'1,j'l_n}^{-m_j}\delta_{i'1,i'l_n}^{-m_i}\s_{i_1i_2,j_1j_2}
\delta_{i'1,i'l_n}^{m_i}\delta_{j'1,j'l_n}^{m_j}.
$$
\end{defi}

\begin{examp}
In case $n=2,i_1=1,j_1=2$ there are the following examples:

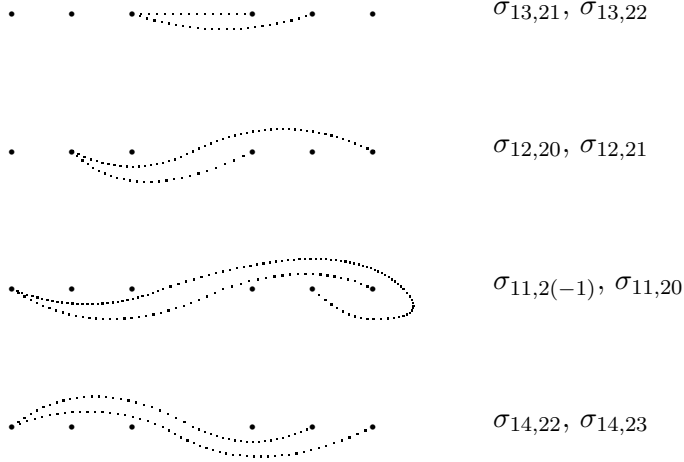
\setlength{\unitlength}{4mm}

\begin{picture}(20,4)

\put(1,2){\circle*{.2}}
\put(3,2){\circle*{.2}}
\put(5,2){\circle*{.2}}
\put(9,2){\circle*{.2}}
\put(11,2){\circle*{.2}}
\put(13,2){\circle*{.2}}

\put(17,2){$\s_{13,21}$, $\s_{13,22}$}

\bezier{15}(5,2)(7,2)(9,2)
\bezier{25}(5,2)(8,1)(11,2)

\end{picture}

\begin{picture}(20,4)

\put(1,2){\circle*{.2}}
\put(3,2){\circle*{.2}}
\put(5,2){\circle*{.2}}
\put(9,2){\circle*{.2}}
\put(11,2){\circle*{.2}}
\put(13,2){\circle*{.2}}

\put(17,2){$\s_{12,20}$, $\s_{12,21}$}

\bezier{25}(3,2)(5,0)(9,2)
\bezier{20}(3,2)(5,1)(7,2)
\bezier{30}(7,2)(10,3.5)(13,2)

\end{picture}

\begin{picture}(20,4)

\put(1,2){\circle*{.2}}
\put(3,2){\circle*{.2}}
\put(5,2){\circle*{.2}}
\put(9,2){\circle*{.2}}
\put(11,2){\circle*{.2}}
\put(13,2){\circle*{.2}}

\put(17,2){$\s_{11,2(-1)}$, $\s_{11,20}$}

\bezier{30}(1,2)(4,0)(8,2)
\bezier{20}(8,2)(11,3)(13,2)
\bezier{30}(1,2)(3.5,1)(6,2)
\bezier{45}(6,2)(12,4)(14,2)
\bezier{15}(14,2)(15,1)(13,1)
\bezier{10}(13,1)(12,1)(11,2)

\end{picture}

\begin{picture}(20,4)

\put(1,2){\circle*{.2}}
\put(3,2){\circle*{.2}}
\put(5,2){\circle*{.2}}
\put(9,2){\circle*{.2}}
\put(11,2){\circle*{.2}}
\put(13,2){\circle*{.2}}

\put(17,2){$\s_{14,22}$, $\s_{14,23}$}

\bezier{30}(13,2)(9,0)(6,2)
\bezier{20}(6,2)(3.5,3)(1,2)
\bezier{20}(11,2)(9,1)(7,2)
\bezier{25}(7,2)(3.5,4)(1,2)

\end{picture}
\end{examp}

\begin{lemma}
\labell{cable-orb}
Conjugation by the cable twist $\cable$ induces a level preserving bijection
of braids $\tau^2_{i_1i_2,j_1j_2}$ with $1\leq i_1<j_1\leq l_1,1\leq
j_2-i_2<l_2$ and of braids $\s^2_{i_1i_2,j_1j_2}$ with $1\leq i_1<j_1\leq
l_1,1\leq j_2-i_2<l_2$ such that:
\begin{eqnarray*}
\tau^2_{i_1i_2,j_1j_2} & \mapsto & \tau^2_{i_1i_2^+,j_1j_2^+},\\
\s^2_{i_1i_2,j_1j_2} & \mapsto & \s^2_{i_1i_2^+,j_1j_2^+}
\end{eqnarray*}
\end{lemma}

\begin{nota}
By $A(l_1)$ denote the singular function $x^{l_1+1}$.
By $BP(l_1,...,l_n)$ denote the singular function
$x_1^{l_1+1}+\cdots+x_n^{l_n+1}$.

With respect to the Hefez Lazzeri base these functions determine by braid
monodromy well defined subgroups of the braid group $\br_\mu$, which
are denoted
by
$$
\br_{A(l_1)}\quad\text{ resp.}\quad
\br_{BP(l_1,...,l_n)}.
$$
\end{nota}

\subsection*{braid monodromy results}

We have to review our main results for braid monodromy of Brieskorn Pham
polynomials up to now. The first records a slight reformulation of theorem
\ref{mainprop}. The second is the previously noted result for
singularities of type
$A_l$.
\\

We remark once and for all that the constants $\e_1=1,\e_i>0$ are
chosen in such a
way that $\e_{i+1}\ll\e_i$.

\begin{lemma}
\labell{sigma}
The versal braid monodromy of the family
$$
x_1^{l_1+1}-\a(l_1+1)x_1+x_2^{l_2+1}-\e_2(l_2+1)x_2
$$
is generated by twists
\begin{eqnarray*}
\tau^2_{i_1i_2,j_1j_2}, & & 1\leq i_1<j_1\leq l_1,1\leq j_2-i_2<l_2,\\
\s^3_{i_1i_2,i^+_1i_2}, & & 1\leq i_1<i^+_1\leq l_1, 1\leq i_2\leq l_2,\\
\s^2_{i_1i_2,j_1i_2}, & & 1<i_1^+<j_1\leq l_1, 1\leq i_2\leq l_2.
\end{eqnarray*}
\end{lemma}

\proof
By definition $\tau^3_{i_1i_2',i^+_1i_2}=\s^3_{i_1i_2,i^+_1i_2}$, $1\leq
i_1<i^+_1\leq l_1, 1\leq i_2\leq l_2$.\\
We then notice that $\psi_{i_2}(\br(A_{l_1}))$ is generated by
$$
\s^3_{i_1i_2,i^+_1i_2},\s^2_{i_1i_2,j_1i_2},\quad 1\leq i_1<i^+_1<j_1\leq l_1.
$$
By lemma \ref{iso/tri/gen} then $\psi_{i_2}(\br(A_{l_1}))$ is also generated by
$$
\s^3_{i_1i_2,i^+_1i_2},\psi_{i_2}
(\s_{i,i^+}\invv\check\s^2_{i_1,j_1}\s_{i,i^+}^2),
\quad 1\leq i_1<i^+_1<j_1\leq l_1.
$$
Since the latter elements coincide with
$\tau_{i_1i_2',i_1^+i_2}\inv\tau_{i^+_1i_2',j_1i_2}^2\tau_{i_1i_2',i_1^+i_2}$
we
are done.
\qed

In this proof we used the braid monodromy of $A_l$-singularities, which we
recall from \ref{brmon-An} for convenience.

\begin{lemma}
\labell{An}
The braid monodromy of the function $x^{l+1}$ with respect to the Hefez Lazzeri
system of paths is given by
$$
\s_{i,i^+}^3,1\leq i<l,\quad \s_{i,j}^2, 1< i^+<j<l.
$$
\end{lemma}


\section{families of type $g_\a$}

The topic of this section are the versal braid monodromies of families of
functions $g_\a$. We refer by $G_2,G_3$ and $G_n$ to the groups associated to
polynomials of two, three resp.\ $n$ variables.

Let us remark that by lemma \ref{ZB} the groups
determined in this section are in a sense the smaller complement to the the
groups for the families of type $f$, which have been investigated in the plane
case in the last chapter and will also be the topic of the next section.\\

The choice of the constants $\e_i$ is tacitly assumed to be made in such a
way, that the unit parameter disc of the families $g_\a$ contains no singular
parameter apart from the origin.

\begin{lemma}
\labell{ggfam}
The versal braid monodromy $G_2$ of the family $g_\a(x_1,x_2)$
restricted to the
unit disc is generated by the homomorphic images in $\br_{l_1l_2}$ of the
braid monodromy groups $\br_{A(l_2)}\subset\br_{l_2}$ under the primary
homomorphisms
$\phi_{i_1},1\leq i_1\leq l_1$.
\end{lemma}

\proof
The only critical parameter in the disc $|\a|\leq1$ is $\a=0$.
The corresponding critical function is $x_1^{l_1+1}-(l_1+1)x_1+x_2^{l_2+1}$,
which has $l_1$ critical point of type $A(l_2)$ with distinct critical values.

The bifurcation divisors of the families of functions parameterized by $\a$,
$$
-l_1\xi+x_2^{l_2+1}+\e_2\a(l_2+1)x_2,\quad \xi^{l_1}=1,
$$
embed into the bifurcation divisor of $g_\a$, and the corresponding
embeddings of
punctured discs induce embeddings of mapping class groups which correspond to
the embeddings $\phi_{i_1}$ under the standard identifications
with the braid groups $\br_{l_1l_2}$ and
$\br_{l_2}$ by the Hefez Lazzeri choice of a strongly distinguished system of
paths.
The versal braid monodromies of
the families above can be identified with the braid monodromy of lemma \ref{An}
and yield then the versal braid monodromy of the family $g_\a$ as claimed.
\qed

\begin{lemma}
\labell{gg}
The group $G_2\subset\br_{l_1l_2}$ is generated by the elements
\begin{eqnarray*}
& \s_{i_1i_2,i_1j_2}^3, & 1\leq i_2=j_2-1<l_2,1\leq i_1\leq l_1,\\
& \s_{i_1i_2,i_1j_2}^2, & 1\leq i_2<j_2-1<l_2,1\leq i_1\leq l_1.
\end{eqnarray*}
\end{lemma}

\proof
The group $\br_{A(l_2)}$ is generated by elements
$$
\s_{i_2,i^+_2}^3,1\leq i_2<l_2,\quad \s_{i_2,j_2}^2, 1\leq i_2<j_2-1<l_2,
$$
so the claim holds by lemma \ref{ggfam} since their images under the primary
homomorphisms $\phi_{i_1}$ are the elements of the assertion.
\qed

\begin{lemma}
\labell{gggfam}
The versal braid monodromy $G_3$ of the family
$g_\a(x_1,x_2,x_3)$ restricted to the unit disc is the subgroup of
$\br_{l_1l_2l_3}$ generated by the homomorphic images of the braid
monodromy groups
$\br_{BP(l_2,l_3)}\subset\br_{l_2l_3}$ under the primary homomorphisms
$\phi_{i_1},1\leq i_1\leq l_1$.
\end{lemma}

\proof
We give the proof stressing the analogy to the case $n=2$ of
\ref{ggfam}:\\
The only critical parameter in the disc $|\a|\leq1$ is $\a=0$.
The corresponding critical function is
$x_1^{l_1+1}-(l_1+1)x_1+x_2^{l_2+1}+x_3^{l_3+1}$,
which has $l_1$ critical points of type $BP(l_2,l_3)$ with distinct critical
values.

The bifurcation divisors of the families of functions
$$
-l_1\xi+x_2^{l_2+1}+\e_2\a(l_2+1)x_2+x_3^{l_3+1}+\e_3\a(l_3+1)x_3,\quad
\xi^{l_1}=1,
$$
embed into the bifurcation divisor of $g_\a$, and the corresponding
embeddings of
punctured discs induce embeddings of mapping class groups which correspond to
the embeddings $\phi_{i_1}$ under the Hefez Lazzeri identifications with the
braid groups $\br_{l_1l_2l_3}$ and
$\br_{l_2l_3}$. The versal braid monodromies of
the families above, which are given by the braid monodromy of polynomials of
type $BP(l_2,l_3)$, then yield the claim.
\qed

For the next lemma we have to resort for the first time to the forthcoming
sections on $l$-companion families.

\begin{lemma}
\labell{gggband}
The group $G_3$ contains the image of $G_2$ under the $l_3$-band homomorphism.
\end{lemma}

\proof
The group $G_2$ is generated by the braid monodromy of singular polynomials of
type $A_{l_1}$ under the primary homomorphisms. The $l_3$-companion family of the
discriminant family of $g_\a(x_1,x_2)$ is the discriminant family of:
$$
x_1^{l_1+1}-(l_1+1)x_1+x_2^{l_2+1}+\e_2\a(l_2+1)x_2+x_3^{l_3+1}-\e_3(l_3+1)x_3.
$$
Similarly its versal braid monodromy group is generated by the images under the
primary homomorphisms of the versal braid monodromy for the
$l_3$-companion of the
discriminant family of the families of functions
$$
-l_1\xi+x_2^{l_2+1}+\e_2\a(l_2+1)x_2,\quad \xi^{l_1}=1.
$$
Their versal braid monodromy coincides with the braid monodromy of the versal
unfolding of a singular function of type $A_{l_2}$. So by lemma
\ref{versal-band}
we conclude that the versal braid monodromy of the family of functions
$$
x_1^{l_1+1}-(l_1+1)x_1+x_2^{l_2+1}+\e_2\a(l_2+1)x_2+x_3^{l_3+1}-\e_3(l_3+1)x_3
$$
contains the images under the $\phi_{i_1}$ and $\eta_{l_3}$ of the braid
monodromy of the function
$$
-l_1\xi+x_2^{l_2+1}+\e_2\a(l_2+1)x_2,\quad \xi^{l_1}=1.
$$
Hence also $\eta_{l_3}(G_2)$ is contained. But then it must be a
subgroup of the versal braid monodromy of the family of functions
$$
x_1^{l_1+1}-(l_1+1)x_1+x_2^{l_2+1}+\e_2\a(l_2+1)x_2+x_3^{l_3+1}-\e_3\a(l_3+1)x_3
$$
at $\a=0$, which is $G_3$.
\qed

Analogous to the cases $n=2,3$ dealt with in \ref{ggfam},\ref{gggfam} we
obtain the generalisation to arbitrary $n$.

\begin{lemma}
\labell{ggnfam}
The versal braid monodromy $G_n$ of the family $g_\a(x_1,...,x_n)$
restricted to
the unit disc is generated by the homomorphic images of the braid
monodromy groups $\br_{BP(l_2,...,l_n)}\subset\br_{l_2\cdots l_n}$
under the primary homomorphisms $\phi_{i_1},1\leq i_1\leq l_1$.
\end{lemma}

Similarly we can extend the assertion of lemma \ref{gggband} for the pair $2,3$
to arbitrary pairs $n-1,n$.

\begin{lemma}
\labell{ggnband}
The group $G_n$ contains the image of $G_{n-1}$ under the $l_n$-band
homomorphism.
\end{lemma}

Since the generators given by Catanese and Wajnryb \cite{CW}, cf.\
\ref{perm/stab}, are cable twists we can guess that there are many cable twists
in the groups considered here.

\begin{lemma}
\labell{subcable}
Suppose $i'i_n<i'j_n$ is a pair of indices, then
$\delta_{i,j}^{j_n-i_n+2}\in G_n$, especially
$\cable^{l_n+1}\in G_n$.
\end{lemma}

\proof
If we consider the family of functions
$$
x_1^{l_1+1}-(l_1+1)x_1+(\sum_{\nu=2}^{n-1}x_\nu^{l_\nu+1}-\e_\nu(l_\nu+1)x_\nu)
+x_n^{l_n+1}-\a\e_n(l_n+1)x_n,
$$
we see that its monodromy at $\a=0$ is contained in the monodromy at $\a=0$ of
the family $g_\a(x_1,...,x_n)$.
Similarly to \ref{ggfam} the family can be shown to have
$l_1...l_{n-1}$ critical
points of type $A(l_n)$. Their local monodromies embed via compositions
$\phi_{i'}:=\phi_{i_1}\circ...\circ\phi_{i_{n-1}}$ of primary
homomorphisms into
$\br_{BP(l_1...l_n)}$, so $\d_{i,j}^{j_n-i_n+2}$ is in $G_n$ if
$\d_{i_n,j_n}^{j_n-i_n+2}$ is in $\br_{A(l_n)}$.\\[-4mm]

Since the classical geometric monodromy of the $A_n$ singularity is
\begin{eqnarray*}
\cable^{n+1} & = &
\s_{1,2}^3\s_{1,3}^2...\s_{1,n}^2\s_{2,3}^3\s_{2,4}^2
...\s_{2,n}^2...\s_{n-1,n}^3,
\end{eqnarray*}
we can deduce in fact
\begin{eqnarray*}
\delta_{i_n,j_n}^{j_2-i_2+2} & = &
\prod_{i\leq k< j} (\s_{k,k^+}^3\prod_{k^+<k'\leq j}\s_{k,k'}^2)
\quad\in\quad \br_{A(l_n)}.
\end{eqnarray*}
The additional claim follows from
$\cable^{l_2+1}=\prod_{i'}\delta_{i'1,i'l_n}^{l_n+1}$.
\qed

\begin{lemma}
\labell{coil/ggn}
For given $i'<j'$ the same braid subgroup is generated by the elements
$$
\s_{i'i_n,j'j_n}^2,\quad1\leq i_n-j_n< l_n,
$$
and by suitably chosen $G_n$-conjugates of elements
$$
\begin{array}{cl}
\s^2_{i'i_n,j'j_n} & 1\leq i_n,j_n\leq l_n, i_n,i_n^+\neq j_n,\\[1mm]
\s^2_{i'i_n,j'i_n}\s^2_{i'i_n,j'j_n}\s^{-2}_{i'i_n,j'i_n}
& 1<i_n^+=j_n\leq l_n,\\[1mm]
\s^2_{i'i_n,i'j_n}\s^2_{i'i_n,j'j_n}\s^{-2}_{i'i_n,i'j_n}
& 1<i_n^+=j_n\leq l_n.
\end{array}
$$
\end{lemma}

\proof
We introduce filtrations $T=T_3\supset T_2$ and
$S=S_3\supset S_2\supset S_1$ on the two sets of elements by
\begin{eqnarray*}
S_1 & = &
\{\s^2_{i'i_n,j'j_n}|1\leq j_n<i_n\leq l_n\}\\
S_2 & = &
S_1 \cup\{\s^2_{i'i_n,j'j_n}|1\leq i_n<i_n^+<j_n\leq l_n\}\\
T_2 & = &
T_1\cup\{\s^2_{i'i_n,j'j_n}|i_n,j_n\not\equiv 0\mod (l_n+1)\}
\end{eqnarray*}
Then it suffices to show that for each $s\in S_2$ there is a $t\in T_2$ with
$t$ equal to some $G_n$-conjugate of $s$ and that for each $s\in S_3-S_2$ there
is a $\in T_3-T_2$ with a $G_n$-conjugate of $s$ equal to $t$ conjugated by some
elements of $T_2$ and vice versa, cf.\ lemma \ref{filt}.

Since $\cable^{l_n+1}\in G_n$ by \ref{subcable}, we may conjugate the
elements of
$S_1$ by all powers $\cable^{m(l_n+1)}$.
Similarly $\d_{j'1,j'l_n}^{l_n+1}$ is an element of $G_n$ by
\ref{subcable}, hence all elements
\begin{eqnarray*}
\cable^{m(l_n+1)}(\d_{j'1,j'l_n}^{l_n+1}
\s^2_{i'i_n,j'j_n}\d_{j'1,j'l_n}^{-l_n-1})\cable^{-m(l_n+1)}
& = &
\cable^{m(l_n+1)}\s^2_{i'i_n,j'j'_n}\cable^{-m(l_n+1)},
\end{eqnarray*}
where $j_n'=j_n-l_n-1$, are $G_n$-conjugates of elements in $S_2-S_1$. In fact
the elements thus obtained are just all the elements in $T_2$.

Finally we observe that the braids \zspace$\delta_{j'i_n^+,j'l_n}^{l_n-i_n+1}$,
resp.\ $\delta_{i'1,i'j_n^-}^{j_n}$ are elements of $G_n$ due to \ref{subcable}
again. Hence we may invoke \ref{tau/sigma} to show that the elements
in $S_3-S_2$
have $G_n$-conjugates which are equal up to conjugation by elements
in $T_2$ to
elements
$$
\s_{i'i_n,j'0}^2,\s^2_{i'l_n^+,j'j_n},\quad  1\leq i_n<l_n,1<j_n\leq l_n,
$$
which in turn are contained in $T_3-T_2$.\\[-4mm]

Because conjugation by powers $\cable^{m(l_n+1)}$ yields all elements
of $T_3-T_2$
and preserves the set $T_2$, all elements in $T_3-T_2$ up to conjugation by
elements in $T_2$ are $G_n$-conjugates of elements in $S_3-S_2$, so
we are done.
\qed



\section{families of type $f_\a$}

We turn our attention now to the families of type $f_\a$, for which the versal
braid monodromy has to be computed not only locally. But in fact all geometric
insight is in the case $n=2$ dealt with in chapter \ref{Plane} and the notion
of $l$-companion families, which will be exploited in the next two section. So
here we mainly have to translate between various results and to organize them
in such a way we need for the induction.\\

\begin{lemma}
\labell{ff}
The versal braid monodromy of a family of functions $f_\a(x_1,x_2)$
is generated by the elements
$$
\begin{array}{rl}
\s_{i_1i_2,j_1j_2}^3, & i_2=j_2, \, i_1<j_1 \text{ correlated},\\
\s_{i_1i_2,j_1j_2}^2, & i_2=j_2, \, i_1<j_1 \text{ not correlated},\\
\s^2_{i_1i_2,j_1j_2} & \text{with }1\leq i_1<j_1\leq l_1,1\leq i_2-j_2<l_2.
\end{array}
$$
\end{lemma}

\proof
We have to show that the elements in \ref{ff} and those in \ref{sigma} generate
the same subgroup of $\br_{l_1l_2}$.
Since both generator sets have the elements with equal second index
component in
common, it suffices to prove that the remaining elements of each set
generate the
same braid subgroup.

Notice that both sets are filtered by level -- which is underlined -- with
\begin{eqnarray*}
S_1 & = &
\{\s^2_{i_1i_2,j_1j_2}\big|1\leq i_1<j_1\leq l_1,\underline{j_1-i_1=1},
1\leq i_2-j_2<l_2\}\\
S_2 & = &
S_1\cup\{\s^2_{i_1i_2,j_1j_2}\big|1\leq i_1<j_1\leq l_1,\underline{j_1-i_1=2},
1\leq i_2-j_2<l_2\}\\
& \vdots \\
S_{l_1} & = &
\{\s^2_{i_1i_2,j_1j_2}\big|1\leq i_1<j_1\leq l_1,1\leq i_2-j_2<l_2\}\\[2mm]
T_1 & = &
\{\tau^2_{i_1i_2,j_1j_2}\big|1\leq i_1<j_1\leq l_1,\underline{j_1-i_1=1},
1< j_2-i_2\leq l_2\}\\
T_2 & = &
T_1\cup\{\tau^2_{i_1i_2,j_1j_2}\big|1\leq i_1<j_1\leq
l_1,\underline{j_1-i_1=2},
1< j_2-i_2\leq l_2\}\\
& \vdots\\
T_{l_1} & = &
\{\tau^2_{i_1i_2,j_1j_2}\big|1\leq i_1<j_1\leq l_1, 1< j_2-i_2\leq l_2\}.
\end{eqnarray*}

For the proof we need therefore to check the hypotheses of \ref{filt} only:
The first, $S_1=T_1$, is immediate, since elements of level one coincide almost
by definition
$$
\s_{i_1i_2,i_1^+j_2}=\tau_{i_1i_2',i_1^+j_2},\quad i_2'+l_2=i_2.
$$

For the inductive hypothesis lemma \ref{ff(st)} yields, that elements
$\tau^2_{i_10,k_1k_2}$ and $\s^2_{i_1l_2,k_1k_2}$ with $1\leq
i_1<k_1\leq l_1,1\leq k_2<l_2$ are equal up to conjugation by elements
in $S_{k_1-i_1-1}\cup T_{k_1-i_1-1}$, i.e.\ by elements of smaller level.

To extend this result to the remaining elements we consider the
action of overall
conjugation by $\cable$.
Since this conjugation is level preserving, we get, that
$$
\s^2_{i_1i_2,k_1j_2}=\cable^{i_2'}\s^2_{i_1l_2,k_1k_2}\cable^{-i_2'},
\tau^2_{i_1i_2',k_1j_2}=\cable^{i_2'}\tau^2_{i_10,k_1k_2}\cable^{-i_2'}
$$
are equal up to conjugation by elements of smaller level since
$\s^2_{i_1l_2,k_1k_2},\tau^2_{i_10,k_1k_2}$ are.
The hypotheses are hence met, for each generator $\s^2_{i_1i_2,j_1j_2}$ or
$\tau^2_{i_1i_2,j_1j_2}$ is in the conjugation orbit of a
$\s^2_{i_1l_2,k_1k_2}$
resp.\ $\tau^2_{i_10,k_1k_2}$ by $\cable$.
\qed[4mm]

The generator set has still the draw back that it is not finite.
Though this could
be amended we even take a step further and proceed to a generator set which
generates the same group only up to the subgroup $G_2$.

\begin{lemma}
\labell{ff/gg}
The versal braid monodromy of $f_\a(x_1,x_2)$ is
generated by elements $G_2$-conjugate to
$$
\begin{array}{ccl}
\s_{i_1i_2,j_1j_2}^3&: & i_2=j_2, \, i_1<j_1 \text{ correlated},\\[1mm]
\s_{i_1i_2,j_1j_2}^2&: & i_2=j_2, \, i_1<j_1 \text{ not correlated},\\[1mm]
\s^2_{i_1i_2,j_1j_2} &: & i_1<j_1,i_2\neq j_2, \, i<j \text{ not
correlated},\\[1mm]
\s^2_{i_1i_2,j_1j_2}\s^2_{i_1i_2,k_1k_2}\s^{-2}_{i_1i_2,j_1j_2}
&: & i_2+1= k_2, \, i<j<k \text{ correlated}.
\end{array}
$$
\end{lemma}

\proof
We have to show, that we can assign $G_2$ conjugates to the given
elements, such
that these conjugates generate the same braid subgroup as the elements of
\ref{ff}.

The first two rows of elements obviously coincide. If we apply the case $n=2$
of lemma \ref{coil/ggn} to the last row of \ref{ff}, then we get the bottom
rows here.
\qed

\begin{lemma}
\labell{ff(cor)}
\sloppy
The versal braid monodromy of the family of functions
$$
x_1^3-3\a x_1+x_2^{l_2+1}-\e_2(l_2+1)x_2
$$
is generated by elements
$$
\begin{array}{rl}
\s_{1i_2,2j_2}^3, & i_2=j_2,\\
\s^2_{1i_2,2j_2} & \text{with }1\leq i_2-j_2<l_2.
\end{array}
$$
\end{lemma}

\proof
This is the special case $l_1=2$ of lemma \ref{ff}.
\qed

Again we have to use the results of
the next two sections.

\begin{lemma}
\labell{fff/gg}
The versal braid monodromy of $f_\a(x_1,x_2,x_3)$
is generated by elements $\eta(G_2)$-conjugate to
$$
\begin{array}{ccl}
\s_{i_1i_2i_3,j_1j_2j_3}^3 &:& i_1=j_1-1,i_2=j_2,i_3=j_3,
\\[1mm]
\s_{i_1i_2i_3,j_1j_2j_3}^2 &:& i_1^+=j_1,i_2=j_2,1\leq i_3-j_3<l_3\\[1mm]
\s_{i_1i_2i_3,j_1j_2j_3}^2 &:& i_1< j_1-1, \,i_2= j_2,
\\[1mm]
\s_{i_1i_2i_3,j_1j_2j_3}^2 &:& i_1< j_1, \,i_2\neq j_2,
i_1i_2<j_1j_2 \text{ not correlated},
\\[1mm]
\eta(\s_{i,j}^2)\s_{ii_3,kk_3}^2\eta(\s_{i,j}\invv)
&:& i_1i_2<j_1j_2 <k_1k_2 \text{ correlated}.
\end{array}
$$
\end{lemma}

\proof
First note that the discriminant family associated to
$f_\a(x_1,x_2,x_3)$ is the
$l_3$-companion family of the discriminant family associated to
$f_\a(x_1,x_2)$.
So by lemma \ref{comp/gen} the versal braid monodromy is generated by the
$l_3$-companions of the generators given in lemma \ref{ff}.

\begin{sloppypar}
But instead of the $l_3$-companions of the generators in \ref{ff} we take the
$l_3$-companions of the elements in \ref{ff/gg}. The definition implies that up
to conjugation by elements of $\eta(G_2)$ they generate the same group, which
suffices for our claim.
\end{sloppypar}

So we need only to run the list of \ref{ff/gg} through the procedure given in
\ref{std-comp} to get the list of the claim.
\qed

\begin{lemma}
\labell{fff/ggg}
\sloppy
The versal braid monodromy of the family $f_\a(x_1,x_2,x_3)$
is generated by a set of braids $G_3$-conjugate to elements of
$$
\begin{array}{ccl}
\s_{i,j}^3 &:& i_1^+=j_1,i_2=j_2,i_3=j_3,
\\[1mm]
\s_{i,k}^2 &:& i_1< k_1,
i<k \text{ not correlated},\\[1mm]
\s_{i,j}^2\s_{i,k}^2\s_{i,j}\invv &:& i_1<k_1,i<k \text{ correlated, }
j\in\{i_1^+i_2i_3,i_1k_2k_3\}.
\end{array}
$$
\end{lemma}

\proof
Since $\eta(G_2)\subset G_3$ by \ref{gggband} it suffices to show that the
elements in \ref{fff/gg} are $G_3$-conjugates of those given here.
Our method of proof consists in replacing elements of a row by others obtained
through conjugation with elements of previous rows and elements of $G_3$, since
under such transformations the group generated by $G_3$-conjugates does not
change.

Let us start with the fifth row. We have there the
$2(l_1-1)(l_2-1)l_3^2$ elements
\begin{eqnarray*}
\eta(\s_{i_1i_2,j_1j_2}^2)\s_{i_1i_2i_3,k_1k_2k_3}^2\,
\eta(\s_{i_1i_2,j_1j_2}\invv)
& &
1<i_1^+=k_1\leq l_1,1<i_2^+=k_2\leq l_2,\\
& &
1\leq i_3,k_3\leq l_3, j_1j_2\in\{i_1k_2,k_1^+i_2\}.
\end{eqnarray*}
By \ref{uncoil-i} those with $j_1j_2=i_1^+i_2$ equal \zspace[-3mm]
$\s_{i,j_1j_2i_3}^2\s_{i,k}^2\s_{i,j_1j_2i_3}\invv=
\s_{j_1j_2i_3,k}\invv\s_{i,k}^2\s_{j_1j_2i_3,k}^2$ up to conjugation
by elements of
the second row. Since $j_1=k_1$, the twist $\s_{j_1j_2i_3,k}^2$ is in $G_3$, if
$\s_{j_2i_3,k_2k_3}^2$ is in the braid monodromy group
$\br_{BP(l_2,l_3)}$, cf.\
\ref{gggfam}.

If $i<k$ are not correlated, then $j_2i_3<k_2k_3$ aren't either, with
$j_2<k_2$.
We deduce with \ref{ff/gg} that $\s_{j_1j_2i_3,k}^2$ then is in $G_3$.
Hence the elements of the fifth row with $j_1j_2=i_1^+i_2$ are equal to
\begin{enumerate}
\item
$\s_{i_1i_2i_3,j_1j_2i_3}^2\s_{i_1i_2i_3,k_1k_2k_3}^2
\s_{i_1i_2i_3,j_1j_2i_3}\invv$ if $i<k$ are
correlated,
\item
$\s_{i_1i_2i_3,k_1k_2k_3}$ if $i<k$ are not correlated,
\end{enumerate}
up to conjugation by elements of the second row and of $G_3$.

Similarly by \ref{uncoil-k} those with $j_1j_2=i_1k_2$ are equal to
\begin{enumerate}
\item[i')]
$\s_{i_1i_2i_3,j_1j_2k_3}^2\s_{i_1i_2i_3,k_1k_2k_3}^2
\s_{i_1i_2i_3,j_1j_2k_3}\invv$ if $i<k$ are
correlated,
\item[ii')]
$\s_{i_1i_2i_3,k_1k_2k_3}$ if $i<k$ are not correlated,
\end{enumerate}
up to conjugation by elements of the second row and of $G_3$.

Note that $ii)$ and $ii')$ yield the same elements and that we have
the following
restrictions on the indices, $1\leq i_\nu,j_\nu,k_\nu\leq l_\nu,\nu=1,2,3$,
assumed in all cases:
\begin{enumerate}
\item
$i_1^+i_2^+i_3,i_1^+i_2^+i_3^+=k_1k_2k_3, j_1j_2j_3=i_1^+i_2i_3$,
\item
$i_1^+i_2^+=k_1k_2,i_3,i_3^+\neq k_3$,
\item[i')]
$i_1^+i_2^+i_3,i_1^+i_2^+i_3^+=k_1k_2k_3, j_1j_2j_3=i_1k_2k_3$.
\end{enumerate}

Second we want to replace the elements of the second row. Lemma
\ref{coil/ggn} implies in case $i'=i_1i_2,j'=j_1j_2,n=3$ that the elements
$$
\s^2_{i_1i_2i_3,j_1j_2j_3},\quad i_1^+i_2=j_1j_2,1\leq i_n-j_n<l_n,
$$
are $G_3$-conjugates of elements
\begin{center}
\begin{enumerate}
\setcounter{enumi}{2}
\item
$\s^2_{i_1i_2i_3,i_1^+i_2k_3}$, $i_3,i_3^+\neq k_3$,
\item
$\s^2_{i_1i_2i_3,i_1^+i_2i_3}\s^2_{i_1i_2i_3,i_1^+i_2i_3^+}
\s_{i_1i_2i_3,i_1^+i_2i_3}\invv$,
\item
$\s^2_{i_1i_2i_3,i_1i_2i_3^+}\s^2_{i_1i_2i_3,i_1^+i_2i_3^+}
\s_{i_1i_2i_3,i_1i_2i_3^+}\invv$.
\end{enumerate}
\end{center}

We reassemble the elements thus obtained according to the form in
which they are
given, twists of exponent $3$, twists of exponent $2$ given as $\s_{i,k}^2$ and
twists of exponent $2$ given as $\s_{i,j}^2\s_{i,k}^2\s_{i,j}\invv$. The two
final steps then consists in collecting the corresponding sets of indices.

For the elements which are full twists of the form $\s_{i,k}^2$ we get:
$$
\begin{array}{cll}
& \{i,k|i_1^+i_2=k_1k_2,i_3,i_3^+\neq k_3\} & iii)\\[1mm]
\cup & \{i,k|i_1^+< k_1, \,i_2= k_2\} & \text{$3^{rd}$ row of
\ref{fff/gg}}\\[1mm]
\cup & \{i,k|i_1< k_1, \,i_2\neq k_2,i_1^+i_2^+\neq k_1k_2 \} &
\text{$4^{th}$ row of \ref{fff/gg}}\\[1mm]
\cup & \{i,k|i_1^+i_2^+=k_1k_2, i_3,i_3^+\neq k_3\} & ii)\\[2mm]
= & \{i,k|i_1<k_1, i<k \text{ not correlated}\} & 2^{nd}\text{ row of
\ref{fff/ggg}}
\end{array}
$$

For the elements which are full twists of the form
$\s_{i,j}^2\s_{i,k}^2\s_{i,j}\invv$ we get:
$$
\begin{array}{clr}
&
\{i,j,k|j_1j_2j_3=i_1^+i_2i_3,k_1k_2k_3\in\{i_1^+i_2^+i_3,i_1^+i_2^+i_3^+\}
\}
& i)\\[1mm]
\cup &
\{i,j,k|k_1k_2k_3\in\{i_1^+i_2^+i_3,i_1^+i_2^+i_3^+\},j_1j_2j_3=i_1k_2k_3\}
& i')\\[1mm]
\cup & \{i,j,k|j_1j_2j_3=i_1^+i_2i_3,k_1k_2k_3=i_1^+i_2i_3^+\}
& iv)\\[1mm]
\cup &\{i,j,k|j_1j_2j_3=i_1i_2i_3^+,k_1k_2k_3=i_1^+i_2i_3^+\}
& v)\\[2mm]
= & \{i,j,k|i_1<k_1,i<k \text{ correlated, }
j\in\{i_1^+i_2i_3,i_1k_2k_3\}\}
\end{array}
$$

So the set we obtained replacing the elements of \ref{fff/gg} by
$G_3$-conjugates
of another generating set is exactly that of our claim, and we are done.
\qed

\begin{lemma}
\labell{ffn/ggn}
The versal braid monodromy of the family $f_\a(x_1,...,x_n)$
is generated by a set of $G_n$-conjugates of braids
$$
\begin{array}{ccl}
\s_{i,j}^3 &:& i_1^+=j_1,i_2...i_n=j_2...j_n,
\\[1mm]
\s_{i,k}^2 &:& i_1< k_1,
i<k \text{ not correlated},\\[1mm]
\s_{i,j}^2\s_{i,k}^2\s_{i,j}\invv &:& i_1^+=k_1,i<k \text{ correlated, }
j\in\{i^+_1,i_2...i_n,i_1k_2...k_n\}.
\end{array}
$$
\end{lemma}

\proof
This can be proved by induction on $n$, with the cases $n=2,3$ already done,
\ref{ff/gg},\ref{fff/ggg}.
So let us assume the claim holds for $n-1$.

By \ref{comp/gen}
generators are given by the $l_n$-companions of
generators for the braid monodromy group of $f_\a(x_1,...,x_{n-1})$.

If such generators are given up to conjugation by elements of
$G_{n-1}$ only, then
their $l_n$-companions are known up to conjugation by elements of
$\eta(G_{n-1})$ due to
the definition. This is fine here, since we are interested in generators up
to conjugation by elements of $G_n$ which is even coarser because
$\eta(G_{n-1})\subset G_n$ by \ref{ggnband}.

Hence we conclude that up to conjugation by elements of $G_n$
generators for the braid monodromy of $f_\a(x_1,...x_n)$ are obtained from the
list of $G_{n-1}$-conjugates of generators for
$f_\a(x_1,...,x_{n-1})$ by taking the
$l_n$-companions according to the rule \ref{std-comp}:
$$
\begin{array}{cccl}
\s_{i',k'}^3 & : & \s_{i,k}^3 & i_1^+i_2...i_n=k_1k_2...k_n,
\\[1mm]
& & \s_{i,k}^2, & i_1^+i_2...i_{n-1}=k_1k_2...k_{n-1},\\
& & &1\leq i_n-k_n<l_n
\\[1mm]
\s_{i',k'}^2 & : & \s_{i,k}^2, & i_1< k_1,
i'<k' \text{ not correlated},\\[1mm]
\s_{i',j'}^2\s_{i',k'}^2\s_{i',j'}\invv & : &
\eta(\s_{i',j'}^2)\s_{i,k}^2\eta(\s_{i',j'}\invv) & i_1<k_1, i'<k'
\text{ correlated},\\
&&& j'\in\{i^+_1i_2...i_{n-1},i_1k_2...k_{n-1}\}.
\end{array}
$$
We proceed in strict analogy to the proof of \ref{fff/ggg}. Without
changing the
group which $G_n$-conjugates generate, we replace the elements of
the fourth row by
\begin{enumerate}
\item
$\s_{i,j'i_n}^2\s_{i,k}^2\s_{i,j'i_n}\invv$, if $i<k$ correlated,
$j'=i_1^+i_2...i_{n-1}$,
\item[i')]
$\s_{i,j'k_n}^2\s_{i,k}^2\s_{i,j'k_n}\invv$, if $i<k$ correlated,
$j'=i_1k_2...k_{n-1}$,
\item[ii)]
$\s_{i,k}^2$ if $i<k$ not correlated, $i_1^+=k_1, i'<k'$ correlated.
\end{enumerate}
We go on and replace the elements of the second row by elements
\begin{enumerate}
\setcounter{enumi}{2}
\item
$\s^2_{i_1...i_n,i_1^+i_2...i_{n-1}k_n}$, $i_n,i_n^+\neq k_n$,
\item
$\s^2_{i,j}\s^2_{i,k}\s_{i,j}\invv$,
$i_2...i_{n-1}=j_2...j_{n-1}=k_2...k_{n-1}$,
$i_1^+i_n^+=j_1j_n^+=k_1k_n$,
\item
$\s^2_{i,j}\s^2_{i,k}\s_{i,j}\invv$,
$i_2...i_{n-1}=j_2...j_{n-1}=k_2...k_{n-1}$,
$i_1^+i_n^+=j_1^+j_n=k_1k_n$.
\end{enumerate}
A final check that the lists of index pairs and triples for the generators thus
obtained and the elements of the claim coincide completes the proof.
\qed

\begin{thm}
\labell{ffnggn}
The braid monodromy of the function $x_1^{l_1+1}+\cdots+x_n^{l_n+1}$ is
generated by
\begin{eqnarray*}
S & = & \left\{
\begin{array}{ccl}
\s_{i,k}^3 & i<k & \text{correlated},\\
\s_{i,k}^2 & i<k & \text{not correlated},\\
\s_{i,j}^2\s_{i,k}^2\s_{i,j}^{-2} & i<j<k &
\text{correlated}
\end{array}\right\}
\end{eqnarray*}
\end{thm}

\proof
The proof can be obtained by induction on $n$. Since the claim has been proved
already for $n=1,2$ we may suppose $n>2$ and that the case $n-1$ is
already known.

So we deduce from \ref{ggnfam} that the braid monodromy group $G_n$ is
generated by
\begin{eqnarray*}
\s_{i,k}^3 && i_1=k_1,i<k  \text{ correlated},\\
\s_{i,k}^2 && i_1=k_1,i<k  \text{ not correlated},\\
\s_{i,j}^2\s_{i,k}^2
\s_{i,j}^{-2} && i_1=k_1,i<j<k  \text{ correlated}
\end{eqnarray*}
To get generators for the total braid monodromy we have -- by
\ref{ZB} -- to add
the elements of \ref{ffn/ggn}.
$$
\begin{array}{cl}
\s_{i,j}^3 & i_1^+=j_1,i_2...i_n=j_2...j_n,
\\[1mm]
\s_{i,k}^2, & i_1< k_1,
i<k \text{ not correlated},\\[1mm]
\s_{i,j}^2\s_{i,k}^2\s_{i,j}\invv & i_1^+=k_1,i<k \text{ correlated, }
j\in\{i^+_1,i_2...i_n,i_1k_2...k_n\}.
\end{array}
$$
By a check on the indices occurring in these two sets, we see, that
in order to get
the claim we have to add elements
$$
\begin{array}{cl}
\s_{i,j}^3 & i_1^+=j_1,i_2...i_n\neq j_2...j_n,
\\[1mm]
\s_{i,j}^2\s_{i,k}^2\s_{i,j}\invv & i_1^+=k_1,i<j<k \text{ correlated, }
j\not\in\{i^+_1i_2...i_n,i_1k_2...k_n\}.
\end{array}
$$
Of course they may be added without harm if and only if they are
elements of the
braid monodromy.

So let us first consider triples of correlated indices $i<j<k$ with
$j_1=i_1^+$.
Then $i<i_1^+i_2...i_n<j<k$ is a correlated quadruple
such that the
full twists associated to the correlated triples of $i<k$, $i<j$
and $j<k$ with $i_1^+i_2...i_n$
are among the given generators of the braid monodromy.
Hence we may conclude with \ref{braid/aa} that all elements
$\s_{i,j}^2\s_{i,k}^2\s_{i,j}\invv$ with $j_1=i_1^+$, $i<j<k$
correlated, are in
the braid monodromy.

Similarly we argue for correlated index triples $i<j<k$ with $j_1=i_1$.
Since $i<j<i_1k_2...k_n<k$ is a correlated quadruple
then, we may conclude as above that all elements
$\s_{i,j}^2\s_{i,k}^2\s_{i,j}\invv$ with $j_1=i_1$, $i<j<k$ correlated, are in
the braid monodromy.

Finally let $i<k$ be any pair of correlated indices.
If $\s_{i,k}^3$ is not among the generators given for the braid monodromy, set
$j:=i_1^+i_2...i_n$.
Then $i<j<k$ is a correlated triple such that $\s_{i,j}^3$ and $\s_{j,k}^3$ are
among the given generators of the braid monodromy.
Since so is $\s_{i,j}^2\s_{i,k}^2\s_{i,j}\invv$, we conclude with \ref{braid/a}
that also $\s_{i,k}^3$ is an element of the braid monodromy.
\qed


\newcommand{\cfam}{\tilde\efami\ups{l}} 
\newcommand{\xii}{\xi_{\nu_1}}
\newcommand{\xij}{\xi_{\nu_2}}
\newcommand{\Bd}{B_{\d_0}(z_i)}
\newcommand{\Bz}{B_{\d_0}(z_i)}
\newcommand{\efamic}{\efami\ups{l}}


\section{$l$-companion models}

We introduce in this section the notion of a companion family, which
is associated to a discriminant family by replacing the singular value
divisor by a number of parallel copies.
If the discriminant family arises from a family of functions, so does the
companion family, hence we may try to relate not only the braid monodromy but
also the versal braid monodromies of both families.

In this section
we have a closer look at the companion families of the simplest discriminant
families, the ordinary node and the ordinary cusp, and provide the necessary
definitions and arguments to get first results on these companion families and
their monodromy.

\begin{defi}
Given a discriminant family $\efami$ defined by a polynomial $d(z,\a)$, a
discriminant family $\efami\ups{l}$ defined by a polynomial
$d\ups{l}$ is called {\it $l$-companion family} if
\begin{eqnarray*}
d\ups{l}(z,\a) & = &
\prod_{\xi^l=1}d(z-\xi\e,\a)
\end{eqnarray*}
for some $0<\e\ll1$.
\end{defi}

So for example we get an $l$-companion family of the ordinary node $z^2-\a^2$
defined by the equation
\begin{eqnarray*}
\prod_{\xi^l=1}\big((z-\xi\e)^2-\a^2\big) & = & 0.
\end{eqnarray*}

The close relation between the two equations should -- and shall -- result in
more ties between the two families $\efami$ and $\efami\ups{l}$.\\

First note that there is a natural way to pass from a distinguished system of
paths for a regular fibre of $\efami$ at a parameter $\a$ to a system of paths
for the fibre of $\efamic$ at $\a$.

Of course we have to impose $\e$ to be small enough to get again a
regular fibre.
For a suitable choice of $\e$ then the given system of paths for $\efami$ meets
the boundary of the $2\e$-discs at the punctures once only. By changing the
system in its isotopy class only, we may even assume that this boundary point
corresponds to puncture translated by $2\e$.

We then split each path into a bunch of $l$ paths and connect it with the
punctures of the $\efamic$ fibre in the order of increasing angle $arg(\xi)$.

The system of paths thus obtained we call the {\it Hefez-Lazzeri
refinement}, as
the procedure mimics the iteration step given in the Hefez-Lazzeri
article for a
special system of paths.

\begin{examp}
Cut discs of radius $2\e$ off the original fibre and replace the
truncated paths
by $l$-bunches.

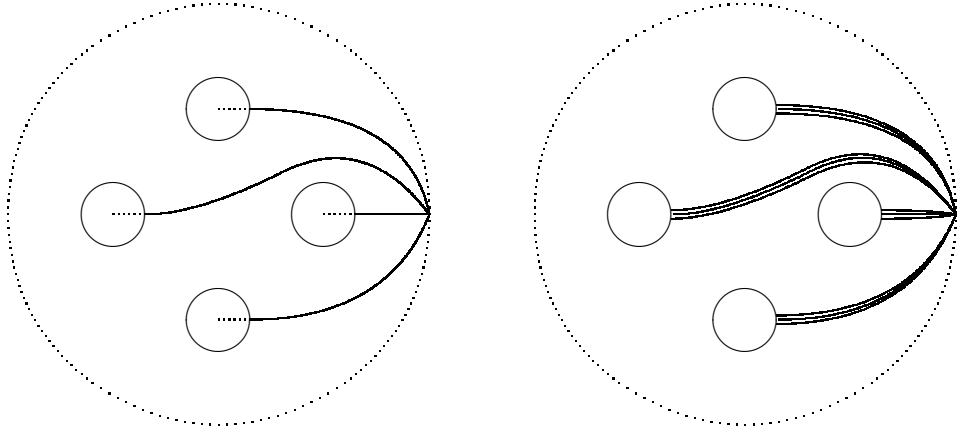
\begin{center}
\setlength{\unitlength}{2.8mm}
\begin{picture}(26,20)(0,-10)

\put(-10,-10){\setlength{\unitlength}{5.6mm}\begin{picture}(10,10)

\bezier{30}(0,5)(0,7.5)(2,9)
\bezier{30}(2,9)(5,11)(8,9)
\bezier{30}(8,9)(10,7.5)(10,5)
\bezier{30}(10,5)(10,2.5)(8,1)
\bezier{30}(8,1)(5,-1)(2,1)
\bezier{30}(2,1)(0,2.5)(0,5)
\end{picture}}

\put(15,-10){\setlength{\unitlength}{5.6mm}\begin{picture}(10,10)

\bezier{30}(0,5)(0,7.5)(2,9)
\bezier{30}(2,9)(5,11)(8,9)
\bezier{30}(8,9)(10,7.5)(10,5)
\bezier{30}(10,5)(10,2.5)(8,1)
\bezier{30}(8,1)(5,-1)(2,1)
\bezier{30}(2,1)(0,2.5)(0,5)
\end{picture}}

\put(5,0){\circle{3}}
\put(-5,0){\circle{3}}
\put(0,-5){\circle{3}}
\put(0,5){\circle{3}}

\bezier{6}(5,0)(5.75,0)(6.5,0)
\bezier{6}(-5,0)(-4.25,0)(-3.5,0)
\bezier{6}(0,-5)(.75,-5)(1.5,-5)
\bezier{6}(0,5)(.75,5)(1.5,5)

\bezier{360}(6.5,0)(8,0)(10,0)
\bezier{360}(1.5,5)(9,5)(10,0)
\bezier{360}(1.5,-5)(8,-5)(10,0)
\bezier{360}(-3.5,0)(-1,0)(3,2)
\bezier{360}(3,2)(7,4)(10,0)

\put(25,0){\setlength{\unitlength}{.28mm}\begin{picture}(10,10)

\bezier{135}(65,0)(90,0)(100,0)
\bezier{135}(65,2)(90,2)(100,0)
\bezier{135}(65,-2)(90,-2)(100,0)
\bezier{180}(16,50)(90,50)(100,0)
\bezier{180}(15,52)(90,52)(100,0)
\bezier{180}(15,48)(90,48)(100,0)
\bezier{180}(15,-48)(80,-48)(100,0)
\bezier{180}(16,-50)(80,-50)(100,0)
\bezier{180}(15,-52)(80,-52)(100,0)
\bezier{140}(-35,2)(-10,2)(30,22.3)
\bezier{140}(-34,0)(-10,0)(31,20.5)
\bezier{140}(-35,-2)(-10,-2)(30,17.7)
\bezier{140}(30,22.3)(70,42)(100,0)
\bezier{140}(31,20.5)(70,40)(100,0)
\bezier{140}(30,17.7)(70,38)(100,0)

\end{picture}}

\put(30,0){\circle{3}}
\put(20,0){\circle{3}}
\put(25,-5){\circle{3}}
\put(25,5){\circle{3}}

\end{picture}
\end{center}

Glue in copies of a $2\e$ disc punctured at points of absolute value $\e$ and
argument $\xi,\xi^l=1$, and provided with a standard system of paths.

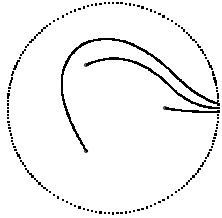
\begin{center}
\setlength{\unitlength}{1.4mm}
\begin{picture}(10,20)(-4,-10)

\put(-10,-10){\setlength{\unitlength}{2.8mm}\begin{picture}(10,10)

\bezier{30}(0,5)(0,7.5)(2,9)
\bezier{30}(2,9)(5,11)(8,9)
\bezier{30}(8,9)(10,7.5)(10,5)
\bezier{30}(10,5)(10,2.5)(8,1)
\bezier{30}(8,1)(5,-1)(2,1)
\bezier{30}(2,1)(0,2.5)(0,5)
\end{picture}}

\put(5,0){\circle{.3}}
\put(-2.5,4.1){\circle{.3}}
\put(-2.5,-4.1){\circle{.3}}

\bezier{60}(5,0)(7,-.4)(10,-.4)
\bezier{90}(-2.5,4.1)(1.5,6)(5.5,2)
\bezier{40}(5.5,2)(7.5,0)(10,0)

\bezier{120}(-2.5,-4.1)(-7,3)(-3,6)
\bezier{110}(-3,6)(1,8)(6,3)
\bezier{60}(6,3)(8,1)(10,.4)

\end{picture}
\end{center}
\end{examp}

So from now on given a system of paths for a fibre of $\efami$, we may tacitly
assume that a system of paths in the fibre of $\efamic$ is given by the
Hefez-Lazzeri refinement. We will still sometimes say so explicitly but even if
we won't this should be understood without mentioning.\\

In particular this means that given the braid monodromy group of a family as a
subgroup of an abstract braid group $\br_{l'}$, which involves an
implicit choice
of a system of paths, we may consider the braid monodromy group of a
$l$-companion family to be a subgroup of $\br_{l'l}$ since we may make the
necessary identification by implicit use of the Hefez-Lazzeri refined system of
paths.

Now let us have some exercise with the simplest discriminant families:

\begin{lemma}
\labell{node/comp}
The versal braid monodromy of the $l$-companion
family associated to the discriminant family defined by $z^2-\a^2$
with defining
equation
\begin{eqnarray*}
\prod_{\xi^l=1}\big((z-\xi\e)^2-\a^2\big) & = & 0
\end{eqnarray*}
is generated by
\begin{eqnarray*}
\s_{1i,2j}^2, & & 1\leq i,j\leq l.
\end{eqnarray*}
and is isomorphic to $\ker(\pbr_{2l}\to\pbr_l\times\pbr_{l_2})$.
\end{lemma}

\proof
The $l$-companion family in this case coincides with the model discriminant
family for the family of function $f_\a(x_1,x_2)$ in case
$l_1=2, l_2=l$. Hence the versal braid monodromy is almost that of the
discriminant family of $f_\a(x_1,x_2)$. We need only remark that the group
locally assigned to the singular fibre at the origin is generated by
full twists
on v-arcs as in the case of all other fibres. Then the same methods show then
that
$$
\tau_{1i,2j},\quad1\leq j-i\leq l_2
$$
generate the braid monodromy.
But this is just the claim up to inner conjugations.
\qed

\begin{lemma}
\labell{cusp/comp}
\sloppy
The versal braid monodromy of the $l$-companion
family associated to the discriminant family defined by $z^2-\a^3$ is equal to
the subgroup of
$
\br_{BP(2,l)}
$
generated by
$$
\begin{array}{rl}
\s_{1i_2,2j_2}^3, & i_2=j_2,\\
\s^2_{1i_2,2j_2} & \text{with }1\leq i_2-j_2<l_2.
\end{array}
$$
\end{lemma}

\proof
The $l$-companion family associated to the discriminant family defined by
$z^2-\a^3$ coincides with the discriminant family of the family of functions
$$
x_1^3-3\a x_1+x_2^{l_2+1}-\e_2(l_2+1)x_2.
$$
The versal braid monodromy has be computed in
\ref{ff(cor)}.
\qed[3mm]

In order to avoid giving lists of elements as in the preceding lemmas, we
introduce some new notions.

\begin{defi}
An element $g\in\br_n$ is called {\it positive twist element} of {\it exponent}
$k\in\NN$, if there is
$h\in\br_n$ such that $hgh\inv=\s_1^k$.
\end{defi}

\begin{defi}
The {\it $l$-companions} in $\br_{l'l}$ of $\s_1^2\in\br_{l'}$ are the elements
\begin{eqnarray*}
\s_{1i,2j}^2, & & 1\leq i,j\leq l,
\end{eqnarray*}
the {\it $l$-companions} in $\br_{l'l}$ of $\s_1^3\in\br_l$ are the elements
\begin{eqnarray*}
\s_{1i,2i}^3, & & 1\leq i\leq l,\\
\s^2_{1i,2j}, & & 1\leq i-j<l.
\end{eqnarray*}
\end{defi}

\begin{defi}
If $g\in\br_{l'}$ is a twist of exponent $2$ and $hgh\inv=\s_1^2$,
then the {\it
$l$-companions} in
$\br_{l'l}$ of $g\in\br_{l'}$ are the elements
\begin{eqnarray*}
\eta(h\inv)\s_{1i,2j}^2\eta(h), & & 1\leq i,j\leq l.
\end{eqnarray*}
\end{defi}

\begin{defi}
If $g\in\br_{l'}$ is a twist of exponent $3$ and $hgh\inv=\s_1^3$,
then the {\it
$l$-companions} in
$\br_{l'l}$ of $g\in\br_{l'}$ are the elements
\begin{eqnarray*}
\eta(h\inv)\s_{1i,2i}^3\eta(h), & & 1\leq i\leq l,\\
\eta(h\inv)\s_{1i,2j}^2\eta(h), & & 1\leq i-j< l.
\end{eqnarray*}
\end{defi}

Of course in the last two definitions there arises the question of
well-definedness. But this is easily taken care of:

\begin{lemma}
\labell{stab-indp}
The definition of $l$-companions is independent on the \mbox{choice of $h$}.
\end{lemma}

\proof
If we replace $h$ by $h'h$ with $h'\s_1=\s_1h'$, i.e.\ $h'$ in the stabilizer
of $\s_1$, then $\eta_l(h')$ commutes with $\s_{1i,2j}$.
\qed[3mm]

With these definitions we can easily determine $l$-companions of
various elements:

\begin{lemma}
\labell{std-comp}
The following table lists the $l$-companions of the elements given in the left
column:
$$
\begin{array}{rccl}
\s_{i',j'}^2 &:& \s_{i'i_n,j'j_n}^2, & i_n,j_n\leq l,\\[4mm]
\s_{i',j'}^3 &:& \s_{i'i_n,j'j_n}^3, & i_n=j_n< l,\\
&& \s_{i'i_n,j'j_n}^2, & 1\leq i_n-j_n< l,\\[4mm]
\s_{i',j'}^2\s_{i',k'}^2\s_{i',j'}\invv &:&
\eta(\s_{i',j'}^2)\s_{i'i_n,k'k_n}^2\eta(\s_{i',j'}\invv), & i_n,j_n\leq l.
\end{array}
$$
\end{lemma}

\begin{lemma}
\labell{multinode/comp}
The versal braid monodromy of the $l$-companion
family associated to the discriminant family defined by a homogeneous
polynomial
$d(z,\a)$ of degree $l_1$ with defining equation
\begin{eqnarray*}
\prod_{\xi^{l_2}=1}d(z-\xi\e,\a) & = & 0
\end{eqnarray*}
is generated by
\begin{eqnarray*}
\s_{i_1i_2,j_1j_2}^2, & & 1\leq i_1,j_1\leq l,1\leq i_2,j_2\leq 2.
\end{eqnarray*}
and is isomorphic to $\ker(\pbr_{l_1l_2}\to\times_{l_1}\pbr_l)$.
\end{lemma}

\proof
Without loss of generality we may assume $d(z,\a)=z^{l_1}-\a^{l_1}$.
The $l_2$-companion family in this case coincides with the model discriminant
family for the family of function $f_\a(x_1,x_2)$.
Hence the versal braid monodromy
is almost that of the discriminant family of $f_\a(x_1,x_2)$. We need
only remark
that the group locally assigned to the singular fibre at the origin
is generated
by full twists on v-arcs as in the case of all other fibres. Then the
same methods
show then that
$$
\tau_{i_1i_2,j_1j_2},\quad1\leq i_1,i_2\leq l_1,\,1\leq j_2-i_2\leq l_2
$$
generate the braid monodromy.
But this is just the claim up to inner conjugations.
\qed

\begin{lemma}
\labell{An/comp}
\sloppy
The versal braid monodromy of the $l_2$-companion
family associated to the discriminant family of the families of functions
$$
x^{l_1+1}-\a(l_1+1)x
$$
is equal to the subgroup of $\br_{l_1l_2}$
generated by the $l_2$-companions of $\s_{i_1,i_1^+}^3$ and $\s_{i_1,j_1}^2$,
$1<i_1^+<j_1\leq l_1$.
\end{lemma}

\proof
The $l_2$-companion family coincides with the discriminant family of the family
of functions
$$
x_1^{l_1+1}-\a(l_1+1)x_1+x_2^{l_2+1}-\e_2(l_2+1)x_2.
$$
A set of generators for the versal braid monodromy has been given in lemma
\ref{ff}. By close inspection these are the $l_2$-companions of the given
elements.
\qed[3mm]


\section{$l$-companion monodromy}

The aim of this section is to compute the versal braid monodromy
of the $l$-companion family $\efami\ups{l}$ associated to a
family $\efami$ in terms of the braid monodromy data given for $\efami$.

So we suppose $\efami$ is a discriminant family with critical parameters
confined to the interior of the unit disc $\EE$.
\\

Moreover we suppose that
suitable strongly distinguished systems of paths are given in the
parameter disc
with base point at $\a=1$ and in the corresponding fibre.

The first provides us with a basis for the fundamental group of the base of the
punctured disc bundle associated to $\efami$, the second is used to provide an
isomorphism between the mapping class group of the punctured fibre at
$\a=1$ with
the abstract braid group $\br_{l'}$ on $l'$ strands in bijection to
the punctures.

With these data given, the versal braid monodromy is the subgroup of
$\br_{l'}$ generated by the parallel transport of generators of the
groups locally
assigned to regular fibres close to singular ones.
\\

 From now on we pursue the strategy to build up the monodromy for
$\efami\ups{l}$
from small pieces related to the critical parameters of $\efami$ in order to
obtain a close relation between the monodromies over these pieces in
terms of the
results of the previous section.
To do so we actually need only to choose $\e$ sufficiently small, but to make
this precise we start with the following technical tool:

\begin{defi}
Given a monic polynomial $p$ in $\CC[z]$, define the {\it root distance}
by
$$
\delta(p)=\min\{|z_1-z_2|\big| (z-z_1)(z-z_2) \text{ divides }p\}.
$$
\end{defi}

\begin{defi}
Given a discriminant polynomial $p\in\CC[z,\a]$ which is monic in the variable
$z$. Define the {\it discriminant root distance} to be the function
\begin{eqnarray*}
\delta(\a) & := & \delta(p_\a),
\end{eqnarray*}
which is a continuous non-negative function.
\end{defi}

\begin{lemma}['minimum principle']
\labell{min/princ}
If $\delta$ has a strictly positive minimum on a bounded domain, then
it attains
the minimum on the boundary.
\end{lemma}

\proof
Assume, that $\delta$ attains the minimum at an interior point
$\a_0$. Since by assumption $\delta(\a_0)>0$, there is a neighbourhood $U$ of
$\a_0$ such that $p(z,\a)$ is given over $U$ as the union of graphs of
holomorphic functions $\phi_i:U\to\CC$. Hence $\delta$ is given over $U$ by
$\min_{i\neq j}|\phi_i(\a)-\phi_j(\a)|$. By the minimum principle for
non-vanishing holomorphic functions, we conclude, that $\delta$ is constant on
$U$. But we can extend this argument to the whole closure of our domain, so
$\delta$ attains the minimum on the boundary, too.
\qed[3mm]

Suppose the critical points of $\efami$ are $z_{i,\nu},\a_i$, which are the
critical points for the projection to the parameter space of the zero
set of the
defining polynomial $p=p(z,\a)$.

We then want to introduce small discs $U_i:= B_\eta(\a_i)$ in the
parameter base
subjected to a lot of properties which actually all hold if we choose $\eta$
appropriately small enough.
\begin{enumerate}
\item
Since the $\a_i$ are distinct and in the
interior of $\EE$, we may impose $U_i\subset\EE$, and $U_i\cap U_{i'}=\{\,\}$
if $i\neq i'$.
\item
The distinguished system of paths intersects the boundary of each $U_i$ in a
unique point $\a_i'$. So its restriction to the complement of the $U_i$ is a
distinguished system of paths for the $U_i$.
\item
At a critical parameter $\a_i$ the polynomial $p_{\a_i}$ of a generic
discriminant family factors as
$$
\prod_\nu (z-z_{i,\nu})^{m_{i,\nu}}
$$
according to the multiplicities $m_{i,\nu}$ of the punctures $z_{i,\nu}$.
So Hensel's lemma yields -- $\eta$ again assumed small enough -- a
factorization
over $U_i$
$$
p=\prod q_{i,\nu}, \text{ with }q_{i,\nu}(z,\a_i)=(z-z_{i,\nu})^{m_{i,\nu}}.
$$
\item
There is a positive constant $\d_0$ such that for
$\a\in U_i$:
\begin{eqnarray*}
q_{i,\nu}(z,\a)=0 & \implies & z\in B_{\frac12\d_0}(z_{i,\nu})
\end{eqnarray*}
and roots of different factors are at least $2\d_0$ apart from each other.
\end{enumerate}

Next we introduce $\EE':=\EE-\bigcup U_i$ and show that the discriminant root
distance is bounded away from $0$ on $\EE'$.

\begin{lemma}
\labell{dist/bound}
There is a positive constant $\d_\eta$ such that for the defining
polynomial $p$
of the discriminant family $\efami$
$$
\d(\a) >\d_\eta,\quad\forall\,\a\in\EE'.
$$
\end{lemma}

\proof
The set $\bar\EE'$ is a domain on which $\d$ is strictly positive. Moreover by
the minimum principle for $\d$ the function attains its minimum on $\del \EE'$
which we may take as $2\d_\eta$ to get the claim.
\qed[3mm]

We can now prove a kind of continuity of critical parameters.

\begin{lemma}
\labell{crit/conti}
There is a locally generic perturbation $\efami\ups{l}$ of an
$l$-companion family
of $\efami$, such that all critical parameters in $\EE$ are confined
to $\bigcup
U_i$.
\end{lemma}

\proof
First we prove the claim for the $l$-companion family itself, which
is defined by
a polynomial $p\ups{l}:=\prod p(z-\xi\e,\a)$.
Since a parameter is critical if and only if it is in the zero locus of the
discriminant root distance $\d,$ it suffices to show that $\d$ does not vanish
on $\EE'$.

For $\a\in\EE'$ by lemma \ref{dist/bound} the distance of roots of $p_\a$ is
bounded from below by some positive number $\d_\eta$.
If $z_\k$ are the roots of $p_\a$, then $z_\k+\xi_\nu\e$ are those of
$p\ups{l}_\a$. The distance of a pair of roots is then
\begin{eqnarray*}
|z_{\k_1}+\xii\e-z_{\k_2}-\xij\e| & = &
|\xii\e-\xij\e|\,\neq\,0,\quad\text{in case }\k_1=\k_2,\nu_1\neq\nu_2,\\
|z_{\k_1}+\xii\e-z_{\k_2}-\xij\e| & \geq &
|z_{\k_1}-z_{\k_2}|-|\xii\e-\xij\e|\\
& > & \d_\eta-2\e,\quad\text{in case }\k_1\neq\k_2.
\end{eqnarray*}
So if we impose $\e<\frac12\d_\eta$, we get the claim for $\efami\ups{l}$.
\qed[3mm]

Now that we got little discs $U_i$ carrying all the degenerations and hence the
local contributions to the monodromy, we may take the next step and
compute the monodromy over $U_i$ for a local generic perturbation
of $\efami\ups{l}|_{U_i}$.

Let us first have a look at the monodromy of $\efami|_{U_i}$ in the fibre over
$\a_i'\in\del U_i$. The fibre contains the discs $B_{\d_0}(z_{i,\nu})$.
Outside polydiscs $U_i\times B_{\d_0}(z_{i,\nu})$ the fibration is locally
trivial, therefore all mapping classes of the monodromy can be realized by
diffeomorphisms supported on $B_{\d_0}(z_{i,\nu})$.
In fact they are given by the locally assigned group.\\

Let us choose now a strongly distinguished system of paths in the fibre of
$\efami$ at $\a_i'$. We refine this system of paths to a system
for the fibre of $\efami\ups{l}$ following the Hefez-Lazzeri construction.
Note that the local triviality in the complement of $U_i\times\Bd$
carries over to $\efami\ups{l}$.\\

In the final step we want to obtain generators for the braid monodromy group of
the $l$-companion family, which we identify with a subgroup of
$\br_{l'l}$ by the
choice of the Hefez-Lazzeri refined system of paths in the fibre at $\a=1$.

Since we have got a system of paths for the $U_i$, it suffices to understand
parallel transport to $\a=1$ to get generators for the global monodromy group
from local ones.\\

The chosen system of paths in the fibre at $\a_i'$ for $\efami$ and its
Hefez-Lazzeri refinement for $\efami\ups{l}$ yield system of paths in
the fibres
at $\a=1$ under parallel transport along the truncated system of
paths in $\EE'$.
This has to be subjected to a thorough investigation in order to determine
properties which are preserved:

\begin{remark}
Recall the existence of smooth {\it 'bump'} functions
$\chi,\chi_\e:\CC\to\RR$ for
any real $\e>0$:
\begin{eqnarray*}
  \chi &:& 0\leq\chi(z)=\chi(|z|)\leq 1,\chi(z)=0\text{ if }|z|\geq1,
\chi(z)=1\text{ if }|z|\leq\frac12,\\
  \chi_\e &:& \chi_\e(z)=\chi(z/\e),
\end{eqnarray*}
with support contained in the unit disc, resp.\ the disc of radius $\e$.
\end{remark}

\begin{lemma}
\labell{vec-tube}
Suppose $\cg$ parameterizes a path in $\EE'$.
Given a vector field on $\cg^*\efami$ such that the punctures form
integral curves,
then for $\e$ sufficiently small a vector field can be found such that all
parallels of distance bounded by $\e$ form integral curves.
\end{lemma}

\proof
Let $v(z,t)$ be the given field, $z_i(t)$ the punctures, then define
\begin{eqnarray*}
v_\e(z,t) & := &
v(z-\sum_i\chi_{2\e}(z-z_i(t))(z-z_i(t)),t).
\end{eqnarray*}
The parallels of punctures of distance bounded by $\e$ are given by functions
$z_{i,c}(t):=z_i(t)+c$ with $|c|\leq\e$.

Since $|z_{i,c}(t)-z_j(t)|\geq2\e$ for $i\neq j$ (given that
$|z_i-z_j|\geq3\e$)
all summands except for the $i^{th}$ vanish at $z_{i,c}$, hence
$$
v_\e(z_{i,c}(t),t)\,=\,
v(z_{i,c}-\chi_{2\e}(c)c,t)\,=\,
v(z_i,t).
$$
The claim is therefore proved, because $\frac\del{\del t}z_{i,c}=\frac\del{\del
t}z_i=v(z_i,t)$ by hypothesis.
\qed

\begin{lemma}
If parallel transport in $\efami$ over a path in $\EE'$ is realized by a
diffeomorphism
$\phi$, then parallel transport in $\efami$ and $\efami\ups{l}$ is
realized by a
diffeomorphism
defined by
\begin{eqnarray*}
\phi\ups{l}(z) & := &
\phi(z)+\sum_i\chi_{2\e}(z-z_i)(\phi(z_i)-z_i-\phi(z)+z).
\end{eqnarray*}
\end{lemma}

\proof
If $\phi$ is the integrated flow of a vector field $v$, just take $\phi\ups{l}$
to be the integrated flow of the vector field provided by \ref{vec-tube}.
\qed[3mm]

We conclude that parallel transport commutes with Hefez-Lazzeri refinement: The
parallel transport of a Hefez-Lazzeri refinement is the Hefez-Lazzeri
refinement
of the parallel transport.\\

The change of path system from the transported one to the given one 
on the $\a=1$
fibre is realized by the action of a suitable diffeomorphism. This induces
a map
on the braid group given by conjugation with the corresponding mapping class.

The same considerations apply to the family $\efami\ups{l}$ and the
corresponding
mapping class is the image under the band homomorphism $\eta_l$ of the former.
\\

So now we can determine the generators transported from the fibre at $\a_i'$.
We know the transported generator of $\efami$. This determines the conjugating
mapping class up to stabilizers. The image under $\eta_l$ is the conjugating
mapping class up to stabilizers for the local generators of $\efami\ups{l}$, if
these are only among the $l$-companions of the generators for $\efami$.

As in the proof of \ref{stab-indp} the stabilizers do not influence the
conjugation
of $l$-com\-panions, and we may conclude, that the braid monodromy of
$\efami\ups{l}$ is generated by the $l$-companions of generators of
the monodromy
of $\efami$.

\begin{lemma}
\label{fn/gen}
\label{comp/gen}
For $n\geq2$ the versal braid monodromy of a family of functions
$f_a(x_1,...,x_n)$,
is generated by the $l_n$-companions of generators of the versal braid
monodromy of the family of functions $f_\a'$ given by
$$
x_1^{l_1+1}-\a(l_1+1)x_1+\sum_{i=2}^{n-1}\left(x_i^{l_i+1}-\e_i(l_i+1)x_i\right)
$$
\end{lemma}

\proof
In the family $f_\a'$ there are only ordinary multiple points and
multiple points
which are the images of critical points of type $A_{l_n}$. In fact the local
pieces are isomorphic to the cases given in lemma
\ref{multinode/comp} and lemma \ref{An/comp}. So versal braid
monodromy over the
$U_i$ is given by the $l_n$-companions of the generators of the
locally assigned
groups for $f'_\a$. Hence the claim follows by the discussion above.
\qed

\begin{lemma}
\label{versal-band}
The versal braid monodromy of the $l$-companion of the discriminant family
$\efami$ contains the image under the $l$-band
homomorphism $\eta_l$ of the braid monodromy of $\efami$.
\end{lemma}

\proof
For any given element of the braid monodromy of $\efami$ there is a path in the
base to which it is assigned. This path can be taken outside the
$U_i$ so we can
conclude that in the associated $l$-companion family the image under
$\eta_l$ is
assigned to this path.
\qed



\chapter{bifurcation braid monodromy of elliptic fibrations}

Given a family over a base $T$ of smooth regular elliptic surfaces with an
elliptic fibration induced by a global map to $\P$. Suppose all surfaces have
smooth fibres only except for $k$ fibres of type $I_1$, $l$ of type
$I^*_0$, the
divisor of critical values defined in $T\times\P$ is a finite cover of $T$ of
degree $k+l$.

The associated monodromy homomorphism takes values in the braid group of the
sphere. We show that its image is contained up to conjugacy in a
subgroup associated to a family $\xfami_{k,l}$ of elliptic fibrations.

On the other hand a fibration preserving
topological automorphism of an elliptic fibration induces an mapping class of
the base $\P$ punctured at the base points of the singular fibres.
We give a topological characterization of a subgroup of induced mapping classes
which we show to contain the image of the braid monodromy homomorphism and to
coincide with the image in case of the families $\xfami_{6,l}$.

\section{introduction}

The monodromy problems we want to discuss fit quite nicely into the following
general scheme:
Given an algebraic object $X$ consider an algebraic family $g:\xfami\to T$ such
that a fibre $g\inv(t_0)$ is isomorphic to $X$ and such that the
restriction to
a connected subfamily $g|:\xfami'\to T'$ containing $X$ is a locally trivial
$C^\infty$ fibre bundle. If $G$ is the structure group of this bundle, the
geometric monodromy is the natural homomorphism $\rho:\pi_1(T',t_0)\to G$. A
monodromy map with values in a group $A$ is obtained by composition with some
representation $G\to A$.\\

In the standard setting $X$ is a complex manifold, e.g.\ a smooth complex
projective curve. In this case $\xfami$ is a flat family of compact curves
containing
$X$, the subfamily $\xfami'$ contains only the smooth curves and is a locally
trivial bundle of $C^\infty$ surfaces with structure group the
mapping class group $M\!\:\!ap(X)$. From the geometric monodromy one can
obtain the
algebraic monodromy by means of the natural representation
$M\!\:\!ap(X)\to Aut(H_1(X))$.\\

In the present paper we investigate the monodromy of regular elliptic
fibrations.
So $X$ is an elliptic surface with a map $f:X\to\P$ onto the projective line.
We consider families $g:\xfami\to T$ of elliptic surface containing $X$
with a map
$f_T:\xfami\to\P$ which extends $f$ and induces an elliptic fibration on each
surface. Subfamilies $\xfami'$ are to be chosen as $C^\infty$ fibre bundles
with
structure group $\diff_f(X)$, the group of isotopy classes of diffeomorphism
which commute with the fibration map up to a diffeomorphism of the base.

This
structure group has a natural representation in the mapping class group of the
base $\P$ punctured at the singular values of the fibration map $f$. The
corresponding monodromy homomorphism is called the braid monodromy of the
family
$\xfami'$, since the mapping class group of the punctured base is
isomorphic to
a braid group of the sphere.
\begin{eqnarray*}
\sbr_n & = & \left\langle \s_1,...,\s_{n-1}\,\left|\,
\begin{array}{l}
\s_i\s_{i+1}\s_i=\s_{i+1}\s_i\s_{i+1},\\
\s_i\s_j=\s_j\s_i \text{ if }|i-j|\geq2,\\
(\s_1\cdot...\cdot\s_{n-1})^n=1
\end{array}
\right.\right\rangle.
\end{eqnarray*}

 From all possible $C^\infty$-types of elliptic fibrations we choose the subset
only of those represented by elliptic fibrations with singular fibres of types
$I_1,I^*_0$ only, which is significant for any
elliptic fibrations without multiple fibres deforms to such a fibration.
We call a subgroup $E$ of a spherical braid group the {\sl braid monodromy
group}
of a fibration $X$, if $E$ is the smallest subgroup (w.r.t.\ inclusion)
such that
for all admissible $\xfami$ the image of the braid monodromy is a subgroup
of $E$
up to conjugation and prove:

\begin{main}
\labell{mainE}
The braid monodromy group of a regular elliptic fibration $X$
with no singular fibres except $6k$ fibres of type $I_1$ and $l$ fibres of type
$I^*_0$ is a subgroup of $\sbr_{6k+l}$ representing the conjugacy class of
$$
E^s_{6k,l}:=\left\langle \s_{ij}^{m_{ij}},\,i<j\,\left|\quad m_{ij}=\left\{
\begin{array}{ll}
1 & \text{ if }i,j\leq l \,\,\vee\,\, i\equiv j\,(2),i,j>l\\
2 & \text{ if } i\leq l<j\\
3 & \text{ if } i,j>l, i\not\equiv j \,(2)
\end{array}
\right.\right.\right\rangle.
$$
(Here $\s_{i,i+1}:=\s_i$, while for $j>i+1$ we define
$\s_{i,j}:=\s_{j-1}\cdots\s_{i+1}\s_i\s_{i+1}\inv\cdots\s_{j-1}\inv$.)
\end{main}

Along with the proof of the theorem we will notice that each mapping class in
the braid monodromy group $E^s_n(X)$ is represented by a diffeomorphism which
can be lifted to a diffeomorphism of $X$ inducing the trivial mapping class on
some generic fibre. Hence we ask for the converse:
\begin{quote}
Does every diffeomorphism of $X$, isotopic to the identity mapping on some
generic fibre, induce a mapping class of the punctured base which is in the
monodromy group of $X$?
\end{quote}

A positive answer would yield a topological characterisation of the braid
monodromy group!\\

In fact we show that the group of mapping classes induced in the said way
coincides with the stabiliser group of an appropriate Hurwitz action. Then we
use \cite{hurwitz} to give an affirmative answer to the question above in case
the number of fibres of type $I_1$ does not exceed $6$.

\begin{thm}
\label{1}
The braid monodromy group $E^s_n$ of a regular elliptic fibration $X$ with no
singular fibres except $6$ fibres of type $I_1$ and $l$ fibres of type $I^*_0$
coincides with the group of mapping classes of the punctured base of $X$ which
are induced by diffeomorphisms of $X$ which respect the fibration, preserve a
fibre $F$ and induce a map on $F$ isotopic to the identity.
\end{thm}

\section{bifurcation braid monodromy}

With each locally trivial bundle one can associate the structure homomorphism
defined on the fundamental group of the base with respect to any base point. It
takes values in the mapping classes of the fibre over the base point.

Given a curve $C$ in the affine plane we can take a projection to the
affine line
which restricts to a finite covering $C\to\CC$. The complement of the curve
and its
vertical tangents is the total space of a punctured disc bundle over the
complement of the branch points in the affine line.

The structure homomorphism of this bundle is called the braid monodromy of the
plane curve with respect to the projection, and it can be naturally
regarded as a
homomorphism from the fundamental group of the branch point complement to the
braid group, since the latter is naturally isomorphic to the mapping class
group
of the punctured disc.

This definition is readily generalized to the case of a divisor in the
Cartesian product of the affine or projective line with an irreducible base
$T$.
Then the structure homomorphism takes values in a braid group, resp.\ in a
mapping
class group of a punctured sphere which is naturally isomorphic to the
spherical
braid group $\sbr_n$.

In situations as we are interested in, such a divisor is defined as the
locus of
critical values of a family of algebraic functions of
constant bifurcation degree
with values in $\LL\cong\P$ or $\CC$.
Thus we give the relevant definitions:

\begin{description}

\item[Definition:]
A flat family $\xfami\to T$ with an algebraic morphism $f:\xfami\to \LL$ is
called a
{\it framed family of functions} $(\xfami,T,f)$.

\item[Definition:]
The {\it bifurcation set} of a framed family of functions
over $T$ is the smallest Zariski closed subset $\bfami$ in $T\times \LL$ such
that the diagonal map $\xfami\to T\times \LL$ is smooth over the complement of
$\bfami$.

\item[Definition:]
The {\it discriminant set} of a framed family of functions over $T$ is the
divisor in $T$ such that the bifurcation set $\bfami$ is an unbranched cover
over its complement by the restriction of the natural projection
$T\times \LL\to
T$.

\item[Definition:]
A framed family of functions is called {\it of constant
bifurcation degree} if the bifurcation set is a finite cover of $T$.

\item[Definition:]
The {\it bifurcation braid monodromy} of a framed family of functions with
constant bifurcation degree over an irreducible base $T$ is defined to be the
braid monodromy of $\bfami$ in $T\times\LL$ over $T$.

\end{description}

Note that this definition of braid monodromy differs slightly from the
definition
given in the introduction but that the resulting objects are the same.

\section{families of divisors in Hirzebruch surfaces}

Given a Hirzebruch surface $\fk$ with a unique section $\isec$ of
selfintersection $-k$, we consider families of divisors on which the ruling of
the Hirzebruch surface defines families of functions with constant bifurcation
type.\\
We can pull back divisors form the base along the ruling to get divisors on
$\fk$ which we call vertical, among others the fibre divisor $L$.\\
Consider now the family of divisors on $\fk$ which consist of a vertical part
of degree $l$ and a divisor in the complete linear system of
$\ofami_{\fk}(4\isec+3kL)$ called the horizontal part. It is a family
parameterized by $T=\ph(\ofami_{\P}(l))\times\ph(\ofami_{\fk}(4\isec+3kL))$
with
total space
$$
\dfami_{k,l}=\left\{(x,t)\in\fk\times T\,|\,x\in D_t\subset\fk\right\}.
$$
Let $T'$ be the Zariski open subset of $T$ which is the base of the family
$\dfami'_{k,l}$ of divisors in $\dfami_{k,l}$ with reduced horizontal part.

\begin{lemma}
\labell{dfami}
The ruling on $\fk$ defines a morphisms $\dfami_{k,l}'\to\P$ by which it
becomes a framed family of functions of constant bifurcation degree.
\end{lemma}

\proof
The critical value set of the vertical part of a divisor is the divisor of
which it is the pull back, thus it is constant of degree $l$.\\
The assumption on reducedness forces the horizontal part to be without fibre
components. We may even conclude that a reduced horizontal part consists of
$\isec$ and a disjoint divisor which is a branched cover of the base of degree
$3$. The critical values set is therefore the branch set which is of constant
degree $6k$, and we are done.
\qed

Note that we can consider the abstract braid group presented as
\begin{eqnarray*}
\br_n & = & \langle \s_1,...,\s_{n-1}\,|\,
\s_i\s_{i+1}\s_i=\s_{i+1}\s_i\s_{i+1},
\s_i\s_j=\s_j\s_i \text{ if } |i-j|\geq2\rangle.
\end{eqnarray*}
resp.\ $\sbr_n$ as presented in the introduction to be realised by the mapping
class group of the punctured disc, resp.\ sphere. Such an identification is
given if each
$\s_i$ is realised by the half-twist on an embedded arc $a_i$ connecting two
punctures provided that
$a_i\cap a_{i+1}$ is a single puncture and $a_i\cap a_j$ is empty if
$|i-j|\geq2$.

\begin{prop}
\labell{dmono}
The image of the bifurcation braid monodromy homomorphism of the family
$\dfami'_{k,l}$ is conjugated to the subgroup of $\sbr_{6k+l}$:
$$
E^s_{6k,l}:=\left\langle \s_{ij}^{m_{ij}},\,i<j\,\left|\quad m_{ij}=\left\{
\begin{array}{ll}
1 & \text{ if }i,j\leq l \,\,\vee\,\, i\equiv j\,(2),i,j>l\\
2 & \text{ if } i\leq l<j\\
3 & \text{ if } i,j>l, i\not\equiv j \,(2)
\end{array}
\right.\right.\right\rangle.
$$
\end{prop}

The proof of this proposition and a couple of preparatory results will take the
rest of the section.

First note that our whole concern lies in the understanding of the bifurcation
set $\bfami$ in $T'\times\P$ with its projection to $T'$. As an approximation
we will consider families of affine plane curves given by families of
polynomials in affine coordinates $x,y$ with the regular map induced by the
affine projection $(x,y)\mapsto x$.

Their bifurcation sets are contained in the Cartesian product of the family
bases with the affine line $\CC$, and it will soon be shown that this pair can
be induced from $(T'\times\P,\bfami)$.
Eventually we can extract all necessary information from such families to
prove our claim.

\begin{lemma}
\labell{div/discr}
Consider $y^3-3\,p(x)\,y\,+\,2\,q(x)$ as a family of polynomial functions
\mbox{$\CC^2\times T\to\CC$} parametrised by a base $T$ of pairs $p,q$ of
univariate polynomials. Then the bifurcation set is the zero set of
$g(x):=p^3(x)-q^2(x)$, the discriminant set is the zero set of the discriminant
of $g$ with respect to
$x$.
\end{lemma}

\proof
The bifurcation divisor is cut out by the discriminant polynomial of
$y^3-3p(x)y+2q(x)$ with respect to $y$. The first claim is then immediate since
$g$ is proportional to the corresponding Sylvester determinant:
$$
\begin{vmatrix}
1 & 0 & -3p & 2q \\
  & 1 & 0 & -3p & 2q \\
3 & 0 & -3p\\
  & 3 & 0 & -3p\\
  & & 3 & 0 & -3p
\end{vmatrix}
$$
For the second claim we only note that a pair $p,q$ belongs to the
discriminant set if and only if $p^3-q^2$ has a
multiple root hence this locus is cut out by the discriminant of
$g$ with respect to $x$.
\qed

\begin{lemma}
\labell{dis-comp}
The discriminant locus of a family $y^3+3r(x)y^2-3p(x)y+2q(x)$ is the union of
the degeneration component of triples $p,q,r$ defining singular curves and the
cuspidal component of triples defining polynomial maps with a
degenerate critical
point.
\end{lemma}

\proof
In general a branched cover of $\CC$ has not the maximal number of
branch points
only if the cover is singular, or the number of preimages of a branch point
differs by more than one from the degree of the branching. The second
alternative
occurs only if there is a degenerate critical point in the preimage or if
there are
two critical points. Since the last case can not occur in a cover of degree
only
three we are done.
\qed

\begin{lemma}
\labell{cusp-mult}
Given the family $y^3+3r(x)y^2-3p(x)y+2q(x)$ the cuspidal component of the
discriminant is the zero set of the resultant of $p(x)+r^2(x)$ and
$2q(x)-r^3(x)$ with respect to $x$.\\
Its equation - considered a polynomial in the variable $\l_0$ - is of
degree $n$
with coprime coefficients if\\[-2mm]
$$
p(x)=\sum_{i=0}^d\l_ix^i,\quad q(x)=x^n+\sum_{i=0}^{n-1}\x_ix^i,\quad
r(x)=\sum_{i=0}^{\lfloor n/3\rfloor}\zeta_ix^i.
$$
\end{lemma}

\proof
The cuspidal discriminant is the locus of all parameters for which there is a
common zero of $f,\del_yf,\del^2_yf$. Since $\del^2_yf=0$ is linear in $y$,
we can eliminate $y$ and get the resultant of $p(x)+r^2(x)$ and $2q(x)-r^3(x)$
with respect to $x$.\\
By the degree bound on $q$ and $r$ the discriminant equation is the resultant
of a matrix in which the variable $\l_0$ occurs exactly $n$ times. Moreover
the diagonal determines the coefficient of $\l_0^n$ to be
$(1-\zeta_{n/3}^3)^{max(d,2n/43)}$ resp.\
$1$ depending on whether $n/3\in\ZZ$ or not. Even in the first case the
coefficients are coprime since the resultant is not divisible by
$(1-\zeta_{n/3}^3)$.
\qed

\begin{lemma}
\labell{dgn-mult}
For the family
$y^3+3r(x)y^2-3p(x)y+2q(x)$
the degeneration component of the discriminant is
the locus of triples for which there is a common zero in
$x,y$ of the polynomial and its two partial derivatives.
Its equation - considered a polynomial in the variable $\xi_0$ - is monic of
degree $2n-2$ if\\[-4mm]
$$
p(x)=\sum_{i=0}^d\l_ix^i,\quad q(x)=x^n+\sum_{i=0}^{n-1}\x_ix^i,\quad
r(x)=\sum_{i=0}^{\lfloor n/3\rfloor}\zeta_ix^i.
$$
\end{lemma}
\proof
The degeneration locus is given by the Jacobian criterion as claimed. The
equation of the discriminant of the subdiagonal unfolding of the
quasihomogeneous singularity $y^3-x^n$ is known to the quasihomogeneous and of
degree $2n-2$ in $\xi_0$. Since the unfolding over the $\xi_0$-parameter is a
Morsification the coefficient of $\xi^{2n-2}_0$ must be constant.
\qed

\begin{lemma}
\labell{special}
The bifurcation braid monodromy of the family $y^3-3\pnull
y+2(x^n+\qnull)$ maps
onto a subgroup of $\br_{2n}$ which is conjugated to the subgroup generated by
$$
(\s_1\cdots \s_{2n-1})^3,(\s_{2n-2,2n}\cdots \s_{2,4})^{n+1},
(\s_{2n-3,2n-1}\cdots \s_{1,3})^{n+1}.
$$
\end{lemma}

\proof
As one can show with the help of the preceding lemmas, the discriminant
locus is the union of the degeneration locus and the cuspidal component
which are
cut out respectively by the polynomials $\pnull^3-\qnull^2$ and $\pnull$.

By Zariski/van Kampen the fundamental group of the complement with base point
$(\pnull,\qnull)=(1,0)$ is generated by the fundamental group of the complement
restricted to the line $\pnull=1$ and the homotopy class of a loop which links
the line $\pnull=0$ once.

For the pair $(1,0)$ the set of regular values of the polynomial consists
of the
affine line punctured at the $(2n)$-th roots of unity, which we number
counterclockwise, $1$ the first puncture.
To express the bifurcation braid group in terms of abstract generators, we
identify the elements $\s_{i}$ with the half twist on the circle segment
between the $i$-th and $i+1$-st puncture.

For the line $\pnull=1$ the bifurcation locus is given by
$(x^n+\qnull-1)(x^n+\qnull+1)$. This locus is smooth but branches of degree $n$
over the base at $\qnull=\pm1$. The corresponding monodromy transformations are
the second and third transformation given in the claim.

Associated to the degeneration path $(\pnull,\qnull)=(1-t,\imag t)$,
$t\in[0,1]$
there is a loop in the complex line $\pnull=1+\imag\qnull$ which links the line
$\pnull=0$. For this degeneration the
bifurcation divisor is regular and contains points of common absolute value
determined by $t$ only, except for $t=1$ where it has $n$ ordinary cusps with
horizontal tangent cone. Since a cusp corresponds to a triple half twist
and the first and second puncture merge in the degeneration, the monodromy
transformation for our loop is the first braid of the claim.
\qed

\begin{lemma}
\labell{semi-gen}
The bifurcation braid monodromy of the family $y^3-3\pnull
y+2(x^n+\qnull+\e x)$,
$\e$~small and
fix, is in the conjugation class of the subgroup of the braid group
generated by
$$
(\s_{1}\s_3\cdots \s_{2n-1})^3,\s_{i,i+2},i=1,...,2n-2.
$$
\end{lemma}

\proof
The discriminant locus in the $\pnull,\qnull$ parameter plane
consists again of the
cuspidal component $\pnull=0$ and the degeneration component.

Since the perturbation $\e$ is arbitrarily small, some features of the
family of
lemma~\ref{special} are preserved. The conclusion of the Zariski/van Kampen
argument still holds,
each braid group generators $\s_{i}$ is now realized as half twist on
segments of a slightly distorted circle, and the loop linking
$\pnull=0$ is only
slightly perturbed. So the monodromy transformation associated to this loop is
formally the same as before, the first braid in the claim.

The dramatic change occurs in the bifurcation curve over the line $\pnull=1$.
Now the
bifurcation locus is the union of two disjoint smooth components each of which
branches simply of degree $n$ with all branch points near $\qnull=1$, resp.\
$\qnull=-1$.
Since the local model $x^n+\e x$ has the full braid group as its monodromy
group,
the monodromy along $\pnull=1$ is generated by the elements
$\s_{i,i+2}$ as claimed.
\qed

\begin{lemma}
\labell{generic}
The bifurcation braid monodromy of the
family
$$
y^3-3(\pnull+\peins x)y+2(x^n+\qnull+\qeins x)
$$
is in the conjugation class of the subgroup generated by
$$
\s_{i}^{3},i\equiv1(2),\s_{i,i+2},0<i<2n-1.
$$
\end{lemma}

\proof
The components of the discriminant are the degeneration component and
the cuspidal component.

The line $\pnull=1,\peins=0,\qeins=\e$ small and fix, is generic for the
degeneration component
and we may conclude from lemma \ref{semi-gen} that there are elements in the
fundamental group of the discriminant complement with respect to
$(\pnull,\peins,\qnull,\qeins)=(1,0,0,0)$ which map to $\s_{i,i+2}$ as in lemma
\ref{semi-gen}.

Since the line $\qnull=i,\qeins=\peins=0$ is transversal for the
cuspidal component
so are
parallel lines with $\peins=\e'$ small and fix. The bifurcation curve is then
given by
$(\pnull+\e' x)^3=(x^n+\imag)^2$. For $\pnull=0$ the critical values are
distributed in pairs along a circle in the affine line which merge pairwise for
$\peins\to0$.

But for $\e'$ sufficiently small and by varying $\pnull$ to an appropriate
small extend such
that the degeneration component is not met, there are $n$ obvious degenerations
when $\pnull$ is of the same modulus as $\e'$ and $\pnull+\e' x$ is a
factor of
$x^n+\imag$. By the local nature of the degeneration a degree
argument shows that
these are all possibilities. Moreover one can easily see that the corresponding
monodromy transformations are triple half twist for each one of the pairs.

Moreover these twists are transformed to the transformations $\s_{i,i+1}^3$,
$i$
odd, when transported along $(1-t,\imag t,0,\e't)$, $t\in[0,1]$ to the chosen
reference point.

The monodromy group is thus completely determined since the fundamental
group is
generated by elements which map to the given group under the monodromy
homomorphism.
\qed

\begin{lemma}
\labell{big+gen}
Let a family of plane polynomials be given which is of the form
\begin{equation*}
\label{pol}
y^3+3\Big(\sum_{i=0}^{d_r}\zeta_ix^i\,\Big)y^2
-3\Big(\sum_{i=0}^{d_p}\l_ix^i\,\Big)y+
2\big(x^n+\sum_{i=0}^{n-1}\xi_ix^i\,\big),
\end{equation*}
\hfill$2\leq n,0<3d_p\leq 2n,3d_r\leq n$.\\
Then the bifurcation braid monodromy group
is in the conjugation class of the subgroup of the
braid group generated by
$$
\s_{i}^{3},i\equiv1(2),\s_{i,i+2},0<i<2n-1.
$$
\end{lemma}

\proof
Since the family considered in the previous lemma is a subfamily now and
has the
claimed monodromy, we have to show that the new family has no additional
monodromy
transformations.

In the proof above we have seen that the cuspidal component is cut in $n$
points
by a line in $\l_0$ direction. The component is reduced since its
multiplicity at the origin is $n$, too, by lemma \ref{cusp-mult}. The
degeneration component is reduced by the analogous argument relying on lemma
\ref{dgn-mult} and transversally cut in $2n-2$ points by a line in $\xi_0$
direction.
Hence by Zariski/van Kampen arguments as proved in \cite{Be} the fundamental
group of the discriminant complement of the subfamily surjects onto the
fundamental group of the family considered now.
\qed[4mm]

\begin{lemma}
\labell{red-gen}
Define a subgroup of the braid group $\br_{2n+l}$:
$$
E_{2n,l}:=\left\langle \s_{i,j}^{m_{ij}},\,i<j\,\left|\quad m_{ij}=\left\{
\begin{array}{ll}
1 & \text{ if }i,j\leq l \vee i\equiv j (2),i>l\\
2 & \text{ if } i\leq l<j\\
3 & \text{ if } i,j>l \wedge i\not\equiv j (2)
\end{array}
\right.\right.\right\rangle.
$$
Then the same subgroup is generated also by
$$
\s_{i},i<l,\s_{i,i+2},l<i,\s^2_{i,j},i\leq l<j,\s^3_{i},i>l,i\not\equiv
l (2).
$$
\end{lemma}

\proof
We have to show that the redundant elements can be expressed in the
elements of the
bottom line. This is immediate from the following relations $(i<j)$:
$$
\begin{array}{lll}
\s_{i,j} & = & \s\inv_{j-1}\cdots \s\inv_{i+1} \s_{i}\low
\s_{i+1}\low\cdots \s_{j-1}\low
,\quad j\leq l,\\
\s_{i,j} & = & \s\inv_{j-2,j}\cdots \s\inv_{i+2,i+4} \s_{i,i+2}\low
\s_{i+2,i+4}\low\cdots \s_{j-2,j}\low
,\quad l< i,i\equiv j(2),\\
\s^3_{i,j} & = & \s\inv_{j-2,2}\cdots \s\inv_{i+1,i+3} \s^3_{i}
\s_{i+1,i+3}\low\cdots
\s_{j-2,j}\low ,\quad l<i,i\not\equiv l,j(2),\\
\s^3_{i,j} & = & \s\inv_{j-2,2}\cdots \s\inv_{i+1,i+3} \s_{i-1,i+1}\low
\s^3_{i-1,i}
\s\inv_{i-1,i+1} \s_{i+1,i+3}\low\cdots \s_{j-2,j}\low
  ,\, l<i,j\not\equiv l,i(2).
\end{array}
$$
\qed

\begin{lemma}
\labell{mono-full}
Consider a family $\left(y^3-3p(x)y+2q(x)\right)a(x)$
parametrised by triples $p,q,a$, with $p$ from the vector space of univariate
polynomials of degree at most $2n/3$, $q,a$ from the affine space of monic
polynomials of degree $n$ and $l$ respectively.
Then the subgroup $E_{2n,l}$ of $\br_{2n+l}$ is conjugate to a subgroup of the
image of the bifurcation braid monodromy.
\end{lemma}

\proof
We choose our reference divisor to be $(y^3-3y+2x^n)\prod_i^l(x-l-2+i)$ with
corresponding bifurcation set $x_i=l+2-i,i\leq l$ on the real axis and
$x_{l+1}=1$ and
the $x_i,i>l+1$ equal to the $2k^{\text{th}}$-roots of unity in
counterclockwise
numbering.\\
We identify the elements $\s_{i,j}$ of the braid group with the half twist
on arcs
between $x_i,x_j$, which are chosen to be
\begin{enumerate}
\item
a circle segment through the lower half plane, if $i,j\leq l$,
\item
a circle secant in the unit disc, if $i,j>l$,
\item
the union of a secant in the unit disc to a point on its boundary between
$x_{2n+l}$
and $1$ with an arc through the lower half plane, if $i\leq l<j$.
\end{enumerate}
(Of each kind we have depicted one in the following figure.)

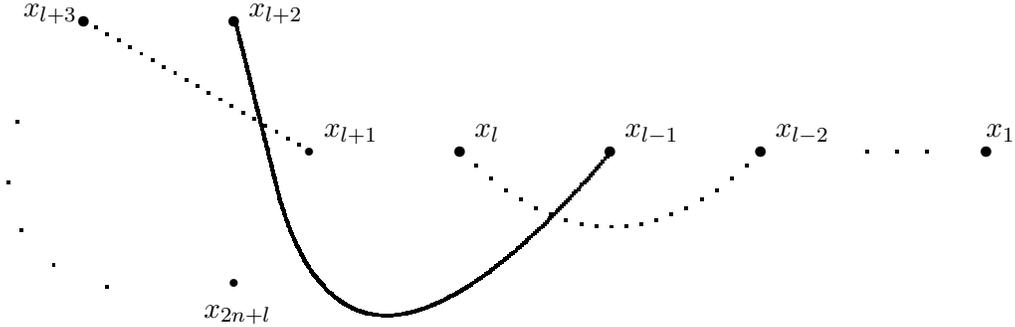
\setlength{\unitlength}{2mm}
\begin{picture}(68,22)
\linethickness{1pt}
\put(15,1.3){\circle*{.6}}
\put(20,10){\circle*{.5}}
\put(30,10){\circle*{.7}}
\put(40,10){\circle*{.7}}
\put(50,10){\circle*{.7}}
\put(65,10){\circle*{.7}}
\put(5,18.7){\circle*{.7}}
\put(15,18.7){\circle*{.7}}

\bezier{4}(.5,12)(-2,3)(6.5,1)
\bezier{2}(57,10)(59,10)(61,10)

\bezier{20}(30,10)(40,0)(50,10)              
\bezier{20}(20,10)(12.5,14.35)(5,18.7)       

\bezier{120}(15,18.7)(16.5,12.7)(18,6.7)     
\bezier{240}(18,6.7)(23,-10)(40,10)

\put(51,11){$x_{l-2}$}
\put(65,11){$x_1$}
\put(16,19){$x_{l+2}$}
\put(1,19){$x_{l+3}$}
\put(21,11){$x_{l+1}$}
\put(31,11){$x_l$}
\put(41,11){$x_{l-1}$}
\put(13,-1){$x_{2n+l}$}

\end{picture}
\vspace*{6mm}\\

Since keeping the horizontal part $y^3-3y+2x^n$ fix, the
bifurcation divisor of the vertical is equivalent to that of the universal
unfolding of the function $x^l$ we have the elements $\s_{i},i<l$ in the
braid monodromy.
These elements are obtained for example in families
$$
a(x)=\left((x-l+i-3/2)^2+\l\right)\prod_{j\neq i,i+1}^l(x-l-2+j).
$$

The second set of elements, $\s_{i,j}^2,i\leq l<j$ is obtained by families
of the
kind
$$
(y^3-3y+2x^n)(x-l-2+i-\l)\prod_{j\neq i}^l(x-l-2+j)
$$
since the zero $l+2-i+\l$ may trace any given path in the range of the
projection, in particular that around an arc on which the full twist
$\s_{i,j}^2$
is performed.\\
Finally varying the horizontal part as in lemma \ref{generic}
while keeping the
$a(x)$ factor fix proves that the braid group elements $\s_{i,i+2},l<i$ and
$\s_{i},l<i,i\not\equiv l(2)$ are in the image of the monodromy. So we may
conclude that this image contains $E_{2n,l}$ up to conjugacy.
\qed

\proof[ of prop.\ \ref{dmono}]
Denote by $S$ the Zariski open subset of $T'$ which parameterizes divisors
of the
family $\dfami_{k,l}'$ which have no singular value at a point $\infty\in\P$.
The corresponding family in $\fk\times S$ may then be restricted to a
family $\ffami_{k,l}$ in $\CC\times\CC\times S$, where $\fk$ is
trivialized as $\CC\times\CC$ in the complement of the negative section
$\isec$ and the fibre over $\infty$.
By construction $\ffami_{k,l}$ has constant bifurcation degree.\\
Consider now the family of polynomials
$$
\left( y^3+3r(x)y^2-3p(x)y+2q(x)\right)a(x),
$$
where $r,p,q,a$ are taken from the family of all quadruples of polynomials in
one variable subject to the conditions that
\begin{enumerate}
\item
$r,p$ are of respective degrees $k$ and $2k$,
\item
$q$ is monic of degree $3k$, $a$ is monic of degree $l$
\item
the discriminant of $y^3+3r(x)y^2-3p(x)y+2q(x)$ is not identically zero.
\end{enumerate}
This family can be naturally identified with $\ffami_{k,l}$.
By lemma \ref{mono-full}, up to conjugacy, $E_{6k,l}$ is contained in the
monodromy image $\rho(\pi_1(S\setminus Discr(\ffami_{k,l})))$.

For the converse we note that the bifurcation set
of the family decomposes into the bifurcation sets $\bif_h$ of the horizontal
part $y^3+3r(x)y^2-3p(x)y+2q(x)$ and $\bif_v$ of the vertical part
$a(x)$. Hence
the monodromy is contained in the subgroup $\br_{(6k,l)}$ of braids
which do not
permute points belonging to different components. $\br_{(6k,l)}$ has natural
maps to $\br_{6k}$ and $\br_l$ which commute with the braid monodromies of both
bifurcation set considered on their own.

The discriminant decomposes into the discriminants of $\bif_h$,
$\bif_v$ and the divisor of parameters for which $\bif_h\cap\bif_v$ not empty.
They give rise in turn to braids which can be considered as elements in
$$
\br_{6k},\,\br_l\text{ resp.\ }\br_{(6k,l)}^{0,0}:=\{\b\in\br_{(6k,l)}|\,
\b\text{ trivial in }\br_{6k}\times\br_l\}.
$$
Now with lemma \ref{red-gen} we can identify $E_{6k.l}$ as
the subgroup of $\br_{6k+l}$ generated by
$E_{6k}\subset\br_{6k}$, $\br_l$ and $\br_{(6k,l)}^{0,0}$
which are generated in turn by the elements
$$
\{\s_{i,i+2},\s^3_{i},l<i\},\{\s_{i},i<l\},\{\s^2_{i,j},i\leq l<j\}\text{
resp.}
$$
And by lemma \ref{big+gen} the image can not contain more elements.
\pagebreak

Since the bifurcation diagram of $\ffami_{k,l}$ embeds in the
bifurcation diagram of $\dfami_{k,l}'$ with complement of codimension
one, there
is a commutative diagram
$$
\begin{array}{ccc}
\pi_1(S\setminus Discr(\ffami_{k,l})) & \tto\!\!\!\!\!\!\!\!\tto
& \pi_1(T'\setminus Discr(\dfami_{k,l}'))\\
\downarrow & & \\[-4.7mm]
\downarrow & & \Big|\\[-.1mm]
E_{6k,l} & & \big|\\[-.5mm]
\rceil & & \Big|\\[-3mm]
\hspace*{.2mm}\downarrow & & \hspace*{.1mm}\big\downarrow\\
\br_{6k+l} & \tto\!\!\!\!\!\!\!\!\tto & \sbr_{6k+l}
\end{array}
$$
from which we read off our claim.
\qed

\begin{cor}
For any element $\b$ in the braid monodromy group of $\dfami_{k,l}'$ there is a
diffeomorphism of the base $\P$ which fixes a neighbourhood of
$\infty\in\P$ and
which represents the mapping class $\b$.
\end{cor}

\proof
The element $\b$ is image of an element $\b'$ in the braid monodromy of the
bifurcation diagram of $\ffami_{k,l}$. The bifurcation set does not meet the
boundary so integration along a suitable vector field yields a realisation of
$\b'$ as a diffeomorphism acting trivially on a neighbourhood of the boundary.
Its trivial extension to the point $\infty$ is the diffeomorphism sought for.
\qed

\section{families of elliptic surfaces}

In this section we start investigating families of regular elliptic surfaces
for which the type of singular fibres is restricted to $I_1$ and $I_0^*$.
We will go back and forth between a family of elliptic fibrations, its
associated
family of fibrations with a section and a corresponding Weierstrass model
of the
latter, so we note some of their properties:

\begin{prop}
\labell{jac}
Given a family of elliptic fibrations with constant bifurcation type over an
irreducible base $T$, there is a family of elliptic fibrations with a section,
such that the bifurcation sets of both families coincide.
\end{prop}

\proof
Given a family as claimed there is the associated family of Jacobian
fibrations,
cf.\ \cite[I.5.30]{fm}. The bifurcation sets of both families coincide.
\qed

In turn, for each family of elliptic fibrations with a section there is a
corresponding family of Weierstrass fibrations, cf.\ Seiler \cite{sei}.

A regular Weierstrass fibration $W$ is defined by an equation
$$
wz^2=4y^3-3Pw^2y+2Qw^3
$$
in the projectivisation of the vector bundle
$\ofami\oplus\ofami(2\chi)\oplus\ofami(3\chi)$  over the projective line $\P$
where $\chi$ is the holomorphic Euler number of the fibration,
$w,y,z$ are 'homogeneous coordinates' of the bundle, and $P,Q$ are sections of
$\ofami(4\chi),\ofami(6\chi)$ respectively.\\
So $W$ is a double cover of the Hirzebruch surface $\hirz_{2\chi}={\mathbf
P}(\ofami\oplus\ofami(2\chi))$ branched along the section $\s_{2\chi}$ and the
divisor in its complement $\ofami(2\chi)$ defined by the equation
$y^3-3Py+2Q=0$.\\

A framed family of Weierstrass fibrations over a parameter space $T$ is a
given by
data as before where now $P,Q$ are sections of the pull backs to $T\times\P$
of $\ofami(4\chi),\ofami(6\chi)$ such that for each parameter $\l\in T$ they
define a Weierstrass fibration. In the sequel $P,Q$ are referred to as the
coefficient data of the Weierstrass family.

\begin{lemma}
\labell{w-fact}
Let $\wfami$ be the Weierstrass family associated to a framed family over
$T$ of
regular elliptic fibrations in which all surfaces have no singular fibres
except
for $l$ of type $I_0^*$ and $6k$ of $I_1$ with coefficient data $P,Q$, then
there
are three families of sections $a,p,q$ of $\ofami(l),\ofami(2k),\ofami(3k)$
respectively, such that $p,q$ have no common zero,
\begin{eqnarray*}
& p\cdot a^2=P, & q\cdot a^3=Q,
\end{eqnarray*}
and the bifurcation set is given by
\begin{eqnarray*}
a\left(p^3-q^2\right)=0 & \subset & T\times\P.
\end{eqnarray*}
\end{lemma}

\proof
By the classification of Kas \cite{kas} at base points of regular fibres the
discriminant $P^3-Q^2$ does not vanish, at base points of fibres of type $I_1$
the discriminant vanishes but neither $P$ nor $Q$ and at base points of fibres
of type $I_0^*$  the vanishing order of $P$ is two, the vanishing order of
$Q$ is
three.\\
Since by hypothesis the locus of base points of singular fibres of type $I_0^*$
form a family of point divisors of degree $l$ there is a section $a$ of
$\ofami(l)$ such that $P$ has a factor $a^2$ and $Q$ a factor $a^3$.\\
With $deg P=2(l+k),
deg Q=3(l+k)$ we get the other degree claims.\\
Finally the discriminant of the Weierstrass fibration is given by $P^3-Q^2$
which
has -- by the above -- the same zero set as $a\left(p^3-q^2\right)$.
\qed

\begin{description}
\item[Remark:] In the situation of the lemma, a family of divisors is given for
$\hirz_{k}$ by the equation $a(y^3-3pw^2y+2qw^3)=0$, $a$ cutting out the
vertical part. The double cover along this divisor is a family of fibrations
obtained from the original family by contracting all smooth rational curves of
selfintersection $-2$, of which there are four for each fibre of type $I_0^*$.
\end{description}

We are now prepared to come back to the main theorem:\\

\proof[ of the main theorem]
Given any framed family of regular elliptic fibrations containing $X$ we
consider
a Weierstrass model $\wfami$ of the associated Jacobian family. Since $\wfami$
is again framed there is an induced family of divisors on a Hirzebruch
surface obtained as before.\\
This family of divisors is a pull back from the space $\dfami_{k,l}$ so the
monodromy is a subgroup of the bifurcation monodromy of Hirzebruch divisors.

On the other hand for the family of triples of polynomials
$p(x),q(x),a(x)$ with $p$ of degree at most $2k$ and $q,a$ monic of degree $3k$
respectively $l$, we can form the family given by
$$
z^2=y^3-3p(x)a^2(x)+2q(x)a^3(x),
$$
which is Weierstrass in the complement of parameters where
$a(x)\left(p^3(x)-q^2(x)\right)$ has a multiple root or vanishes identically.
At least after suitable base change, cf.\
\cite[p.\ 163]{fm}, this Weierstrass family has a simultaneous resolution
yielding
a family $\xfami_{k,l}$ of elliptic surfaces with a section.\\
The Jacobian of $X$ is contained in $\xfami_{k,l}$, since its Weierstrass data
consist of sections $P,Q$ which are factorisable as $a^2p,a^3q$ according to
lemma \ref{w-fact} and after the choice of a suitable $\infty$ this data
can be identified with polynomials in this family.\\
The fibration $X$ is deformation equivalent to its Jacobian with
constant local analytic type, cf.\ \cite[thm.\ I.5.13]{fm} and hence of
constant
fibre type.
The monodromy group therefore contains the bifurcation monodromy group of
divisors on Hirzebruch surfaces $\dfami_{k,l}$ and so the two groups even
coincide.
\qed

Regarding elements in the braid monodromy as mapping classes again they can be
shown to be induced by diffeomorphism of the elliptic fibration, but more
is true
in fact:

\begin{prop}
\labell{lift}
For each braid $\b$ in the framed braid monodromy group there is a
diffeomorphism
of the elliptic fibration which preserves the fibration, induces $\b$ on
the base
and the trivial mapping class on some fibre.
\end{prop}

\proof
As we have seen in the corollary to prop.\ \ref{dmono} we can find a
representative $\bdiff$ for the braid $\b$ by careful integration of a suitable
vector field such that $\bdiff$ is the identity next to a point $\infty$.\\
In \cite[II.1.2]{fm} there is a proof for families of nodal elliptic fibrations
and sufficient hints for more general families of constant singular fibre
types,
that a horizontal vector field on the total family can be found which fails
to be
a lift only in arbitrarily small neighbourhoods of singular points on singular
fibres. Integration of such a vector field yields a diffeomorphism
$\Phisl$ which
is a lift of $\bdiff$.\\
We have seen that the monodromy generators arising from the horizontal part can
be realized over a suitable polydisc parameter space, cf.\ lemma \ref{generic}.
Since the vertical part as in lemma \ref{mono-full} does not have any effect on
the fibre $F_\infty$ over $\infty$ we can conclude that this fibration
family is
the trivial family next to $F_\infty$. So we apply the argument above to get a
lift $\Phisl$ which induces the trivial mapping class on $F_\infty$.
\qed[2mm]

\section{Hurwitz stabilizer groups}

In this section we determine the stabilizers of the action of the braid group
$\br_n$ on homomorphisms defined on the free group $F_n$ generated by elements
$t_1,...,t_n$.
The action is given by precomposition with the Hurwitz automorphism of $F_n$
associated to a braid in $\br_n$:
$$
\br_n\to\ON{Aut} F_n:\quad
\s_{i,i+1}\mapsto\left(t_j\mapsto
\begin{cases}
t_j & j\neq i,i+1\\
t_it_jt_i\inv & j=i\\
t_i & j=i+1
\end{cases}
\right).
$$
We start with a result from \cite{hurwitz}:

\begin{prop}
\labell{l}
Let $F_n:=\langle t_i,1\leq i\leq n\,|\quad\rangle$ be the free group on $n$
generators, define a homomorphism $\phi_n:F_n\to\br_3=\langle
a,b\,|\,aba=bab\rangle$ by
$$
\phi_n(t_i)=\left\{
\begin{array}{ll}
a & i\text{ odd}\\
b & i\text{ even}
\end{array}\right.
$$
and let $\br_n$ act on homomorphisms $F_n\to\br_3$ by Hurwitz
automorphisms of $F_n$. Then the stabilizer group $Stab_{\phi_n}$
contains the braid subgroup
$$
E_n=\langle \s_{i,j}^{m_{ij}}\:|\: m_{ij}=1,3
\text{ if $j\equiv i$, resp.\ $i\not\equiv j \mod 2$}\rangle
$$
with $E_n=Stab_{\phi_n}$, if $n\leq6$.
\end{prop}

Note that the action in \cite{hurwitz} was defined on tuples
$\left(\phi_n(t_1),...,\phi_n(t_n)\right)$ but that it is obviously
equivalent to
the action considered here.\\
This result can now be applied to find stabilizers of similar homomorphisms:

\begin{prop}
\labell{slcor}
Let $F_n:=\langle t_i,1\leq i\leq n\,|\quad\rangle$ be the free group on $n$
generators, define a homomorphism $\psi_n:F_n\to\slz$ by
$$
\psi_n(t_i)=\left\{
\begin{array}{ll}
\left(\begin{smallmatrix}1&1\\0&1\end{smallmatrix}\right) & i\text{ odd}\\[3mm]
\left(\begin{smallmatrix}\phantom{-}1&0\\-1&1\end{smallmatrix}\right)
& i\text{ even}
\end{array}\right.
$$
and let $\br_n$ act on homomorphisms $F_n\to\slz$ by Hurwitz automorphisms of
$F_n$.
Then the stabilizer group $Stab_{\psi_n}$ of
$\psi_n$ is equal to the stabilizer group $Stab_{\phi_n}$ of $\phi_n$.
\end{prop}

\proof
Both groups, $\slz$ and $\br_3$, are central extensions of $\psl$, and both
$\phi_n$ and $\psi_n$ induce the same homomorphism $\chi_n:F_n\to\psl$.
Of course $Stab_\chi$ contains $Stab_\phi$ and $Stab_\psi$ and thus our claim
is proved as soon as we can show the opposite inclusions.\\
First note that the braid action defined on homomorphisms as above is
equivalent to the Hurwitz action on the tuples of images of the specified
generators $t_i\in F_n$, hence the braid action will not change the conjugation
class of these images.\\
Now let $\b$ be a braid in $Stab_\chi$. Then
$\phi\circ\b(t_i)=(ab)^{3k_i}\phi(t_i)$ since $(ab)^3$ is the fundamental
element of $\br_3$ which generates the center of $\br_3$ and thus the
kernel of the extension $\br_3\to\psl$.
The degree homomorphism $d:\br_3\to\ZZ$ is a class function with value one on
all $\phi(t_i)$, hence $d((ab)^{3k_i})=0$. Since $d(ab)=2$ we conclude
$k_i=0$ and $\b\in Stab_\phi$.\\
Similarly we have $\psi\circ\b(t_i)=\pm\psi(t_j)$ for $\b\in Stab_\chi$.
Since the trace is a class function on $\slz$ which has value $2$ on all
$\psi(t_i)$ while it is $-2$ on $-\psi(t_i)$, we also get $\b\in Stab_\psi$.
\qed

\begin{prop}
\labell{slstab}
Let $F_n:=\langle t_i,1\leq i\leq n\,|\quad\rangle$ be the free group on
$n=l+l'$ generators, define a homomorphism $\psi_{l,l'}:F_n\to\slz$ by
$$
\psi_{l,l'}(t_i)=\left\{
\begin{array}{ll}
\left(\begin{smallmatrix}1&1\\0&1\end{smallmatrix}\right) & i>l, i\not\equiv l
\mod 2\\[2mm]
\left(\begin{smallmatrix}\phantom{-}1&0\\-1&1\end{smallmatrix}\right)
& i>l,i\equiv l \mod 2\\[2mm]
\left(\begin{smallmatrix}-1&\phantom{-}0\\\phantom{-}0&-1\end{smallmatrix}
\right) & i\leq l
\end{array}\right.
$$
and let $\br_n$ act on homomorphisms $F_n\to\slz$ by Hurwitz automorphisms of
$F_n$.
Then the stabilizer group $Stab_{\psi_{l,l'}}$ of $\psi_{l,l'}$
is generated by the image of $Stab_{\psi_{l'}}$ under the inclusion
$\br_{l'}\hookrightarrow\br_n$ mapping to braids with only the last $l'$
strands
braided and
$$
E_{l,l'}:=\left\langle \s_{ij}^{m_{ij}},\,1\leq i<j\leq n\,\left|\quad
m_{ij}=\left\{
\begin{array}{lcl}
1 & \text{if} & j\leq l\vee i\equiv j (2),i>l\\
2 & \text{ if } & i\leq l<j\\
3 & \text{ if } & i>l,i\not\equiv j(2)
\end{array}
\right.\right.\right\rangle.
$$
If $l'\leq6$ then even $Stab_{\psi_{l,l'}}=E_{l,l'}$.
\end{prop}

\proof
Again we argue with the equivalent Hurwitz action on images of the generators.
First we consider the induced action on conjugacy classes.
On $n$-tuples of conjugacy classes the Hurwitz action induces an action of
$\br_n$ through the natural homomorphism $\pi$ to the permutation group $S_n$.
Since the tuple induced from $\psi$ consists of $l$ copies of the conjugacy
class of $-id$ followed by $l'$ copies of the distinct conjugacy class of
$\psi(t_1)$, the associated stabilizer group is $\tilde E:=\pi\inv(S_{l}\times
S_{l'})$, and as in \cite{klui} one can check that
\begin{eqnarray*}
\tilde E & = & \langle \s_{ij}, i<j\leq l\text{ or }l<i<j;
\tau_{ij}:=\s_{ij}^2, i\leq l<j\rangle.
\end{eqnarray*}
So as a first step we have $Stab_\psi$ contained in $\tilde E$.\\

Since $-id$ is central it is the only element in its conjugacy class and we may
conclude that the $\tilde E$ orbit of $\psi$ contains only homomorphisms which
map the first $l$ generators onto $-id$. With a short calculation using that
$-id$ is a central involution we can check that the $\tau_{ij}$ act
trivially on
such elements:
$$
\begin{array}{cl}
& \tau_{ij}  (-id,...,-id,M_{l+1},...,M_{n})\\
= & \s_{i+1}\inv\cdots\s_{j-1}\inv\s_j^2\s_{j-1}\cdots\s_{i+1}
(-id,...,-id,M_{l+1},...,M_{n})\\
= & \s_{i+1}\inv\cdots\s_{j-1}\inv\s_j^2
(-id,...,-id,M_{l+1},...,M_{j-1},-id,M_j,...,M_{n})\\
= & \s_{i+1}\inv\cdots\s_{j-1}\inv
(-id,...,-id,M_{l+1},...,M_{j-1},-id,M_j,...,M_{n})\\
= &  (-id,...,-id,M_{l+1},...,M_{n})
\end{array}
$$
Therefore given $\b\in\tilde E$ as a word $w$ in the generators
$\s_{ij},\tau_{ij}$ of $\tilde E$ the action of $\b$ on $\psi$ is the same as
that of $\b'$ where $\b'$ is given by a word $w'$ obtained from $w$ by dropping
all letters $\tau_{ij}$.
By the commutation relations of the $\s_{ij}$ we may collect all letters
$\s_{ij},i,j\leq l$ to the right of letters $\s_{ij},i,j>l$ without changing
$\b'$ and get a factorization $\b'=\b_1'\b_2'$ with
$\b_1'\in\br_{l},\b_2'\in\br_{l'}$.\\

Hence $\b\in\tilde E$ acts trivially on $\psi$ if and only if $\b_1'\b_2'$ does
so if and only if $\b_2'$ acts trivially on $\psi_{l'}$.
Thus $Stab_{\psi_{l,l'}}$ is generated by the $\tau_{ij}$ the
$\s_{ij},i,j\leq l$
and the $\b_2'\in Stab_{\psi_{l'}}$.
Both conclusions of the proposition then follow since $\s_{ij},i,j>l$ are
contained in $Stab_{\psi_{l'}}$ and since they are even generators if
$l'\leq6$,
prop.\ \ref{l}.
\qed

\section{mapping class groups of elliptic fibrations}

We return to elliptic fibrations and obtain some results concerning mapping
classes of elliptic fibrations. We need in fact to enrich the structure a bit:

\begin{description}

\item[Definition:]
A {\it marked elliptic fibration} is an elliptic fibration with a distinguished
regular fibre, $f:X,F\to B,b_0$, which can be thought of as given by a {\it
marking} $F\hookrightarrow E$.

\item[Definition:]
A {\it fibration preserving map} of a marked elliptic surface $f:X,F\to B,b_0$
is a homeomorphism $\Phisl_X$ of $X$ such that
$f\circ\Phisl_X=\bdiff_{B,b_0}\circ f$
for a homeomorphism $\bdiff_{B,b_0}$ of $(B,b_0)$ and such that
$\Phisl_X|_F$ is
isotopic to the identity on $F$.\\
The map $\bdiff_{B,b_0}$ is called the {\it induced base homeomorphism}.

\end{description}

An induced homeomorphism necessarily preserves the set $\Delta(f)$ of singular
values of the fibration map $f$ and therefore can be regarded as a
homeomorphism of
the punctured base $B,\Delta(f)$ preserving the base point.\\

On the other hand with each elliptic fibration $f:X\to B$ we have a torus
bundle
over $B-\Delta(f)$. Its structure homomorphism is the natural map
$$
\psi:\pi_1(B,b_0)\tto \diff(F)
$$
to the group of isotopy classes of diffeomorphisms of the distinguished fibre.

\begin{prop}
\labell{moi-equiv}
Given a marked elliptic fibration and a braid $\b$ representing an isotopy
class
of homeomorphisms of its punctured base, there is a fibration preserving map
inducing $\b$ if and only if $\b$ stabilizes the structure map of the
associated
torus bundle.
\end{prop}

\proof
A fibration preserving homeomorphism $\Phi$ of an unmarked elliptic surface
induces
a map $\bdiff_B$ of the punctured base. By the classification
of torus bundles there exists then a commutative diagram
\begin{eqnarray*}
\pi_1(B-\Delta(f),b_0) & \stackrel{(\bdiff_B)_*}{-\!\!\!-\!\!\!\tto}
& \pi_1(B-\Delta(f),\bdiff_B(b_0))\\
\downarrow\psi_{b_0} & & \downarrow\psi_{\bdiff(b_0)}\\
\diff(F) & \stackrel{(\Phisl|_F)_*}{-\!\!\!-\!\!\!-\!\!\!\tto}
& \diff(\Phisl(F))
\end{eqnarray*}
But the result of Moishezon \cite[p.\ 169]{Mo1} implies that the reverse
implication is true in the absence of multiple fibres.\\
If now $\Phisl$ is a fibration preserving homeomorphism of a marked elliptic
surface
then the bottom map is the identity and the claim is immediate.
\qed

\begin{prop}
\labell{moi-normal}
Given a marked elliptic fibration with only fibres of types $I_1,I^*_0$ there
is a choice of free generators for $\pi_1(B,b_0)$, an isomorphism
$\diff(F)\cong\slz$ and an isomorphism of the abstract braid
group onto the mapping class group such that the structure homomorphism of the
associated bundle is $\psi_{6k,l}$ and such that the action of $\br_n$ on $F_n$
commutes with the action of $\diff(B,\Delta)$ on $\pi_1$.
\end{prop}

\proof
The proof proceeds along the lines of Moishezon's proof, cf.\ \cite{fm}, for the
normal form of an elliptic surface with only fibres of type $I_1$. The same
strategy leads to our claim since fibres of type $I_0^*$ have local
monodromy in
the center of $\slz$.
\qed

By now we have finally got all necessary results to prove theorem 1 as stated
in the introduction.\\

\proof[ of theorem \ref{1}]
The mapping class group of the base punctured at the base points of singular
fibres is isomorphic to the braid group $\sbr_{6+l}$ of the sphere.
We have previously shown that the mapping classes induced by fibration
preserving maps are those acting trivially on the structure homomorphism of the
torus bundle given with the elliptic fibration, prop.\ \ref{moi-equiv}.

By prop.\ \ref{moi-normal} and prop.\ \ref{slstab} the corresponding group is
conjugation equivalent to
$E_{6,l}$. On the other hand so is the monodromy group by the main theorem.
Moreover for each braid of the monodromy group there is by prop.\ \ref{lift} a
fibration preserving diffeomorphism, so we get an inclusion and hence both
groups
coincide.
\qed



\chapter{braid monodromy and fundamental groups}


In this concluding chapter we want to address four topics. First we will deduce
by the method of van Kampen a presentation for the fundamental group of the
discriminant complement of a Brieskorn Pham polynomial. Second we want to
relate the algebraic monodromy and the braid monodromy of Brieskorn Pham
polynomials by means of Dynkin diagrams.

In a third section we give all corollaries for arbitrary singular functions
which can be deduced immediately from our results and we finish in a four
section with some conjectures and speculations.

\section{fundamental groups}

For convenience we restate our result on generators of the braid monodromy.

\begin{thm}
\label{mainref}
The braid monodromy group of a Brieskorn-Pham polynomial\\
$x_1^{l_1+1}+\cdots
x_n^{l_n+1}$ is generated by the following twist powers:
$$
\begin{array}{ccll}
\s_{i,j}^3 &:& i<j & \text{ correlated}\\
\s_{i,j}^2 &:& i<j & \text{ not correlated}\\
\s_{j,k}\s_{i,j}^2\s_{j,k}\inv &:& i<j<k &
\text{ correlated}
\end{array}
$$
\end{thm}

The most important corollary drawn from this theorem is a presentation of the
fundamental group of the discriminant complement which can be computed by the
Zariski van Kampen method.

\begin{thm}
\labell{mainfund}
The fundamental group of the discriminant complement of a versal unfolding of a
Brieskorn-Pham polynomial $x_1^{l_1+1}+\cdots x_n^{l_n+1}$ has a presentation
given with respect to the multiindex set $I=I(l_1,...,l_n)$ of cardinality
$\mu=l_1\cdots l_n$:
$$
\begin{array}{rcl}
\lspan t_i,\, i\in I& | & t_i t_j t_i=t_j t_i t_j,\qquad i,j\in I,\, i<j\text{
correlated,}\\
&& t_i t_j=t_j t_i,\quad\qquad i,j\in I,\, i<j\text{ not correlated,}\\
&& t_it_jt_kt_i=t_jt_kt_it_j,\qquad i,j,k\in I,\, i<j<k\text{ correlated}\,
\rspan
\end{array}
$$
\end{thm}

\proof
A presentation of the fundamental group can be obtained from generators of the
braid monodromy according to \ref{cpl-pres}.

So we have generators $t_i$, $i\in I_n$ in bijection to the critical points.
We obtain the relations from the generators of the braid monodromy group, which
are given in theorem \ref{mainref}. In fact a generator of the first two rows
can be factored as $\b_0\s_1^{2}\b_0\inv$ and $\b_0\s_1^{3}\b_0\inv$
respectively in such a way that
$$
\b_0(t_1)=t_i,\quad\b_0(t_2)=t_j.
$$
and similarly generators of the last row can be conjugated such that
$$
\b_0(t_1)=t_i,\quad\b_0(t_2)=t_jt_kt_j\inv .
$$
We then compute a sufficient set of relators using lemma \ref{few-rela}:
\begin{enumerate}
\item
in case $\s_{i,j}^3,\,i<j$ correlated:
$$
t_i\inv\b(t_i)=t_i\inv t_j\inv t_i\inv t_jt_it_j,t_j\inv\b(t_j)=t_j\inv
t_j\inv t_i\inv t_j\inv t_it_jt_it_j.
$$
\item
in case $\b=\s_{i,j}^2,\,i<j$ not correlated:
$$
t_i\inv\b(t_i)=t_i\inv t_j\inv t_it_j,t_j\inv\b(t_j)=t_j\inv t_i\inv
t_i\inv t_jt_it_j.
$$
\item
in case $\s_{j,k}\s_{i,j}^2\s_{j,k}\inv,\, i<j<k$ correlated:
\begin{eqnarray*}
t_i\inv\b(t_i)&=&t_i\inv (t_jt_kt_j\inv )\inv t_it_jt_kt_j\inv\\
& = & t_i\inv t_jt_k\inv t_j\inv t_it_jt_kt_j\inv,\\
(t_jt_kt_j\inv )\inv\b(t_jt_kt_j\inv )&=&(t_jt_kt_j\inv )\inv
(t_jt_kt_j\inv )\inv t_i\inv (t_jt_kt_j\inv )t_i(t_jt_kt_j\inv )\\
& = & t_jt_k\inv t_j\inv t_jt_k\inv t_j\inv t_i\inv t_jt_kt_j\inv
t_it_jt_kt_j\inv.
\end{eqnarray*}
\end{enumerate}

In all cases the relators are conjugate, so we can do with the following
relations:
\begin{eqnarray*}
& t_it_jt_i \,=\, t_jt_it_j & \text{for } i,j \text{ correlated},\\
& t_it_j \,=\, t_jt_i & \text{for } i,j \text{ not correlated},\\
& t_it_jt_kt_j\inv  \,=\, t_jt_kt_j\inv t_i & \text{for } i,j,k \text{
correlated.}
\end{eqnarray*}
Any relation of the third kind we can multiply by $t_jt_i$ on the
right. Using a
relation of the first kind and cancellation of inverse letters, we arrive at
our claim:
\begin{eqnarray*}
t_it_jt_kt_j\inv  &=& t_jt_kt_j\inv t_i\\
\iff\quad t_it_jt_kt_i &=& t_jt_kt_j\inv t_it_jt_i\\
\iff\quad t_it_jt_kt_i &=& t_jt_kt_j\inv t_jt_it_j\\
\iff\quad t_it_jt_kt_i &=& t_jt_kt_it_j.
\end{eqnarray*}
\qed

\section{Dynkin diagrams}

We want to interpret the theorems of the first section in terms of the Dynkin
diagrams which Pham associated to the functions under consideration.

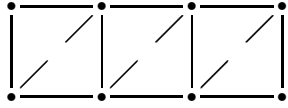
\begin{center}
\setlength{\unitlength}{6mm}
\begin{picture}(10,3)(0,-1.5)

\put(2,1){\circle*{.2}}
\put(4,1){\circle*{.2}}
\put(6,1){\circle*{.2}}
\put(8,1){\circle*{.2}}
\put(2,-1){\circle*{.2}}
\put(4,-1){\circle*{.2}}
\put(6,-1){\circle*{.2}}
\put(8,-1){\circle*{.2}}

\drawline(2.2,1)(3.8,1)
\drawline(4.2,1)(5.8,1)
\drawline(6.2,1)(7.8,1)
\drawline(2.2,-1)(3.8,-1)
\drawline(4.2,-1)(5.8,-1)
\drawline(6.2,-1)(7.8,-1)
\drawline(2.2,-.8)(2.8,-.2)
\drawline(4.2,-.8)(4.8,-.2)
\drawline(6.2,-.8)(6.8,-.2)
\drawline(3.2,.2)(3.8,.8)
\drawline(5.2,.2)(5.8,.8)
\drawline(7.2,.2)(7.8,.8)
\drawline(2,-.8)(2,.8)
\drawline(4,-.8)(4,.8)
\drawline(6,-.8)(6,.8)
\drawline(8,-.8)(8,.8)

\end{picture}

{the Dynkin diagram of $x_1^3+x_2^5$ according to Pham}
\end{center}

So let us recall the situation for the simple singularities. There we have the
distinguished Dynkin diagram with no cycles and the corresponding generalised
Artin groups. To these Brieskorn Artin groups we can associate the braid
stabilizer subgroup, which in fact coincides with the braid subgroup
generated by
elements associated to edges and pairs without edge.\\
We generalize this situation in the following way.

\begin{defi}
The braid subgroup associated to a Dynkin diagram obtained from a distinguished
system of paths is defined as follows:
\begin{itemize}
\item
To each edge of weight $\pm1$ associate the generator $\s_{i,j}^3$,
\item
To each non-connected vertex pair associate the generator $\s_{i,j}^2$,
\item
To each edge triangle of weight product $-1$ associate the generator\\
$\s_{j,k}\s_{i,j}^2\s_{j,k}$, $i<j<k$.
\end{itemize}
In each case the indices are those of the vertices involved.
\end{defi}

One can check that we can describe the Dynkin digram obtained by Pham for the
function
$$
x_1^{l_1+1}+\cdots x_n^{l_n+1}
$$
using the multiindex set $I=I(l_1,...,l_n)$.
The pair of vertices of indices $i=i_1i_2\ldots i_n$ and $j=j_1j_2\ldots j_n$
are connected by an edge if $i<j$ are correlated and we assign the weight
$$
(-1)^{\textstyle 1+\sum_\nu(j_\nu-i_\nu)}.
$$
Then we can show the following relation.

\begin{lemma}
\labell{dynkin}
The braid monodromy group of a Brieskorn Pham polynomial in dimension $n$, $f=x_1^{l_1+1}+\cdots
x_n^{l_n+1}$, is given by the braid subgroup of $\br_{l_1\cdots l_n}$ generated
by the elements associated to the Dynkin diagram of $f$.
$$
\begin{array}{cll}
\s_{i,j}^3 & & \text{there is an edge between vertices $i,j$}\\
\s_{i,j}^2 & & \text{there is no edge between vertices $i,j$}\\
\s_{j,k}\s_{i,j}^2\s_{j,k}\inv & &
\text{there is a triangle of weight product $-1$ on the vertices $i,j,k$}
\end{array}
$$
\end{lemma}

\proof
The first two rows of generators in \ref{mainref} are obviously associated to
the diagram. Now in case there is any triangle of edges the indices $i<j<k$ are
correlated and the weight product is
$$
(-1)^{\textstyle
3+\sum_\nu(j_\nu-i_\nu)+\sum_\nu(k_\nu-i_\nu)+\sum_\nu(k_\nu-j_\nu)} = -1.
$$
So there is a bijection between generators given in the theorem and the
generators associated to the Dynkin diagram.
\qed

In a similar way the presentation of the fundamental group of the discriminant
complement can be expressed in terms of the Dynkin diagram.


\section{other functions}

Any versal family of functions is induced from the versal family of a Brieskorn
Pham function, so we may draw some conclusions from this fact:

\begin{lemma}
Suppose a singular function f is adjacent to a Brieskorn Pham
polynomial $\tilde
f$, then its braid monodromy group is contained in the intersection of the
braid monodromy group of $\tilde f$ with some group in the conjugacy class of
$\br_\mu\subset\br_{\tilde\mu}$.
\end{lemma}

The reverse adjacency relation implies:

\begin{lemma}
Suppose a Brieskorn Pham polynomial $g$ is adjacent to a singular function $f$,
then the braid monodromy group of $f$ contains a subgroup isomorphic to the
braid monodromy group of $g$.
\end{lemma}

The corollary concerning the fundamental group of the discriminant complements
can be formulated as follows.

\begin{lemma}
Suppose a singular function $f$ is adjacent to a Brieskorn Pham
polynomial $\tilde f$ with Milnor numbers $\mu$ and $\tilde\mu$ respectively,
then the fundamental group of its discriminant complement
fits into a commutative diagram
\begin{eqnarray*}
F_\mu & \inj & F_{\tilde\mu}\\
\downarrow & & \downarrow\\[-4mm]
\downarrow & & \downarrow\\
\pi_1(\CC^\mu\smin\dfami_f) & \tto & \pi_1(\CC^{\tilde\mu}\smin\dfami_{\tilde
f}).
\end{eqnarray*}
\end{lemma}


\section{conjectures and speculations}

Now that we have defined and computed the braid monodromy groups of Brieskorn
Pham polynomials we are in a position to make an educated guess what
they are in
case of general functions. We formulate two approaches both based on the fact
that the braid monodromy determines the fundamental group of the discriminant
complement, which in turn is restricted by the algebraic or geometric
monodromy.

Since the braid monodromy acts by the Artin representation on the algebraic
monodromy homomorphism, the immediate guess is:

\begin{conj}
The braid monodromy group is the stabilizer subgroup of the algebraic monodromy
homomorphism.
\end{conj}

The other is much more explicit and relies on the fact that the algebraic
monodromy is encoded into the Dynkin diagram. It presumes that the braid
subgroups associated to Dynkin diagrams in an orbit under the braid group
action belong to a single conjugacy class.

\begin{conj}
The braid monodromy group of a singular function coincides with the braid
subgroup associated to a Dynkin diagram up to conjugacy.
\end{conj}

It leads immediately to a preferred finite presentation of the fundamental
group of the discriminant complement.

\begin{conj}
\labell{dynkin-pres}
The fundamental group of the discriminant complement of a versal unfolding of a
singular function $f$ has a presentation
given with respect to a Dynkin diagram of $f$ of vertex cardinality
$\mu$:
$$
\begin{array}{rcl}
\lspan t_i,\, 1\leq i\leq\mu & | & t_i t_j t_i=t_j t_i t_j,\qquad i,j\in I,\,
i<j\text{ joint by an edge of weight $\pm1$,}\\
&& t_i t_j=t_j t_i,\quad\qquad i,j\in I,\, i<j\text{ not joint by an edge,}\\
&& t_it_jt_kt_i=t_jt_kt_it_j,\qquad i,j,k\in I,\, i<j<k\text{ in a edge
triangle}\\
&& \hspace*{\fill}\text{ of weight product $-1$}\,
\rspan
\end{array}
$$
\end{conj}

In any case the braid monodromy group and the fundamental group would be shown
to be topological invariants despite the observations by Brieskorn in \cite{Br}
which show a dependence of the braid monodromy homomorphism on analytic
invariants.\\

But we are not only interested into the fundamental groups themselves but also
in their relations. Even if we knew that the fundamental groups are determined
from the Dynkin diagrams, which implies that the braid monodromy groups embed
for adjacent singular functions, it is by no means obvious what we speculate
now:

\begin{conj}
Suppose $f$ is adjacent to $g$, so a versal unfolding of $g$ is a versal
unfolding for $f$, then the natural map of fundamental groups of discriminant
complements injects.
\begin{eqnarray*}
F_\mu & \inj & F_{\tilde\mu}\\
\downarrow & & \downarrow\\[-4mm]
\downarrow & & \downarrow\\
\pi_1(\CC^\mu\smin\dfami_f) & \tto & \pi_1(\CC^{\tilde\mu}\smin\dfami_{\tilde
f}).
\end{eqnarray*}
\end{conj}

The final speculation seems even further off. We have speculated on the image
of the braid monodromy homomorphism. Now we consider the domain. It is the
fundamental group of the bifurcation complement. So the question might be, how
much of this fundamental group is captured by the braid monodromy group.
Boldly put:

\begin{conj}
The braid monodromy group is the fundamental group of the complement of the
bifurcation diagram.
\end{conj}

Here again one should be cautious, since the observation allude to above,
implies that the generic plane section yields curves of different multiplicity
in cases of equal topological invariants.



\begin{appendix}
\chapter{braid computations}

\newcommand{\ul}{\underline}

This appendix is designed to serve several purposes. First the progress in the
chapters is eased if some of the computational obstacles are hidden in this
appendix. Second the arguments are often similar and it is easier to get used
to them, if they are used in one place instead of being scattered throughout.

We distinguish the cases that the indices are just natural numbers, or pairs,
or multiindices, and we add a final section, which contains a very helpful
criterion to show, when two set of elements generate the same group.\\

We have to introduce some notation, which partly applies only for a single
section. Since all index sets we use are ordered, we can always denote
by $i^+$ the immediate successor of $i$ in some index set. The same notation
applies also to single components of multiindices.

As before our most general index set is defined in terms of possible exponents
of Brieskorn Pham polynomials.

\begin{nota}
Given a finite sequence $l_1,...,l_n$ of positive integers, define the
multiindex set $\idx=\idx(l_1,...,l_n)$ to be
\begin{eqnarray*}
\idx  & := & \{i_1...i_n\,|\,1\leq i_\nu\leq l_\nu,\,1\leq\nu\leq n\}
\end{eqnarray*}
equipped with the natural lexicographical order.
\end{nota}

Since we mostly deal with conjugates we underline a conjugated element to make
the structure more obvious.

\section*{simple indices}

In this section $i',j'$ etc.\ are used just to denote additional natural
numbers.

\begin{remark}
The twist $\check\s_{i,k}:=(\prod_{i<j<k}\s_{i,j}^2)\ul{\s_{i,k}}
(\prod_{i<j<k}\s_{i,j}^2)\inv$ is the twist on the horizontal arc from $i$ to
$k$ passing behind all intermediate punctures (as opposed to the arc of
$\s_{i,k}$ which passes in front).
\end{remark}

\begin{remark}
\label{sig}
The half twist on the arc from $i$ to $k$ passing in front up
to $j$ and behind from $j+1$ onwards can be given as
\begin{eqnarray*}
(\prod_{j< j'<k}\s^2_{i,j'})\ul{\s_{i,k}}(\prod_{j< j'<k}\s^2_{i,j'})\inv.
\end{eqnarray*}

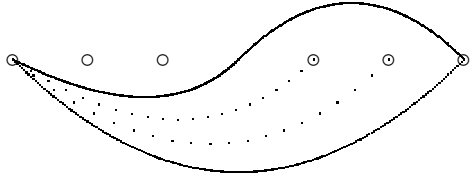
\begin{center}
\setlength{\unitlength}{5mm}
\begin{picture}(20,7)(-4,-5)

\put(0,0){\circle{.3}}
\put(2,0){\circle{.3}}
\put(4,0){\circle{.3}}

\put(8,0){\circle{.3}}
\put(10,0){\circle{.3}}
\put(12,0){\circle{.3}}

\bezier{200}(0,0)(6,-6)(12,0)
\bezier{20}(0,0)(5.5,-4.5)(10,0)
\bezier{20}(0,0)(4.8,-3.2)(8,0)
\bezier{200}(0,0)(4,-2)(6,0)
\bezier{200}(6,0)(9,3)(12,0)

\end{picture}
\end{center}
\end{remark}

\begin{lemma}
\labell{iso/tri/gen}
The braid subgroup $\br(A_n)\subset\br_n$ is generated by elements
\begin{eqnarray*}
& \s_{i,i+1}^3 & 1\le i<n\\
& \s_{i,i^+}\invv\check\s_{i,j}^2\s_{i,i^+}^2 & 1<i,i^+<j\leq n.
\end{eqnarray*}
\end{lemma}

\proof
Consider the following two filtered sets of elements of $\br_n$.
\begin{eqnarray*}
& S_1:=\{\s_{i,i+1}^3\} & S_k:= S_1\cup \{\s^2_{i,j}|1<j-i\leq k\},\\
& T_1:=\{\s_{i,i+1}^3\} & T_k:=
T_1\cup \{\s_{i,i^+}\invv\check\s_{i,j}^2\s_{i,i^+}^2|1<j-i\leq k\}
\end{eqnarray*}
By the first remark we get the relation
\begin{eqnarray*}
\s_{i,i^+}\invv\check\s_{i,j}^2\s_{i,i^+}^2 & = &
(\prod_{i^+< j'<j}\s^2_{i,j'})\ul{\s_{i,j}}(\prod_{i^+< j'<j}\s^2_{i,j'})\inv.
\end{eqnarray*}
so $S_2=T_2$ and
the other hypotheses of lemma \ref{filt} hold as well. Therefore the assertion
is proved, since $S_n$ is known to generate $\br(A_n)$.
\qed

\begin{lemma}
The subgroup $\br_{A(l)}$ of $\br_l$ is generated by
\begin{enumerate}
\item
$\s_{i,i+1}^3$,
\item
$\s_{i,i+2}^2$,
\item
$(\prod\limits_{j=i+2}^{k-1}\s_{i,j}^2)\ul{\s_{i,k}^2}
(\prod\limits_{j=i+2}^{k-1}\s_{i,j}^2)\inv,\quad i+1<k-1$.
\end{enumerate}
\end{lemma}

\proof
This is the same result as \ref{iso/tri/gen} just with different notation.
\qed

\begin{lemma}
\labell{braid/b}
Suppose in a punctured disc two otherwise distinct arcs meet in the punctures
$p,q$ thus defining a inner disc. If there is a system of arcs such that
\begin{enumerate}
\item
each puncture in the inner disc is connected by an arc with either $p$ or $q$,
\item
apart from $p,q$ all arcs have no points in common,
\end{enumerate}
then the twist on the outer arcs are equal up to conjugation by full twists on
the inner arcs.
\end{lemma}

\proof
We may identify the mapping class group of neighbourhood of the inner disc with
the abstract braid group, such that the twists on inner arcs correspond to
$\s_{1,j}^2,1<j\leq m$ and $\s_{j,n}^2, m<j<n$ and the twist on the outer arcs
correspond to $\s_{1,n}$ and $\check\s_{1,n}(1)$.
The claim then follows from
\begin{eqnarray*}
(\prod_{j=2}^m\s_{1,j}^2)\ul{\s_{1,n}} (\prod_{j=2}^m\s_{1,j}^2)\inv & = &
(\prod_{j=m+1}^{n-1}\s_{j,n}^2)\ul{\check\s_{1,n}}
(\prod_{j=m+1}^{n-1}\s_{j,n}^2)\inv.
\end{eqnarray*}

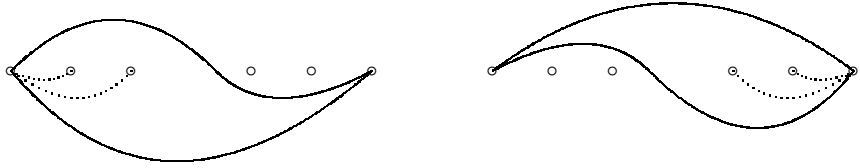
\begin{center}
\setlength{\unitlength}{4mm}
\begin{picture}(25,6)

\put(0,0){\begin{picture}(20,6)(1,-3)

\put(0,0){\circle{.3}}
\put(2,0){\circle{.3}}
\put(4,0){\circle{.3}}

\put(8,0){\circle{.3}}
\put(10,0){\circle{.3}}
\put(12,0){\circle{.3}}

\bezier{200}(0,0)(3.4,3.4)(6.8,0)
\bezier{12}(0,0)(1.2,-.6)(2,0)
\bezier{200}(6.8,0)(8.6,-1.8)(12,0)
\bezier{20}(0,0)(2.4,-1.8)(4,0)
\bezier{300}(0,0)(5,-6)(12,0)

\end{picture}}

\put(10,0){\begin{picture}(20,6)(-5,-3)

\put(0,0){\circle{.3}}
\put(2,0){\circle{.3}}
\put(4,0){\circle{.3}}

\put(8,0){\circle{.3}}
\put(10,0){\circle{.3}}
\put(12,0){\circle{.3}}

\bezier{200}(12,0)(9,-3.8)(5.2,0)
\bezier{12}(12,0)(10.8,-.6)(10,0)
\bezier{200}(5.2,0)(3.4,1.8)(0,0)
\bezier{20}(8,0)(9.6,-1.8)(12,0)
\bezier{300}(0,0)(6,4.5)(12,0)

\end{picture}}
\end{picture}
\end{center}

\qed

\begin{lemma}
\labell{braid/x}
For any $j$, $i<j\leq k$ the twist $\check\s_{i,k}(i)$ can be given as
$\s_{i,k}$ suitably conjugated by braids $\s_{i',j'}$, $i\leq i'<j\leq j'\leq
k$.
\end{lemma}

\proof
First note that $\s_{i,k}$ and $\check\s_{i,k}$ are twists on arc which meet
the initial hypothesis of lemma \ref{braid/b}. By the second remark above, the
full twists on arcs from the puncture of index $i'$ passing in front up to
$j-1$ and behind from $j$ onwards is in the group generated by elements
$\s_{i',k}^2,\s_{i',j'}^2$, $i<i'<j\leq j'<k$.

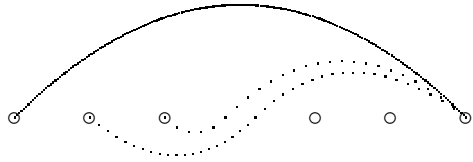
\begin{center}
\setlength{\unitlength}{5mm}
\begin{picture}(20,6)(-4,-3)

\put(0,0){\circle{.3}}
\put(2,0){\circle{.3}}
\put(4,0){\circle{.3}}

\put(8,0){\circle{.3}}
\put(10,0){\circle{.3}}
\put(12,0){\circle{.3}}

\bezier{20}(2,0)(4.4,-2)(6.4,0)
\bezier{5}(4,0)(4.8,-.8)(5.6,0)
\bezier{20}(6.4,0)(8.8,2.4)(12,0)
\bezier{20}(5.6,0)(8.6,3.0)(12,0)
\bezier{200}(0,0)(6,6)(12,0)

\end{picture}
\end{center}

On the other hand these arcs and the arcs to which the $\s_{i,j'}^2$, $j\leq
j'<k$ are associated can be chosen simultaneously to meet the remaining
hypotheses of lemma \ref{braid/b}. So we get our claim.
\qed

\section*{pair indices}

In this section we will not use primed integer variables except for $i_2'$
which always is only a shorthand for $i_2-l_2$.

\begin{lemma}
\labell{iso/zero/conj}
The elements $\tau_{i_1i_2',j_1i_2'}$ and $\tau_{i_1i_2',j_1i_2}$ are equal up
to conjugation by twists $\tau^2_{i_1i_2,j_1j_2}$, $1\leq j_2-i'_2<l_2$ for all
$1\leq i_1<j_1\leq l_1$, $i_2'=i_2-l_2$.
\end{lemma}

\proof
In fact we have
$$
(\prod_{i_2'<j_2<i_2}\tau^2_{i_1i_2',j_1j_2})\inv\ul{\tau_{i_1i_2',j_1i_2'}}
(\prod_{i_2'<j_2<i_2}\tau^2_{i_1i_2',j_1j_2}) =
\tau_{i_1i_2',j_1i_2}.
$$

\begin{center}
\setlength{\unitlength}{4.5mm}
\begin{picture}(24,7)(0,-3)

\put(0,0){\circle{.3}}
\put(2,0){\circle{.3}}
\put(4,0){\circle{.3}}
\put(6,0){\circle{.3}}

\bezier{2}(12,0)(13,0)(14,0)

\put(18,0){\circle{.3}}
\put(22,0){\circle{.3}}
\put(24,0){\circle{.3}}
\put(20,0){\circle{.3}}

\bezier{200}(6.8,-0)(5,-1.8)(2,0)
\bezier{20}(7.6,0)(5.2,-2.4)(2,0)
\bezier{20}(8.4,0)(5.4,-3.0)(2,0)
\bezier{29}(2,0)(5.6,-3.6)(9.2,0)

\bezier{120}(16,0)(18.2,-1.1)(20,0)
\bezier{130}(10,0)(12,2)(16,0)
\bezier{129}(2,0)(5.8,-4.2)(10,0)

\bezier{30}(9.2,0)(12.2,3)(18,0)
\bezier{50}(8.4,0)(12.8,4.4)(24,0)

\bezier{50}(7.6,0)(13,5.4)(24,0.7)
\bezier{530}(6.8,0)(13.2,6.4)(24,1.4)

\bezier{180}(20,0)(23,-3)(25,-1)
\bezier{12}(22,0)(23.5,-1.5)(24.5,-.5)

\bezier{120}(25,-1)(26.3,0.3)(24,1.4)
\bezier{6}(24.5,-.5)(25.2,.2)(24,.7)

\end{picture}
\end{center}

Again the claim can also be proved by checking the hypotheses of lemma
\ref{braid/b}.
\qed

\begin{lemma}
\labell{tri/conj}
Up to conjugation by twists $\tau^2_{i_1i_2,j_1j_2}$, $1\leq i_1<j_1\leq
l_1$, $1\leq j_2-i_2<l_2$,
elements $(i_2':=i_2-l_2)$
\begin{eqnarray*}
\tau_{i_1i_2',i^+_1i_2'}\inv\tau_{i^+_1i_2',j_1i_2'}^2\tau_{i_1i_2',i^+_1i_2'}
& \text{if} & j_1-i_1\leq l_1/2,\\
\tau_{i_1i_2',i^+_1i_2'}\inv\tau_{i^+_1i_2',j_1i_2}^2\tau_{i_1i_2',i^+_1i_2'}
& \text{if} & j_1-i_1\geq l_1/2+1,\\
\tau_{i_1i_2',i^+_1i_2}\inv\tau_{i^+_1i_2',j_1i_2}^2\tau_{i_1i_2',i^+_1i_2}
& \text{if} & j_1-i_1=l_1/2+1/2.
\end{eqnarray*}
correspond bijectively to the elements
\begin{eqnarray*}
\tau_{i_1i_2',i^+_1i_2}\inv\tau_{i^+_1i_2',j_1i_2}^2\tau_{i_1i_2',i^+_1i_2}
&  & 1<i_1+1=i_1^+<j_1\leq l_1, 1\leq i_2\leq l_2.
\end{eqnarray*}
\end{lemma}

\proof
The elements of the third row already have the claimed factorisation.
In the other cases we have to conjugate in such a way that central twist and
conjugating twist are conjugated simultaneously to the claimed twists.\\
For elements of the second row, it is only the conjugating twist which has not
the claimed form.

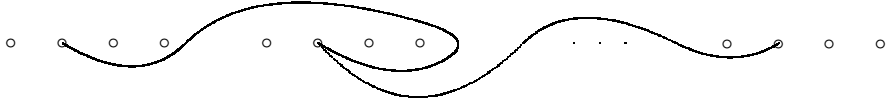
\begin{center}
\setlength{\unitlength}{3.4mm}
\begin{picture}(34,7)(-10,-3)

\put(-10,0){\circle{.3}}
\put(-8,0){\circle{.3}}
\put(-6,0){\circle{.3}}
\put(-4,0){\circle{.3}}

\put(0,0){\circle{.3}}
\put(2,0){\circle{.3}}
\put(4,0){\circle{.3}}
\put(6,0){\circle{.3}}

\bezier{2}(12,0)(13,0)(14,0)

\put(18,0){\circle{.3}}
\put(22,0){\circle{.3}}
\put(24,0){\circle{.3}}
\put(20,0){\circle{.3}}

\bezier{120}(16,0)(18.2,-1.1)(20,0)
\bezier{130}(10,0)(12,2)(16,0)
\bezier{129}(2,0)(5.8,-4.2)(10,0)

\bezier{200}(-3.2,-0)(-5,-1.8)(-8,0)
\bezier{180}(2,0)(5,-1.8)(7,-.6)
\bezier{120}(7,-.6)(8.3,0.18)(6,.84)
\bezier{230}(-3.2,0)(-.5,2.7)(6,.84)

\end{picture}
\end{center}

The proof of \ref{iso/zero/conj} shows, that we can get
$$
(\prod_{i_2'<j_2<i_2}\tau^2_{i_1i_2',j_1j_2})\inv\ul{\tau_{i_1i_2',i^+_1i_2'}}
(\prod_{i_2'<j_2<i_2}\tau^2_{i_1i_2',j_1j_2}) =
\tau_{i_1i_2',i^+_1i_2}.
$$
Moreover we check at once, that the conjugating factor commutes with
$\tau_{i^+_1i_2',j_1i_2}$, hence an overall conjugation of
$\tau_{i_1i_2',i^+_1i_2'}\inv\tau_{i^+_1i_2',j_1i_2}^2\tau_{i_1i_2',i^+_1i_2'}$
with $(\prod_{i_2'<j_2<i_2}\tau^2_{i_1i_2',j_1j_2})$ yields
$\tau_{i_1i_2',i^+_1i_2}\inv\tau_{i^+_1i_2',j_1i_2}^2\tau_{i_1i_2',i^+_1i_2}$.\\

In case of elements of the first row, neither twist is in the right shape:

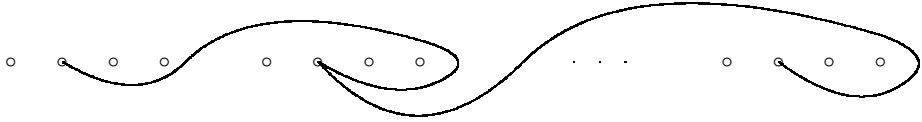
\begin{center}
\setlength{\unitlength}{3.4mm}
\begin{picture}(34,7)(-10,-3)

\put(-10,0){\circle{.3}}
\put(-8,0){\circle{.3}}
\put(-6,0){\circle{.3}}
\put(-4,0){\circle{.3}}

\put(0,0){\circle{.3}}
\put(2,0){\circle{.3}}
\put(4,0){\circle{.3}}
\put(6,0){\circle{.3}}

\bezier{2}(12,0)(13,0)(14,0)

\put(18,0){\circle{.3}}
\put(22,0){\circle{.3}}
\put(24,0){\circle{.3}}
\put(20,0){\circle{.3}}

\bezier{229}(2,0)(5.8,-4.2)(10,0)
\bezier{530}(10,0)(14,4)(24,1.05)
\bezier{120}(25,-.75)(26.3,0.225)(24,1.05)
\bezier{180}(20,0)(23,-2.25)(25,-.75)

\bezier{200}(-3.2,-0)(-5,-1.8)(-8,0)
\bezier{180}(2,0)(5,-1.8)(7,-.6)
\bezier{120}(7,-.6)(8.3,0.18)(6,.84)
\bezier{230}(-3.2,0)(-.5,2.7)(6,.84)

\end{picture}
\end{center}

We first take care of the middle factor $\tau_{i^+_1i_2',j_1i_2'}$, as before
just by conjugation with $(\prod_{i_2'<j_2<i_2}\tau^2_{i^+_1i_2',j_1j_2})$.

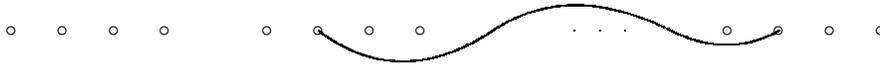
\begin{center}
\setlength{\unitlength}{3.4mm}
\begin{picture}(34,7)(-10,-3)

\put(-10,0){\circle{.3}}
\put(-8,0){\circle{.3}}
\put(-6,0){\circle{.3}}
\put(-4,0){\circle{.3}}

\put(0,0){\circle{.3}}
\put(2,0){\circle{.3}}
\put(4,0){\circle{.3}}
\put(6,0){\circle{.3}}

\bezier{2}(12,0)(13,0)(14,0)

\put(18,0){\circle{.3}}
\put(22,0){\circle{.3}}
\put(24,0){\circle{.3}}
\put(20,0){\circle{.3}}


\bezier{220}(8.8,0)(5.2,-2.4)(2,0)

\bezier{120}(15.8,0)(18,-1.1)(20,0)
\bezier{130}(8.8,0)(11.8,2)(15.8,0)
\end{picture}
\end{center}

But since we have to conjugate overall, also the conjugating factors are
conjugated, and they are not unaffected:

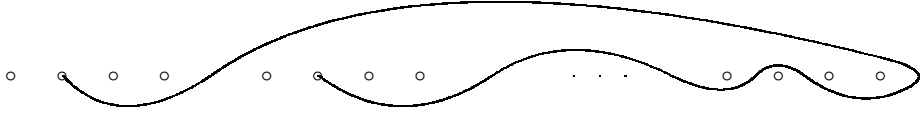
\begin{center}
\setlength{\unitlength}{3.4mm}
\begin{picture}(34,7)(-10,-3)

\put(-10,0){\circle{.3}}
\put(-8,0){\circle{.3}}
\put(-6,0){\circle{.3}}
\put(-4,0){\circle{.3}}

\put(0,0){\circle{.3}}
\put(2,0){\circle{.3}}
\put(4,0){\circle{.3}}
\put(6,0){\circle{.3}}

\bezier{2}(12,0)(13,0)(14,0)

\put(18,0){\circle{.3}}
\put(22,0){\circle{.3}}
\put(24,0){\circle{.3}}
\put(20,0){\circle{.3}}

\bezier{120}(-2.2,0)(-5.8,-2.4)(-8,0)

\bezier{220}(8.8,0)(5.2,-2.4)(2,0)

\bezier{120}(15.8,0)(18,-1.1)(19.1,0)
\bezier{130}(8.8,0)(11.8,2)(15.8,0)

\bezier{100}(19.1,0)(19.9,.8)(21,0)

\bezier{530}(-2.2,0)(5.4,5.4)(24,.7)
\bezier{120}(25,-.5)(26.3,0.15)(24,.7)
\bezier{180}(21,0)(23,-1.5)(25,-.5)

\end{picture}
\end{center}

But we can now conjugate by $(\prod_{i_2'<j_2<i_2}\tau^2_{i_2+i_1',j_1j_2})$
and $(\prod_{i_2'<j_2<i_2}\tau^2_{i_1i_2',i^+_1j_2})$, which are twists on arcs
disjoint to that of $\tau_{i^+_1i_2',j_1i_2}$, hence commutating with it.
\qed

\begin{lemma}
\labell{ff(st)}
Given $1\leq i_1<k_1\leq l_1,1\leq k_2<l_2$ the elements
\begin{eqnarray*}
\tau^2_{i_10,k_1k_2}
& \text{ and } &
\s^2_{i_1l_2,k_1k_2}
\end{eqnarray*}
coincide up to conjugation by elements
\begin{enumerate}
\item
$\s^2_{1l_2,j_1j_2}$, $i_1<j_1<k_1$, $1\leq j_2<l_2$,
\item
$\tau^2_{j_10,k_1k_2}$, $i_1<j_1<k_1$.
\end{enumerate}
\end{lemma}

\proof
It suffices to check that arcs to which the given twists are associated can be
chosen simultaneously in such a way that
\begin{enumerate}
\item
they are confined to the disc with boundary given by the arcs corresponding to
$\tau_{i_10,k_1k_2}$ and $\s_{i_1l_2,k_1k_2}$,
\item
they are distinct outside the punctures of indices $i_1l_2$ and $k_1k_2$,
\item
all punctures in the disc are joint by some arc with either the puncture of
index $i_1l_2$ or that of index $k_1k_2$.
\end{enumerate}
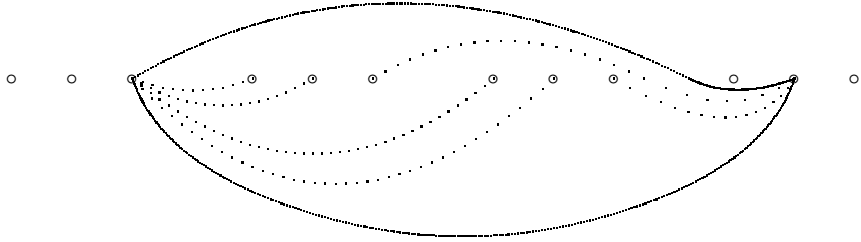
\begin{center}
\setlength{\unitlength}{4mm}
\begin{picture}(25,8)(1,-5)

\put(0,0){\circle{.3}}
\put(2,0){\circle{.3}}
\put(4,0){\circle{.3}}

\put(8,0){\circle{.3}}
\put(10,0){\circle{.3}}
\put(12,0){\circle{.3}}

\put(16,0){\circle{.3}}
\put(18,0){\circle{.3}}
\put(20,0){\circle{.3}}

\put(26,0){\circle{.3}}
\put(24,0){\circle{.3}}
\put(28,0){\circle{.3}}

\bezier{20}(12,0)(16,2.5)(21,0)
\bezier{8}(21,0)(24,-1.5)(26,0)

\bezier{200}(4,0)(12.5,5)(22.5,0)
\bezier{80}(22.5,0)(24,-.75)(26,0)

\bezier{40}(4,0)(10,-5)(16,0)

\bezier{40}(4,0)(10.5,-7)(18,0)

\bezier{12}(4,0)(6,-.8)(8,0)

\bezier{20}(4,0)(7,-1.8)(10,0)

\bezier{16}(20,0)(24.4,-2.6)(26,0)

\bezier{100}(4,0)(5,-3)(10,-4.5)
\bezier{100}(26,0)(25,-3)(20,-4.5)
\bezier{100}(10,-4.5)(15,-6)(20,-4.5)

\end{picture}
\end{center}
So by lemma \ref{braid/b} we may conclude that the assertion holds.
\qed

\section*{multiindices}

In this section we reserve the notation $i'$ etc.\ for the multiindex
$i':=i_1...i_{n-1}$ in $\idx(l_1,...,l_{n-1})$ naturally associated to
$i=i_1...i_n$.

\begin{lemma}
\labell{braid/a}
Suppose the indices $i,j,k$ form a correlated triple, then
\begin{eqnarray*}
\s_{i,k}^3 & = &
(\s_{j,k}^3\s_{i,j}^2\s_{i,k}^2\s_{i,j}^{-2})\ul{\s_{i,j}^3}
(\s_{j,k}^3\s_{i,j}^2\s_{i,k}^2\s_{i,j}^{-2})\inv.
\end{eqnarray*}
\end{lemma}

\proof
Since we claim a conjugation relation it suffices to show that
\begin{eqnarray*}
\s_{i,k} & = &
\s_{j,k}\s_{i,k}\s_{j,k}\s_{i,k}\inv\s_{j,k}\inv\\
& = &
\s_{j,k}\s_{i,k}^2\s_{j,k}\s_{i,k}\s_{j,k}\inv\s_{i,k}^{-2}\s_{j,k}\inv\\
& = &
\s_{j,k}^3\s_{j,k}^{-2}\s_{i,k}^2\s_{j,k}^2\s_{i,j}\s_{j,k}^{-2}
\s_{i,k}^{-2}\s_{j,k}^2\s_{j,k}^{-3}\\
& = &
\s_{j,k}^3(\s_{i,j}^2\s_{i,k}^2\s_{i,j}\invv)\s_{i,j}
(\s_{i,j}^2\s_{i,k}^{-2}\s_{i,j}^{-2})\s_{j,k}^{-3}
\end{eqnarray*}
\qed

\begin{lemma}
\labell{braid/aa}
Suppose the indices $i,j,m,k$ form a correlated quadruple, then
\begin{eqnarray*}
\s_{i,m}^2\s_{i,k}^2\s_{i,m}\invv & = &
(\s_{j,m}^2\s_{j,k}^2\s_{j,m}\invv\s_{i,j}^2\s_{i,m}^2\s_{i,j}\invv)
\s_{i,j}^2\s_{i,k}^2\s_{i,j}\invv
(\s_{j,m}^2\s_{j,k}^2\s_{j,m}\invv\s_{i,j}^2\s_{i,m}^2\s_{i,j}\invv)\inv.
\end{eqnarray*}
\end{lemma}

\proof
It remains to prove that conjugation by
$\s_{j,m}^2\s_{j,k}^2\s_{j,m}\invv\s_{i,j}^2$ acts trivially on
$\s_{i,m}^2\s_{i,k}^2\s_{i,m}\invv$.

Since the twists
$\s_{j,m}^2\s_{j,k}^2\s_{j,m}\invv$, $\s_{i,j}^2$ and
$\s_{i,m}^2\s_{i,k}^2\s_{i,m}\invv$ correspond to three arcs which
form a triangle with no punctures in its
interior, our claim is the homomorphic image of the following claim in $\br_3$:
\begin{eqnarray*}
\s_{1,3}^2 & = & \s_{2,3}^2\s_{1,2}^2\s_{1,3}^2\s_{1,2}\invv\s_{2,3}\invv,
\end{eqnarray*}
which is immediately seen to be true, because $\s_{1,2}^2\s_{1,3}\s_{1,2}\invv
=\s_{2,3}\invv\s_{1,3}\s_{2,3}^2$.
\qed

\begin{lemma}
\labell{band=prod}
Given multiindices $i<k$ and $j'$ such that $i'<j'<k'$, then
\begin{eqnarray*}
\eta(\s_{i',j'}^2)\ul{\s_{i,k}}\,\eta(\s_{i',j'}\invv) & = &
\big(\prod_{1\leq j_n\leq l_n}\s_{i,j'j_n}^2\big)\ul{\s_{i,k}}
\big(\prod_{1\leq j_n\leq l_n}\s_{i,j'j_n}^2\big)\inv\\
& = &
\big(\prod_{j_n=i_n-l_n+1}^{i_n}\s_{i,j'j_n}^2\big)\ul{\s_{i,k}}
\big(\prod_{j_n=i_n-l_n+1}^{i_n}\s_{i,j'j_n}^2\big)\inv\\
& = &
\big(\prod_{j_n=k_n}^{k_n+l_n-1}\s_{i,j'j_n}^2\big)\inv\ul{\s_{i,k}}
\big(\prod_{j_n=k_n}^{k_n+l_n-1}\s_{i,j'j_n}^2\big).
\end{eqnarray*}
\end{lemma}

\proof
We start with the observation that the arcs corresponding to the
twists $\s_{i,k}$
and $\eta(\s_{i',j'}^2)\ul{\s_{i,k}}\,\eta(\s_{i',j'}\invv)$ bound a disc which
contains the punctures with indices $j'j_n$, $1\leq j_n\leq l_n$ and
only these.

The claim can now be deduced from \ref{braid/b}, since the conjugating braids
correspond to arc systems which meet the hypotheses of \ref{braid/b}, cf.\ the
following illustrations.

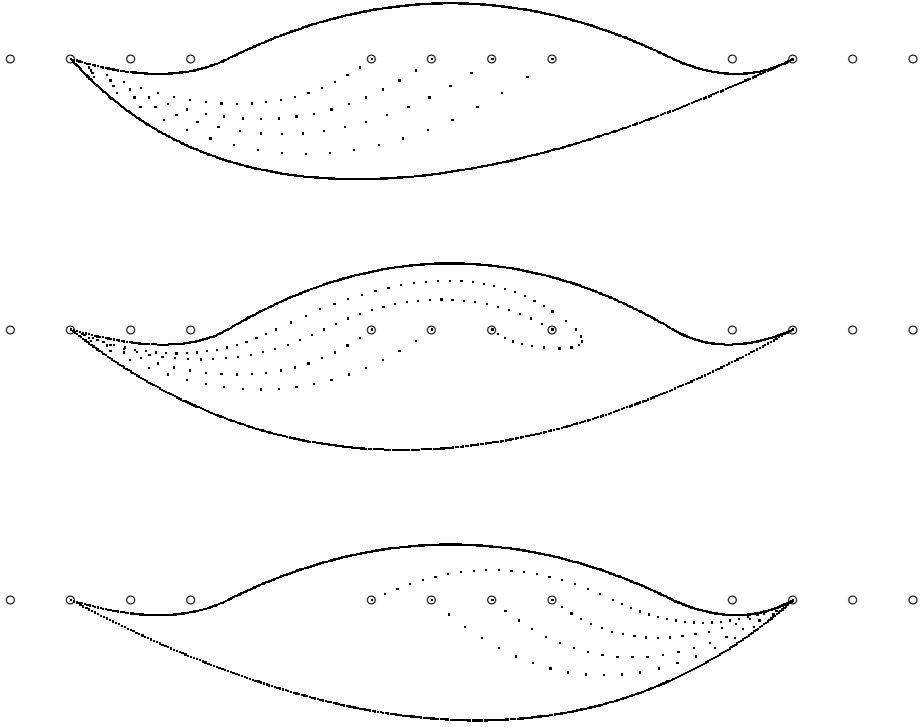
\begin{center}
\setlength{\unitlength}{4mm}
\begin{picture}(30,26)

\put(0,23){\begin{picture}(10,5)

\put(0,0){\circle{.3}}
\put(2,0){\circle{.3}}
\put(4,0){\circle{.3}}
\put(6,0){\circle{.3}}

\put(18,0){\circle{.3}}
\put(12,0){\circle{.3}}
\put(14,0){\circle{.3}}
\put(16,0){\circle{.3}}

\put(30,0){\circle{.3}}
\put(28,0){\circle{.3}}
\put(24,0){\circle{.3}}
\put(26,0){\circle{.3}}

\bezier{400}(2,0)(9,-8)(26,0)

\bezier{20}(2,0)(8.6,-4)(14,0)
\bezier{20}(2,0)(8,-3)(12,0)

\bezier{20}(2,0)(9,-5)(16,0)
\bezier{20}(2,0)(9.6,-6.3)(18,0)

\bezier{70}(2,0)(5.2,-1)(7.2,0)
\bezier{300}(7.2,0)(14.6,3.7)(22,0)
\bezier{100}(26,0)(24.2,-1)(22,0)

\end{picture}}

\put(0,14){\begin{picture}(10,5)

\put(0,0){\circle{.3}}
\put(2,0){\circle{.3}}
\put(4,0){\circle{.3}}
\put(6,0){\circle{.3}}

\put(18,0){\circle{.3}}
\put(12,0){\circle{.3}}
\put(14,0){\circle{.3}}
\put(16,0){\circle{.3}}

\put(30,0){\circle{.3}}
\put(28,0){\circle{.3}}
\put(24,0){\circle{.3}}
\put(26,0){\circle{.3}}

\bezier{300}(2,0)(12,-8)(26,0)

\bezier{20}(2,0)(8.6,-4)(14,0)
\bezier{20}(2,0)(8,-3)(12,0)

\bezier{20}(2,0)(6.4,-2)(10.4,0)
\bezier{20}(10.4,0)(14.4,2)(18,0)

\bezier{6}(16,0)(16.5,-.5)(17.7,-.6)
\bezier{8}(17.7,-.6)(20.1,-.8)(18,.6)
\bezier{23}(18,.6)(14.6,2.9)(8.8,0)
\bezier{20}(2,0)(5.6,-1.6)(8.8,0)

\bezier{70}(2,0)(5.5,-1)(7.2,0)
\bezier{300}(7.2,0)(14.6,4.4)(22,0)
\bezier{100}(26,0)(23.7,-1)(22,0)

\end{picture}}

\put(0,5){\begin{picture}(10,5)

\put(0,0){\circle{.3}}
\put(2,0){\circle{.3}}
\put(4,0){\circle{.3}}
\put(6,0){\circle{.3}}

\put(18,0){\circle{.3}}
\put(12,0){\circle{.3}}
\put(14,0){\circle{.3}}
\put(16,0){\circle{.3}}

\put(30,0){\circle{.3}}
\put(28,0){\circle{.3}}
\put(24,0){\circle{.3}}
\put(26,0){\circle{.3}}

\bezier{300}(2,0)(17,-8)(26,0)

\bezier{20}(14,0)(19.4,-5)(26,0)
\bezier{20}(16,0)(20.3,-3.8)(26,0)

\bezier{20}(18,0)(21,-2.5)(26,0)

\bezier{20}(20,0)(23,-1.5)(26,0)
\bezier{19}(20,0)(16,2)(12,0)

\bezier{70}(2,0)(5.2,-1)(7.2,0)
\bezier{300}(7.2,0)(14.6,3.7)(22,0)
\bezier{100}(26,0)(24.2,-1)(22,0)

\end{picture}}

\end{picture}
\end{center}
\qed

\begin{lemma}
\labell{uncoil-i}
Suppose $i<k$ are multiindices in $\idx =\idx (l_1,...,l_n)$ with $i'<k'$
correlated and assume $j'=i_1^+i_2...i_{n-1}$.
Then up to conjugation by elements $\s^2_{i'i_n,j'j_n}$, $1\leq i_n-j_n<l_n$,
$$
\eta(\s_{i',j'}^2)\ul{\s_{i,k}}\,\eta(\s_{i',j'}\invv)
$$
is equal to
$\s^2_{i,j'i_n}\ul{\s_{i,k}}\s_{i,j'i_n}\invv=
\s_{j'i_n,k}\invv\ul{\s_{i,k}}\s_{j'i_n,k}^2$.
\end{lemma}

\proof
By \ref{band=prod} the given elements are equal to
$$
\big(\prod_{j_n=i_n-l_n+1}^{i_n}\s_{i,j'j_n}^2\big)\ul{\s_{i,k}}
\big(\prod_{j_n=i_n-l_n+1}^{i_n}\s_{i,j'j_n}^2\big)\inv,
$$
hence they are equal to $\s^2_{i,j'i_n}\ul{\s_{i,k}}\s_{i,j'i_n}\invv$ up to
conjugation by elements $\s^2_{i'i_n,j'j_n}$, $1\leq i_n-j_n<l_n$. The claim
$\s^2_{i,j'i_n}\ul{\s_{i,k}}\s_{i,j'i_n}\invv
=\s_{j'i_n,k}\invv\ul{\s_{i,k}}\s_{j'i_n,k}^2$ is immediate.
\qed

\begin{lemma}
\labell{uncoil-k}
Suppose $i<k$ are multiindices in $\idx =\idx (l_1,...,l_n)$ with $i'<k'$
correlated and assume $j'=i_1k_2...k_{n-1}$.
Then up to conjugation by elements $\s^2_{j'j_n,k'k_n}$, $1\leq j_n-k_n<l_n$,
$$
\eta(\s_{i',j'}^2)\ul{\s_{i,k}}\,\eta(\s_{i',j'}\invv)
$$
is equal to
$\s^2_{i,j'k_n}\ul{\s_{i,k}}\s_{i,j'k_n}\invv=
\s_{j'k_n,k}\invv\ul{\s_{i,k}}\s_{j'k_n,k}^2$.
\end{lemma}

\proof
This result is obtained using the last equation of \ref{band=prod} along the
lines of the proof above, \ref{uncoil-i}.
\qed

\begin{lemma}
\labell{tau/sigma}
If $1\leq i_n<l_n$ and $j_n=0$ then
\begin{eqnarray*}
&(\prod\limits_{k_n=1}^{i_n^-}\s^2_{i'i_n,j'k_n})\inv
(\prod\limits_{k_n=i_n^+-l_n}^{-1}\s_{i'i_n,j'k_n}^2)\ul{\s_{i'i_n,j'j_n}}
(\prod\limits_{k_n=i_n^+-l_n}^{-1}\s_{i'i_n,j'k_n}^2)\inv
(\prod\limits_{k_n=1}^{i_n^-}\s^2_{i'i_n,j'k_n}),\\[1mm]
&(\delta_{j'i_n^+,j'l_n}^{l_n-i_n+1})\inv\s_{i'i_n,j'i_n}^2
\ul{\s_{i'i_n,j'i_n^+}}
\s_{i'i_n,j'i_n}\invv\delta_{j'i_n^+,j'l_n}^{l_n-i_n+1},
\end{eqnarray*}
are equal up to
conjugation by twists $\s_{i'i_n,j'k_n}$ with $k_n\neq j_n$,
$i_n-l_n<k_n<i_n$, and the sub-cable twist
$\delta_{j'i_n^+,j'l_n}^{l_n-i_n+1}$.\\

If $i_n=l_n+1$ and $1< j_n\leq l_n$ then
\begin{eqnarray*}
&(\prod\limits_{k_n=j_n^+}^{l_n}\s_{i'k_n,j'j_n}^2)
(\prod\limits_{k_n=i_n^+}^{j_n^-+l_n}\s^2_{i'i_n,j'k_n})\inv
\ul{\s_{i'i_n,j'j_n}}
(\prod\limits_{k_n=i_n^+}^{j_n^-+l_n}\s^2_{i'i_n,j'k_n})
(\prod\limits_{k_n=j_n^+}^{l_n}\s_{i'k_n,j'j_n}^2)\inv, \\[1mm]
&\delta_{i'1,i'j_n^-}^{j_n}\s_{i'j^-_n,i'j_n}^2\ul{\s_{i'j^-_n,j'j_n}}
\s_{i'j^-_n,i'j_n}\invv\delta_{i'1,i'j_n^-}^{-j_n},
\end{eqnarray*}
are equal up to conjugation by twists $\s_{i'k_n,j'j_n}$ with $k_n\neq
i_n$, $j_n<k_n<j_n+l_n$, and the sub-cable twist
$\delta_{i'1,i'j_n^-}^{j_n}$.
\end{lemma}

\proof
The claims are symmetric to each other under the symmetry induced by
the exchange
of indices $i'\leftrightarrow j'$, $i_n\rightarrow l_n+1-j_n$, $j_n\rightarrow
l_n+1-i_n$.\\[4mm]
The claim is illustrated in the case of $\s_{i'2,j'0}$ being conjugated to
$\s_{i'2,j'2}^2\ul{\s_{i'2,j'3}}\s_{i'2,j'2}\invv$.

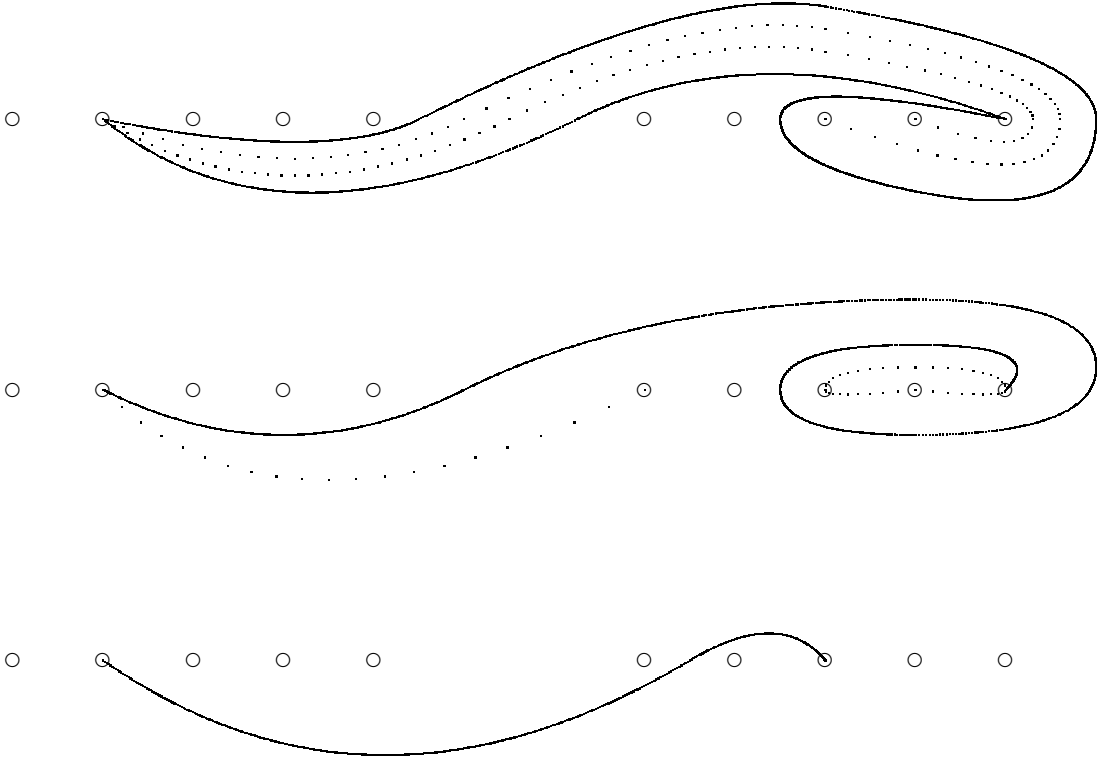
\begin{center}
\setlength{\unitlength}{6mm}
\begin{picture}(22,18)(0,0)

\put(0,15){\begin{picture}(20,4)

\put(0,0){\circle{.3}}
\put(2,0){\circle{.3}}
\put(4,0){\circle{.3}}
\put(8,0){\circle{.3}}
\put(6,0){\circle{.3}}

\put(22,0){\circle{.3}}
\put(14,0){\circle{.3}}
\put(16,0){\circle{.3}}
\put(18,0){\circle{.3}}
\put(20,0){\circle{.3}}

\bezier{200}(22,0)(17,1)(17,0)
\bezier{200}(17,0)(17,-1)(20,-1.6)
\bezier{200}(20,-1.6)(24,-2.4)(24,0)
\bezier{200}(18,2.5)(24,1.5)(24,0)
\bezier{180}(2,0)(7,-1)(9,0) 
\bezier{290}(9,0)(15,3)(18,2.5)

\bezier{240}(22,0)(16.5,2)(12.5,0)
\bezier{240}(12.5,0)(6,-3.25)(2,0)  

\bezier{20}(2,0)(6.5,-1.75)(10,0)
\bezier{20}(10,0)(15,2.5)(18,2)    
\bezier{20}(18,2)(23.2,1.15)(23.2,0)
\bezier{18}(18,0)(23.2,-2)(23.2,0)

\bezier{30}(2,0)(6,-2.5)(11,0)
\bezier{20}(11,0)(15,2)(18,1.5)     
\bezier{10}(20,0)(22.6,-1)(22.6,0)
\bezier{18}(18,1.5)(22.6,.75)(22.6,0)

\end{picture}}

\put(0,9){\begin{picture}(20,4)

\put(0,0){\circle{.3}}
\put(2,0){\circle{.3}}
\put(4,0){\circle{.3}}
\put(8,0){\circle{.3}}
\put(6,0){\circle{.3}}

\put(22,0){\circle{.3}}
\put(14,0){\circle{.3}}
\put(16,0){\circle{.3}}
\put(18,0){\circle{.3}}
\put(20,0){\circle{.3}}

\bezier{100}(22,0)(23,1)(20,1)
\bezier{100}(20,1)(17,1)(17,0)
\bezier{100}(17,0)(17,-1)(20,-1)
\bezier{100}(20,-1)(24,-1)(24,0.5)
\bezier{100}(24,0.5)(24,2)(20,2)
\bezier{200}(20,2)(14,2)(10,0)
\bezier{200}(10,0)(6,-2)(2,0)

\bezier{10}(22,0)(22,.5)(20,.5)
\bezier{10}(20,.5)(18,.5)(18,0)
\bezier{10}(18,0)(18,-.2)(20,0)
\bezier{10}(20,0)(22,-.2)(22,0)

\bezier{20}(2,0)(6,-4)(14,0)


\end{picture}}

\put(0,3){\begin{picture}(20,4)

\put(0,0){\circle{.3}}
\put(2,0){\circle{.3}}
\put(4,0){\circle{.3}}
\put(8,0){\circle{.3}}
\put(6,0){\circle{.3}}

\put(22,0){\circle{.3}}
\put(14,0){\circle{.3}}
\put(16,0){\circle{.3}}
\put(18,0){\circle{.3}}
\put(20,0){\circle{.3}}

\bezier{300}(2,0)(8,-4.2)(15,0)
\bezier{100}(15,0)(17,1.2)(18,0)

\end{picture}}

\end{picture}
\end{center}

The important point to note is the fact, that conjugation by
$(\prod\limits_{k_n=i_n^+-l_n}^{-1}\s_{i'i_n,j'k_n}^2)$ equals conjugation by
$(\delta_{j'i_n^+,j'l_n}^{l_n-i_n})$ and commutes with conjugation by
$(\prod\limits_{k_n=1}^{i_n^-}\s^2_{i'i_n,j'k_n})$ on the twists under
consideration.
Hence we are left to show that
\begin{eqnarray*}
&&(\prod\limits_{k_n=1}^{i_n^-}\s^2_{i'i_n,j'k_n})\inv
(\delta_{j'i_n^+,j'l_n})\ul{\s_{i'i_n,j'j_n}}
(\delta_{j'i_n^+,j'l_n})\inv
(\prod\limits_{k_n=1}^{i_n^-}\s^2_{i'i_n,j'k_n})\\
& = & \s_{i'i_n,j'i_n}^2
\ul{\s_{i'i_n,j'i_n^+}}
\s_{i'i_n,j'i_n}\invv
\end{eqnarray*}
Since the arcs of $(\delta_{j'i_n^+,j'l_n})\ul{\s_{i'i_n,j'j_n}}
(\delta_{j'i_n^+,j'l_n})\inv$ and $\s_{i'i_n,j'i_n^+}$ can be chosen to bound a
disc which contains the puncture of indices $j'1$ to $j'i_n$ we can conclude
that 
\begin{eqnarray*}
(\delta_{j'i_n^+,j'l_n})\ul{\s_{i'i_n,j'j_n}}
(\delta_{j'i_n^+,j'l_n})\inv
& = & (\prod\limits_{k_n=1}^{i_n}\s^2_{i'i_n,j'k_n})
\ul{\s_{i'i_n,j'i_n^+}}
(\prod\limits_{k_n=1}^{i_n}\s^2_{i'i_n,j'k_n})\inv
\end{eqnarray*}
from which we deduce the claim.
\qed

\section*{criterion for change of generators}

\begin{lemma}
\labell{filt}
Given two finitely filtered sets of elements of a group
\begin{eqnarray*}
S=S_n\supset S_{n-1}...\supset S_1,&&T_n\supset T_{n-1}...\supset T_1
\end{eqnarray*}
Then $S$ and $T$ generate the same subgroup if
\begin{enumerate}
\item
$T_1=S_1$,
\item
given $t\in T_k-T_{k-1}$ there is $s\in S_k$, such that $t$ is equal to $s$ up
to conjugation by elements in $\lspan S_{k-1}\rspan$,
\item
given $s\in S_k-S_{k-1}$ there is $t\in T_k$, such that $s$ is equal to $t$ up
to conjugation by elements in $\lspan S_{k-1}\rspan$.
\end{enumerate}
The last hypothesis may be replaced by
\begin{enumerate}
\item[iii')]
given $s\in S_k-S_{k-1}$ there is $t\in T_k$ such that $s$ is equal to $t$ up
to conjugation by elements in $\lspan T_{k-1},S_{k-1}\rspan$.
\end{enumerate}
\end{lemma}

\proof
We show $\lspan T_k\rspan =\lspan S_k\rspan$. So $i)$ starts the induction.
Then $\lspan T_k\rspan \subset\lspan S_k\rspan$ since by induction $\lspan
T_{k-1}\rspan \subset\lspan S_{k-1}\rspan\subset\lspan S_k\rspan$ and by $ii)$
$t\in T_k-T_{k-1}$ implies $t\in\lspan S_k\rspan$.

On the other hand by induction $\lspan S_{k-1}\rspan\subset\lspan
T_{k-1}\rspan$, therefore $s\in S_k-S_{k-1}$ implies $s\in\lspan
T_k,S_{k-1}\rspan\subset\lspan T_k\rspan$ if $iii)$ holds resp.\  if $iii')$
holds. In either case we get $\lspan S_k\rspan\subset\lspan T_k\rspan$.
\qed

\end{appendix}


\end{document}